\documentclass[a4paper]{amsart}

\theoremstyle{plain}
\newtheorem{theorem}{Theorem}[section]

\newtheorem{lemma}[theorem]{Lemma}
\newtheorem{corollary}[theorem]{Corollary}
\newtheorem{definition}{Definition}
\usepackage{amssymb, latexsym}
\usepackage{enumerate}

\begin{document}

\title{Projective sets, intuitionistically}
\author{Wim Veldman}
\address{Institute for Mathematics, Astrophysics and Particle Physics, Faculty of Science,
Radboud University,
Postbus 9010, 6500 GL Nijmegen, the Netherlands}
\email{ Wim.Veldman@ru.nl}

\begin{abstract}
\noindent {\small We try to develop intuitionistic descriptive set theory  and study `definable' subsets of Baire space $\mathcal{N}=\omega^\omega$.  The logic of our arguments is intuitionistic and we also use L.E.J.~Brouwer's  Thesis on bars in $\omega^\omega$ and his continuity axioms. We avoid the operation of taking the complement of a subset of $\mathcal{N}$ as much as possible, as the resulting sets, like negative statements, are not very useful in constructive mathematics.   

 A subset of $ \omega^\omega$ is \textit{(positively) projective} if it  results from a closed or an open subset of $\omega^\omega \times \omega^\omega (= \omega^\omega)$  by a finite number of applications of the two operations of \textit{projection} and \textit{universal projection} or: \textit{co-projection}.
 
A subset of $\omega^\omega$ is $\mathbf{\Sigma}^1_1$ or: \textit{analytic} if it is the  projection of a closed subset of $\omega^\omega$.

We give some examples of $\mathbf{\Sigma}^1_1$ subsets of $\omega^\omega$  like the set of (the codes of) all closed  subsets of $\omega^\omega$  that are \textit{positively uncountable} and also the set of (the codes of) all  closed subsets of $\omega^\omega$ containing an element coding  a (positively) infinite subset of  $\omega^\omega$.

   A subset of $\omega^\omega$ is  called \textit{strictly analytic} if it is the projection of a {\it spread}, i.e. a closed and {\it located} subset of $\omega^\omega$.
   
  Some analytic subsets of $\omega^\omega$ fail to be strictly analytic. 
  
  We will see that  Brouwer's Thesis on bars in $\omega^\omega$  proves separation and boundedness theorems for  strictly analytic
subsets of $\omega^\omega$.  

A subset of $\omega^\omega$ is called $\mathbf{\Pi}^1_1$ or:
\textit{co-analytic} if it is the co-projection of an open subset of $\omega^\omega \times \omega^\omega(=\omega^\omega)$. Most co-analytic sets are {\it not} the complement of an analytic set. There   is no  symmetry between analytic and co-analytic sets as there is in  classical descriptive set theory. 

As an example of a  $\mathbf{\Pi}^1_1$ set we consider  the set of the codes of all closed  subsets of $\omega^\omega$ all of whose members code an \textit{almost-finite} subset of $\omega$.

  We also study the set of the codes of closed and located subsets of $\omega^\omega$  that are \textit{almost-countable}, or, equivalently, \textit{reducible in Cantor's sense}.   This set is probably not $\mathbf{\Pi}^1_1$.

    Finally, we explain the important fact that the (positive) projective hierarchy \textit{collapses}:  every (positively) projective set  is $\mathbf{\Sigma}^1_2$ i.e.  the projection of a co-analytic subset of $\omega^\omega$.}
\end{abstract}

\maketitle

\section{Introduction}

This paper on descriptive set theory is one in a series.  We  explore the field of study laid bare by {\it pre-intuitionists}\footnote{Brouwer uses this term in \cite[p. 140]{brouwer52b} and \cite[p. 1]{brouwer1}.} like R.~Baire, \'E.~Borel, H.~Lebesgue, N.~Lusin and M.~Souslin,  and consider it from L.E.J.~Brouwer's {\it intuitionistic} point of view. 
In \cite{veldman08}, we proved an intuitionistic  Borel hierarchy theorem. In \cite{veldman09}, we discovered the  fine structure of the intuitionistic Borel hierarchy, and, in particular, the fine structure of the class $\mathbf{\Sigma}^0_2$, consisting of the countable unions of closed subsets of $\omega^\omega$. In both \cite{veldman08} and \cite{veldman09}, the argument is far from classical and  essential use is made of Brouwer's Continuity Principle. 

We now are going to treat projective sets. The earlier paper \cite{veldman05} already contains some  surprising results on apparently simple analytic and co-analytic subsets of $\omega^\omega$.

This introductory Section is divided into three parts.
In the first part, we briefly present the basic assumptions of intuitionistic analysis and we agree on a number of  notations.
In the second part, we introduce {\it intuitionistic descriptive set theory}. The reader may decide to skip these first two parts and use them only if further reading makes it necessary to consult them.
In the third part, we describe the further contents of the paper.

\subsection{The language and axioms of intuitionistic analysis}\hfill

The logical constants are used in their intuitionistic sense. A statement $P \vee Q$ is considered proven only if one either has a proof of $P$ or a proof of $Q$. A statement $\exists x \in V[P(x)]$ is considered proven only if one is able to produce an element $x$  of $V$ with a proof of the fact that $x$ has the property $P$.

  Brouwer not only refined the language of mathematics but also introduced a number of assumptions one should call {\it axiomatic}. He was of course the first to use them, see  \cite{brouwer18,brouwer0,brouwer1,brouwer2,brouwer3}. The question how to state and defend them has been further discussed 
by others, see \cite{heyting,howard,kleene,myhill,troelstra,troelstra0,
veldman1,veldman6,veldman08,veldman10, veldman21}. One finds them below in Subsubsections \ref{SSS:countablechoice}, \ref{SSS:bcpcontchoice}, \ref{SSS:fantheorem}, \ref{SSS:stumps}, \ref{SSS:barinduction} and \ref{SSS:creasubj}.

\subsubsection{Finite sequences of natural numbers}\label{SSS:finseq} \hfill

$\omega$ is the set of the natural numbers. 
We use $m,n,\ldots,s,t \ldots$ as variables over  $\omega$. 

$S:\omega\rightarrow \omega$ is the successor function: $\forall  n[S(n) = n+1]$.

\smallskip $p:\omega\rightarrow \omega$ is the function enumerating the primes: $p(0)= 2, p(1)=3, p(2)=5,\ldots$.

\smallskip We code finite sequences of natural numbers by natural numbers: 
$\langle\;\rangle :=0$ is the (code number of) the {\it empty sequence}, and, for all $k>0$, for all $m_0, m_1, \ldots m_{k-1}$,
 $\langle m_0, m_1, \ldots, m_{k-1}\rangle :=\prod_{i<k}p(i)^{m_i}\cdot p(k-1) -1$.

 \smallskip
 
$length(0):=0$ and, for each $s>0$,  $\mathit{length}(s) := 1\;+$
\textit{the largest $k$ such that  $p(k)$ divides $s+1$}.

\smallskip For each $s$, for each $i$, if $i<\mathit{length}(s)-1$, then
 $s(i):=$  \textit{the largest m such that  $p(i)^m$ divides  $s+1$}, and, if $i=\mathit{length}(s)-1$, then $s(i):=$ \textit{the largest m such that  $p(i)^{m+1}$ divides  $s+1$}, and, if $i\ge\mathit{length}(s)$, then $s(i):=0$.
 
 \smallskip Observe: for each $s,k$, if $\mathit{length}(s) = k$, then $s =\langle s(0), s(1), \ldots, s(k-1)\rangle$. 
 
 \smallskip For each $n$, $\omega^n:=\{s\mid length(s)=n\}$ and \\$[\omega]^n:=\{s\in \omega^n\mid\forall i[i+1< n\rightarrow s(i)<s(i+1)]\}$. 
 
 \smallskip $[\omega]^{<\omega}:=\bigcup_n [\omega]^n$.

 \smallskip
For all $s,t$, $s\ast t$ is the number $u$ satisfying: $\mathit{length}(u) =\mathit{length}(s) +\mathit{length}(t)$ and $\forall i<\mathit{length}(s)[u(i) =s(i]$ and $\forall j<\mathit{length}(t)[u(\mathit{length}(s) +j)=t(j)]$.

 \smallskip For all $s,n$ such that $n\le \mathit{length}(s)$, $\overline s(n) :=\overline s n:=\langle s(0), s(1), \ldots,s(n-1)\rangle$.
  
 \smallskip
 For all $s,t$: $s\sqsubseteq t\leftrightarrow \exists u[t=s\ast u]$ and: $s\sqsubset t \leftrightarrow (s\sqsubseteq t \;\wedge\;s\neq t)$ and:
 $s \sqsupset t \leftrightarrow t \sqsubset s$ and $s<_{lex}t\leftrightarrow \exists n[n<\mathit{length}(s) \;\wedge \;\overline s n \sqsubset t \;\wedge\; s(n) <t(n)]$ and: $s\perp t \leftrightarrow s\;\#\;t\leftrightarrow (s<_{lex}t\;\vee\;t<_{lex} s)$ and: $s<_{KB}t \leftrightarrow (t\sqsubset s\;\vee\;s<_{lex}t)$.
 
 $<_{KB}$ is a linear ordering of $\omega$, the \textit{Kleene-Brouwer-ordering}, also called the \textit{Lusin-Sierpinski-ordering}, see \cite[Section 2.G, p. 11]{kechris}.
 
 \smallskip For all $s,i$, $s^i$ is the number $u$ satisfying:\\ $\mathit{length}(u)=$ \textit{the least $k$ such that $\langle i \rangle \ast k \ge \mathit{length}(s)$} and $\forall j<\mathit{length}(u)[u(j) =s(\langle i \rangle \ast j)]$.
 
 \smallskip Note that, for each $i$,  $\langle\;\rangle^i= \langle \;\rangle$. 
 
 \emph{Note that also, for each $p$, for each $i$, $\langle p \rangle^i=\langle\;\rangle$}. 
 
 \smallskip For all $n,m$, $J(n,m):=(\langle n \rangle \ast m) -1$. 
 
 For each $n$, $K(n), L(n)$ are the numbers satisfying $n=J\bigl(K(n),L(n)\bigr)$.

\smallskip For all $s,t$ such that $\mathit{length}(s)=\mathit{length}(t)$ , $\ulcorner s, t \urcorner$ is the number $u$ satisfying \\$\mathit{length}(u)=\mathit{length}(s)$ and 
  $\forall i<\mathit{length}(s)[u(i)=J\bigl(s(i),t(i)\bigr)]$.
 
 For each $u$, $u_{I},u_{II}$ are the elements $s,t$ of $\omega$ such that  $u=\ulcorner s,t\urcorner$, i.e.\\ $\mathit{length}(u_I)=\mathit{length}(u_{II})=\mathit{length}(u)$ and
 \\$\forall i<\mathit{length}(u)[ u_I(i)=K\bigl(u(i)\bigr)  \;\wedge\; u_{II}(i)=L\bigl(u(i)\bigr)]$.

 For each $u$, \\$u_{I,I}:=(u_I)_I$ and: $u_{I,II}:=(u_I)_{II}$ and: 
$u_{II,I}:=(u_{II})_I$ and: $u_{II,II}:=(u_{II})_{II}$.

 \smallskip
 $Bin:=\{s\mid\forall i<\mathit{length}(s)[s(i)=0\;\vee\;s(i)=1]\}$ is the set of the codes of \textit{finite binary sequences}.

 For each $m$, $Bin_m:=\{s\in Bin\mid\mathit{length}(s) =m\}$.
 
 \smallskip For all $R\subseteq \omega$, $\forall m\forall n[mRn\leftrightarrow J(m,n) \in R]$.
 
 \smallskip For all $A,B\subseteq \omega$, $A\times B:=\{J(m,n)\mid m \in A, n \in B\}$. 
 
 \smallskip For all $A\subseteq \omega$, $n=\mu p[A(p)]$ if and only if $A(n)$ and $\forall p<n[\neg A(p)]$.

 \subsubsection{Infinite sequences of natural numbers}\label{SSS:infseq} \hfill
 
 \textit{Baire space} $\omega^\omega$ is the set of all infinite sequences of natural numbers. \\We use $\alpha,\beta, \ldots,\varphi, \psi, \ldots \sigma, \tau, \ldots$ as variables over $\omega^\omega$.
 
 An element of $\omega^\omega$ is a function from $\omega$ to $\omega$, and, given $\alpha,n$ we denote the result of applying $\alpha$ to $n$ by $\alpha(n)$. 
 
 \smallskip
 $[\omega]^\omega:=\{\zeta\mid\forall n[\zeta(n)<\zeta(n+1)]\}$.
 
 \smallskip For every $X\subseteq\omega$, $X^\omega :=\{\alpha\mid\forall n[\alpha(n)\in X]\}$.

 \smallskip
 For all $\alpha,\beta$, $\alpha\circ\beta$ is the element $\gamma$ of $\omega^\omega$ satisfying:  $\forall n[\gamma(n)=\alpha\bigl(\beta(n)\bigr)]$.

\smallskip For all $\alpha, t$, $\alpha\circ t$ is the number $u$ satisfying: 
 
 $\mathit{length}(u) =length(t)$ and $\forall n<length(t)[u(n) = \alpha\bigl(t(n)\bigr)]$.
 
 In particular, for each $t$, $S\circ t$ is the number $u$ satisfying: 
 
 $\mathit{length}(u) =length(t)$ and $\forall n<length(t)[u(n) = t(n)+1]$.

\smallskip

For all $\alpha, \beta$:  $\alpha \;\#\; \beta\leftrightarrow \alpha\perp\beta\leftrightarrow\exists n[\alpha(n) \neq \beta(n)]$, and: $\alpha = \beta\leftrightarrow \forall n[\alpha(n) = \beta(n)]$.

It is a well-known fact that the relation $\#$, called \textit{apartness}, is {\it co-transitive}, i.e.: \\for all $\alpha, \beta, \gamma$, if $\alpha\;\#\;\beta$, then either $\alpha\;\#\;\gamma$ or $\gamma\;\#\;\beta$. 

\smallskip
For each $s$, for each $\alpha$, $s\ast\alpha$ is the element $\gamma $ of $\omega^\omega$ such that  $\forall i<length(s)[\gamma(i) =s(i)]$ and $\forall i[\gamma\bigl(length(s) +i\bigr) =\alpha(i)]$. 

\smallskip For each $s$, for each $\mathcal{X}\subseteq\omega^\omega$, $s\ast\mathcal{X}:=\{s\ast\alpha\mid\alpha\in\mathcal{X}\}$.

\smallskip
For each $\alpha$, for each $n$, $\overline \alpha(n) :=\overline \alpha n:=\langle \alpha(0), \alpha(1), \ldots,\alpha(n-1)\rangle$.

$\overline \alpha (0):=\overline \alpha 0 :=\langle \;\rangle=0$.

For all $s,\alpha$: $s\sqsubset\alpha \leftrightarrow \exists n[s =\overline \alpha n]$ and: $s \perp\alpha\leftrightarrow \alpha \perp s \leftrightarrow \neg(s\sqsubset \alpha)$.

\smallskip Note: for all $a,b$, for all $\gamma$, if $a\perp b$, then either $a\perp \gamma $ or $\gamma\perp b$.

\smallskip For all $s$, $\omega^\omega\cap s:= \{\alpha\mid s\sqsubset\alpha\}$.

\smallskip
For each $m$, $\underline m$ is the element $\gamma$ of $\omega^\omega$ such that $\forall n[\gamma(n) = m]$.

For all $\alpha, i$, $\alpha^i$ is the element $\gamma$ of $\omega^\omega$ such that $\forall n[\gamma(n) =\alpha(\langle i\rangle \ast n)]$.

\smallskip For all $\alpha, m, n$, $\alpha^{m,n}:=(\alpha^m)^n$. 

Note: for all $m,n,p$, $\alpha^{m,n}(p)=\alpha(\langle m, n\rangle\ast p)$. 

\smallskip For all $\alpha$, for all $s$, $^s\alpha$ is the element $\gamma$ of $\omega^\omega$ such that $\forall n[\gamma(n)=\alpha(s\ast n)]$.

Note: $^{\langle m \rangle} \alpha = \alpha^m$. 

\smallskip
For every $\mathcal{X}\subseteq \omega^\omega$, $\mathcal{X}^\omega:=\{\alpha\mid\forall n[\alpha^n \in \mathcal{X}]\}$.

\smallskip
For all $\alpha,\beta$,  $\ulcorner \alpha,\beta\urcorner$ is the element $\gamma$ of $\omega^\omega$ such that $\forall n[\gamma(n) =J\bigl(\alpha(n),\beta(n)\bigr)]$.

For each $\gamma$, $\gamma_I, \gamma_{II}$ are the elements $\alpha,\beta$ of $\omega^\omega$ such that $\gamma =\ulcorner \alpha,\beta\urcorner$, that is: $\forall n[\gamma_I(n)=K\bigl(\gamma(n)\bigr) \;\wedge\;\gamma_{II}(n)=L\bigl(\gamma(n)\bigr)]$.

For each $\alpha$, $\alpha_{I,I}:=(\alpha_I)_I$ and: $\alpha_{I,II}:=(\alpha_I)_{II}$ and: 

$\alpha_{II,I}:=(\alpha_{II})_I$ and: $\alpha_{II,II}:=(\alpha_{II})_{II}$.

\smallskip For all $\mathcal{R}\subseteq \omega^\omega$, $\forall \alpha\forall \beta[\alpha\mathcal{R}\beta\leftrightarrow\ulcorner\alpha,\beta\urcorner \in\mathcal{R}]$.

 \smallskip For all $\mathcal{R}\subseteq \omega^\omega$,  $\forall \alpha\forall n[\alpha\mathcal{R}n \leftrightarrow n\mathcal{R}\alpha\leftrightarrow \langle n\rangle\ast\alpha \in\mathcal{R}]$.
 
\smallskip For all $\mathcal{A}\subseteq\omega^\omega, B\subseteq\omega$, $\mathcal{A}\times B:=B\times\mathcal{A} :=\{\langle n \rangle\ast\alpha\mid \alpha \in\mathcal{A}, n\in B\}$.
 
\smallskip For all $\mathcal{A},\mathcal{B}\subseteq \omega^\omega$, $\mathcal{A}\times \mathcal{B};=\{\ulcorner \alpha,\beta\urcorner\mid \alpha \in \mathcal{A}, \beta \in \mathcal{B}\}$.
 
\smallskip For all $\mathcal{A}\subseteq \omega^\omega$, for all $n$, $\mathcal{A}\upharpoonright n :=\{\alpha\mid\langle n \rangle \ast \alpha \in \mathcal{A}\}$. 

\smallskip For all $\mathcal{X}\subseteq \omega^\omega$, for all $n$, $\mathcal{X}_n:=\{\alpha\mid\langle n \rangle \ast \alpha \in \mathcal{X}\}$. \\An infinite sequence $\mathcal{X}_0, \mathcal{X}_1, \ldots$ of subsets of $\omega^\omega$ is the \textit{same} as the set \\$\mathcal{X}=\{\langle n\rangle\ast\alpha\mid n\in \omega, \alpha \in \mathcal{X}_n\}$. 

\smallskip For all $\mathcal{A}, \mathcal{B}\subseteq \omega^\omega$, \\$\mathcal{A} \subseteq \mathcal{B} \leftrightarrow \forall \alpha[\alpha\in \mathcal{A} \rightarrow \alpha \in \mathcal{B}]$, and: \\$\mathcal{A} \subsetneq \mathcal{B}\leftrightarrow \bigl(\mathcal{A}\subseteq \mathcal{B} \;\wedge\;\neg(\mathcal{B}\subseteq \mathcal{A})\bigr)$ and: \\$\mathcal{A}=\mathcal{B}\leftrightarrow (\mathcal{A}\subseteq \mathcal{B}\;\wedge\;\mathcal{B}\subseteq \mathcal{A})$, and: \\$\mathcal{A} \neq \mathcal{B}\leftrightarrow \neg(\mathcal{A} =\mathcal{B})$.

\smallskip For all  $\mathcal{X}_0,\mathcal{X}_1\subseteq \omega^\omega$, $\mathcal{X}_0 \;\#\;\mathcal{X}_1 \leftrightarrow \forall \alpha[\forall i<2[\alpha^i \in \mathcal{X}_i]\rightarrow \alpha ^0 \;\#\;\alpha^1]$. 

If $\mathcal{X}_0\;\#\;\mathcal{X}_1$, then $\mathcal{X}_0\cap\mathcal{X}_1=\emptyset$, but the converse may fail to be true. 

\smallskip For every infinite sequence $\mathcal{X}_0, \mathcal{X}_1, \ldots$ of subsets of $\omega^\omega$, we define:\\  $\#_n(\mathcal{X}_n) \leftrightarrow \forall \alpha[\forall n[\alpha^n \in \mathcal{X}_n]\rightarrow \exists i\exists j[\alpha ^i \;\#\;\alpha^j]]$. 

If $\#_n(\mathcal{X}_n)$,  then $\bigcap_n \mathcal{X}_n=\emptyset$, but the converse may fail to be true.

\smallskip
\textit{Cantor space} $2^\omega:=\{\alpha\mid\forall n[\alpha(n)<2]\}$. 

\smallskip
For each $\alpha$, \\$D_\alpha:=\{n\mid\alpha(n)\neq 0\}$ is the \textit{subset of $\omega$ decided by $\alpha$}, and: \\$E_\alpha:= \{m\mid\exists n[\alpha(n) =m+1]\}$ is the \textit{subset of $\omega$ enumerated by $\alpha$}.

\smallskip For each $s$, \\$D_s:=\{n<length(s)\mid s(n) \neq 0\}$ and \\$E_s:= \{m\mid\exists n<length(s)[s(n) =m+1]\}$.

Note: for each $\alpha$, $D_\alpha = \bigcup_n D_{\overline{\alpha}n}$ and: $E_\alpha = \bigcup_n E_{\overline{\alpha}n}$.

\smallskip
For each $X\subseteq \omega$, \\ $X$ is \textit{inhabited} if and only if $\exists n[n\in X]$ and: \\$X$ is  \textit{decidable} if and only if $\exists\alpha[X=D_\alpha]$ and: \\$X$ is \textit{enumerable} if and only if $\exists \alpha[X=E_\alpha]$. 

\smallskip
For each $\alpha$, $T_\alpha:=\{s\mid \forall t\sqsubset s[\alpha(t)=0]\}$. \\$T_\alpha$ is called the \textit{tree determined by $\alpha$}. Note: $\forall \alpha[0=\langle\;\rangle\in T_\alpha]$.

\smallskip
For all $\alpha,\beta$, for all $\gamma$, we define:\\
$\gamma:\alpha \le^\ast \beta \leftrightarrow$
$\bigl(\forall s[s\in T_\alpha\rightarrow \gamma(s) \in T_\beta]\;\wedge\;\forall s\forall t[s\sqsubset t \rightarrow \gamma(s) \sqsubset \gamma(t)]\bigr)$, and: \\
$\gamma:\alpha <^\ast \beta \leftrightarrow$
$\bigl(\forall s[s\in T_\alpha\rightarrow \gamma(s) \in T_\beta]\;\wedge\;\forall s\forall t[s\sqsubset t \rightarrow \gamma(s) \sqsubset \gamma(t)]\;\wedge\;\gamma(\langle\;\rangle)\neq\langle\;\rangle\bigr)$. 

For all $\alpha, \beta$,  we define:
$\alpha<^\ast\beta \leftrightarrow \exists \gamma[\gamma:\alpha<^\ast \beta]$, and: $\alpha\le^\ast\beta \leftrightarrow \exists \gamma[\gamma:\alpha\le^\ast \beta]$,

\smallskip
For  each $\delta$, $En_\delta:=\{\delta^n\mid n \in \omega\}$ is the subset of $\omega^\omega$ {\it enumerated by $\delta$}.

\subsubsection{Axioms of Countable Choice}\label{SSS:countablechoice} \hfill

First Axiom of Countable Choice:

\begin{quote}
$\mathbf{AC}_{0,0}$: {\it For all $R\subseteq \omega\times \omega$,  if $\forall m\exists n[mRn]$, then $ \exists \alpha\forall m[mR\alpha(m)]$.}
\end{quote}

\smallskip
 Second Axiom of Countable Choice:

\begin{quote}
$\mathbf{AC}_{0,1}$: {\it For all $\mathcal{R}\subseteq \omega^\omega\times \omega$,  if $\forall m\exists \alpha[m\mathcal{R}\alpha]$, then $\exists \alpha\forall m[mR\alpha^m]$.}
\end{quote}
\subsubsection{Open and closed subsets of $\omega^\omega$, and spreads}\label{SSS:openclosedspreads} \hfill

For each $\beta$, $\mathcal{G}_\beta:=\{\alpha\mid\exists n[\beta(\overline \alpha n)\neq 0]\}$ and $\mathcal{F}_\beta:=\{\alpha\mid\forall n[\beta(\overline \alpha n)= 0]\}$.

The pair of sets $(\mathcal{G}_\beta, \mathcal{F}_\beta)$ is called a \textit{complementary pair of rank 1}. 

For each  $\mathcal{X}\subseteq\omega^\omega$, 

\smallskip $\mathcal{X}$ is \textit{open} or $\mathbf{\Sigma}^0_1$ if  and only if $\exists \beta[\mathcal{X} =\mathcal{G}_\beta]$ and: 

\smallskip$\mathcal{X}$ is \textit{closed} or $\mathbf{\Pi}^0_1$ if  and only if $\exists \beta[\mathcal{X} =\mathcal{F}_\beta]$, and:
 
 \smallskip $\mathcal{X}$ is \textit{inhabited} if and only if $\exists \gamma[\gamma \in \mathcal{X}]$, and:

\smallskip $\mathcal{X}$ is \textit{located} if and only if $\exists \gamma[D_\gamma=\{ s\mid \exists \alpha \in \mathcal{X}[s\sqsubset \alpha]\}]$, and:

 \smallskip $\mathcal{X}$ is \textit{semi-located} if and only if $\exists \gamma[E_\gamma=\{ s\mid \exists \alpha \in \mathcal{X}[s\sqsubset \alpha]\}]$. 

\smallskip For every $\mathcal{X}\subseteq \omega^\omega$, $cl( \mathcal{X}):=\{\alpha\mid\forall n\exists \gamma \in \mathcal{X}[\overline \alpha n \sqsubset \gamma]\}$.  \\$cl(\mathcal{X})$ is called \textit{the closure of $\mathcal{X}$}.  $cl(\mathcal{X})$ is not necessarily $\mathbf{\Pi}^0_1$.\footnote{One may see this as follows. For every $\alpha$, define $\mathcal{Y}_\alpha:=\{\gamma\mid\gamma =\underline 0 \;\wedge\;\alpha \;\#\;\underline 0\}$ and note:  $\mathcal{Y}_\alpha:=cl(\mathcal{Y}_\alpha)$. Assume: every $\mathcal{Y}_\alpha$ is $\mathbf{\Pi}^0_1$. Then $\forall \alpha\exists \beta\forall \gamma[\gamma \in \mathcal{Y}_\alpha\rightarrow \gamma \in \mathcal{F}_\beta]$, and, therefore,\\ $\forall \alpha\exists \beta[\alpha\;\#\;\underline 0\leftrightarrow \forall n[\beta(\overline{\underline 0}n)=0]]$. Using Axiom $\mathbf{AC}_{1,1}$, see Subsubsection \ref{SSS:bcpcontchoice}, one may derive a contradiction.} 

\smallskip One easily proves: for every $\mathcal{X}\subseteq\omega^\omega$, $cl\bigl(cl(\mathcal{X})\bigr) =cl(\mathcal{X})$ and:  \\$\mathcal{X}$ is (semi-)located if and only if $cl(\mathcal{X})$ is (semi-)located.

\smallskip $\mathcal{F}\subseteq \omega^\omega$ is a \textit{spread} if and only if $cl(\mathcal{F})=\mathcal{F}$ and $\mathcal{F}$ is located.

\smallskip
For each $\beta$, we define:  {\it $\beta$ is a spread-law}, $Spr(\beta)$ if and only if \\$\beta \in 2^\omega$ and $ \forall s[ \beta(s)=0 \leftrightarrow \exists n[\beta(s\ast\langle n \rangle)=0]$.

One easily proves that  $\mathcal{F}\subseteq \omega^\omega$ is a spread if and only if $\exists \beta[Spr(\beta) \;\wedge\; \mathcal{F}=\mathcal{F}_\beta]$. 

Note: for all $\beta$, if $Spr(\beta)$, then $\mathcal{F}_\beta =\emptyset$ if and only if $\beta(0)=1$ if and only if $\beta =\underline 1$, and $\exists \gamma[\gamma \in \mathcal{F}_\beta]$ ($\mathcal{F}_\beta$ is \emph{inhabited}) if and only if $\beta(0)=0$. 

The empty set $\emptyset$ thus is a spread, and one may decide, for every spread $\mathcal{F}$, \\{\it either} $\mathcal{F}=\emptyset$ {\it or} $\exists\gamma[\gamma \in \mathcal{F}]$.

\smallskip Assume $Spr(\beta)$ and $\beta(c)=0$. We define: $\mathcal{F}_\beta\cap c:=\{\gamma\in\mathcal{F}_\beta\mid c\sqsubset \gamma\}$. \\Note that $\mathcal{F}_\beta\cap c$ itself is a spread.

\smallskip For each $\beta$, we define: {\it $\beta$ is a  perfect-spread-law}, $Pfspr(\beta)$,  if and only if \\$Spr(\beta)\;\wedge\; \beta(0)=0\;\wedge\;\forall s[\beta(s)=0\rightarrow\exists t\exists u[s\sqsubset t\;\wedge\;s\sqsubset u\;\wedge t\perp u\;\wedge\;\beta(t)=\beta(u)=0]]$. 

$\mathcal{F}\subseteq \omega^\omega$ is a \textit{perfect spread} if and only if $\exists \beta[Pfspr(\beta)\;\wedge\;\mathcal{F} =\mathcal{F}_\beta]$.

\subsubsection{Continuous functions}\label{SSS:continuousfunctions}\hfill

For all $\varphi,\alpha,m$, we define: {\it$ \varphi$ maps $\alpha$ onto $m$},  $\varphi:\alpha\mapsto m$, if and only if  \\$\exists n[\varphi(\overline{\alpha} n) =m +1 \;\wedge\;\forall i<n[\varphi(\overline \alpha i)=0]]$.

If $\exists m[\varphi:\alpha\mapsto m]$, we let $\varphi(\alpha)$ denote the unique $m$ such that $\varphi:\alpha\mapsto m$. 

For every $\mathcal{X}\subseteq \omega^\omega$, for all $\varphi$, we define: {\it $\varphi$ codes a function from $\mathcal{X}$ to $\omega$}, $\varphi:\mathcal{X}\rightarrow \omega$, if and only if  $ \forall \alpha \in \mathcal{X}\exists m[\varphi:\alpha \mapsto m]$.

\smallskip  $\varphi(\mathcal{X}):=\{m\mid\exists \alpha \in \mathcal{X}[\varphi:\alpha\mapsto m]\}=\{\varphi(\alpha)\mid\alpha\in\mathcal{X}\}$.  

\smallskip
For every $\mathcal{X}\subseteq \omega^\omega$, $\omega^\mathcal{X}:=\{\varphi\mid\varphi:\mathcal{X}\rightarrow \omega\}$.

\smallskip

For all $\varphi,\alpha, \beta$, we define: {\it $\varphi$ maps $\alpha$ onto $\beta$}, $\varphi:\alpha\mapsto \beta$, if and only if  \\$\forall n[\varphi^n:\alpha\mapsto \beta(n)]$.

If $\exists \beta[\varphi:\alpha\mapsto \beta]$, we let $\varphi|\alpha$ denote the unique $\beta$ such that $\varphi:\alpha\mapsto \beta$. 

For every $\mathcal{X},\mathcal{Y}\subseteq \omega^\omega$, for all $\varphi$, we define: {\it $\varphi$ maps $\mathcal{X}$ into $\mathcal{Y}$}, $\varphi:\mathcal{X}\rightarrow \mathcal{Y}$, if and only if  $\forall \alpha \in \mathcal{X}\exists \beta\in\mathcal{Y}[\varphi:\alpha \mapsto \beta]$.

 $\varphi|\mathcal{X}:=\{\beta\mid\exists \alpha \in \mathcal{X}[\varphi:\alpha\mapsto \beta]\}=\{\varphi|\alpha\mid\alpha\in\mathcal{X}\}$. 

\smallskip For all $\mathcal{X},\mathcal{Y}\subseteq \omega^\omega$, for all $\varphi$, we define: \textit{$\varphi$ embeds $\mathcal{X}$ into $\mathcal{Y}$}, $\varphi:\mathcal{X}\rightarrowtail\mathcal{Y}$,  if and only if $\varphi:\mathcal{X}\rightarrow \mathcal{Y}$ and $\forall \alpha\in\mathcal{X}\forall\beta\in\mathcal{X}[\alpha\;\#\;\beta\rightarrow \varphi|\alpha\;\#\;\varphi|\beta]$. 

$Emb(\mathcal{X},\mathcal{Y}):=\{\varphi\mid\varphi:\mathcal{X}\rightarrowtail\mathcal{Y}\}$.

For all $\mathcal{X},\mathcal{Y}\subseteq \omega^\omega$, $\mathcal{X}$ \textit{embeds into} $\mathcal{Y}$ if and only if $\exists\varphi[\varphi:\mathcal{X} \rightarrowtail \mathcal{Y}]$.

\smallskip For all $\mathcal{X},\mathcal{Y}\subseteq \omega^\omega$, for all $\varphi$, we define: \textit{$\varphi$ is a surjective mapping from $\mathcal{X}$ onto $\mathcal{Y}$}, $\varphi:\mathcal{X}\twoheadrightarrow\mathcal{Y}$, if and only if $\varphi:\mathcal{X}\rightarrow \mathcal{Y}$ and $\forall \beta\in\mathcal{Y}\exists \alpha\in\mathcal{X}[\varphi|\alpha=\beta]$.

$\mathcal{X}$ {\it maps onto} $\mathcal{Y}$ if and only if there exists a surjective mapping from $\mathcal{X}$ onto $\mathcal{Y}$. 

\smallskip For all $\mathcal{X}\subseteq\omega^\omega$, $(\omega^\omega)^\mathcal{X}:=\{\varphi\mid\varphi:\mathcal{X}\rightarrow\omega^\omega\}$.

Note: $(\omega^\omega)^{(\omega^\omega)}=\{\varphi\mid\varphi:\omega^\omega\rightarrow \omega^\omega\}=\{\varphi \in \omega^{(\omega^\omega)}\mid\varphi(0)=0\}$.

\smallskip For all $\varphi,s$ we let $\varphi|s$ be the largest number $t$ such that $length(t)\le length(s)$ and 

$\forall j <length(t)\exists p\le length(s)[\varphi^j(\overline s p) = t(j)+1 \;\wedge\;\forall i <p[\varphi^j(\overline s i)=0]]$. 

\smallskip Note: $\forall\varphi\forall s[length(\varphi|s)\le length(s)]$.

Note: $\forall \varphi\forall \alpha \forall \beta[\varphi:\alpha\mapsto \beta \leftrightarrow \forall n \exists m[\overline{\beta}n \sqsubseteq \varphi|\overline{\alpha}m]]$.

\smallskip For all $\varphi,\psi$ in $(\omega^\omega)^{(\omega^\omega)}$, we define  $\varphi\star\psi$ in $(\omega^\omega)^{(\omega^\omega)}$ such that, for all $n$, for all $s$, for all $p$, $\varphi^n(s)=p+1$ if and only if $n<length\bigl(\varphi|(\psi|s)\bigr)$ and $\bigl(\varphi|(\psi|s)\bigr)(n)=p+1$.\\Note:   $\forall \alpha[(\varphi\star\psi)|\alpha=\varphi|(\psi|\alpha)]$.

\medskip Let $\mathcal{F}\subseteq \omega^\omega$ be an inhabited spread. Find $\beta$ such that $Spr(\beta)$ and $\mathcal{F}=\mathcal{F}_\beta$.  \\Now define $\rho:\omega^\omega\rightarrow\omega^\omega$ such that, for all $\alpha$, for all $m$, \\ 
 if $\beta(\overline{\rho|\alpha} m \ast\langle \alpha(m)\rangle) = 0$, then $  (\rho|\alpha)(m) =\alpha(m)$, and,\\ if $\beta(\overline{\rho|\alpha} m \ast\langle \alpha(m)\rangle) \neq 0$, then  $(\rho|\alpha)(m) =\mu k[\beta(\overline{\rho|\alpha}m\ast\langle k \rangle)=0]$. 
\\$\rho$ is called the \textit{canonical retraction} of $\omega^\omega$ onto $\mathcal{F}$. 
\\Note: $\forall \alpha[\rho|\alpha \in \mathcal{F}]$, and: $\forall \alpha[\rho|\alpha\;\#\;\alpha\leftrightarrow \exists m[\beta(\overline{\alpha}m)\neq 0]]$, and: $\forall \alpha \in \mathcal{F}[\rho|\alpha=\alpha]$.  

\smallskip Assume:  $Spr(\beta)$ and:  $B\subseteq \omega$ is a {\it bar} in $\mathcal{F}_\beta$, i.e.:  $\forall \gamma \in \mathcal{F}_\beta\exists n[\overline \gamma n\in B]$. \\Define $B':=B\cup\{s\mid \beta(s)\neq 0\}$. Then: $B'$ is a bar in $\omega^\omega$, i.e.: $\forall \gamma \exists n[\overline \gamma n\in B']$. 
\\In order to see this, we use  the canonical retraction $\rho$ of $\omega^\omega$ onto $\mathcal{F}_\beta$. \\Let  $\gamma$ be given. Find $n$ such that $\overline{\rho|\gamma}n\in B$. \\{\it Either}: $\overline{\rho|\gamma}n=\overline \gamma n$ and $\overline \gamma n\in B$, {\it or}: $\overline{\rho|\gamma}n\neq\overline \gamma n$ and $\exists m\le n[\beta(\overline \gamma m)\neq 0]$. \\In both cases, $\overline \gamma n \in B'$. 

\subsubsection{Brouwer's Continuity Principle and the Axioms of Continuous Choice}\label{SSS:bcpcontchoice} \hfill

Brouwer's Continuity Principle:
\begin{quote}
$\mathbf{BCP}$: {\it For every spread $\mathcal{F}$, for every $\mathcal{R}\subseteq \mathcal{F}\times \omega$,

if $\forall \alpha \in \mathcal{F}\exists n[\alpha\mathcal{R}n] $, then $\forall \alpha \in \mathcal{F}\exists m \exists n\forall \beta\in \mathcal{F}[\overline{\alpha}m\sqsubset \beta \rightarrow \beta\mathcal{R}n]]$.} 
\end{quote}

\smallskip
First Axiom of Continuous Choice:
\begin{quote}
$\mathbf{AC}_{1,0}:$  {\it For every spread $\mathcal{F}$, for all $\mathcal{R} \subseteq \mathcal{F} \times \omega$,

if $\forall\alpha\in\mathcal{F}\exists n[\alpha\mathcal{R}n] $, then  $\exists \varphi[\varphi:\mathcal{F}\rightarrow \omega \;\wedge\;\forall \alpha \in \mathcal{F}[\alpha\mathcal{R}\varphi(\alpha)]]$.} 

\end{quote}

\smallskip Second Axiom of Continuous Choice:
\begin{quote}
$\mathbf{AC}_{1,1}:$ {\it For every spread $\mathcal{F}$, for all $\mathcal{R} \subseteq \mathcal{F} \times \omega^\omega$,

if $\forall\alpha\in\mathcal{F}\exists \beta[\alpha\mathcal{R}\beta]$, then $ \exists \varphi[\varphi:\mathcal{F}\rightarrow \omega^\omega \;\wedge\;\forall \alpha \in \mathcal{F}[\alpha\mathcal{R}\varphi|\alpha]]$.} 
\end{quote}
\subsubsection{The Fan Theorem}\label{SSS:fantheorem}\hfill

For all $\mathcal{X}\subseteq \omega^\omega$, for all $B\subseteq \omega$, we define: $Bar_\mathcal{X}(B) \leftrightarrow \forall \gamma \in \mathcal{X}\exists n[\overline{\gamma}n \in B]$.

\smallskip
For each $\beta$, we define: $Fan(\beta)\leftrightarrow\bigl(Spr(\beta)\;\wedge\;\forall s \exists n\forall m >n[\beta(s\ast\langle m \rangle)\neq 0]\bigr)$. 

\smallskip If $Fan(\beta)$, one says: \textit{$\beta$ is a fan-law}.

\smallskip $\mathcal{F}\subseteq \omega^\omega$ is a \textit{fan} if and only if $\exists \beta[Fan(\beta)\;\wedge\;\mathcal{F}=\mathcal{F}_\beta]$.

\smallskip The Fan Theorem:

\begin{quote}
{\it For every fan $\mathcal{F}\subseteq \omega^\omega$, for every $B\subseteq \omega$,
\\if $Bar_\mathcal{F}(B)$, then  $\exists s[D_s\subseteq B \;\wedge\; Bar_\mathcal{F}(D_s)]$.}

\end{quote}

  \smallskip
  The  \textit{restricted} Fan Theorem:
 
 \begin{quote}$\mathbf{FT}$: 
  {\it For each fan $\mathcal{F}\subseteq \omega^\omega$,  for every $\delta$,  if $Bar_\mathcal{F}(D_\delta)$, then  $ \exists n[ Bar_\mathcal{F}(D_{\overline{\delta}n})]$.}
  \end{quote}

\newpage \subsubsection{Stumps}\label{SSS:stumps}\hfill

Axiom on the existence of the set of stumps: 
\begin{quote}
$\mathbf{STP}$: {\it  $\mathcal{STP}$ is a subset of $2^\omega$ such that:\footnote{There is a small difference between the set $\mathcal{STP}$ as it is introduced here and the sets called $\mathbf{Stp}$ in \cite{veldman08}, \cite{veldman09}, respectively.} } {\it \begin{enumerate}[\upshape (i)]

 \item $1^\ast:=\underline 1\in \mathcal{STP}$, and,  

\item for all $\sigma$ in $2^\omega$, if \begin{enumerate} \item $\sigma(0) =0$ and, \item for all $n$,  $\sigma ^n \in \mathcal{STP}$, then $\sigma \in \mathcal{STP}$, \end{enumerate} and, 

\item  for all $\mathcal{Q}\subseteq \mathcal{STP}$, if \begin{enumerate} \item $1^\ast \in \mathcal{Q}$ and, \item for all $\sigma$ in $\mathcal{STP}$, \\if $\sigma(0)=0$ and, for all $n$,  $\sigma ^n \in \mathcal{Q}$, then $\sigma \in \mathcal{Q}$,\end{enumerate} then $\mathcal{STP}=\mathcal{Q}$. 

\end{enumerate}}
\end{quote}

The elements of $\mathcal{STP}$ are called \textit{stumps}.

 For each $\beta$ in $\omega^\omega$, we define $\beta^\ast$ in $2^\omega$ by: \\for all $s$, $\beta^\ast = 1$ if  $\beta(s) = 0$ and $\beta^\ast(s)=0$ if $\beta(s)\neq 0$.

\smallskip
$1^\ast$  is/codes {\it the empty stump}. 
 \\For each $\sigma$ in $\mathcal{STP}$, , $\sigma = 1^\ast$ if and only if $\sigma(0)=1$. 

\smallskip\noindent
For each $\sigma\neq 1^\ast$ in $\mathcal{STP}$,  for each $n$, $\sigma^n$ is a stump: \textit{the $n$-th immediate substump of $\sigma$}.

\smallskip (Also, for each $n$, $(1^\ast)^n=1^\ast$ is a stump.)

\smallskip\noindent
We define relations $<, \le$ on  $\mathcal{STP}$ by simultaneous transfinite induction: \\
{\it for all $\sigma, \tau$} in $\mathcal{STP}$,\begin{enumerate}[\upshape (i)]
\item  $\sigma \le \tau\leftrightarrow \bigl(\sigma\neq 1^\ast\rightarrow\forall n[\sigma^n < \tau]\bigr)$,  {\it and} \item $\sigma < \tau\leftrightarrow\bigl(\tau\neq 1^\ast \;\wedge\;\exists n[\sigma \le \tau^n]\bigr)$. \end{enumerate}

\smallskip
Using the axiom $\mathbf{STP}$ one proves the following

\smallskip

Principle of Induction on $\mathcal{STP}$: \begin{quote}{\it For all $\mathcal{Q}\subseteq \mathcal{STP}$, \\if $\forall \sigma \in \mathcal{STP}[\forall \tau \in \mathcal{STP}[\tau<\sigma\rightarrow \tau \in \mathcal{Q}]\rightarrow \sigma \in \mathcal{Q}]$, then $\mathcal{STP}=\mathcal{Q}$.}\end{quote} 

\smallskip One may  prove:\footnote{The relation $\le^\ast$ has been defined at the end of Subsubsection \ref{SSS:infseq}.}  for all $\sigma, \tau$ in $ \mathcal{STP}$, $\sigma \le\tau$ if and only if  $\sigma \le^\ast\tau$.

\smallskip For all $\alpha$, we let $S^\ast(\alpha)$ be the element $\beta$ of $\omega^\omega$ such that $\beta(0)=0$ and $\forall n[\beta^n=\alpha]$.  $S^\ast(\alpha)$ is called the \textit{successor of $\alpha$}.

Note: $\forall \alpha \in \mathcal{STP}[S^\ast(\alpha)\in\mathcal{STP}]$. 

\subsubsection{Bar Induction}\label{SSS:barinduction}\hfill

Brouwer's\textit{ Thesis on bars in $\omega^\omega$}:

\begin{quote}$\mathbf{BT}$: {\it For each $B\subseteq \omega$, if $Bar_{\omega^\omega}(B)$, then $\exists \sigma \in \mathcal{STP}[Bar_{\omega^\omega}(B\cap T_\sigma)]$.}
\end{quote}

Recall, from Subsubsection \ref{SSS:infseq}, that $T_\sigma=\{s\mid \forall t\sqsubset s[\sigma(t)=0]\}$.

\smallskip
$B\subseteq \omega$ is \textit{monotone} if and only if $\forall s\forall n[s \in B\rightarrow s\ast\langle n\rangle \in B]$.

$C\subseteq \omega$ is \textit{inductive} if and only if $\forall s[\forall n[s\ast\langle n \rangle \in C]\rightarrow s \in C]$.

\smallskip  $\mathbf{BT}$   proves the following

Principle of Bar Induction: \begin{quote} $\mathbf{BI}$:
{\it For all $B,C\subseteq \omega$, \\if $Bar_{\omega^\omega}(B)$, and $B\subseteq C$, and $C$ is monotone and  inductive, \\then $0=\langle \; \rangle \in C$}.
\end{quote}

\smallskip

Assume $Spr(\beta)$. 
We define: \\$B\subseteq \omega$ is \textit{monotone within $\{s\mid \beta(s)=0\}$} if and only if \\$\forall s[\bigl(\beta(s)=0\;\wedge\;s\in B\bigr)\rightarrow\forall n[\beta(s\ast\langle n \rangle)=0\rightarrow s\ast\langle n \rangle \in B]]$, and:
\\$C\subseteq \omega$ is \textit{inductive within $\{s\mid \beta(s)=0\}$} if and only if \\$\forall s[\bigl(\beta(s) =0 \;\wedge\;\forall n[\beta(s\ast\langle n \rangle)=0\rightarrow s\ast\langle n \rangle \in C]\bigr)\rightarrow s\in C]$.

$\mathbf{BI}$ admits of the following extension:

\begin{quote} $\mathbf{BI}$, extended to spreads:\\
{\it For all $\beta$ such that $Spr(\beta)$ and $\beta(0)=0$, for all $B,C\subseteq \omega$, \\if $Bar_{\mathcal{F}_\beta}(B)$, and $B\subseteq C$, \\and $C$ is monotone and  inductive within $\{s\mid \beta(s)=0\}$,  \\then $0=\langle \; \rangle \in C$}.
\end{quote}

Using $\mathbf{BI}$  and calling to aid the canonical retraction $\rho$ of $\omega^\omega$ onto $\mathcal{F}_\beta$, one  easily proves this extended form of $\mathbf{BI}$ from $\mathbf{BI}$ itself.

\subsubsection{The creating subject}\label{SSS:creasubj} \hfill

The \textit{Brouwer-Kripke axiom}, also called: \textit{Kripke's scheme}\footnote{Kripke's scheme  plays a role in the proof of Theorem \ref{T:kripkesan} and it is mentioned in Section \ref{S:a11e11}.}  is the following statement: \begin{quote}  $\mathbf{KS}$: {\it Given any \emph{definite} mathematical proposition $P$, \\one may build $\alpha$ such that $P\leftrightarrow \exists n[\alpha(n)\neq 0]$}.\end{quote}  

The idea underlying the axiom is that, once  $P$ is given, I may, identifying myself with the creating subject,  start thinking upon it, and the truth of $P$ will consist in my finding a proof of $P$,  at some point of time. Time is supposed to be divided into stages that are numbered by natural numbers. For each $n$, $\alpha(n) \neq 0$ if and only if, at stage $n$, I possess a proof of $P$.

This is a rather wild idea,  actually too wild, if we allow $P$ to be a statement about an object that is itself unfinished, like an infinite sequence $\beta=\beta(0), \beta(1), \ldots$ of natural numbers I  am creating step by step, freely choosing each one of its values. At any stage, only finitely many values will have been determined, and the statement: $\forall n[\beta(n)=0]$, provided it has not been violated already, is unprovable at any stage, although possibly true \textit{`in the end'}.

We therefore require $P$ to be \textit{definite}\footnote{The term \textit{`definite'} will also be applied to (other) mathematical objects. The infinite sequence $\underline 0$, for instance, deserves to be called definite.}: $P$ should not be about unfinished objects. In a formal contect, one forbids that the formula corresponding to the proposition contain a free variable over elements of Baire space. 

If one do not take this precaution,  $\mathbf{KS}$ leads to a contradiction with $\mathbf{AC}_{1,1}$, as was first observed by J.~Myhill, see \cite{myhill}:

Assume $\forall \beta \exists \alpha[\beta = \underline 0 \leftrightarrow \exists n[\alpha(n) \neq 0]]$. Applying $\mathbf{AC}_{1,1}$, find $\varphi:\omega^\omega\rightarrow\omega^\omega$ such that $\forall \beta[\beta = \underline 0 \leftrightarrow \exists n[(\varphi|\beta)(n)\neq 0]]$. Then find $n$ such that $(\varphi|\underline{0})(n) \neq 0$. Finally, find $m$ such that $\forall \beta[\overline{\underline{0}}m\sqsubset \beta\rightarrow (\varphi|\beta)(n) = (\varphi|\underline{0})(n)]$ and conclude: $\forall \beta[ \underline{\overline{0} }m \sqsubset \beta \rightarrow \beta = \underline 0]$, a contradiction. 

Myhill wanted to give up $\mathbf{AC}_{1,1}$ because of this argument.  Johan de Iongh proposed the restriction of $\mathbf{KS}$ to definite propositions, see \cite[\S 3]{gielen}.

\begin{theorem}[Consequences of $\mathbf{KS}$]\label{T:kripkescheme}\hfill

\begin{enumerate}[\upshape (i)]\item If $X\subseteq \omega$ is definite, then $\exists \delta[X=E_\delta]$, i.e.: $X$ is enumerable. \item If $\mathcal{X}\subseteq \omega^\omega$ is definite, then $\exists \delta[E_\delta = \{s\mid\exists\gamma\in\mathcal{X}[s\sqsubset \gamma]\}]$, i.e.: $\mathcal{X}$ is semi-located. \end{enumerate}\end{theorem}

\begin{proof} (i) Let $X\subseteq \omega$ be definite. By $\mathbf{KS}$, $\forall n \exists \alpha[n\in X\leftrightarrow \exists m[\alpha(m)\neq 0]]$.\\ Using $\mathbf{AC}_{0,1}$, find $\alpha$ such that $\forall n[n\in X \leftrightarrow \exists m[\alpha ^n(m) \neq 0]]$.\\Now define $\delta$ such that $\delta(0)=0$, and, for all $n,m$, \textit{if} $\alpha^n(m) \neq 0$, then $\delta(\langle n \rangle \ast m)= n+1$, and, \textit{if not}, then $\delta(\langle n \rangle \ast m)=0$, and  
note: $X=E_\delta$.

\smallskip (ii) Let $\mathcal{X}\subseteq \omega^\omega$ be definite. The set $\{s\mid\exists\gamma\in\mathcal{X}[s\sqsubset \gamma]\}$ also is definite, and one may apply (i). \end{proof}

\subsubsection{Semi-classical principles}\label{SSS:lpomp}\hfill

The Limited Principle of Omniscience: \begin{quote}$\mathbf{LPO}$: 
$\forall \alpha[\exists n[\alpha(n)\neq 0]\;\vee\;\forall n[\alpha(n)=0]]$.\end{quote}

\smallskip 
The Lesser Limited Principle of Omniscience\footnote{$\mathbf{LPO}$ and $\mathbf{LLPO}$ were introduced by E.~Bishop, as special cases of the principle of the the excluded middle $X\vee\neg X$.  If one reads well-known theorems constructively, many of them turn out to be equivalent to one of these `principles'. From a constructive point of view, these `principles' are, of course, totally wrong. } and : \begin{quote}$\mathbf{LLPO}$: 
$\forall \alpha[\forall m[2m\neq \mu p[\alpha(p)\neq 0 ]\;\vee\;\forall m[2m+1\neq\mu p[\alpha(p)\neq 0]]$. 
\end{quote}
 Note: $\mathbf{LPO}\rightarrow \mathbf{LLPO}$: \\
Let $\alpha$ be given. Define $\beta$ such that $\forall n[\beta(n) \neq 0\leftrightarrow 2n+1=\mu p[\alpha(p)\neq 0]]$. \\If $\exists n[\beta(n)\neq 0]]$, then $\forall m[2m\neq\mu p[\alpha(p)\neq 0 ]$, and, \\if $\forall n[\beta(n)=0]$, then $\forall m[2m+1\neq\mu p[\alpha(p)\neq 0 ]]$.

\smallskip $\mathbf{LLPO}$ and $\mathbf{BCP}$ together give a contradiction: assuming both, find $p$ such that $\forall \alpha[\underline{\overline 0}p\sqsubset \alpha\rightarrow \forall m[2m\neq\mu p[\alpha(p)\neq 0 ]]$ or $\forall \alpha[\underline{\overline 0}p\sqsubset \alpha\rightarrow \forall m[2m+1\neq\mu p[\alpha(p)\neq 0 ]]$. \\The sequences $\underline{\overline 0}(2p)\ast\underline 1$ and $\underline{\overline 0}(2p+1)\ast\underline 1$ show that both alternatives are false.

\medskip Markov's Principle: \begin{quote}$\mathbf{MP}$: $\forall \alpha[\neg\neg\exists n[\alpha(n)\neq 0]\rightarrow\exists n[\alpha(n)\neq 0]]$\end{quote} has been defended by Markov for algorithmically computable $\alpha$ only.
\subsection{Descriptive set theory}
\hfill

Information on classical descriptive set theory may be found in \cite{lusin}, \cite{moschovakis}, \cite{kechris} and \cite{srivastava}. Some results on the borderline of classical and intuitionistic descriptive set theory may be found in \cite{moschovakis2} and \cite{moschovakis18}.
\subsubsection{Some basic notions}\hfill

For all $\mathcal{X},\mathcal{Y}\subseteq \omega^\omega$, 
 for all $\varphi:\omega^\omega\rightarrow\omega^\omega$, we define: \\{\it $\varphi$ reduces $\mathcal{X}$  to $ \mathcal{Y}$} if and only if $\forall \alpha[\alpha \in \mathcal{X}\leftrightarrow \varphi|\alpha \in \mathcal{Y}]$.

We define: {\it $\mathcal{X}$ reduces to  $\mathcal{Y}$},   $\mathcal{X}\preceq\mathcal{Y}$, if and only if \\there exists $\varphi:\omega^\omega\rightarrow\omega^\omega$ reducing $\mathcal{X}$ to $\mathcal{Y}$.

\smallskip For all $\mathcal{X}_0,\mathcal{X}_1,\mathcal{Y}_0,\mathcal{Y}_1\subseteq \omega^\omega$, we define:\\ $(\mathcal{X}_0,\mathcal{X}_1)$ {\it simultaneously reduces to $(\mathcal{Y}_0, \mathcal{Y}_1)$}, $(\mathcal{X}_0, \mathcal{X}_1)\preceq(\mathcal{Y}_0, \mathcal{Y}_1)$, if and only if \\there exists $\varphi:\omega^\omega\rightarrow \omega^\omega$ reducing $\mathcal{X}_0$ to $\mathcal{Y}_0$ and also $\mathcal{X}_1$ to $\mathcal{Y}_1$.

\smallskip Let $\mathfrak{K}$ be a class of subsets of $\omega^\omega$.

Assume $\mathcal{X}\subseteq \omega^\omega$. We often say: \textit{`$\mathcal{X}$ is $\mathfrak{K}$'} for: \textit{`the set $\mathcal{X}$ belongs to the class $\mathfrak{K}$'} .

\smallskip We define: $\mathcal{X}\subseteq \omega^\omega$ is $\mathfrak{K}$-\textit{complete} if and only if \\$\mathfrak{K}$ is the class of all $\mathcal{Y}\subseteq \omega^\omega$ reducing to $\mathcal{X}$, and 

We define: $\mathcal{X}\subseteq\omega^\omega$ is $\mathfrak{K}$-\textit{universal} if and only if \\$\mathfrak{K}$ is the class of all sets of the form $\mathcal{X}\upharpoonright\alpha$, for some $\alpha$ in $\omega^\omega$.

Note: if $\mathcal{X}$ is $\mathfrak{K}$-universal, then $\mathcal{X}$ is $\mathfrak{K}$-complete.

\subsubsection{Open sets and closed sets}\label{SSS:openandclosed}\hfill

$\mathbf{\Sigma}^0_1:=\{\mathcal{G}_\beta\mid\beta \in \omega^\omega\}$ and: $\mathbf{\Pi}^0_1:=\{\mathcal{F}_\beta\mid\beta \in \omega^\omega\}$.

\smallskip
$\mathcal{E}_1:=\{\alpha\mid\exists n[\alpha(n)\neq 0]\}=\{\alpha\mid\alpha\;\#\;\underline{0}\}$ and $\mathcal{A}_1:=\{\alpha\mid\forall n[\alpha(n)= 0]\}=\{\underline 0\}$.

$\mathcal{E}_1$ is $\mathbf{\Sigma}^0_1$-complete and $\mathcal{A}_1$ is $\mathbf{\Pi}^0_1$-complete.

\smallskip $\mathcal{US}_1:=\{\alpha\mid\alpha_{II} \in \mathcal{G}_{\alpha_I}\} =\{\alpha\mid\exists n[\alpha_I(\overline{\alpha_{II}}n)\neq 0]\}$ and

 $\mathcal{UP}_1:=\{\alpha\mid\alpha_{II} \in \mathcal{F}_{\alpha_I} \}=\{\alpha\mid\forall n[\alpha_I(\overline{\alpha_{II}}n)= 0]\}$.
 
 $\mathcal{US}_1$ is $\mathbf{\Sigma}^0_1$-universal and $\mathcal{UP}_1$ is $\mathbf{\Pi}^0_1$-universal.

\subsubsection{Borel sets of finite rank}\label{SSS:borelfiniterank} \hfill

For each $m>0$, for each $\beta$, we define $\mathcal{G}^m_\beta,\mathcal{F}^m_\beta\subseteq \omega^\omega$,  by induction.

$\mathcal{G}^1_\beta:=\mathcal{G}_\beta$ and $\mathcal{F}^1_\beta:=\mathcal{F}_\beta$, and, for each $m>0$, $\mathcal{G}^{m+1}_\beta =\bigcup_n\mathcal{F}^m_{\beta^n}$ and $\mathcal{F}^{m+1}_\beta =\bigcap_n\mathcal{G}^m_{\beta^n}$.

\smallskip
For each $m>0$, for each $\beta$, the pair of sets $(\mathcal{G}^m_\beta,\mathcal{F}^m_\beta)$ is called a \textit{complementary pair of (positively) Borel sets of rank $m$}.

\smallskip
For each $m>0$, $\mathbf{\Sigma}^0_m:=\{\mathcal{G}^m_\beta\mid\beta\in\omega^\omega\}$ and $\mathbf{\Pi}^0_m:=\{\mathcal{F}^m_\beta\mid\beta\in\omega^\omega\}$.

\smallskip For each $m>0$, we define $\mathcal{E}_m, \mathcal{A}_m\subseteq \omega^\omega$, by induction.

$\mathcal{E}_1, \mathcal{A}_1$ were defined in Subsubsection \ref{SSS:openandclosed}. 

For each $m>0$, $\mathcal{E}_{m+1}:=\{\alpha\mid\exists n[\alpha^n \in \mathcal{A}_m]\}$, and $\mathcal{A}_{m+1}:=\{\alpha\mid\forall n[\alpha^n \in \mathcal{E}_m]\}$.

For each $m>0$, $\mathcal{E}_m$ is $\mathbf{\Sigma}^0_m$-complete and $\mathcal{A}_m$ is $\mathbf{\Pi}^0_m$-complete.

For each $m>0$, $(\mathcal{E}_m, \mathcal{A}_m)$ is a complementary pair of rank $m$. 

\smallskip For each $m>0$, $\mathcal{US}_m:=\{\alpha\mid\alpha_{II}\in \mathcal{G}^m_{\alpha_I}\}$ and $\mathcal{UP}_m:=\{\alpha\mid\alpha_{II}\in \mathcal{F}^m_{\alpha_I}\}$.

For each $m>0$, $\mathcal{US}_m$ is $\mathbf{\Sigma}^0_m$-universal and $\mathcal{UP}_m$ is $\mathbf{\Pi}^0_m$-universal.

For each $m>0$, $(\mathcal{US}_m, \mathcal{UP}_m)$ is a complementary pair of rank $m$. 
\subsubsection{Borel sets in general}\label{SSS:borelgeneral}\label{SSS:borelsetsgeneral}\hfill

The set $\mathcal{HRS}$ of the \textit{hereditarily repetitive stumps} is defined inductively: for each stump $\sigma$: $\sigma \in \mathcal{HRS}\leftrightarrow \bigl(\sigma(0)=0\rightarrow (\forall n[\sigma^n \in \mathcal{HRS} \;\wedge\;\forall m\exists n>m[\sigma^n = \sigma^m]\bigr)$.

\smallskip For each $\sigma$ in $\mathcal{HRS}$, for each $\beta$, we define $\mathcal{G}^\sigma_\beta, \mathcal{F}^\sigma_\beta\subseteq \omega^\omega$, by induction:

if $\sigma=1^\ast$, then $\mathcal{G}^\sigma_\beta = \mathcal{G}_\beta$ and $\mathcal{F}^\sigma_\beta = \mathcal{F}_\beta$, and,

if $\sigma\neq 1^\ast$, then $\mathcal{G}^\sigma_\beta:=\bigcup_n\mathcal{F}^{\sigma^n}_{\beta^n}$ and $\mathcal{F}^\sigma_\beta:=\bigcap_n\mathcal{G}^{\sigma^n}_{\beta^n}$.

\smallskip Note: for each $\sigma$ in $\mathcal{HRS}$, for each $\beta$, $\mathcal{G}_\beta^\sigma\;\#\;\mathcal{F}_\beta^\sigma$. 

The pair of sets $(\mathcal{G}_\beta^\sigma,\mathcal{F}_\beta^\sigma)$ is called a \textit{complementary pair of (positively) Borel sets of rank $\sigma$}. 

\smallskip
For each $\sigma$ in $\mathcal{HRS}$, $\mathbf{\Sigma}^0_\sigma:=\{\mathcal{G}^\sigma_\beta\mid\beta\in \omega^\omega\}$ and $\mathbf{\Pi}^0_\sigma:=\{\mathcal{F}^\sigma_\beta\mid\beta\in \omega^\omega\}$.

\smallskip For each $\sigma$ in $\mathcal{HRS}$, we define $\mathcal{E}_\sigma, \mathcal{A}_\sigma \subseteq \omega^\omega$, by induction:

if $\sigma=1^\ast$, then $\mathcal{E}_\sigma :=\mathcal{E}_1$ and $\mathcal{A}_\sigma :=\mathcal{A}_1$, and

if $\sigma\neq 1^\ast$, then $\mathcal{E}_\sigma:=\{\alpha\mid\exists n[\alpha^n \in \mathcal{A}_{\sigma^n}]\}$ and $\mathcal{A}_\sigma:=\{\alpha\mid\forall n[\alpha^n \in \mathcal{E}_{\sigma^n}]\}$.

\smallskip For each $\sigma$ in $\mathcal{HRS}$, $\mathcal{E}_\sigma$ is $\mathbf{\Sigma}^0_\sigma$-complete and $\mathcal{A}_\sigma$ is $\mathbf{\Pi}^0_\sigma$-complete and $(\mathcal{E}_\sigma,\mathcal{A}_\sigma)$ is a complementary pair of rank $\sigma$.

\smallskip For each $\sigma$ in $\mathcal{HRS}$, $\mathcal{US}_\sigma := \{\alpha\mid\alpha_{II} \in \mathcal{G}^\sigma_{\alpha_I}\}$ and $\mathcal{UP}_\sigma := \{\alpha\mid\alpha_{II} \in \mathcal{F}^\sigma_{\alpha_I}\}$.

For each $\sigma$ in $\mathcal{HRS}$, $\mathcal{US}_\sigma$ is $\mathbf{\Sigma}^0_\sigma$-universal and $\mathcal{UP}_\sigma$ is $\mathbf{\Pi}^0_\sigma$-universal and $(\mathcal{US}_\sigma, \mathcal{UP}_\sigma)$ is a complementary pair of rank $\sigma$. 

\smallskip The function $S^\ast:\omega^\omega\rightarrow\omega^\omega$ has been defined in Subsubsection \ref{SSS:stumps}.
Note:\\ $\forall \sigma \in \mathcal{HRS}[S^\ast(\sigma) \in  \mathcal{HRS}]$. 

\smallskip

Define: $1^\ast:=\underline 1$ and, for all $m$, $(m+1)^\ast=S^\ast(m^\ast)$. Note: for all $m>0$, $\mathbf{\Sigma}^0_m=\mathbf{\Sigma}^0_{m^\ast}$ and $\mathcal{E}_m=\mathcal{E}_{m^\ast}$ and $\mathbf{\Pi}^0_m=\mathbf{\Pi}^0_{m^\ast}$ and $\mathcal{A}_m=\mathcal{A}_{m^\ast}$, $\ldots$

\smallskip $\mathfrak{Borel}:=\{\mathcal{G}^\sigma_\beta\mid\sigma \in \mathcal{HRS}, \beta \in \omega^\omega\}$.

\smallskip The following is proven in \cite[Theorems 4.9, 7.9, 7.10]{veldman08}.

\begin{theorem}[Borel Hierarchy Theorem]\label{T:borelhiertheorem}\hfill

 \begin{enumerate}[\upshape (i)]\item For all $\sigma,\tau$ in $\mathcal{HRS}$, if $\sigma < \tau$, then $\mathcal{E}_\sigma, \mathcal{A}_\sigma$ reduce to both $\mathcal{E}_\tau$ and $\mathcal{A}_\tau$.

\item \emph{(Not using $\mathbf{BCP}$):} For all $\sigma$ in $\mathcal{HRS}$, 

$\forall \varphi:\omega^\omega\rightarrow \omega^\omega\exists \alpha[(\alpha \in \mathcal{E}_\sigma \leftrightarrow \varphi|\alpha \in \mathcal{E}_\sigma)\;\wedge\;(\alpha \in \mathcal{A}_\sigma \leftrightarrow \varphi|\alpha \in \mathcal{A}_\sigma)]$.

\item \emph{(Using $\mathbf{BCP}$)}: For all $\sigma$ in $\mathcal{HRS}$: 

$\forall \varphi:\omega^\omega\rightarrow \omega^\omega[\varphi|\mathcal{E}_\sigma \subseteq \mathcal{A}_\sigma \rightarrow \exists \alpha[\alpha \in \mathcal{A}_\sigma \;\wedge\;\varphi|\alpha \in \mathcal{A}_\sigma]]$ and:

$\forall \varphi:\omega^\omega\rightarrow \omega^\omega[\varphi|\mathcal{A}_\sigma \subseteq \mathcal{E}_\sigma \rightarrow \exists \alpha[\alpha \in \mathcal{E}_\sigma \;\wedge\;\varphi|\alpha \in \mathcal{E}_\sigma]]$, or, equivalently:

for all $\mathcal{X}$ in $\mathbf{\Pi}^0_\sigma$, if $\mathcal{E}_\sigma\subseteq \mathcal{X}$, then $\exists\alpha\in\mathcal{A}_\sigma[\alpha \in \mathcal{X}]$, and:

for all $\mathcal{X}$ in $\mathbf{\Sigma}^0_\sigma$, if $\mathcal{A}_\sigma\subseteq \mathcal{X}$, then $\exists\alpha\in\mathcal{E}_\sigma[\alpha \in \mathcal{X}]$. \end{enumerate}\end{theorem}

Theorem \ref{T:borelhiertheorem}(iii) implies: $\mathcal{E}_\sigma$ positively fails to be $\mathbf{\Pi}^0_\sigma$  and $\mathcal{A}_\sigma$ positively fails to be $\mathbf{\Sigma}^0_\sigma$.  For the intuitionistic mathematician, Theorem \ref{T:borelhiertheorem}(ii) does \textit{not} establish the hierarchy, as, for almost every $\sigma$ in $\mathcal{HRS}$,  he is unable to prove: $\neg \exists \alpha[\alpha \notin \mathcal{E}_\sigma \;\wedge\;\alpha \notin \mathcal{A}_\sigma]$.

\subsubsection{On disjunction}\label{SSS:disjunction} \hfill

For every infinite sequence $\mathcal{X}_0,\mathcal{X}_1,\mathcal{X}_2,\ldots$ of subsets of $\omega^\omega$, we define:\begin{quote}$\mathbb{D}_n(\mathcal{X}_n):=\{\alpha\mid\exists n[\alpha^n\in \mathcal{X}_n]\}$ and  $\mathbb{C}_n(\mathcal{X}_n):=\{\alpha\mid\forall n[\alpha^n\in \mathcal{X}_n]\}$. \end{quote}

$\mathbb{D}_n(\mathcal{X}_n),\mathbb{C}_n(\mathcal{X}_n)$ are the \textit{disjunction} and the \textit{conjunction} of the infinite sequence $\mathcal{X}_0,\mathcal{X}_1,\mathcal{X}_2,\ldots$, respectively.

Note that, for each $\sigma$ in $\mathcal{HRS}$, if $\sigma\neq 1^\ast$,  then $\mathcal{E}_\sigma = \mathbb{D}_n(\mathcal{A}_{\sigma^n})$ and $\mathcal{A}_\sigma = \mathbb{C}_n(\mathcal{E}_{\sigma^n})$.

\smallskip

For all $\mathcal{X}_0,\mathcal{X}_1\subseteq \omega^\omega$, we define: \\$\mathbb{D}(\mathcal{X}_0,\mathcal{X}_1):=\{\alpha\mid\exists i<2[\alpha^i\in \mathcal{X}_i]\}$ and $\mathbb{D}^2(\mathcal{X}_0):=\mathbb{D}(\mathcal{X}_0, \mathcal{X}_0)$.

$\mathbb{D}(\mathcal{X}_0,\mathcal{X}_1)$ is called the \textit{disjunction of $\mathcal{X}_0$ and $\mathcal{X}_1$}. 

Note that $\mathcal{Z}\subseteq\omega^\omega$ reduces to $\mathbb{D}(\mathcal{X}_0,\mathcal{X}_1)$ if and only if there exist $\mathcal{Z}_0,\mathcal{Z}_1$ such that $\mathcal{Z}=\mathcal{Z}_0\cup\mathcal{Z}_1$ and $\forall i <2[\mathcal{Z}_i\preceq \mathcal{X}_i]$. 

The fiollowing result is not difficult but very important. 

\begin{theorem}\label{T:disja1} $\neg \bigl(\overline{\mathbb{D}^2(\mathcal{A}_1)}\subseteq \mathbb{D}^2(\mathcal{A}_1)\bigr)$. \end{theorem}

\begin{proof} Assume $\overline{\mathbb{D}^2(\mathcal{A}_1)}\subseteq \mathbb{D}^2(\mathcal{A}_1)=\{\alpha\mid\alpha^0=\underline 0\;\vee\;\alpha^1 =\underline 0\}$. Note: $\overline{\mathbb{D}^2(\mathcal{A}_1)}$ is a spread containing $\underline 0$. Applying $\mathbf{BCP}$, find $m$ such that \textit{either}: $\forall \alpha \in \overline{\mathbb{D}^2(\mathcal{A}_1)}[\underline{\overline{0}}m\sqsubset \alpha \rightarrow \alpha^0 = \underline 0]$ \textit{or}: $\forall \alpha \in \overline{\mathbb{D}^2(\mathcal{A}_1)}[\underline{\overline{0}}m\sqsubset \alpha \rightarrow \alpha^ 1 = \underline 0]$. Both alternatives are false.\end{proof}

Theorem \ref{T:disja1} shows that the union of two $\mathbf{\Pi}^0_1$-sets is not always $\mathbf{\Pi}^0_1$: $\mathbb{D}(\mathcal{A}_1, \mathcal{A}_1)$ does not reduce to $\mathcal{A}_1$. This result admits of a vast extension:

\smallskip

Assume: $\sigma \in \mathcal{HRS}$.  Define, as in \cite[p. 39]{veldman08}:

$\sigma$   \textit{is  weakly comparative} $\leftrightarrow \bigl(\sigma(0)=0\rightarrow \forall m \forall n\exists p[\sigma^m \le \sigma^p\;\wedge \;\sigma^n\le \sigma^p]\bigr)$,

The following result  is \cite[Theorem 8.8]{veldman08}.

\begin{theorem}[The persisting difficulty of disjunction]\label{T:persisdiff}
 For each $\sigma$ in $\mathcal{HRS}$, if $\sigma$ is  weakly comparative, then $\mathbb{D}(\mathcal{A}_1, \mathcal{A}_\sigma)$ does not reduce to  $\mathcal{A}_{S^\ast(\sigma)}$.
\end{theorem}

\subsubsection{Perhaps}\label{SSS:perhaps}\hfill

For every $\mathcal{X}\subseteq \omega^\omega$, $\mathsf{Perhaps}(\mathcal{X}):=\{\alpha\mid\exists\beta\in\mathcal{X}[\alpha\;\#\;\beta \rightarrow \alpha \in \mathcal{X}]\}$. 

If $\mathcal{X}$ is inhabited, then $\mathcal{X}\subseteq \mathsf{Perhaps}(\mathcal{X})$.

$\mathcal{X}\subseteq \omega^\omega$ is \textit{perhapsive} if and only if $\mathcal{X}= \mathsf{Perhaps}(\mathcal{X})$.

In \cite{waaldijk}, perhapsive subsets of $\omega^\omega$ are called \textit{weakly stable}. \cite{waaldijk} is the birthplace of the notion of `perhapsity'. The notion has been studied further in \cite{veldman5}, \cite{veldman05} and \cite{veldman09}.

\begin{theorem}\label{T:perhaps}\hfill

\begin{enumerate}[\upshape (i)] \item For all $\mathcal{X},\mathcal{Y}\subseteq \omega^\omega$, if $\mathcal{X}\preceq \mathcal{Y}$ and $\mathcal{Y}$ is perhapsive, then $\mathcal{X}$ is perhapsive.

\item $\mathbb{D}^2(\mathcal{A}_1)$ and $\mathcal{E}_2$ are not perhapsive.

\item $\mathcal{A}_2$ is perhapsive and $\neg\bigl(\mathbb{D}^2(\mathcal{A}_1)\preceq \mathcal{A}_2\bigr)$.

\end{enumerate} \end{theorem}

\begin{proof} (i) Let $\mathcal{X},\mathcal{Y},\varphi$ be given 
such that \\$\varphi:\omega^\omega\rightarrow \omega^\omega$ reduces $\mathcal{X}$ to $\mathcal{Y}$ and $\mathcal{Y}$ is perhapsive.
\\Let $\alpha,\beta$ be given such that $\beta \in \mathcal{X}$ and $\alpha \;\#\;\beta \rightarrow \alpha\in \mathcal{X}$.\\ Note: $\varphi|\beta \in \mathcal{Y}$, and, if $\varphi|\alpha\;\#\;\varphi|\beta$, then:  $\alpha\;\#\;\beta$,  and: $\alpha \in \mathcal{X}$, and: $\varphi|\alpha \in \mathcal{Y}$. \\As $\mathcal{Y}$ is perhapsive, we conclude: $\varphi|\alpha \in \mathcal{Y}$ and: $\alpha \in \mathcal{X}$.  \\We thus see: $\forall \alpha[\exists \beta \in \mathcal{X}[\alpha\;\#\;\beta\rightarrow \alpha \in \mathcal{X}]\rightarrow \alpha \in \mathcal{X}]$, that is: $\mathcal{X}$ is perhapsive.

\smallskip (ii) Let $\alpha$ in $\overline{\mathbb{D}^2(\mathcal{A}_1)}$ be given. \\Define $\alpha_0$ such that $(\alpha_0)^0=\underline 0$ and $\forall j[\neg\exists n[j=\langle 0\rangle\ast n]\rightarrow \alpha_0(j)=\alpha(j)]$. 
\\Note: $\alpha_0 \in \mathbb{D}^2(\mathcal{A}_1)$ and, if $\alpha\;\#\;\alpha_0$, then: $\alpha^1=\underline 0$, and: $\alpha \in \mathbb{D}^2(\mathcal{A}_1)$. \\We thus see: $\forall \alpha \in \overline{\mathbb{D}^2(\mathcal{A}_1)}[\alpha \in \mathsf{Perhaps}\bigl(\mathbb{D}^2(\mathcal{A}_1)\bigr)]$. \\Using Theorem \ref{T:disja1}, we conclude: $\mathbb{D}^2(\mathcal{A}_1)\neq \mathsf{Perhaps}\bigl(\mathbb{D}^2(\mathcal{A}_1)\bigr)$, that is: $\mathbb{D}^2(\mathcal{A}_1)$ is not perhapsive.

As $\mathbb{D}^2(\mathcal{A}_1)$ is $\mathbf{\Sigma}^0_2$ and reduces to $\mathcal{E}_2$, also $\mathcal{E}_2$ is not perhapsive, by (i). 

\smallskip (iii) Let $\alpha,\beta$ be given such that $\beta\in \mathcal{A}_2$ and $\alpha\;\#\;\beta\rightarrow \alpha \in \mathcal{A}_2$. \\Let $m$ be given. Find $n$ such that $\beta^m(n)\neq 0$. \\\textit{Either}: $\alpha^m(n)=\beta^m(n)\neq 0$, \textit{or}: $\alpha \;\#\;\beta$, and: $\alpha \in \mathcal{A}_2$, and: $\exists p[\alpha^m(p)\neq 0]$. \\We thus see: $\forall m\exists p[\alpha^m(p)\neq 0]$, that is: $\alpha \in \mathcal{A}_2$. \\Conclude: $\forall \alpha[\exists \beta \in \mathcal{A}_2[\alpha\;\#\;\beta\rightarrow \alpha \in \mathcal{A}_2]\rightarrow \alpha \in \mathcal{A}_2]$, i.e. $\mathcal{A}_2$ is perhapsive.

It follows that $\mathbb{D}^2(\mathcal{A}_1)$ does not reduce to $\mathcal{A}_2$, by (ii) and (i).
\end{proof} 

\subsubsection{Projective sets}\label{SSS:projective} \hfill

For each $\mathcal{X}\subseteq \omega^\omega$, $Ex(\mathcal{X}):=\{\alpha\mid\exists \beta[\ulcorner\alpha,\beta\urcorner\in\mathcal{X}]\}=\{\alpha_I\mid\alpha \in \mathcal{X}\}$ and

 $Un(\mathcal{X}):=\{\alpha\mid\forall\beta[\ulcorner\alpha,\beta\urcorner\in \mathcal{X}]\}$. 
 
 $Ex(\mathcal{X})$ is called the \emph{projection} of $\mathcal{X}$, and $Un(\mathcal{X})$ is called the \emph{co-projection} of $\mathcal{X}$. 
 
\smallskip For each $\beta$, $\mathcal{EF}_\beta:=Ex(\mathcal{F}_\beta)$ and $\mathcal{UG}_\beta:=Un(\mathcal{G}_\beta)$.

\smallskip $\mathbf{\Sigma}^1_1:=\{\mathcal{EF}_\beta\mid\beta \in \omega^\omega\}$ is the class of the \textit{analytic} sets and 

$\mathbf{\Pi}^1_1:=\{\mathcal{UG}_\beta\mid\beta \in \omega^\omega\}$ is the class of the \textit{co-analytic} sets.

\smallskip $\mathbf{\Sigma}^1_1$ thus consists of the \emph{projections of the closed subsets of $\omega^\omega$} and  \\$\mathbf{\Pi}^1_1$  consists of the \emph{co-projections of the open subsets of $\omega^\omega$}.

\smallskip For each $\beta$, ($\mathcal{EF}_\beta,\mathcal{UG}_\beta$) is called a \textit{complementary ($\mathbf{\Sigma}^1_1, \mathbf{\Pi}^1_1$)-pair}.

\smallskip $\mathcal{US}^1_1:=\{\alpha\mid\alpha_{II}\in \mathcal{EF}_{\alpha_I}\}$ and $\mathcal{UP}^1_1:=\{\alpha\mid\alpha_{II}\in \mathcal{UG}_{\alpha_I}\}$.

We shall prove that $\mathcal{US}^1_1$ is $\mathbf{\Sigma}^1_1$-universal, see Theorem \ref{T:analytic}(i).

We shall prove that $\mathcal{UP}^1_1$ is $\mathbf{\Pi}^1_1$-universal, see Theorem \ref{T:coan}(i).

\smallskip $\mathcal{E}_1^1:=\{\alpha\mid\exists \gamma\forall n[\alpha(\overline \gamma n)=0]\}$ and $\mathcal{A}_1^1:=\{\alpha\mid\forall \gamma\exists n[\alpha(\overline \gamma n)\neq 0]\}$. 

We shall prove that $\mathcal{E}^1_1$ is $\mathbf{\Sigma}^1_1$-complete, see Theorem \ref{T:analytic}(ii).

We shall prove that $\mathcal{A}^1_1$ is $\mathbf{\Pi}^1_1$-complete, see Theorem \ref{T:coan}(ii).

\smallskip For certain purposes, the class $\mathbf{\Sigma}^1_1$ is too wide. We therefore introduce the class

$\mathbf{\Sigma}^{1,\ast}_1 :=\{\mathcal{EF}_\beta\mid Spr(\beta)\}$  of the \textit{strictly analytic} sets.

$\mathbf{\Sigma}^{1,\ast}_1$  consists of the projections of the  subsets of $\omega^\omega$ that are both closed and {\it located}.

\smallskip  For certain purposes, the class $\mathbf{\Pi}^1_1$ is too narrow. We therefore introduce the class

$\mathbf{\Pi}^{1+}_1:=\{Un(\mathcal{X})\mid\mathcal{X}\in \mathfrak{Borel}\}$.

$\mathbf{\Pi}^{1+}_1$  is the class of the \textit{broadly co-analytic} sets.

\smallskip
For each $\beta$, $\mathcal{UEF}_\beta :=Un(\mathcal{EF}_\beta)$ and  $\mathcal{EUG}_\beta :=Ex(\mathcal{UG}_\beta)$.

\smallskip $\mathbf{\Pi}^1_2:=\{\mathcal{UEF}_\beta\mid\beta\in\omega^\omega\}$ and $\mathbf{\Sigma}^1_2:=\{\mathcal{EUG}_\beta\mid\beta\in\omega^\omega\}$.

\smallskip For each $\beta$, $(\mathcal{EUG}_\beta,\; \mathcal{UEF}_\beta)$ is a \textit{complementary $(\mathbf{\Sigma}^1_2, \;\mathbf{\Pi}^1_2)$-pair.}

\smallskip $\mathcal{E}_2^1:=\{\alpha\mid\exists \delta\forall \gamma\forall n[\alpha(\overline{\ulcorner  \gamma, \delta\urcorner} n)=0]\}$ and 
$\mathcal{A}_2^1:=\{\alpha\mid\forall\delta\exists\gamma\exists n[\alpha(\overline{\ulcorner \gamma, \delta\urcorner} n)\neq 0]\}$. 

We shall prove that $\mathcal{E}^1_2$ is $\mathbf{\Sigma}^1_2$-complete, and that $\mathcal{A}^1_2$ is $\mathbf{\Pi}^1_2$-complete, \\see Theorem \ref{T:secondlevel}(ii).

\smallskip $\mathcal{US}^1_2:=\{\alpha\mid\alpha_{II}\in \mathcal{EUG}_{\alpha_I}\}$ and $\mathcal{UP}^1_1:=\{\alpha\mid\alpha_{II}\in \mathcal{UEF}_{\alpha_I}\}$.

We shall prove that $\mathcal{US}^1_2$ is $\mathbf{\Sigma}^1_2$-universal, and  
 that $\mathcal{UP}^1_2$ is $\mathbf{\Pi}^1_2$-universal, \\see Theorem \ref{T:secondlevel}(i).

\subsection{The main results of this paper}
Apart from this introductory Section, the paper contains Sections numbered 2 to 7.

\smallskip
In Section 2, we first establish some properties of the class $\mathbf{\Sigma}^1_1$. 

We then prove that the set \begin{quote}$\mathcal{IF}:= \{\alpha\mid \exists \beta \in (T_\alpha)^\omega\forall n[\beta(n+1)<_{KB}\beta(n)]\}$,\end{quote} i.e.:
 the set of all $\alpha$ such that the tree $T_\alpha:=\{s\mid\forall t\sqsubset s[\alpha(t)=0]\}$ is (positively) \textit{ill-founded with respect to the Kleene-Brouwer-ordering $<_{KB}$}, is $\mathbf{\Sigma}^1_1$ but not $\mathbf{\Sigma}^1_1$-complete.   We also prove  that the set \begin{quote}$\mathcal{UNC}:=\{\beta\mid \forall \alpha\exists\gamma\in \mathcal{F}_\beta\forall n[\gamma\;\#\;\alpha^n]\}$ \end{quote}of codes of the {\it positively uncountable} closed subsets of $\omega^\omega$ is $\mathbf{\Sigma}^1_1$-complete, and that the same holds for the set \begin{quote}$\mathsf{Share}^\ast(\mathcal{INF}):=\{\beta\mid Spr(\beta)\;\wedge\;\exists \alpha \in \mathcal{F}_\beta\forall m\exists n>m[\alpha(n)\neq0]\}$\end{quote} of codes of the spreads that contain an element $\alpha$ such that $D_\alpha=\{n\mid\alpha(n)\neq 0\}$ is an infinite subset of $\omega$. 
 
 The final Subsection of Section 2 is devoted to the class $\mathbf{\Sigma}^{1\ast}_1$ of the \textit{strictly analytic} subsets of $\omega^\omega$. $\mathbf{\Sigma}^{1\ast}_1$ is a proper subclass of $\mathbf{\Sigma}^1_1$ and is lacking some of the useful closure properties of $\mathbf{\Sigma}^1_1$.

\smallskip
In Section 3, we give intuitionistic proofs of the  Separation Theorems due to Lusin and  Novikov.  Novikov's Theorem is the stronger one and says that, given any infinite sequence  $\mathcal{X}_0,\mathcal{X}_1,\dots$ of $\mathbf{\Sigma}^{1\ast}_1$ subsets of $\omega^\omega$ such that  $\#_n(\mathcal{X}_n)$, (that is, in a constructively strong sense: $\bigcap_n(\mathcal{X}_n)=\emptyset$),  one may find an infinite sequence  $\mathcal{B}_0,\mathcal{B}_1, \dots$  of Borel subsets of $ \omega^\omega$ such that $\forall n[\mathcal{X}_n\subseteq \mathcal{B}_n]$ and $\#_n(\mathcal{B}_n)$.  The proofs use Brouwer's Thesis on bars in $\omega^\omega$. 

We  give an intuitionistic proof of Lusin's result that the  range of a strongly one-to-one function from a spread into $\omega^\omega$ is (positively) Borel. It is  shown that the positively Borel set $\mathbb{D}^2(\mathcal{A}_1):=\{\alpha\mid \alpha^0=\underline 0\;\vee\;\alpha^1=\underline 0\}$ positively fails to be the range of a strongly one-to-one function from a spread into  $\omega^\omega$. 

\smallskip
In Section 4,  we  establish some properties of the class $\mathbf{\Pi}^1_1$ of the co-analytic subsets of $\omega^\omega$.   We prove that the set \begin{quote} $\mathcal{WF}:=\{\alpha\mid\forall \beta\in (T_\alpha)^\omega\exists n[ \beta(n)\le_{KB}\beta(n+1)]\}$,\end{quote} i.e.: the set of all $\alpha$ such that the tree $T_\alpha$ is \textit{well-founded with respect to $<_{KB}$}, coincides with $\mathcal{A}^1_1$ and thus is $\mathbf{\Pi}^1_1$-complete.  The proof uses Brouwer's Thesis on bars in $\omega^\omega$.  We also show that the set \begin{quote}$\mathsf{Sink}^\ast(\mathcal{ALMOST^\ast FIN}):=\\\{\beta\mid Spr(\beta)\;\wedge\; \forall \alpha \in \mathcal{F}_\beta\forall \zeta\in [\omega]^\omega\exists n[\alpha\circ\zeta(n)=0]\}$,\end{quote} consisting of the codes of all spreads all of whose elements $\alpha$ have the property that $D_\alpha$ is an \textit{almost-finite} subset of $\omega$,  is  $\mathbf{\Pi}^1_1$-complete. We then prove that the set \begin{quote}$\mathcal{E}^1_1!:=\{\alpha\mid \exists \gamma[\forall n[\alpha(\overline\gamma n)=0]\;\wedge\forall \delta[\delta\;\#\;\gamma\rightarrow\exists n[\alpha(\overline \delta n)\neq 0]]\}$\end{quote} consisting  of those $\alpha$ that admit exactly one path $\gamma$  is not $\mathbf{\Pi}^1_1$ although, in classical descriptive set theory, $\mathcal{E}^1_1!$ is $\mathbf{\Pi}^1_1$-complete. It remains true that every $\mathbf{\Pi}^1_1$ set reduces to  $\mathcal{E}^1_1!$.

\smallskip
In Section 5, we  prove that there exist $\mathbf{\Sigma}^1_1$ sets that positively fail to be $\mathbf{\Pi}^1_1$ and  $\mathbf{\Pi}^1_1$ sets that positively fail to be $\mathbf{\Sigma}^{1\ast}_1$. We use Kripke's scheme $\mathbf{KS}$ in order to prove  that there are $\mathbf{\Pi}^1_1$ sets that are not  $\mathbf{\Sigma}^{1}_1$.  We also see that some $\mathbf{\Sigma}^1_1$ sets positively fail to be (positively) Borel and that some $\mathbf{\Pi}^1_1$ sets are not (positively) Borel. Using  Brouwer's Thesis on bars in $\omega^\omega$, we prove one half of Souslin's Theorem: $\mathbf{\Sigma}^{1\ast}_1\cap\mathbf{\Pi}^1_1\subseteq \mathfrak{Borel}$.  The converse statement fails intuitionistically.  

\smallskip 

In Section 6, we study the set \begin{quote}$\mathcal{ALMOST^\ast COUNT}:=\\\{\beta\mid Spr(\beta)\;\wedge\;\exists \delta \forall \gamma \in \mathcal{F}_\beta\forall \alpha\exists n[\overline \gamma \alpha(n)=\overline{\delta^n}\alpha(n)]\}$ \end{quote} of codes of \textit{almost-countable spreads}. This set is $\mathbf{\Sigma}^1_2$ and probably not $\mathbf{\Pi}^1_1$, although we have no proof of the latter fact. We prove, again using Brouwer's Thesis on bars in $\omega^\omega$, that the almost-countable spreads are just the spreads that are \textit{reducible} in Cantor's sense and that they form a hierarchy in various senses, the so-called \textit{Cantor-Bendixson hierarchy}.

\smallskip 
In Section 7, we study the class $\mathbf{\Pi}^1_2$ of the co-projections of analytic sets and the class $\mathbf{\Sigma}^1_2 $ of the projections of co-analytic sets. We prove that the Second Axiom of Continuous Choice, $\mathbf{AC}_{1,1}$, implies: $\mathbf{\Pi}^1_2\subseteq\mathbf{\Sigma}^1_2$ and  thus causes {\it the collapse of the (positive) projective hierarchy}. We draw a parallel with arithmetic, where \textit{Church's Thesis} causes the collapse of the (positive) arithmetical hierarchy.

\section{Analytic sets}\label{S:analytic}\subsection{The class $\mathbf{\Sigma}^1_1$}\label{SS:sigma11}\hfill

Some relevant definitions may be found in Subsubsection \ref{SSS:projective}.

\begin{definition}\label{D:analytic} $\mathcal{X}\subseteq \omega^\omega$ is \emph{ analytic} or $\mathbf{\Sigma}^1_1$ if and only if  there exists $\beta$ such that $\mathcal{X}=\mathcal{EF}_\beta:=Ex(\mathcal{F}_\beta)=\{\alpha\mid \exists \gamma[\ulcorner \alpha, \gamma\urcorner \in \mathcal{F}_\beta\}$. \end{definition} 

$\mathcal{X}\subseteq \omega^\omega$ thus is analytic if $\mathcal{X}$ is the projection of a closed subset of $\omega^\omega$. 

\begin{definition}\label{D:souslinsystem} A \emph{Souslin system} is a mapping $s\mapsto \mathcal{P}_s$ that associates to every $s$ a subset $\mathcal{P}_s$ of $\omega^\omega$. The \emph{Souslin operation} applied to such a system produces the set $\mathbb{A}_s\mathcal{P}_s:=\bigcup_\alpha\bigcap_n\mathcal{P}_{\overline\alpha n}$. \end{definition}

The next Theorem shows that the class $\mathbf{\Sigma}^1_1$ behaves nicely. The class is closed under the operations of countable union and countable intersection and contains all (positively) Borel subsets of $\omega^\omega$. Every set reducing to an analytic set is itself analytic. The class $\mathbf{\Sigma}^1_1$ is also closed under projection and under the Souslin operation.

\begin{theorem}\label{T:analytic} \begin{enumerate}[\upshape (i)]\hfill

 \item $\mathcal{US}^1_1:=\{\alpha\mid\alpha_{II}\in \mathcal{EF}_{\alpha_I}\}$ is $\mathbf{\Sigma}^1_1$-universal.
\item $\mathcal{E}^1_1:=\{\alpha\mid\exists \gamma\forall n[\alpha(\overline \gamma n)=0]\}$ is $\mathbf{\Sigma}^1_1$-complete.
\item For every infinite sequence $\mathcal{X}_0,\mathcal{X}_1, \ldots $ in $\mathbf{\Sigma}^1_1$, $\bigcup_n\mathcal{X}_n \in \mathbf{\Sigma}^1_1$ and $\bigcap_n\mathcal{X}_n \in \mathbf{\Sigma}^1_1$, i.e.

$\forall\beta\exists\gamma\exists\delta[\bigcup_n\mathcal{EF}_{\beta^n} =\mathcal{EF}_\gamma \;\wedge\; 
\bigcap_n\mathcal{EF}_{\beta^n} =\mathcal{EF}_\delta]$.

\item $\mathfrak{Borel}\subseteq \mathbf{\Sigma}^1_1$, i.e.      $\forall \sigma \in \mathcal{HRS}\forall \beta \exists\gamma\exists\delta[\mathcal{G}^\sigma_\beta=\mathcal{EF}_\gamma\;\wedge\;\mathcal{F}^\sigma_\beta=\mathcal{EF}_\delta]$.

\item For all $\mathcal{X}\subseteq \omega^\omega$, if $\mathcal{X}\in \mathbf{\Sigma}^1_1$, then $Ex(\mathcal{X})\in \mathbf{\Sigma}^1_1$, i.e.  $\forall\beta\exists\gamma[Ex(\mathcal{EF}_\beta)=\mathcal{EF}_\gamma]$.

\item For all $\mathcal{X},\mathcal{Y}\subseteq\omega^\omega$, if $\mathcal{X}\preceq\mathcal{Y}\in \mathbf{\Sigma}^1_1$, then $\mathcal{X} \in \mathbf{\Sigma}^1_1$, i.e.

$\forall \beta\forall \varphi:\omega^\omega\rightarrow \omega^\omega\exists\gamma[\{\alpha\mid\varphi|\alpha\in\mathcal{EF}_\beta\}=\mathcal{EF}_\gamma]$.

\item For each $\beta$, $\mathbb{A}_s\mathcal{EF}_{\beta^s} \in \mathbf{\Sigma}^1_1$.

\end{enumerate}\end{theorem}

\begin{proof} (i) For each $\alpha$, $\alpha \in \mathcal{US}^1_1\leftrightarrow \alpha_{II}\in \mathcal{EF}_{\alpha_I}\leftrightarrow\exists\gamma[\ulcorner \alpha_{II}, \gamma\urcorner \in \mathcal{F}_{\alpha_I}]\leftrightarrow\\ \exists \gamma \forall n [\alpha_I(\overline{\ulcorner \alpha_{II}, \gamma\urcorner}n)=0]$.

Define $\beta$ such that, for all $n$, for all $a,c$ in $\omega^n$,  \\$\beta(\ulcorner  a,c\urcorner)) \neq 0$ if and only if, for some $m<n$, $\overline{\ulcorner a_{II}, c\urcorner}m<n$ and  $a_{I}(\overline{\ulcorner a_{II}, c\urcorner}m)\neq 0$.

Then, for each $\alpha$, $\alpha \in \mathcal{EF}_\beta$ if and only if $\exists \gamma[\ulcorner \alpha, \gamma\urcorner  \in \mathcal{F}_\beta]$ \\if and only if $\exists \gamma\forall n[\beta(\overline{\ulcorner \alpha, \gamma\urcorner} n)=0]$ if and only if $\exists \gamma \forall n[\alpha_I(\overline{\ulcorner\alpha_{II}, \gamma\urcorner}n)=0]$ \\if and only if $\alpha_{II} \in \mathcal{EF}_{\alpha_I}$ if and only if $\alpha \in \mathcal{US}^1_1$.

 Conclude: $\mathcal{US}^1_1=\mathcal{EF}_\beta\in\mathbf{\Sigma}^1_1$.

Also: for each $\varepsilon$, $\mathcal{EF}_\varepsilon=\mathcal{US}^1_1\upharpoonright \varepsilon$. Conclude: $\mathcal{US}^1_1$ is $\mathbf{\Sigma}^1_1$-universal. 

\smallskip (ii) For each $\alpha$, $\alpha\in \mathcal{E}_1^1\leftrightarrow \exists \gamma\forall n[\alpha(\overline \gamma n)=0]$.

Define $\mathcal{F}:=\{\alpha\mid\forall n[\alpha_I(\overline{\alpha_{II}}n)=0]\}$ and note $\mathcal{E}_1^1=Ex(\mathcal{F})$. 

Define $\beta$ such that $\forall a[\beta(a)=0\leftrightarrow \forall n[\overline{a_{II}}n<length(a_I)\rightarrow a_I(\overline{a_{II}}n)=0]]$ and note: $\mathcal{F}=\mathcal{F}_\beta$. We thus see: $\mathcal{E}_1^1\in\mathbf{\Sigma}^1_1$.

Let $\varepsilon$ be given. Note: $\forall \alpha[\alpha \in \mathcal{EF}_\varepsilon \leftrightarrow \exists \gamma \forall n [\varepsilon(\overline{\ulcorner \alpha, \gamma\urcorner}n)=0]]$. 

Define $\varphi:\omega^\omega\rightarrow\omega^\omega$ such that 
$\forall \alpha\forall k\forall c\in \omega^k[ (\varphi|\alpha)(c)=\varepsilon(\ulcorner\overline \alpha k ,c\urcorner)]$. \\ Note: $\varphi$ reduces $\mathcal{EF}_\varepsilon$ to $\mathcal{E}^1_1$.
Conclude: $\mathcal{E}_1^1$ is $\mathbf{\Sigma}^1_1$-complete.

\smallskip (iii) Let $\mathcal{X}_0,\mathcal{X}_1, \ldots$ be an infinite sequence of analytic subsets of $\omega^\omega$. \\Using $\mathbf{AC}_{0,1}$, find $\beta$ such that $\forall n[\mathcal{X}_n=\mathcal{EF}_{\beta^n}]$. 

Note: for all $\alpha$, $\alpha \in \bigcup_n\mathcal{X}_n \leftrightarrow \exists n \exists \gamma[\ulcorner \alpha, \gamma\urcorner\in \mathcal{F}_{\beta^n}]\leftrightarrow \exists \gamma[\ulcorner\alpha,\gamma\circ S \urcorner\in \mathcal{F}_{\beta^{\gamma(0)}}]$. 
\\Define $\mathcal{Z}_0:=\{\ulcorner\alpha, \gamma\urcorner\mid \forall k [\beta^{\gamma(0)}(\overline{\ulcorner \alpha, \gamma\circ S\urcorner}k)=0]\}$ and note: \\$\mathcal{Z}_0\in\mathbf{\Pi}^0_1$ and $\bigcup_n\mathcal{X}_n = Ex(\mathcal{Z}_0)\in\mathbf{\Sigma}^1_1$.

Note, using $\mathbf{AC}_{0,1}$: for all $\alpha$, $\alpha \in \bigcap_n \mathcal{X}_n\leftrightarrow \forall n \exists \gamma[\ulcorner\alpha,\gamma\urcorner\in \mathcal{F}_{\beta^n}]\leftrightarrow \\\exists\gamma\forall n[\ulcorner\alpha,\gamma^n \urcorner\in \mathcal{F}_{\beta^n}]$. \\Define $\mathcal{Z}_1:=\{\ulcorner\alpha, \gamma\urcorner\mid \forall n\forall m[\beta^n(\overline{\ulcorner\alpha,\gamma^n\urcorner}m)=0]\}$ and note: \\$\mathcal{Z}_1\in\mathbf{\Pi}^0_1$ and $\bigcap_n\mathcal{X}_n = Ex(\mathcal{Z}_1)\in\mathbf{\Sigma}^1_1$.

\smallskip (iv) follows from (iii) by induction on the class of positively Borel sets.

\smallskip (v) Let $\beta$ be given. Note: for every $\alpha$, $\alpha \in Ex(\mathcal{EF}_\beta)\leftrightarrow \exists \gamma[\ulcorner\alpha,\gamma\urcorner\in \mathcal{EF}_\beta\leftrightarrow \exists\gamma\exists\delta[\ulcorner\ulcorner\alpha,\gamma\urcorner,\delta\urcorner \in \mathcal{F}_\beta]\leftrightarrow \exists\gamma[\ulcorner\ulcorner\alpha,\gamma_I\urcorner,\gamma_{II}\urcorner \in \mathcal{F}_\beta]$. 
\\Define $\mathcal{Z}:=\{\ulcorner\alpha, \gamma\urcorner\mid\forall n[\beta[\overline{\ulcorner\ulcorner\alpha,\gamma_I\urcorner,\gamma_{II}\urcorner}n=0] \}$ and  note: \\$\mathcal{Z}\in\mathbf{\Pi}^0_1$ and $Ex(\mathcal{EF}_\beta) = Ex(\mathcal{Z})\in\mathbf{\Sigma}^1_1$.

\smallskip (vi) Let $\varphi:\omega^\omega\rightarrow \omega^\omega$ and $\beta$ be given. For every $\alpha$, $\varphi|\alpha \in \mathcal{EF}_\beta \leftrightarrow \exists \gamma[\ulcorner\varphi|\alpha,\gamma\urcorner \in \mathcal{F}_\beta]$. \\Define $\mathcal{Z}:=\{\ulcorner \alpha, \gamma\urcorner \mid\forall n[\beta(\overline{\ulcorner\varphi|\alpha, \gamma\urcorner}n)=0]\}$ and note: \\$\mathcal{Z}\in\mathbf{\Pi}^0_1$ and $\{\alpha\mid\varphi|\alpha \in \mathcal{EF}_\beta\} = Ex(\mathcal{Z})\in\mathbf{\Sigma}^1_1$.

\smallskip (vii) Let $\beta$ be given. Note, using $\mathbf{AC}_{0,1}$: for each $\alpha$, $\alpha \in \mathbb{A}_s\mathcal{EF}_{\beta^s} \leftrightarrow \\\exists \gamma \forall n[\alpha \in \mathcal{EF}_{\beta^{\overline \gamma n}}] \leftrightarrow \exists \gamma\forall n \exists \delta[\ulcorner \alpha, \delta\urcorner \in \mathcal{F}_{\beta^{\overline \gamma n}}]\leftrightarrow\exists \gamma \exists \delta\forall n[\ulcorner \alpha, \delta^n\urcorner \in \mathcal{F}_{\beta^{\overline \gamma n}}]\leftrightarrow\\ \exists \gamma \forall n[\ulcorner \alpha, (\gamma_{II})^n\urcorner \in \mathcal{F}_{\beta^{\overline {\gamma_I} n}}]$.
\\Define $\mathcal{Z}:=\{\ulcorner, \alpha, \gamma\urcorner\mid \forall n\forall m[\beta^{\overline {\gamma_I} n}(\overline{\ulcorner \alpha_I, (\gamma_{II})^n\urcorner }m)=0]\}$ and note: \\$\mathcal{Z}\in\mathbf{\Pi}^0_1$ and $\mathbb{A}_s\mathcal{EF}_{\beta^s}= Ex(\mathcal{Z})\in\mathbf{\Sigma}^1_1$.

\end{proof}

\subsection{The set $\mathcal{IF}$ }\label{SS:if}
\begin{definition}\label{D:kleenebrouwer} For all $s,t$ in $\omega$ one defines: $s<_{KB} t$ if and only if\\ either $t\sqsubset s$  or $\exists i[i<length(s) \;\wedge i<length(t) \; \overline s i =\overline t i \;\wedge s(i)<t(i)]$. \end{definition}

$<_{KB}$ is a linear ordering on $\omega$. We define, for all $s,t$,  \\$\max_{KB}(s,t):= s$ if $t\le_{KB}s$, and $\max_{KB}(s,t):= t$ otherwise.

$<_{KB}$ is called the \emph{Kleene-Brouwer ordering of $\omega$.}

\begin{definition} We define $\mathcal{IF}:= \{\alpha\mid \exists \beta \in (T_\alpha)^\omega\forall n[\beta(n+1)<_{KB}\beta(n)]\}$.\end{definition}

$\mathcal{IF}$ 
is the set of all $\alpha$ such that the tree $T_\alpha:=\{s\mid\forall t\sqsubset s[\alpha(t)=0]\}$ is (positively) \textit{ill-founded} with respect to the Kleene-Brouwer-ordering $<_{KB}$.

In classical mathematics, $\mathcal{IF}=\mathcal{E}^1_1$, see also Theorem \ref{T:wf}. In our intuitionstic context, the two sets are different. The reason is that the class of all sets reducing to $\mathcal{IF}$ is not closed under the operation of finite union:

\begin{theorem}\label{T:pif}\hfill

\begin{enumerate}[\upshape (i)]
  \item The set $\mathbb{D}^2(\mathcal{A}_1)$ does not reduce to the set $\mathcal{IF}$: $\mathbb{D}^2(\mathcal{A}_1) \npreceq \mathcal{IF}$.
  
  \item The set $\mathcal{E}_1^1$ is a proper subset of the set $\mathcal{IF}$: $\mathcal{E}^1_1 \subsetneqq \mathcal{IF}$.
  
  \item The set $\mathcal{IF}$ is $\mathbf{\Sigma}^1_1$ but  not $\mathbf{\Sigma}^1_1$-complete. 
	    \end{enumerate}
\end{theorem}	    
	    
\begin{proof}
	    Assume: $\varphi:\omega^\omega\rightarrow\omega^\omega$ reduces $\mathbb{D}^2(\mathcal{A}_1)=\{\alpha\mid\alpha^0=\underline 0 \;\vee\;\alpha^1=\underline 0\}$ to $\mathcal{IF}$.
	    
	    Assume: $\alpha\in\overline{\mathbb{D}^2(\mathcal{A}_1)}$.\\ Define $\alpha_0,\alpha_1$ such that $\forall i<2[(\alpha_i)^i=\underline 0 \;\wedge\;\forall j[\neg \exists n[j=\langle i\rangle\ast n]\rightarrow \alpha_i(j)=\alpha(j)]]$.
	    \\Note: $\forall i<2[\alpha_i\in\mathbb{D}^2(\mathcal{A}_1) \;\wedge\; \varphi|\alpha_i \in \mathcal{IF}]$.
	    \\Find $\delta_0,\delta_1$ such that $\forall i<2\forall n[\delta_i(n)\in T_{\varphi|\alpha_i}\;\wedge\; \delta_i(n+1)<_{KB}\delta_i(n)]$.
	   \\Define $\zeta$ such that, for each $n$, \begin{enumerate}[\upshape (1)]\item if $\forall i<2\forall j\le n[\delta_i(j) \in T_{\varphi|\alpha}]$, then $\zeta(n)= \max_{KB}\bigl(\delta_0(n), \delta_1(n)\bigr)$, and, \item for all $i<2$, if $\exists j\le n[\delta_i(j) \notin T_{\varphi|\alpha}]$, then $\zeta(n) =\delta_{1-i}(n)$.\end{enumerate}
	    
	 This is a good definition: if, for some $i<2$,  for some $j$, $\delta_i(j) \notin T_{\varphi|\alpha}$, then $\alpha\;\#\;\alpha_i$, and, therefore, $\alpha =\alpha_{1-i}$, and, for each $j$, $\delta_{1-i}(j) \in T_{\varphi|\alpha}$.   
	\\Note: $\forall n[\zeta(n)\in T_{\varphi|\alpha} \;\wedge\;\zeta(n+1)<_{KB}\zeta(n)]$ and 
	conclude: $\varphi|\alpha\in\mathcal{IF}$, and: $\alpha \in \mathbb{D}^2(\mathcal{A}_1)$.
	 
	 We thus see: $\forall \alpha \in \overline{\mathbb{D}^2(\mathcal{A}_1)}[\alpha \in \mathbb{D}^2(\mathcal{A}_1)]$.\\  According to Theorem \ref{T:disja1}, we have a contradiction.
	 
	 Conclude: $\mathbb{D}^2(\mathcal{A}_1)\npreceq\mathcal{IF}$.
	 
	 \smallskip (ii) Assume: $\alpha\in\mathcal{E}^1_1$. Find $\gamma$ such that $\forall n[\alpha(\overline \gamma n)=0]$. \\Note: $\forall n[\overline \gamma n \in  T_\alpha\;\wedge\;\overline \gamma(n+1)<_{KB}\overline \gamma n]$ and: $\alpha \in \mathcal{IF}$. \\We thus see: $\mathcal{E}^1_1\subseteq \mathcal{IF}$.
	 
	 According to Theorem \ref{T:analytic}, $\mathbb{D}^2(\mathcal{A}_1)\preceq \mathcal{E}_1^1$, but, as we saw in (i), $\mathbb{D}^2(\mathcal{A}_1)\npreceq \mathcal{IF}$. Conclude: $\mathcal{E}_1^1 \neq \mathcal{IF}$ and: $\mathcal{E}_1^1 \subsetneqq \mathcal{IF}$.
	 
	 \smallskip (iii) Define $\mathcal{Z}:=\{\ulcorner\alpha, \gamma\urcorner\mid\forall n[\gamma(n)\in T_{\alpha}\;\wedge\;\gamma(n+1)<_{KB}\gamma(n)]\}$ and note: \\$\mathcal{Z}\in\mathbf{\Pi}^0_1$ and $\mathcal{IF}=Ex(\mathcal{Z})$. Conclude: $\mathcal{IF}$ is  
$\mathbf{\Sigma}^1_1$.\\ As, according to (i), the analytic set $\mathbb{D}^2(\mathcal{A}_1)$ does not reduce to $\mathcal{IF}$, \\$\mathcal{IF}$ is not $\mathbf{\Sigma}^1_1$-complete.	 
    \end{proof}
      
      \subsection{The sets $\mathcal{UNC}$, $\mathcal{UNC}'$ and $\mathcal{UNC}''$ } \hfill

  \begin{definition}    $\mathcal{X}\subseteq\omega^\omega$ is \emph{ (positively) uncountable} if and only if $\forall \alpha\exists \beta \in \mathcal{X}\forall n[\beta\;\#\;\alpha^n]$.
  
  $\mathcal{X}\subseteq \omega^\omega$ is \emph{weakly (positively) uncountable} if and only if \\$\exists \alpha[\alpha\in \mathcal{X}]$ and $\forall \alpha \in \mathcal{X}^\omega\exists\beta\in\mathcal{X}\forall n[\beta\;\#\;\alpha^n]$. \end{definition}
      
      Clearly, every  uncountable subset of $\omega^\omega$ is weakly uncountable. For spreads, the two notions coincide:
      
      \begin{theorem} If $\mathcal{F}\subseteq \omega^\omega$ is a spread and weakly (positively) uncountable, then $\mathcal{F}$ is  (positively) uncountable. \end{theorem}
      
      \begin{proof} Let $\beta$ be given such that $Spr(\beta)$ and $\mathcal{F}:=\mathcal{F}_\beta$ is weakly uncountable. \\Let $\rho$ be the canonical retraction of $\omega^\omega$ onto $\mathcal{F}$. \\Note: $\forall \alpha[\rho|\alpha \in \mathcal{F} \;\wedge\; (\alpha\;\#\;\rho|\alpha \rightarrow \exists n[\beta(\overline \alpha n)\neq 0])]$. 
      
      Let $\alpha$ be given. Find $\delta$ in $\mathcal{F}$ such that $\forall n[\delta \;\#\;\rho|(\alpha^n)]$. 
      
      Let $n$ be given. As the apartness relation $\#$ is co-transitive and $\delta\;\#\; \rho|(\alpha^n)$,  \\\textit{either}: $\delta \;\#\;\alpha^n$, \textit{or}: $\alpha ^n\;\#\;\rho|(\alpha^n)$.
      
      In the latter case, find   $m$  such that $\beta(\overline{ \alpha^n} m)\neq 0$. \\Note  $\beta(\overline \delta m)=0$ and conclude:  $\overline{\alpha^n}m\neq \overline \delta m$, and: $\delta \;\#\;\alpha ^n$. 
      
      Conclude: $\forall n[\delta \;\#\;\alpha ^n]$.
      
      We thus see: $\forall \alpha\exists\delta\in \mathcal{F}\forall n[\delta\;\#\;\alpha^n]$, i.e. $\mathcal{F}$ is  uncountable. \end{proof}

      The following intuitionistic theorem is the same as \cite[Theorem 2.1]{gielen}, see also \cite[Section 8]{veldman20b}, and was first proven by W. Gielen.  G. Cantor's (classical) famous {\it Perfect Set Theorem} states that $2^\omega$ embeds continuously in every uncountable  $\mathbf{\Pi}^0_1$ subset of $\omega^\omega$.  P.S.  Alexandrov and F. Hausdorff, independently,  extended the result to Borel subsets of $\omega^\omega$ and  M. Souslin showed that it also holds for  $\mathbf{\Sigma}^1_1$ subsets of $\omega^\omega$. In our intuitionistic context the Theorem holds for {\it every} subset of $\omega^\omega$. This is  due to the Second Axiom of Continuous Choice, $\mathbf{AC}_{1,1}$, see Subsubsection \ref{SSS:bcpcontchoice}.

      \begin{theorem}\label{T:gielench} $\mathcal{X}\subseteq \omega^\omega$ is  (positively) uncountable if and only if $2^\omega$ embeds into $\mathcal{X}$. \end{theorem}
      
      \begin{proof} (i) First,  assume: $\mathcal{X}\subseteq\omega^\omega$ and $2^\omega$ embeds into $\mathcal{X}$. Find $\varphi:2^\omega\rightarrowtail \mathcal{X}$. 
      
      We now prove that $\mathcal{X}$ is positively uncountable.
      
       Let $\alpha$ be given. 
       Using induction,  define $\delta$ such that, for each $n$,   \\$\delta(n)\in Bin$ and $\delta(n)\sqsubset \delta(n+1)$ and $\varphi|\bigl(\delta(n)\bigr)\perp\alpha^n],$
      as follows. \\ Define $\delta(0)=0=\langle\;\rangle$. \\ Suppose  $n$ is given such that $\delta (n)$ has been defined.\\ Find $p$ such that $\varphi|(\delta(n)\ast\overline{\underline 0}p)\perp\varphi|(\delta(n)\ast\overline{\underline 1}p)$.\\ \textit{If} $\alpha^n\perp \varphi|(\delta(n)\ast\overline{\underline 0}p)$, define $\delta(n+1):=\delta(n)\ast\overline{\underline 0}p$, and, \textit{if not}, define $\delta(n+1):=\delta(n)\ast\overline{\underline 1}p$.
      
      It will be clear that $\alpha$ satisfies the requirements.
      
   \smallskip   Now find $\varepsilon$ in $2^\omega$ such that $\forall n[\delta(n)\sqsubset \varepsilon]$ and define: $\beta:=\varphi|\varepsilon$. \\Note: $\beta\in \mathcal{X}$ and  $\forall n[\alpha^n \;\#\; \varphi|\varepsilon = \beta]$. 
      
      We thus see: $\forall \alpha\exists \beta \in \mathcal{X}\forall n[\alpha^n\;\#\;\beta]$, i.e. $\mathcal{X}$ is (positively) uncountable.
      
      \smallskip (ii) Next, assume $\mathcal{X}\subseteq \omega^\omega$ is (positively) 
      uncountable.
      
      We want to prove that $2^\omega$ embeds into $\mathcal{X}$. 
       
       Using the Second Axiom of Continuous Choice $\mathbf{AC}_{1,1}$, see Subsubsection \ref{SSS:bcpcontchoice}, \\find $\varphi:\omega^\omega\rightarrow\omega^\omega$ such that $\forall \alpha[\varphi|\alpha\in \mathcal{X}\;\wedge\;\forall n[\varphi|\alpha \;\#\;\alpha^n]]$. 
       
       \smallskip We first prove: $\forall s\exists t \exists u[s\sqsubset t\;\wedge\;s\sqsubset u\;\wedge\;\varphi|t \perp\varphi|u]$. \\Let $s$ be given. \\Define $\delta:=\varphi|(s\ast\underline 0)$ and define $\varepsilon$ such that $s\sqsubset \varepsilon\;\wedge\;\varepsilon^s=\varphi|\delta $. \\Note: $\varphi|\varepsilon\;\#\;\varepsilon^s=\varphi|\delta$. Find $m$ such that $\varphi|\overline \delta m\perp \varphi|\overline \varepsilon m$ and define $t:=\overline \delta m$ and $u:=\overline \varepsilon m$. Clearly, $t,u$ satisfy the requirements.
       
       \smallskip Now define $\zeta$ such that $\zeta(0)=0$ and, for each $s$ in $Bin$, $\zeta(s\ast\langle 0\rangle) =u'$ and \\$\zeta(s\ast\langle 1 \rangle) =u''$, where $u$ is the least $v$ such that $\zeta(s)\sqsubset v'\;\wedge\zeta(s)\sqsubset v''\;\wedge\;\varphi|v'\perp\varphi|v''$. \\Note: $\forall s  \in  Bin\forall t\in Bin[s\sqsubset t\rightarrow \zeta(s) \sqsubset \zeta(t) ]$.\\ Find $\rho:2^\omega\rightarrow \omega^\omega$ such that $\forall \gamma \in 2^\omega\forall n[\zeta(\overline \gamma n)\sqsubset \rho|\gamma]$.
       \\Find $\psi:2^\omega\rightarrow \omega^\omega$ such that $\forall \gamma \in 2^\omega\forall n[\psi|\gamma = \varphi|(\rho|\gamma)]$. Note: $\psi:2^\omega\rightarrow \mathcal{X}$. \\Also note: $\forall s\in Bin\forall t \in Bin[s\perp t\rightarrow \varphi|(\zeta(s))\perp\varphi|(\zeta(t))]$. \\Conclude:
 $ \psi:2^\omega\rightarrowtail \mathcal{X}$   and: $2^\omega$ embeds into $\mathcal{X}$.  \end{proof}

      \begin{theorem}\label{T:helpunc}\begin{enumerate}[\upshape (i)]\item The set $\omega^{(2^\omega)}$ is $\mathbf{\Sigma}^0_1$-complete.
      \item The set $(\omega^\omega)^{(2^\omega)}$ is $\mathbf{\Pi}^0_2$-complete.
       \item The set $Emb(2^\omega,\omega^\omega)$ is $\mathbf{\Pi}^0_2$-complete. \end{enumerate}\end{theorem}
       
   \begin{proof} (i) Using the Fan Theorem $\mathbf{FT}$, see Subsubsection \ref{SSS:fantheorem},  note: 
    for all $\varphi$, \\$\varphi \in \omega^{(2^\omega)}\leftrightarrow \forall \gamma \in \mathcal{C}\exists n[\varphi(\overline\gamma n)\neq 0]\leftrightarrow \exists m\forall s \in Bin_m\exists n\le m[\varphi(\overline s n)\neq 0]$. \\Conclude: $\omega^{(2^\omega)}$ is $\mathbf{\Sigma}^0_1$.
    
\smallskip  We now want to prove that the set $\mathcal{E}_1$ reduces to the set $\omega^{(2^\omega)}$. \\  Define $\varphi:\omega^\omega\rightarrow\omega^\omega$ such that 
    $\forall \alpha\forall n\forall s \in Bin_n[(\varphi|\alpha)(s)= \alpha(n)]$. \\Note that, for each $\alpha$, for each $n$, if $n=\mu p[\alpha(p)\neq 0]$ then  \\$\varphi|\alpha:2^\omega\rightarrow \omega$ and $\forall \alpha \in 2^\omega[\varphi(\alpha)=\alpha(n)-1]$. \\Clearly, $\varphi$ reduces $\mathcal{E}_1=\{\alpha\mid\exists n[\alpha(n)\neq 0]\}$ to $\omega^{(2^\omega)}$. \\As $\mathcal{E}_1$ is $\mathbf{\Sigma}^0_1$-complete, so is $\omega^{(2^\omega)}$. 
    
    \smallskip (ii) and (iii). \\We first prove that the two sets $(\omega^\omega)^{(2^\omega)}$ and $Emb(2^\omega, \omega^\omega)$ both belong to $\mathbf{\Pi}^0_2$.
   
   First note: for all $\varphi$, $\varphi \in (\omega^\omega)^{(2^\omega)}$ if and only if $ \forall n[\varphi^n\in \omega^{(2^\omega)}]$. \\Using (i), conclude: $(\omega^\omega)^{(2^\omega)}\in\mathbf{\Pi}^0_2$.
    
    \smallskip Then note,  using the Fan Theorem $\mathbf{FT}$: for all $\varphi$, \\$\varphi \in Emb(2^\omega,\omega^\omega)$ if and only if 
    \\$\varphi \in (\omega^\omega)^{(2^\omega)}$ and $\forall s \in Bin \forall \alpha \in2^\omega\forall\beta\in2^\omega\exists n[\varphi|s\ast\langle 0\rangle \ast \overline \alpha n \perp \varphi|s\ast\langle 1\rangle \ast \overline \beta n]$ \\if and only if  
    \\$\varphi\in(\omega^\omega)^{(2^\omega)}$ and $ \forall s \in Bin \exists n\forall t \in Bin_n\forall u\in Bin_n [\varphi|s\ast\langle 0\rangle \ast t \perp \varphi|s\ast\langle 1\rangle \ast u]$.

Conclude: $Emb(2^\omega,\omega^\omega)\in\mathbf{\Pi}^0_2$.

  \smallskip   We now prove that the set $\mathcal{A}_2$ reduces to both \\the set  $(\omega^\omega)^{(2^\omega)}$ and the set $Emb(2^\omega, \omega^\omega)$.

\smallskip
Define $\psi:\omega^\omega\rightarrow\omega^\omega$ such that, for all $m$,  for all $\alpha$, for all $s$ in $2^{<\omega}$,  \\{\it if} $m<length(s)$  and $\exists n<length(s)[\alpha^m(n)\neq 0]$, then $(\psi|\alpha)^m(s)= s(m)+1$, and,\\ {\it if not}, then $(\psi|\alpha)^m(s)=0$. 

 Note: for all $\alpha$, for all $m$, \\ (i) if $\alpha^m \in \mathcal{E}_1$ then $(\psi|\alpha)^m: 2^\omega\rightarrow \omega$   and, for all $\beta$ in $2^\omega$, $(\psi|\alpha)^m(\beta)=\beta(m)$, and \\(ii)  if $\exists n[(\psi|\alpha)^m(\overline \beta n)\neq 0]$ then $\alpha^m \in \mathcal{E}_1$.

 Conclude: for all $\alpha$, for all $m$, $\alpha^m \in \mathcal{E}_1$ if and only if $(\psi|\alpha)^m:2^\omega\rightarrow \omega$.

 Conclude: for all $\alpha, \alpha \in \mathcal{A}_2$ if and only if $\psi|\alpha: 2^\omega\rightarrow\omega^\omega$.
 
 We thus see that $\psi$ reduces $\mathcal{A}_2$ to $(\omega^\omega)^{(2^\omega)}$. 
 
 As $\mathcal{A}_2$ is $\mathbf{\Pi}^0_2$-complete, also $(\omega^\omega)^{(2^\omega)}$ is $\mathbf{\Pi}^0_2$-complete.
 
 Note: for all $\alpha$, if $\alpha \in \mathcal{A}_2$, then  $\psi|\alpha: 2^\omega\rightarrow \omega^\omega$ and $\forall \beta \in 2^\omega[(\psi|\alpha)|\beta=\beta]$.
 
 Conclude: for all $\alpha$, $\alpha \in \mathcal{A}_2$ if and only if $\psi|\alpha \in Emb(2^\omega,\omega^\omega)$.

We thus see that $\psi$ reduces $\mathcal{A}_2$ to $Emb(2^\omega,\omega^\omega)$. 

As $\mathcal{A}_2$ is $\mathbf{\Pi}^0_2$-complete, also  $Emb(2^\omega,\omega^\omega)$ is $\mathbf{\Pi}^0_2$-complete.
\end{proof}  We will need the next Lemma, Lemma \ref{L:logicalhelp}, in the proof of Theorem \ref{T:perfspr}(iii).
  \begin{lemma}\label{L:logicalhelp}\begin{enumerate}[\upshape (i)]\item For all finite  $A\subseteq\omega$, for every $P\subseteq A$, for every  proposition $Q$,  \\if $\forall m\in A[m\in P\;\vee\; Q]$, then  $\forall m\in A[m\in P]\;\vee\;Q$. \item For all finite sets $A,B\subseteq \omega$, for all $P\subseteq A$, for all $Q\subseteq B$, \\if $\forall m \in A\forall n\in B[m\in P\;\vee\;n\in Q]$, then  $\forall m\in A[m\in P]\;\vee\;\forall n\in B[n\in Q]$.\end{enumerate}\end{lemma}\begin{proof} (i) Use induction on the number of elements of $A$. \\If $A=\emptyset$, the statement is true. \\Now assume the statement has been proven for $A$, and $q \in \omega\setminus A$. \\We prove that the statement is true for $A\cup\{q\}$. \\Assume $P\subseteq A\cup\{q\}$ and   $\forall m\in A\cup\{q\}[m\in P\;\vee\; Q]$.\\ Then, by the induction hypothesis: $\forall m\in A[m\in P] \;\vee\; Q$ but also: $q\in P\;\vee\; Q$. Conclude: $\forall m\in A\cup\{q\}[m\in P] \;\vee\;  Q$.

\smallskip (ii) Assume:  $A,B$ are finite subsets of $\omega$, and $\forall m \in A\forall n\in B[m\in P\;\vee\;n\in Q]$. Using (i), conclude: $\forall n \in B[\forall m\in A[m\in P]\;\vee\;n\in Q]$. \\Using (i) once more, conclude: $\forall m\in A[m\in P]\;\vee\;\forall n \in B[n\in Q]$. \end{proof} 
\begin{definition}\label{D:perfspr}  For each $\beta$, we define:  \emph{$\beta$ is a perfect-spread-law},  $Pfspr(\beta)$, if and only if $Spr(\beta)$ and $\beta(0)=0$ and, for all $s$, if   $\beta(s)=0$, \\then $\exists t\exists u[s\sqsubset t\;\wedge\;s\sqsubset u\;\wedge t\perp u\;\wedge\;\beta(t)=\beta(u)=0]]$.

 If $Pfspr(\beta)$, then $\mathcal{F}_\beta=\{\alpha\mid\forall n[\beta(\overline \alpha n)=0]\}$ is called a  \emph{perfect spread}.  \end{definition}
   
   In intuitionistic real analysis it is not true that the  image of the closed interval $[0,1]$ under a continuous function is itself a closed subset of $\mathcal{R}$. One may see this from the failure of the Intermediate Value Theorem and the failure of the theorem that a continuous function from $[0,1]$ to $\mathcal{R}$ always attains its greatest value. The next Theorem brings to light related facts. The image of Cantor space $2^\omega$ under a continuous function from $2^\omega$ to $\omega^\omega$ is always a located subset of $\omega^\omega$ but not always a closed subset of $\omega^\omega$. The latter remains true, however, if the function is \emph{strongly injective}. \\ $\mathcal{F}\subseteq\omega^\omega$ is a spread if and only if $\mathcal{F}$ is both located and closed, see Subsubsection \ref{SSS:openandclosed}.
     \begin{theorem}\label{T:perfspr}\begin{enumerate}[\upshape(i)] \item Cantor space $2^\omega$  embeds into every perfect spread. \item For each $\varphi:2^\omega\rightarrow\omega^\omega$, $\varphi|2^\omega$ is a located subset of $\omega^\omega$. \item For each $\varphi:2^\omega\rightarrowtail\omega^\omega$, $\varphi|2^\omega$ is a perfect spread and a fan. \item $\neg\forall \varphi \in (\omega^\omega)^{(2^\omega)}\exists\beta[Spr(\beta)\;\wedge\;\varphi|2^\omega=\mathcal{F}_\beta]$.\end{enumerate}\end{theorem}
     
     \begin{proof}  (i) Let $\mathcal{F}\subseteq\omega^\omega$ be a perfect spread. Find $\beta$ such that $Pfspr(\beta)$ and $\mathcal{F}=\mathcal{F}_\beta$. Define $\zeta$ such that $\zeta(0)=0$ and, for all $s$ in $Bin$, $\zeta(s\ast\langle 0\rangle):= u'$ and $\zeta(s\ast\langle 1 \rangle)=u''$ where $u$ is the least $v$ such that $v'\perp v''$ and $\zeta(s)\sqsubset v'$ and $\zeta(s)\sqsubset v''$ and $\beta(v')=\beta(v'')=0$. Define $\varphi:2^\omega\rightarrow\omega^\omega$ such that $\forall \alpha \in2^\omega\forall n[\zeta(\overline \alpha n)\sqsubset \varphi|\alpha]$. \\Note: $\forall \alpha \in 2^\omega[\varphi|\alpha\in \mathcal{F}_\beta]$. \\Also note: for all $\alpha,\beta$ in $2^\omega$, if $\alpha\;\#\;\beta$, then, for some $n$, $\overline \alpha n \perp \overline \beta n$ and: $\zeta(\overline \alpha n)\perp \zeta(\overline \beta n)$, and: $\varphi|\alpha\;\#\;\varphi|\beta$. Conclude: $\varphi:2^\omega\rightarrowtail \mathcal{F}$.
     
     \smallskip (ii) Let $\varphi:2^\omega\rightarrow \omega^\omega$ be given.  We define $\delta$ as follows. Let $s$ be given. \\Note: $\forall \alpha \in 2^\omega\exists m[s\sqsubset\varphi|\overline\alpha m\;\vee\;s\perp\varphi|\overline\alpha m]$. \\Using  $\mathbf{FT}$, find $m$ such that $\forall \alpha \in 2^\omega[s\sqsubset\varphi|\overline\alpha m\;\vee\;s\perp\varphi|\overline\alpha m]$, i.e. \\$\forall t \in Bin_m[ s\sqsubset \varphi|t\;\vee\; s\perp\varphi|t]$. \\Define  $\delta(s): = 0$ if $\exists t\in Bin_m[s\sqsubset\varphi|t]$ and $\delta (s):=1$ if $\forall t \in Bin_m[s\perp\varphi|t]$. \\Conclude: $\forall s[\delta(s)=0\leftrightarrow \exists\alpha\in2^\omega[s\sqsubset\varphi|\alpha]]$ and: $\varphi|2^\omega$ is a located subset of $\omega^\omega$. \\Also note: $Fan(\delta)$ and: $\varphi|2^\omega\subseteq \mathcal{F}_\delta$.  
     
     \smallskip (iii) Let $\varphi:2^\omega\rightarrowtail \omega^\omega$ be given. 
     Using (ii), find $\delta$ such that \\$\forall s[\delta(s)= 0\leftrightarrow \exists \alpha \in 2^\omega[s\sqsubset \varphi|\alpha]]$  and: $Fan(\delta)$ and: $\varphi|2^\omega\subseteq \mathcal{F}_\delta$. 
     
     We first prove: $Pfspr(\delta)$.
     \\Let $s$ be given such that $\delta(s) =0$. Find $\alpha$ in $2^\omega$ such that $s\sqsubset\varphi|\alpha$. Find $m$  such that $s\sqsubset\varphi|\overline \alpha m$. Find $n$ such that $\varphi|(\overline\alpha m\ast\underline{\overline 0}n)\perp\varphi|(\overline\alpha m\ast\underline{\overline 1}n)$ and define: $t:=\varphi|(\overline\alpha m\ast\underline{\overline 0}n)$ and $u:=\varphi|(\overline\alpha m\ast\underline{\overline 1}n)$. Note $\delta(t)=\delta(u)=0$ and $s\sqsubset t$ and $s\sqsubset u$ and $t\perp u$. 
     
     \smallskip
     Assume $s\in Bin$. \\Note:  $\forall \alpha \in 2^\omega[\varphi|(s\ast\langle 0\rangle\ast\alpha_I )\;\#\; \varphi|(s\ast\langle 1\rangle \ast\alpha_{II})]$ and: \\$\forall \alpha \in 2^\omega\exists m[\varphi|\bigl(s\ast\langle 0\rangle\ast\overline{\alpha_I}m \bigr)\perp \varphi|\bigl(s\ast\langle 1\rangle \ast\overline{\alpha_{II}}m\bigr)]$ and, using $\mathbf{FT}$: \\$\exists m\forall \alpha \in 2^\omega[\varphi|\bigl(s\ast\langle 0\rangle\ast\overline{\alpha_I}m \bigr)\perp \varphi|\bigl(s\ast\langle 1\rangle \ast\overline{\alpha_{II}}m\bigr)]$,i.e.  \\$\exists m\forall a \in Bin_{m}\forall b\in Bin_m[\varphi|(s\ast\langle 0\rangle\ast a )\perp \varphi|(s\ast\langle 1\rangle \ast b)]$. 
     
     \smallskip Define $\zeta$ such that, for each $s$ in $Bin$,\\ $\zeta(s)$ is the least $m$ such that $\forall a \in Bin_{m}\forall b\in Bin_m[\varphi|(s\ast\langle 0\rangle\ast a)\perp \varphi|(s\ast\langle 1\rangle \ast b)]$.
     
     \smallskip
     
     We now prove: $\mathcal{F}_\delta\subseteq \varphi|2^\omega$.
    Let $\gamma\in\mathcal{F}_\delta$ be given. 
    Assume: $s\in Bin$.
   \\ Note: $\forall a \in Bin_{\zeta(s)}\forall b\in Bin_{\zeta(s)}[\varphi|(s\ast\langle 0\rangle\ast a)\perp \gamma\;\vee\;\gamma \perp\varphi|(s\ast\langle 1\rangle \ast b)]$.
   \\Conclude, using Lemma \ref{L:logicalhelp}: \\$\forall a \in Bin_{\zeta(s)}[\varphi|(s\ast\langle 0\rangle\ast a)\perp \gamma]\;\vee\;\forall a \in Bin_{\zeta(s)}[\gamma\perp \varphi|(s\ast\langle 1\rangle \ast a)]$.
   
   \smallskip
    Define $\eta$ in $2^\omega$ such that $\forall s \in Bin[\eta(s)= 1\leftrightarrow \forall a \in Bin_{\zeta(s)}[\varphi|(s\ast\langle 0\rangle\ast a)\perp \gamma]]$.
    
    Define $\alpha$ in $2^\omega$ such that $\forall n[\alpha(n) =\eta(\overline \alpha n)]$.
    
   Note that, for all $\beta$ in $2^\omega$, for all $n$, if $n=\mu p[\alpha(p)\neq \beta(p)$, then $\varphi|\beta \perp \gamma$,  i.e.,  \\for all $\beta $ in $2^\omega$, if $\beta\perp \alpha$, then $\varphi|\beta\perp \gamma$.

    We now prove: $\varphi|\alpha=\gamma$.
    \\Assume: 
 $\varphi|\alpha\perp\gamma$. Find $n$ such that $\varphi|\overline \alpha n \perp \gamma $. \\Define $m=n+\zeta(\overline\alpha n)$ and note: $\forall d\in Bim_m[ d \perp\overline \alpha n\rightarrow \varphi|d\perp \gamma]$. \\Conclude: $\forall d \in Bin_m[\varphi|d\perp\gamma]$. \\Note: $\forall d \in Bin_m[length(\varphi|d)\le m]$. Conclude: $\delta(\overline \gamma m)\neq 0$.  Contradiction. \\We thus see: $\neg(\varphi|\alpha\perp\gamma)$, and: $\varphi|\alpha=\gamma$.
    
     Conclude: $\forall \gamma \in \mathcal{F}_\delta\exists \alpha\in 2^\omega[\varphi|\alpha =\gamma]$, and: $\varphi|2^\omega=\mathcal{F}_\delta$.
   
\smallskip       (iv) Assume: $\forall \varphi \in (\omega^\omega)^{(2^\omega)}\exists\beta[Spr(\beta)\;\wedge\;\varphi|2^\omega=\mathcal{F}_\beta]$.

Using Brouwer's Continuity Principle $\mathbf{BCP}$, see Subsubsection \ref{SSS:bcpcontchoice},  we prove that this assumption leads to a contradiction as it implies $\mathbf{LPO}$, see Subsubsection \ref{SSS:lpomp}.
        
         Let $\alpha$ be given. We intend to prove: $\alpha = \underline 0 \;\vee\;\alpha\;\#\;\underline 0$.
         
          Define $\varphi:2^\omega\rightarrow \omega^\omega$ such that $\forall \gamma \in 2^\omega[\varphi|(\langle 0\rangle \ast \gamma) = \alpha \;\wedge\;\varphi|(\langle 1\rangle\ast\gamma)=\underline 0]$. \\Note: $\varphi|2^\omega=\{\alpha,\underline 0\}$.
        Find $\beta$ such that $Spr(\beta)$ and $\{\alpha,\underline 0\}=\mathcal{F}_\beta$. \\Note:  $\forall s[\beta(s) =0\leftrightarrow (s\sqsubset \alpha\;\vee\; s\sqsubset \underline 0)]$.   \\Note: $\forall \gamma \in \mathcal{F}_\beta[\gamma=\alpha\;\vee\;\gamma =\underline 0]$. \\Applying $\mathbf{BCP}$, find $m$ such that \\\textit{either}: $\forall \gamma \in \mathcal{F}_\beta[\underline{\overline 0}m\sqsubset \gamma\rightarrow \gamma = \underline 0]$, and: $ \underline{\overline 0}m\perp \alpha\;\vee\;\alpha=\underline 0$, \\\textit{or}: $\forall \gamma \in \mathcal{F}_\beta[\underline{\overline 0}m\sqsubset \gamma\rightarrow \gamma =\alpha]$, and: $\alpha =\underline 0$. \\Conclude: $\alpha =\underline 0\;\vee\;\alpha\;\#\;\underline 0$. 
        
        We thus see: $\forall \alpha[\alpha = \underline 0 \;\vee\;\alpha\;\#\;\underline 0]$, that is: $\mathbf{LPO}$, a contradiction.
        \end{proof}
     
 \begin{definition}  We introduce three subsets of $\omega^\omega$:

       $\mathcal{UNC}:=\{\beta\mid \forall \alpha\exists\gamma\in \mathcal{F}_\beta\forall n[\gamma\;\#\;\alpha^n]\}$, and $\mathcal{UNC}':=\{\beta \in \mathcal{UNC}\mid Spr(\beta)\}$ and  \\$\mathcal{UNC}'':=\{\beta\mid \forall \alpha\exists\gamma\in \mathcal{EF}_\beta\forall n[\gamma\;\#\;\alpha^n]\}$. \end{definition}
       
       $\mathcal{UNC}$,  $\mathcal{UNC}'$ and $\mathcal{UNC}''$ are the sets of the codes of (positively) uncountable \emph{closed sets},  (positively) uncountable \emph{located closed sets} and (positively) uncountable \emph{analytic sets}, respectively.

       The classical result corresponding to  the following theorem is due to W.~Hurewicz, see \cite[Theorem 27.5]{kechris}. The proof 
       in \cite{kechris} is very different from ours and not constructive.

     \begin{theorem}\label{T:uncsigma11complete}

  $\mathcal{UNC}$,  $\mathcal{UNC'}$ and $\mathcal{UNC''}$ are $\mathbf{\Sigma}^1_1$-complete.
 
    \end{theorem}
   
   \begin{proof} We first prove that $\mathcal{UNC}$ is $\mathbf{\Sigma}^1_1$.
   
   Using Theorem \ref{T:gielench}, note that, for each $\beta$, \\$\beta \in \mathcal{UNC}$ if and only if there exists $\varphi: \mathcal{C}\rightarrowtail \mathcal{F}_\beta$.
   
  Now define  $\mathcal{A}:=\{\ulcorner \beta, \varphi \urcorner\mid \varphi:2^\omega\rightarrowtail\omega^\omega\;\wedge\;\forall s\in 2^{<\omega}\forall t[t\sqsubseteq \varphi|s\rightarrow \beta(t)=0]\}$. 
  \\Then  $\mathcal{UNC}=Ex(\mathcal{A})$. Note, using Theorem \ref{T:helpunc}: $\mathcal{A} \in \mathbf{\Pi}^0_2$. \\Conclude, using Theorem \ref{T:analytic}: $\mathcal{UNC}\in \mathbf{\Sigma}^1_1$. 
   
  \smallskip We now prove that $\mathcal{UNC}$ is $\mathbf{\Sigma}^1_1$-complete.
  
   Define $\varphi:\omega^\omega\rightarrow\omega^\omega$ such that, for all $\alpha$, for all $s$,  $(\varphi|\alpha)(s) =0$ if and only if \\there exists  $u$  such that $\forall t\sqsubseteq u[\alpha(t)=0]$ and $length(u) =length(s)$ and \\$\forall i<length(s)[s(i)=2u(i)+1\;\vee\;s(i)=2u(i)+2]]$. 
   
   We prove that $\varphi$ reduces $\mathcal{E}^1_1$ to $\mathcal{UNC}$. 
   
   First, assume: $\alpha \in \mathcal{E}_1^1$. Find $\gamma$ such that $\forall n[\alpha(\overline \gamma n)=0]$. Define $\beta$ such that, for all $s$,  \\$\beta(s)=0$ if and only if $\forall i<length(s)[s(i)=2\gamma(i)+1\;\vee\;s(i) =2\gamma(i)+2]$. \\Note: $Pfspr(\beta)$ and $\mathcal{F}_\beta \subseteq \mathcal{F}_{\varphi|\alpha}$. Conclude, using Theorems \ref{T:helpunc}(i) and \ref{T:gielench}: $\varphi|\alpha \in\mathcal{UNC}$.
   
   Now let $\alpha$ be given such that  $\varphi|\alpha\in \mathcal{UNC}$. \\Using Theorem \ref{T:helpunc}, find $\beta$ such that $Pfspr(\beta)$ and $\mathcal{F}_\beta\subseteq \mathcal{F}_{\varphi|\alpha}$. \\Find $\delta$ in $\mathcal{F}_\beta$. Find $\gamma$ such that $\forall n[\delta(n) = 2\gamma(n)+1\;\vee\;\delta(n)=2\gamma(n)+2]$. \\ Conclude: $\forall n[\alpha(\overline \gamma n)=0]$ and: $\alpha \in \mathcal{E}_1^1$. 
   
   We thus see: $\mathcal{E}^1_1$ reduces to $\mathcal{UNC}$. As $\mathcal{E}_1^1$ is $\mathbf{\Sigma}^1_1$-complete, see Theorem \ref{T:analytic}, so is $\mathcal{UNC}$.
   
  \medskip We now consider $\mathcal{UNC}'$.
  
  Define $\mathcal{A}':=\{\ulcorner\beta, \varphi\urcorner \in \mathcal{A}\mid Spr(\beta)\}$. \\Note: $\mathcal{A}'\in  \mathbf{\Pi}^0_2$ and $\mathcal{UNC}'=Ex(\mathcal{A}')$. Conclude: $\mathcal{UNC}'\in \mathbf{\Sigma}^1_1$.
  
 \smallskip We now want to prove that $\mathcal{UNC}'$ is $\mathbf{\Sigma}^1_1$-complete.
 
 We would like to use again the function $\varphi$ we used in the previous  paragraph, but, unfortunately, not: for every $\alpha$, $\varphi|\alpha$ is a spread-law. 
 
  We therefore define $\psi:\omega^\omega\rightarrow\omega^\omega$  such that, for all $\alpha$, for all $s$,  \\$(\psi|\alpha)(s) =0$ if and only if there exist $k,t$ such that $(\varphi|\alpha)(t)=0$ and $s= t\ast\underline{\overline 0}k$.  
  
  Observe that, for every $\alpha$, $\psi|\alpha$ is a spread-law and $\mathcal{F}_{\varphi|\alpha}\subseteq \mathcal{F}_{\psi|\alpha}$. 
  
 We prove that $\psi$ reduces $\mathcal{E}_1^1$ to $\mathcal{UNC}'$.  
 
 First, assume: $\alpha \in \mathcal{E}_1^1$. Then $\mathcal{F}_{\varphi|\alpha} \in \mathcal{UNC}$. Note:   $\mathcal{F}_{\varphi|\alpha}\subseteq \mathcal{F}_{\psi|\alpha}$, so also $\mathcal{F}_{\psi|\alpha}$ is (positively) uncountable, and, as $\psi|\alpha$ is a spread-law, $\psi|\alpha \in \mathcal{UNC}'$.

   Now let $\alpha$ be given such that  $\psi|\alpha\in \mathcal{UNC}'$. \\Find $\beta$ such that $Pfspr(\beta)$ and $\mathcal{F}_\beta\subseteq \mathcal{F}_{\psi|\alpha}$. \\Note: for all $s$, if $\exists \gamma \in \mathcal{F}_\beta[s\sqsubset \gamma]$, then  $\forall i<length(s)[s(i)>0]]$, and $(\varphi|\alpha)(s)=0$. \\Conclude: $\mathcal{F}_\beta \subseteq \mathcal{F}_{\varphi|\alpha}$ and $\alpha \in \mathcal{E}_1^1$.

   \medskip We now consider $\mathcal{UNC}''$.
   
   Define $\mathcal{A}'':=\{\ulcorner \beta, \varphi\urcorner \mid \varphi:2^\omega\rightarrowtail \mathcal{EF}_{\beta}\}$. 
   \\Note, using the Second Axiom of Continuous Choice $\mathbf{AC}_{1,1}$, see Subsubsection \ref{SSS:bcpcontchoice}:  for every $\beta$, for every $\varphi$,
   $\varphi:2^\omega\rightarrow \mathcal{EF}_{\beta} $ if and only if $\exists \psi:2^\omega\rightarrow \omega^\omega\forall \gamma \in 2^\omega[\ulcorner \varphi|\gamma, \psi|\gamma\urcorner\in\mathcal{F}_{\beta}]$. \\
   Define $\mathcal{A}^\ast:=\{\ulcorner \beta, \varphi\urcorner \mid \varphi_I:2^\omega\rightarrowtail\omega^\omega\;\wedge\;\varphi_{II}:2^\omega\rightarrow\omega^\omega\;\wedge\;\forall \gamma \in 2^\omega[\ulcorner \varphi_I|\gamma, \varphi_{II}|\gamma\urcorner \in \mathcal{F}_\beta]\}$. Note: $\mathcal{UNC}''=Ex(\mathcal{A}'')=Ex(\mathcal{A}^\ast )$, and, using Theorem \ref{T:helpunc}: $\mathcal{A}^\ast \in \mathbf{\Pi}^0_2$. \\Conclude: $\mathcal{UNC}'' \in \mathbf{\Sigma}^1_1$.
   
   In order to see that $\mathcal{UNC}''$ is $\mathbf{\Sigma}^1_1$-complete, we remind ourselves of the fact: $\mathbf{\Pi}^0_1\subseteq\mathbf{\Sigma}^1_1$. Define $\tau:\omega^\omega\rightarrow\omega^\omega$ such that $\forall \beta\forall s[(\tau|\beta)(s) =\beta(s_I)]$ and note: $\forall \beta[\mathcal{EF}_{\tau|\beta}=\mathcal{F}_\beta]$. Conclude: $\tau$ reduces $\mathcal{UNC}$ to $\mathcal{UNC}''$, and, as $\mathcal{UNC}$ is $\mathbf{\Sigma}^1_1$-complete, so is $\mathcal{UNC}''$.
   \end{proof}

        \subsection{  $\mathsf{Share}(\mathcal{INF})$ and $\mathsf{Share}(\mathcal{INF}\cap 2^\omega)$}\hfill

     The following definition occurs already  in \cite{veldman5}.
     
     \begin{definition}\label{D:share}  For each $\mathcal{X}\subseteq \omega^\omega$,  we define \\$\mathsf{Share}(\mathcal{X}):=\{\beta\mid Spr(\beta)\;\wedge\;\exists \gamma \in \mathcal{F}_\beta[\gamma \in \mathcal{X}]\}$.\end{definition}

  If $\beta \in \mathsf{Share}(\mathcal{X})$, one says: `\textit{The spread} $\mathcal{F}_\beta$ \textit{shares an element with the set} $\mathcal{X}$'. 
  
  \begin{definition}\label{D:infinite}
   $\mathcal{INF}:=\{\alpha\mid \forall m\exists n>m[\alpha(n)\neq 0]\}$. \end{definition}If $\alpha\in\mathcal{INF}$, then $D_\alpha:=\{n\mid\alpha(n)\neq 0\}$ is a decidable and infinite subset of $\omega$. 
  
  \smallskip
  The next result corresponds to a well-known fact in classical descriptive set theory, see \cite[p. 209, Exercise 27]{kechris},  or \cite[p. 137, Exercise 4.2.3]{srivastava}.
  
  \begin{theorem}\label{T:shareinf}
  $\mathsf{Share}(\mathcal{INF})$ and $\mathsf{Share}(\mathcal{INF}\cap2^\omega)$ are $\mathbf{\Sigma}^1_1$-complete.
\end{theorem}

\begin{proof} We first observe that these two sets are indeed $\mathbf{\Sigma}^1_1$. \\Note: $\{\beta\mid Spr(\beta)\}$ is $\mathbf{\Pi}^0_2$. \\For each $\beta$, $\beta \in \mathsf{Share}(\mathcal{INF})$ if and only if \\$Spr(\beta)$ and $\exists \alpha\exists \zeta \in [\omega]^\omega\forall n[\beta(\overline \alpha n)=0\;\wedge \alpha\circ\zeta(n)\neq 0]$. \\Conclude, using Theorem \ref{T:analytic}: $\mathsf{Share}(\mathcal{INF})$ is $\mathbf{\Sigma}^1_1$.

For each $\beta$, $\beta \in \mathsf{Share}(\mathcal{INF}\cap2^\omega)$ if and only if \\$Spr(\beta)$ and $\exists \alpha\in 2^\omega\exists \zeta \in [\omega]^\omega\forall n[\beta(\overline \alpha n)=0\;\wedge \alpha\circ\zeta(n)\neq 0]$. \\ Conclude: $\mathsf{Share}(\mathcal{INF}\cap2^\omega)$ is  $\mathbf{\Sigma}^1_1$.

 \medskip We now  prove that $\mathsf{Share}(\mathcal{INF})$ and $\mathsf{Share}(\mathcal{INF}\cap 2^\omega)$ are $\mathbf{\Sigma}^1_1$-complete. 
 
 First define $\delta$ such that $\delta(0)=0$ and $\forall s \forall n[\delta(s\ast\langle n \rangle)=\delta(s)\ast\underline{\overline 0}n\ast\langle 1 \rangle]$. \\Then define $\varphi:\omega^\omega\rightarrow \omega^\omega$ such that \\$\forall \alpha\forall s[(\varphi|\alpha)(s)=0\leftrightarrow\exists n\exists t[s=\delta(t)\ast\underline{\overline 0}n \;\wedge\;\forall u\sqsubseteq t[\alpha(u)=0]]]$. \\Note that, for each $\alpha$, $Spr(\varphi|\alpha)$, i.e. $\varphi|\alpha$ is a spread-law, and $\mathcal{F}_{\varphi|\alpha}\subseteq 2^\omega$. \\We show that $\varphi$ reduces $\mathcal{E}^1_1$ to both  $\mathsf{Share}(\mathcal{INF}\cap 2^\omega)$ and $\mathsf{Share}(\mathcal{INF})$.
 
 First, assume: $\alpha\in\mathcal{E}^1_1$. Find $\gamma$ such that $\forall n[\alpha(\overline \gamma n)=0]$.
 \\Note: $\forall n\forall t[t\sqsubseteq \delta(\overline\gamma n) \rightarrow (\varphi|\alpha)(t)=0]$. \\Find $\varepsilon$ in $2^\omega$ such that $\forall n[\delta(\overline \gamma n)\sqsubset \varepsilon]$. \\Note: $\varepsilon \in \mathcal{F}_{\varphi|\alpha}$ and, as $\forall n[\varepsilon\bigl(n+\sum_{i=0}^{i=n}\gamma(i)\bigr)=1]$, also $\varepsilon \in \mathcal{INF}$. \\Conclude: $\varphi|\alpha \in \mathsf{Share}(\mathcal{INF}\cap 2^\omega)\subseteq \mathsf{Share}(\mathcal{INF})$.
 
 Now assume: $\varphi|\alpha \in \mathsf{Share}(\mathcal{INF})$.
Find $\varepsilon$ in $\mathcal{INF}\cap\mathcal{F}_{\varphi|\alpha}$. \\Define $\gamma$ such that $\gamma(0):=\mu i[\varepsilon(i) \neq 0]$ and $\forall n[\gamma(n+1)=\mu i[\varepsilon\bigl(\gamma(n)+i+1\bigr) \neq 0]$. \\Note: $\forall n[\delta(\overline \gamma n)\sqsubset \varepsilon]$ and: $\forall n[\alpha(\overline \gamma n)=0]$ and: $\alpha \in \mathcal{E}^1_1$.

\smallskip We thus see that $\varphi$ reduces $\mathcal{E}_1^1$ to both
   $\mathsf{Share}(\mathcal{INF}\cap2^\omega)$ and $\mathsf{Share}(\mathcal{INF})$. \\It follows that these sets, like $\mathcal{E}^1_1$, are $\mathbf{\Sigma}^1_1$-complete. 
     \end{proof}

\subsection{Strictly analytic subsets of $\omega^\omega$}\label{SS:strictlyanalytic}	\hfill

\begin{definition}\label{D:strictlyanalytic} $\mathcal{X}\subseteq \omega^\omega$ is \emph{ strictly analytic} or $\mathbf{\Sigma}^{1\ast}_1$ if and only if  there exists $\beta$ such that $Spr(\beta)$ and $\mathcal{X}=\mathcal{EF}_\beta:=Ex(\mathcal{F}_\beta)=\{\alpha\mid \exists \gamma[\ulcorner \alpha, \gamma\urcorner \in \mathcal{F}_\beta]\}$. \end{definition} 

$\mathcal{X}\subseteq \omega^\omega$ thus is strictly analytic if it its the projection of a closed {\it and located} subset of $\omega^\omega$, see Subsubsection \ref{SSS:openclosedspreads}.

Recall that $\mathcal{X}\subseteq \omega^\omega$ is {\it located} if and only if $\exists \alpha[\{s\mid\exists \gamma \in \mathcal{X}[s\sqsubset \gamma]\}=D_\alpha]$, i.e. \\the set $\{s\mid\exists \gamma \in \mathcal{X}[s\sqsubset \gamma]\}$ is a {\it decidable} subset of $\omega$, and  
\\ $\mathcal{X}\subseteq \omega^\omega$ is {\it semi-located} if and only if $\exists \alpha[\{s\mid\exists \gamma \in \mathcal{X}[s\sqsubset \gamma]\}=E_\alpha]$, i.e. \\the set $\{s\mid\exists \gamma \in \mathcal{X}[s\sqsubset \gamma]\}$ is an {\it enumerable} subset of $\omega$.

Also recall that, for every infinite sequence $\mathcal{X}_0, \mathcal{X}_1,  \ldots$ of subsets of $\omega^\omega$, \\$\mathbb{D}_n(\mathcal{X}_n)=\{ \gamma \mid \exists n[\gamma^n \in \mathcal{X}_n]\}$ and  $\mathbb{C}_n(\mathcal{X}_n)=\{ \gamma \mid \forall n[\gamma^n \in \mathcal{X}_n]\}$, see Subsubsection \ref{SSS:disjunction}. 

The following theorem shows that $\mathbf{\Sigma}^{1\ast}_1$ is a proper subclass of $\mathbf{\Sigma}^1_1$ and behaves less nicely.

Note that, as a consequence of  the first item of the theorem, every strictly analytic subset of $\omega^\omega$ is either empty or inhabited. 

\begin{theorem}\label{T:ansan}\hfill

\begin{enumerate}[\upshape (i)]
\item For every $\mathcal{X}\subseteq \omega^\omega$, $\mathcal{X} \in \mathbf{\Sigma}^{1\ast}_1 \leftrightarrow (\mathcal{X}=\emptyset\;\vee\;\exists\varphi:\omega^\omega\rightarrow\omega^\omega[\mathcal{X}=\varphi|\omega^\omega])$. 
\item For every $\mathcal{X} \subseteq \omega^\omega$, if $\mathcal{X} \in \mathbf{\Sigma}^{1\ast}_1$,  then $\mathcal{X}$ is semi-located. \item For every $\mathcal{X}\subseteq\omega^\omega$,   if $\mathcal{X}$ is inhabited and semi-located, then $\overline{\mathcal{X}}\in \mathbf{\Sigma}^{1\ast}_1$.

\item Not every inhabited and closed subset of $\omega^\omega$ is semi-located, i.e. \\$\neg\forall\beta[\exists \gamma[\gamma \in \mathcal{F}_\beta]\rightarrow \mathcal{F}_\beta$ is semi-located $]$. 

\item Every spread is strictly analytic but not every closed subset of $\omega^\omega$ is strictly analytic, i.e. $\forall \beta[Spr(\beta)\rightarrow \mathcal{F}_\beta \in \mathbf{\Sigma}^{1\ast}_1]$ but $\neg\forall \beta[\mathcal{F}_\beta \in \mathbf{\Sigma^{1\ast}_1}]$, i.e.  $\neg(\mathbf{\Pi}^0_1\subseteq \mathbf{\Sigma}^{1\ast}_1)$. 

\item Semi-located  and closed subsets of $\omega^\omega$ are not always located subsets of $\omega^\omega$, i.e.  \\$\neg\forall \beta[\mathcal{F}_\beta$ is semi-located $\rightarrow\mathcal{F}_\beta$ is located$\;]$.

\item The closure of an open subset of $\omega^\omega$ is not always a closed subset of $\omega^\omega$, i.e.\\ $\neg\forall \beta\exists \gamma[\mathcal{F}_\gamma =\overline{\mathcal{G}_\beta}]$. 
\item $\mathbf{\Sigma}^{1\ast}_1$ is  closed under the operation of (finite) union but $\mathbf{\Sigma}^{1\ast}_1$ is not closed under the operation of (finite) intersection, because: \\$\neg\forall \beta[\{\beta\}\cap\{\underline0\} \in\mathbf{\Sigma}^{1\ast}_1]$ and: $\neg\forall \beta[\{\beta, \underline 1\}\cap\{\underline 0, \underline 1\} \in\mathbf{\Sigma}^{1\ast}_1]$.
\item $\mathbf{\Sigma}^{1\ast}_1$ is not closed under the operation of countable union, because:\\ $\neg\forall \alpha[\bigcup_n\{\beta\mid\beta=\underline 0 \;\wedge\;\alpha(n)\neq 0\}\in\mathbf{\Sigma}_1^{1\ast}]$.
\item For every infinite sequence $\mathcal{X}_0,\mathcal{X}_1, \mathcal{X}_2, \ldots$ of strictly analytic and \emph{inhabited} subsets of $\omega^\omega$, the sets $\bigcup_n\mathcal{X}_n$, $\mathbb{D}_n(\mathcal{X}_n)$  and $\mathbb{C}_n(\mathcal{X}_n)$ are strictly analytic. 
\item For every strictly analytic $\mathcal{X}\subseteq\omega^\omega$, $Ex(\mathcal{X})$ is strictly analytic.

\end{enumerate}
\end{theorem}

\begin{proof} (i) First, assume: $\mathcal{X} \in \mathbf{\Sigma}^{1\ast}_1$. Find $\beta$ such that $Spr(\beta)$ and $\mathcal{X}=Ex(\mathcal{F}_\beta)$. \\There are two cases: $\beta(0)\neq 0$ and $\beta(0)=0$. \\In the first case: $\mathcal{X} =\mathcal{F}_\beta =\emptyset$. \\In the second case, let $\rho:\omega^\omega\rightarrow \mathcal{F}_\beta$ be the canonical retraction\footnote{see Subsubsection \ref{SSS:continuousfunctions}} of $\omega^\omega$ onto $\mathcal{F}_\beta$. \\Define $\varphi:\omega^\omega\rightarrow \omega^\omega$ such that $\forall \alpha[\varphi|\alpha=(\rho|\alpha)_I]$ and note: $\mathcal{X}=\varphi|\omega^\omega$. 

Conversely, let $\mathcal{X}\subseteq \omega^\omega$ and $\varphi:\omega^\omega\rightarrow\omega^\omega$ be given such that $\mathcal{X}=\varphi|\omega^\omega$. 
 \\Define $\beta$ in $\mathcal{C}$ such that $\beta(\langle\;\rangle)=0$ and, for each $n>0$, for each  $s$ in $\omega^n$,    \\ $\beta(s)=0$ if and only if $\forall i<n-1[ s_{II}(i)\sqsubset s_{II}(i+1)]]$ and $\forall i<n[\overline{s_I}(i+1)\sqsubseteq \varphi|\bigl(s_{II}(i)\bigr)]$. \\Note: $Spr(\beta)$ and: $\mathcal{Y}=\mathcal{F}_\beta$ and: $\varphi|\omega^\omega=\mathcal{Y}$. 
  
  \smallskip (ii) Assume: $\mathcal{X}\in \mathbf{\Sigma}^{1\ast}_1$, that is, by (i): either $\mathcal{X}=\emptyset$ or $\exists\varphi:\omega^\omega\rightarrow\omega^\omega[\mathcal{X}=\varphi|\omega^\omega]$. Note: $\emptyset$ is semi-located. Now assume:   $\mathcal{X}$ is inhabited. \\Find $\varphi:\omega^\omega\rightarrow\omega^\omega$ such that $\mathcal{X}=\varphi|\omega^\omega$. Note: $\forall s[\exists \gamma[s\sqsubset\varphi|\gamma]\leftrightarrow \exists t[s\sqsubset\varphi|t]]$. \\Define $\delta$ such that $\forall n[(n_I\sqsubseteq \varphi|n_{II}\rightarrow\delta(n)=n_I+1)\;\wedge\; (\neg(n_I\sqsubseteq \varphi|n_{II})\rightarrow\delta(n)=0)]$. Note: $E_\delta=\{s\mid\exists \gamma[s\sqsubset\varphi|\gamma]\}$ and conclude: $\mathcal{X}=\varphi|\omega^\omega$ is semi-located.

  \smallskip (iii)
  Assume: $\mathcal{X}\subseteq \omega^\omega$ is inhabited and semi-located. \\Find $\delta$ such that $E_\delta=\{s\mid\exists  \gamma \in\mathcal{X}[s\sqsubset\gamma]\}$. \\Note: $\exists n[\delta(n)=\langle\;\rangle +1=1]$ and: $\forall s\in E_\delta\exists n \exists p[\delta(n)=s\ast \langle p\rangle +1]$. \\Define $\varepsilon$ such that $\varepsilon(0)=0$ and, for all $s,n$, \\\textit{if} $\exists p[\delta(n)=\varepsilon(s)\ast\langle p \rangle +1]$, then $\varepsilon(s\ast\langle n \rangle)=\delta(n)-1$, and, \\ \textit{if not}, then $\varepsilon(s\ast\langle n\rangle)=\delta(m)-1$, where $m=\mu q[\exists p[\delta(q) =\varepsilon(s)\ast\langle p\rangle+1]]$. \\Now define $\varphi:\omega^\omega\rightarrow\omega^\omega$ such that $\forall \alpha\forall n[\varepsilon(\overline \alpha n)\sqsubset \varphi|\alpha]$ and note: $\overline{\mathcal{X}}=\varphi|\omega^\omega$. 
  
  \smallskip
  (iv) Define $\varphi:\omega^\omega\rightarrow\omega^\omega$ such that $\forall \alpha\forall s[(\varphi|\alpha)(s) =0\leftrightarrow \bigl(s\sqsubset\underline 0\;\vee\;(s\sqsubset\underline 1 \;\wedge\; \underline{\overline 0}s\sqsubset \alpha )\bigr)]$. Note: $\forall \alpha\forall \gamma[\gamma \in \mathcal{F}_{\varphi|\alpha}\leftrightarrow \bigl(\gamma = \underline 0\;\vee\;(\gamma = \underline 1 \;\wedge\;\alpha = \underline 0)\bigr)]$. \\Assume: $\forall \alpha[\mathcal{F}_{\varphi|\alpha}$ \textit{is semi-located}$]$.  \\Using $\mathbf{AC}_{1,1}$ , find $\psi:\omega^\omega\rightarrow\omega^\omega$ such that $\forall \alpha[E_{\psi|\alpha}=\{s\mid\exists \gamma \in  \mathcal{F}_{\varphi|\alpha}[s\sqsubset \gamma]\}]$. \\Note: $\langle 1 \rangle \in E_{\psi|\underline 0}$. Find $p$ such that $(\psi|\underline 0)(p)=\langle 1 \rangle +1$. \\Find $q$ such that $\psi^p(\overline{\underline 0}q)=\langle 1 \rangle +2$  and $\forall i<q[\psi^p(\underline{\overline 0}i)=0]$. \\ Note: $\forall \alpha[\overline{\underline 0}q\sqsubset \alpha \rightarrow \langle 1 \rangle \in E_{\psi|\alpha}]$. \\Conclude: $\forall \alpha[\overline{\underline 0}q\sqsubset \alpha \rightarrow \alpha = \underline 0]$. Contradiction. 
  
  \smallskip (v) Let $\beta$ be given such that $Spr(\beta)$. Define $\gamma$ such that $\forall s[\gamma(s) =0 \leftrightarrow \beta(s_I)=0]$. Note: $Spr(\gamma)$ and $\mathcal{F}_\beta = Ex(\mathcal{F}_\gamma)$. Conclude: $\mathcal{F}_\beta \in \mathbf{\Sigma}^{1\ast}_1$. 
  
   Assume: $\mathbf{\Pi}^0_1\subseteq \mathbf{\Sigma}^{1\ast}_1$. Then, according to (ii): $\forall \beta[\mathcal{F}_\beta$ is semi-located$]$. \\This conclusion contradicts (iv).
   
   \smallskip (vi) Assume: $\forall \beta[\mathcal{F}_\beta$ is semi-located $\rightarrow \mathcal{F}_\beta$ is located$]$. 
   \\Let $\alpha$ be given. \\Define $\beta$ such that $\forall s[\beta(s) =0 \leftrightarrow \bigl(length(s) \ge 1 \rightarrow \alpha\circ s(0)\neq 0\bigr)]$. \\Note: $\mathcal{F}_\beta=\{\gamma\mid \alpha\circ\gamma(0)\neq 0]\}$. \\Define $\delta$ such that for each $n$, \textit{if} either: $length(n_I) \ge 1$ and $\alpha\circ n_I(0)\neq 0$ or: $n_I =0=\langle\;\rangle$ and  $\alpha(n_{II})\neq 0$, then: $\delta(n)=n_I+1$, and, \textit{if not}, then $\delta(n) =0$.  \\Note: $E_\delta =\{s\mid\exists \gamma \in \mathcal{F}_\beta[s\sqsubset \gamma]\}$. 
   Conclude: $\mathcal{F}_\beta$ is semi-located. \\Using the above assumption, conclude: $\mathcal{F}_\beta$ is located.\\ Find $\varepsilon$ such that $E_\delta=D_\varepsilon$.  Note: if $\varepsilon(0)=0$, then $0\notin D_\varepsilon=E_\delta$ and $\forall n[\alpha(n)=0]$ and, if $\varepsilon(0)\neq 0$, then $0\in D_\varepsilon=E_\delta$ and $\exists n[\alpha(n)\neq 0]$. \\Conclude:  $\forall n[\alpha(n)=0]\;\vee\;\exists n[\alpha(n)\neq 0]$.
   \\We thus see that our assumption implies $\mathbf{LPO}$ and  is contradictory, see Subsubsection \ref{SSS:lpomp}.

  \smallskip (vii) Assume: $\forall \beta\exists \gamma[\mathcal{F}_\gamma=\overline{\mathcal{G}_\beta}]$. 
  
  Define $\varphi:\omega^\omega\rightarrow\omega^\omega$ such that $\forall \alpha\forall s[(\varphi|\alpha)(s) =0 \leftrightarrow (s\perp\underline 0\;\wedge\;\overline\alpha s\perp \underline 0)]$. \\Note: $\mathcal{G}_{\varphi|\underline 0} = \emptyset$, and, for every $\alpha$, if $\alpha\;\#\;\underline 0$, then $\mathcal{G}_{\varphi|\alpha}=\{\delta\mid\delta\;\#\;\underline 0\}$. \\By our assumption: $\forall \alpha \exists \gamma[\mathcal{F}_\gamma = \overline{\mathcal{G}_{\varphi|\alpha}}]$. \\Using $\mathbf{AC}_{1,1}$, find $\rho:\omega^\omega\rightarrow\omega^\omega$ such that $\forall \alpha[\mathcal{F}_{\rho|\alpha}=\overline{\mathcal{G}_{\varphi|\alpha}}]$. \\Note: $\mathcal{F}_{\rho|\underline 0} = \emptyset$, and, for every $\alpha$, if $\alpha\;\#\;\underline 0$, then $\mathcal{F}_{\rho|\alpha}=\omega^\omega$. \\Assume: we find $n$ such that $(\rho|\underline 0)(\overline{\underline 0}n)\neq 0$. \\Determine $p$ such that $\forall \alpha[\overline{\underline 0}p\sqsubset \alpha\rightarrow (\rho|\alpha)(\overline{\underline 0}n)\neq 0]$. \\Conclude: $\forall \alpha[\overline{\underline 0}p\sqsubset \alpha\rightarrow \underline 0 \notin \mathcal{F}_{\rho|\alpha}]$. Contradiction. \\Conclude: $\forall n[(\rho|\underline 0)(\overline{\underline 0}n)= 0]$ and: $\underline 0 \in \mathcal{F}_{\rho|\underline 0}$. Contradiction.
  
  \smallskip (viii) Assume $\mathcal{X}_0,\mathcal{X}_1\subseteq \omega^\omega$ are strictly analytic. It suffices to consider the case that both $\mathcal{X}_0,\mathcal{X}_1$ are inhabited. Find $\varphi$ such that $\forall i<2[\varphi^i:\omega^\omega\rightarrow\omega^\omega \;\wedge\;\mathcal{X}_i=\varphi^i|\omega^\omega]$.  Define $\psi:\omega^\omega\rightarrow\omega^\omega$ such that $\forall \alpha\forall n[\psi|(\langle 0\rangle\ast\alpha)=\varphi^0|\alpha\;\wedge\;\psi|(\langle n+1\rangle\ast\alpha)=\varphi^1|\alpha]$ and note: $\mathcal{X}_0\cup\mathcal{X}_1 =\psi|\omega^\omega$.
  
  \smallskip Assume: $\forall \beta[\{\beta\}\cap\{\underline 0\} \in \mathbf{\Sigma}^{1\ast}_1]$. \\Using (i), conclude: $\forall \beta[\{\beta\}\cap \{\underline 0\} =
  \emptyset\;\vee\;\exists \gamma[\gamma \in \{\beta\}\cap\underline 0\}]]$, and: $\forall \beta[\beta \neq \underline 0\;\vee\;\beta =\underline 0]$. Using $\mathbf{BCP}$, find $p$ such that either: $\forall \beta[\underline{\overline 0}p\sqsubset\beta\rightarrow \beta\neq \underline 0]$ or: $\forall \beta[\underline{\overline 0}p\sqsubset\beta\rightarrow \beta = \underline 0]$. Both alternatives are false, so we obtain a contradiction.
  
  \smallskip Now assume: $\forall \beta[\{\beta,\underline 1\}\cap\{\underline 0, \underline 1\} \in \mathbf{\Sigma}^{1\ast}_1]$.
  According to (ii), for each $\beta$, \\$\{\beta,\underline 1\}\cap\{\underline 0, \underline 1\} $ is semi-located, i.e.: $\exists \delta[E_\delta=\{s\mid \exists \gamma \in\{\beta,\underline 1\}\cap\{\underline 0, \underline 1\}[s\sqsubset \gamma]\}]$. \\Using $\mathbf{AC}_{1,1}$, find $\varphi:\omega^\omega\rightarrow\omega^\omega$ such that, for each $\beta$,  \\$E_{\varphi|\beta}=\{s\mid \exists \gamma \in\{\beta,\underline 1\}\cap\{\underline 0, \underline 1\}[s\sqsubset \gamma]\}]$. \\Note: $\langle 0\rangle \in E_{\varphi|\underline 0}$ and find $p$ such that $(\varphi|\underline 0)(p) =\langle 0\rangle +1$. \\Find $m$ such that $\forall\beta[\underline{\overline 0}m\sqsubset \beta\rightarrow (\varphi|\beta)(p)=(\varphi|\underline 0)(p)]$. \\Conclude: $\forall \beta[\underline{\overline 0}m\sqsubset \beta\rightarrow \langle 0\rangle \in E_{\varphi|\beta}]$ and:  $\forall \beta[\underline{\overline 0}m\sqsubset \beta\rightarrow  \underline 0 \in \{\beta,\underline 1\}\cap\{\underline 0, \underline 1\}]$, that is: $\forall \beta[\underline{\overline 0}m\sqsubset \beta\rightarrow \beta = \underline 0]$, a contradiction.
  
  \smallskip (ix) Assume: $\forall \alpha[\bigcup_n\{\beta\mid\beta=\underline 0\;\wedge\;\alpha(n)\neq 0\} \in \mathbf{\Sigma}^{1\ast}_1$. Then, according to (i), $\forall \alpha[\bigcup_n\{\beta\mid\beta=\underline 0\;\wedge\;\alpha(n)\neq 0\}=\emptyset \;\vee\; \exists \gamma[\gamma \in \bigcup_n\{\beta\mid\beta=\underline 0\;\wedge\;\alpha(n)\neq 0\}]]$, and: $\forall \alpha[\forall n[\alpha(n)=0]\;\vee\;\exists n[\alpha(n)\neq 0]]$, that is: $\mathbf{LPO}$, a contradiction, see Subsubsection \ref{SSS:lpomp}.
  
  \smallskip (x)
   Let $\mathcal{X}_0,\mathcal{X}_1,\ldots$ be an infinite sequence of inhabited strictly analytic subsets of $\omega^\omega$. Using (i) and $\mathbf{AC}_{0,1}$, find $\varphi$ such that $\forall n[\varphi^n:\omega^\omega\rightarrow\omega^\omega\;\wedge\;\mathcal{X}_n=\varphi^n|\omega^\omega]$. 
   
   Define $\psi:\omega^\omega\rightarrow\omega^\omega$ such that, for all $n$, for all $\alpha$, $\psi|(\langle n \rangle \ast \alpha) = \varphi^n|\alpha$ and note: $\bigcup_n\mathcal{X}_n = \psi|\omega^\omega$ is strictly analytic.
  
  Define $\rho:\omega^\omega\rightarrow\omega^\omega$ such that, for all $n$, for all $\alpha$, $\bigl(\rho|(\langle n \rangle \ast \alpha)\bigr)^n = \varphi^n|(\alpha^n)$ and, for all $i\neq n$,  $\bigl(\rho|(\langle n\rangle\ast\alpha)\bigr)^i=\alpha^i$ and note: $\mathbb{D}_n\mathcal{X}_n = \rho|\omega^\omega$ is strictly analytic.
  
  Define $\tau:\omega^\omega\rightarrow\omega^\omega$ such that, for all $n$, for all $\alpha$, $(\tau| \alpha)^n = \varphi^n|(\alpha^n)$ and conclude: $\mathbb{C}_n\mathcal{X}_n = \tau|\omega^\omega$ is strictly analytic.
  
  \smallskip (xi) Assume $\mathcal{X}\subseteq \omega^\omega$ is strictly analytic. \\Then, according to (i), one may decide: $\mathcal{X}=\emptyset$ or: $\mathcal{X}$ is inhabited. \\Note: $Ex(\emptyset)=\emptyset$ is strictly analytic. \\If $\mathcal{X}$ is inhabited, find $\varphi:\omega^\omega\rightarrow\omega^\omega$ such that $\mathcal{X} = \varphi|\omega^\omega$. Define $\psi:\omega^\omega\rightarrow \omega^\omega$ such that $\forall \alpha[\psi|\alpha=(\varphi|\alpha)_I]$ and note $Ex(\mathcal{X})= \psi|\omega^\omega$ is strictly analytic.
\end{proof}

Using Theorem \ref{T:ansan}(x), one may prove: for every $\sigma$ in $\mathcal{HRS}$, $\mathcal{E}_\sigma$ and $\mathcal{A}_\sigma$ are strictly analytic. The sets $\mathcal{E}_\sigma, \mathcal{A}_\sigma$, are the leading sets of the intuitionistic Borel hierarchy, see Subsubsection \ref{SSS:borelsetsgeneral}.

\medskip We conclude our discussion of strictly analytic subsets of $\omega^\omega$ by observing that {\it Kripke's scheme} $\mathbf{KS}$, see Subsubsection \ref{SSS:creasubj}, makes the gap between analytic and strictly analytic subsets of $\omega^\omega$ somewhat smaller.
\begin{theorem}\label{T:kripkesan} (Using $\mathbf{KS}$:)\begin{enumerate}[\upshape (i)] \item Every \emph{inhabited} and \emph{definite} closed subset of $\omega^\omega$ is strictly analytic. \item Every \emph{inhabited} and \emph{definite} analytic subset of $\omega^\omega$ is strictly analytic. \end{enumerate} \end{theorem}

\begin{proof} (i) Assume $\mathcal{F}\subseteq \omega^\omega$ is inhabited,  definite and  closed. According to Theorem \ref{T:kripkescheme} in Subsubsection \ref{SSS:creasubj}, $\mathcal{F}$ is semi-located. According to theorem \ref{T:ansan}(iii), $\mathcal{F}$ is strictly analytic.

\smallskip (ii) Assume $\mathcal{X}\subseteq\omega^\omega$ is inhabited, definite and analytic. Find $\mathcal{F}$ in $\mathbf{\Pi}^0_1$ such that $\mathcal{X}=Ex(\mathcal{F})$. Note that $\mathcal{F}$ is inhabited. We assume that also $\mathcal{F}$ is definite. According to (i), $\mathcal{F}$ is strictly analytic. According to Theorem \ref{T:ansan}(xi), also $\mathcal{X}=Ex(\mathcal{F})$ is strictly analytic. 
\end{proof}

  John Burgess, in \cite{Burgess}, also studies strictly analytic subsets of $\omega^\omega$, or, as he called them, using a term of 
of Brouwer's and following \cite{gielen},  ``{\em dressed spreads}''.  Avoiding $\mathbf{AC}_{1,1}$, he does not restrict application of the Brouwer-Kripke scheme to definite propositions and concludes: every inhabited analytic subset of $\omega^\omega$ is strictly analytic. The argument given  for Theorem \ref{T:kripkesan}(ii) is essentially his.

 \section{Separation theorems}
 
\subsection{Results by Lusin and Novikov} \hfill
\begin{definition}\label{D:apartsets}Let $\mathcal{X}, \mathcal{Y}$ be subsets of $ \omega^\omega$. 

We define:    the pair $(\mathcal{X}, \mathcal{Y})$ is \emph{positively disjoint}, notation: $\mathcal{X}\;\#\;\mathcal{Y}$,  if and only if,\\ for all $\alpha$ in $\mathcal{X}$, for all $\beta$ in $\mathcal{Y}$, $\alpha\;\#\;\beta$.

We also define:    the pair $(\mathcal{X}, \mathcal{Y})$ is \emph{Borel-separable}), notation: $\mathcal{X}\;\#^\mathfrak{Borel}\;\mathcal{Y}$, if and only if there exist (positively) Borel sets $\mathcal{A}, \mathcal{B}$ such that  $\mathcal{X}\subseteq \mathcal{A}$, $\mathcal{Y}\subseteq \mathcal{B}$ and $\mathcal{A} \;\#\;\mathcal{B}$. \end{definition}

\begin{lemma}\label{L:borelseparable1} Let $\mathcal{Y}, \mathcal{X}_0, \mathcal{X}_1, \mathcal{X}_2, \ldots$ be an infinite sequence of subsets of $\omega^\omega$.

If, for each $n$, $\mathcal{Y}\;\#^\mathfrak{Borel}\;\mathcal{X}_n$, then $\mathcal{Y}\;\#^\mathfrak{Borel}\;\bigcup_n \mathcal{X}_n$.  \end{lemma}
\begin{proof}Assume: for each $n$, $\mathcal{Y} \;\#^\mathfrak{Borel}\;\mathcal{X}_n$. \\Find\footnote{We are silently applying the Second Axiom of Countable Choice $\mathbf{AC}_{0,1}$, as Borel sets should be thought as given by means of their codes, see Subsubsection \ref{SSS:borelsetsgeneral}. We do so at other occasions too, without further warning.}, for each $n$, Borel sets $\mathcal{A}_n, \mathcal{B}_n$ such that $\mathcal{Y}\subseteq \mathcal{A}_n$ and $\mathcal{X}_n\subseteq \mathcal{B}_n$ and $\mathcal{A}_n \;\#\;\mathcal{B}_n$.\\Define $\mathcal{A}:= \bigcap_n \mathcal{A}_n$ and $\mathcal{B}:=\bigcup_n \mathcal{B}_n$. \\ Note: $\mathcal{A}, \mathcal{B}$ are Borel and $\mathcal{Y}\subseteq \mathcal{A}$ and $\bigcup_n \mathcal{X}_n \subseteq \mathcal{B}$ and $\mathcal{A}\;\#\;\mathcal{B}$.\\Conclude: $\mathcal{Y}\;\#^\mathfrak{Borel}\;\bigcup_n\mathcal{X}_n$. \end{proof}

A version of the next theorem occurs in \cite[Theorem 18.4.1, p. 163]{veldman1}. A related result is proven in \cite{aczel}.

\begin{theorem}[Lusin's Separation Theorem]\label{T:lusinseparation}\hfill

 Let $\mathcal{X},\mathcal{Y}\subseteq\omega^\omega$ be strictly analytic. If  $\mathcal{X}\;\#\; \mathcal{Y}$, then $\mathcal{X}\;
 \#^\mathfrak{Borel}\; \mathcal{Y}$.  \end{theorem}

\begin{proof} Let $\mathcal{X},\mathcal{Y}\subseteq\omega^\omega$ be strictly analytic. Assume: $\mathcal{X}\;\#\;\mathcal{Y}$. 

If  $\mathcal{X}=\emptyset$,  we define $\mathcal{A}:=\emptyset$ and $\mathcal{B}:=\omega^\omega$,  and are done. 

If   $\mathcal{Y}=\emptyset$,   we define  $\mathcal{A}:=\omega^\omega$ and $\mathcal{B}:=\emptyset$,  and are done.

We thus may assume that $\mathcal{X},\mathcal{Y}$ are inhabited. \\Find $\varphi, \psi:\omega^\omega\rightarrow \omega^\omega$ such that $\mathcal{X} =\varphi|\omega^\omega$ and $\mathcal{Y}=\psi|\omega^\omega$. 

\smallskip Define $B:=\{s\mid\varphi|s^0\perp\psi|s^1\}$. 

We first prove that $B$ is a bar in $\omega^\omega$. 
 
 Let $\alpha$ be given. Find $n$ such that $\overline{\varphi|\alpha^0}n\perp\overline{\psi|\alpha^1}n$. \\Then  find $m$ such that $\overline{\varphi|\alpha^0}n\sqsubseteq \varphi|\overline{\alpha^0}m$ and $\overline{\psi|\alpha^1}n\sqsubseteq \psi|\overline{\alpha^1}m$.  \\Find $p$ such that $\overline{\alpha^0}m\sqsubseteq (\overline \alpha p)^0$ and $\overline{\alpha^1}m\sqsubseteq (\overline \alpha p)^1$ and note: $\overline \alpha p \in B$. 
 
 We thus see: $\forall \alpha\exists p[\overline \alpha p \in B]$, i.e. $B$ is a bar in $\omega^\omega$.

 \smallskip
 Now define $C:=\{s\mid \varphi|(\omega^\omega\cap s^0) \;\#^\mathfrak{Borel}\;\psi|(\omega^\omega\cap s^1)  \}$.

 \smallskip We first  prove: $B\subseteq C$.
 
 Let $s$ in $B$ be given. Then $\varphi|s^0\perp\psi|s^1$. 
 
 Define $\mathcal{A}:= \omega^\omega \cap (\varphi|s^0)$ and $\mathcal{B}:= \omega^\omega \cap (\psi|s^0)$. \\Note: $(\mathcal{A}, \mathcal{B})$  is a (positively) disjoint pair of basic open sets and  \\$\varphi|(\omega^\omega\cap s^0)\subseteq \mathcal{A}$ and $\psi|(\omega^\omega\cap s^1)\subseteq \mathcal{B}$. \\Conclude: $s \in C$.
 
 We thus see: $\forall s \in B[s\in C]$, i.e. $B\subseteq C$.
 
 \smallskip
 Note that $C$ is monotone: for each $s$, for each $n$, $s^0\sqsubseteq (s\ast\langle n \rangle)^0$ and $s^1\sqsubseteq (s\ast\langle n \rangle)^1$, and, therefore, if $s\in C$, also, $s\ast\langle n \rangle \in C$.
 
 \smallskip
 We finally prove that $C$ is inductive.
 
 Let $s$ be given such that $\forall n[s\ast\langle n\rangle \in C]$. We want to prove: $s\in C$. 
 
 Consider $k:=length(s)$ 
 and  distinguish three cases.
 
 {\it Case (a)}. $\neg \exists i<2\exists t[k=\langle i \rangle\ast t]$. Then, for each $n$, 
    $(s\ast\langle n\rangle)^0=s^0$ and $(s\ast\langle n\rangle)^1=s^1$. \\Note: $s\ast\langle 0 \rangle \in C$, and, therefore, also $s \in C$.

{\it Case (b)}.   $\exists t[k=\langle 0\rangle \ast t]$. Then, for all $n$, $(s\ast\langle n \rangle)^0=s^0\ast\langle n \rangle$  and $(s\ast\langle n \rangle)^1= s^1$.
\\Conclude: for all $n$, $\varphi|(\omega^\omega\cap s^0\ast\langle n \rangle)\;\#^\mathfrak{Borel}\;\psi|(\omega^\omega\cap s^1)$. \\Note: $\varphi|(\omega^\omega\cap s^0)=\bigcup_n\varphi|(\omega^\omega\cap s^0\ast\langle n \rangle)$. \\Conclude, using Lemma \ref{L:borelseparable1}, $\varphi|(\omega^\omega\cap s^0)\;\#^\mathfrak{Borel}\;\psi|(\omega^\omega\cap  s^1)$, i.e. $s\in C$.

{\it Case (c)}. $\exists t[k=\langle 1 \rangle\ast t]$. Then, for all $n$, $(s\ast\langle n \rangle)^0=s^0$  and $(s\ast\langle n \rangle)^1= s^1\ast \langle n \rangle$.
\\Conclude: for all $n$, $\varphi|(\omega^\omega\cap s^0)\;\#^\mathfrak{Borel}\;\psi|(\omega^\omega\cap s^1\ast\langle n \rangle))$. \\Note: $\psi|(\omega^\omega\cap s^1)=\bigcup_n\psi|(\omega^\omega\cap s^1\ast\langle n \rangle)$. \\Conclude, using Lemma \ref{L:borelseparable1}, $\varphi|(\omega^\omega\cap s^0)\;\#^\mathfrak{Borel}\;\psi|(\omega^\omega\cap  s^1)$, i.e. $s\in C$.

 \smallskip
 
 Using the Principle of Bar Induction $\mathbf{BI}$, see Subsubsection \ref{SSS:barinduction},  we conclude: \\$\langle\;\rangle \in C$, i.e. $\varphi|\omega^\omega\;\#^\mathfrak{Borel} \;\psi|\omega^\omega$.
\end{proof}

\begin{definition}\label{D:nopointincommon} Let $\mathcal{X}_0,\mathcal{X}_1, \ldots$ be an infinite sequence of subsets of $\omega^\omega$. \\We define: the infinite sequence $\mathcal{X}_0, \mathcal{X}_1, \ldots$ \emph{positively refuses to have a common point}, or: \emph{is $\omega$-separate},  notation: $\#_n\mathcal{X}_n$, if and only if, for every $\alpha$, \\if $\forall n[\alpha^n \in \mathcal{X}_n]$, then  $\exists i\exists j[\alpha^i\perp\alpha^j].$

We also define:  the infinite sequence $\mathcal{X}_0, \mathcal{X}_1, \ldots$ is \emph{$\omega$-Borel separable}, \\notation: $\#_n^\mathfrak{Borel}\mathcal{X}_n$,  if and only if there exists an infinite sequence $\mathcal{B}_0, \mathcal{B}_1$ $\ldots$\\ of (positively) Borel sets such that $\forall n[\mathcal{X}_n\subseteq \mathcal{B}_n]$ and $\#_n \mathcal{B}_n$. \end{definition}

\begin{lemma}\label{L:borelnovikovseparable} Let $\mathcal{Y}_0, \mathcal{Y}_1, \ldots$ and $\mathcal{X}_0, \mathcal{X}_1, \ldots$ be infinite sequences of subsets of $\omega^\omega$.

If, for each $n$, the infinite sequence $\mathcal{Y}_n, \mathcal{X}_0, \mathcal{X}_1, \ldots$ is $\omega$-Borel separable, \\then also the infinite sequence $\bigcup_n \mathcal{Y}_n, \mathcal{X}_0, \mathcal{X}_1, \ldots$ is $\omega$-Borel-separable. \end{lemma}

\begin{proof}Assume: for each $n$, the infinite sequence $\mathcal{Y}_n, \mathcal{X}_0, \mathcal{X}_1, \ldots$ is $\omega$-Borel separable, \\Find, for each $n$, an infinite and $\omega$-separate sequence  $\mathcal{B}_n, \mathcal{C}_{n,0}, \mathcal{C}_{n,1}, \ldots $ of (positively) Borel sets such that  $\mathcal{Y}\subseteq \mathcal{B}_n$ and, for each $i$,  $\mathcal{X}_i\subseteq \mathcal{C}_{n,i}$.\\Define $\mathcal{B}:= \bigcup_n \mathcal{B}_n$ and, for each $i$, $\mathcal{C}_i:=\bigcap_n \mathcal{C}_{n,i}$. \\ Note: $\mathcal{B}$ is Borel, and for each $i$, $\mathcal{C}_i$ is Borel and $\bigcup_n\mathcal{Y}_n\subseteq \mathcal{B}$ and, for each $i$,  $\mathcal{X}_i \subseteq \mathcal{C}_i$ and the infinite sequence $\mathcal{B},\mathcal{C}_0, \mathcal{C}_1, \ldots$ is $\omega$-separate.\\Conclude: the infinite sequence $\bigcup_n \mathcal{Y}_n, \mathcal{X}_0, \mathcal{X}_1, \ldots$ is $\omega$-Borel-separable. \end{proof}
 \begin{theorem}[Novikov's Separation Theorem]\label{T:novikov}
\hfill

 Let $\mathcal{X}_0,\mathcal{X}_1,\ldots$ be an infinite sequence  of  \emph{inhabited} strictly analytic subsets of $\omega^\omega$. \\If $\#_n(\mathcal{X}_n)$, then $\#^\mathfrak{Borel}_n(\mathcal{X}_n)$. 
\end{theorem}

\begin{proof} Let $\mathcal{X}_0,\mathcal{X}_1, \ldots$ be an infinite sequence of \emph{inhabited} strictly analytic subsets of $\omega^\omega$ such that $\#_n(\mathcal{X}_n)$. Using $\mathbf{AC}_{0,1}$, find $\varphi$ such that $\forall n[\varphi^n:\omega^\omega\rightarrow\omega^\omega\;\wedge\;\mathcal{X}_n=\varphi^n|\omega^\omega]$. 

\smallskip
Define $B:=\{s\mid\exists i\exists j[\varphi^i|s^i\perp\varphi^j|s^j]\}$. 

We first prove that $B$ is a bar in $\omega^\omega$.  
 
 Let $\alpha$ be given. Find $i,j,n$ such that $\overline{\varphi^i|\alpha^i}n\perp\overline{\varphi^j|\alpha^j}n$. \\Then find $m$ such that   $\overline{\varphi^i|\alpha^i}n\sqsubseteq \varphi^i|(\overline{\alpha^i}m)$ and $\overline{\varphi^j|\alpha^j}n\sqsubseteq \varphi^j|(\overline{\alpha^j}m)$. \\ Find $p$ such that $\overline{\alpha^i}m\sqsubseteq (\overline \alpha p)^i$ and $\overline{\alpha^j}m\sqsubseteq (\overline \alpha p)^j$ and note: $\overline \alpha p \in B$. 
 
 We thus see: $\forall\alpha\exists p[\overline \alpha p \in B]$.

\smallskip Define $C:=\{s\mid \#_n^\mathfrak{Borel} \varphi^n|(\omega^\omega\cap s^n)\}$. 
 
 \smallskip \emph{Note that, for each $p$, $\langle p \rangle \in C$ if and only if $\langle \;\rangle \in C$}, as, for each $p$, $\langle p \rangle^0=\langle p\rangle ^1=\langle\;\rangle$, see Subsubsection \ref{SSS:finseq}.

  We prove: $B\subseteq C$.
 
 Let $s$ in $B$ be given. Find $i,j$ such that $\varphi^i|s^i\perp\varphi^j|s^j$. \\Define an infinite sequence $\mathcal{B}_0, \mathcal{B}_1,  \ldots$ of subsets of $\omega^\omega$ such that $\mathcal{B}_i = \omega^\omega \cap \varphi^i|s_i$ and $\mathcal{B}_j = \omega^\omega \cap \varphi^j|s_j$, and, for all $k$, if $k\neq i$ and $k\neq j$, then $\mathcal{B}_k=\omega^\omega$. \\ Note that, for all $n$, $\mathcal{B}_n$ is Borel and $\varphi^n|(\omega^\omega \cap s^n)\subseteq \mathcal{B}_n$. Also note: $\#_n\mathcal{B}_n$. \\Conclude  $s \in C$.
 
 We thus see: $\forall s  \in B[s \in C]$, i.e. $B\subseteq C$.

 \smallskip Note that $C$ is monotone as, \\for all $s,t$, for all $\psi:\omega^\omega\rightarrow \omega^\omega$, if $s\sqsubseteq t$, then $\psi|(\omega^\omega\cap t)\subseteq \psi|(\omega^\omega\cap s)$.
 
 \smallskip

 We finally prove that $C$ is inductive.
 
 Let $s$ be given such that $\forall n[s\ast\langle n\rangle \in C]$. We want to prove: $s\in C$.
 
 Consider $k:=length(s)$. 
 
 \textit{Cases (a)}. $k=0$.  Then $s=\langle\;\rangle$ and $s\ast\langle 0\rangle =\langle 0\rangle$ and $s\ast\langle 0  \rangle \in C$ and, therefore,  $s\in C$. 
 
 \textit{Case (b)}. $k\neq 0$. Find $i$ such that $k=\langle i\rangle\ast t$. \\Note: for each $n$, $(s\ast\langle n \rangle)^i =s^i\ast\langle n\rangle$, and, for all $j\neq i$. $(s\ast\langle n \rangle)^{j}=s^{j}]]$. \\Conclude: for each $n$, the infinite sequence of sets $$\varphi_0|(\omega^\omega\cap s_0), \varphi_1|(\omega^\omega\cap s_1), \ldots \varphi_{i-1}|(\omega^\omega\cap s_{i-1}), \varphi_i|(\omega^\omega\cap s_i\ast\langle n \rangle), \varphi_{i+1}|(\omega^\omega\cap s_{i+1}), \ldots$$ is $\omega$-Borel-separable.  Note: $\varphi_i|(\omega^\omega\cap s_i)=\bigcup_n \varphi_i|(\omega^\omega\cap s_i\ast\langle n \rangle)$. \\Conclude, using Lemma \ref{L:borelnovikovseparable} that the infinite sequence of sets: $$\varphi_0|(\omega^\omega\cap s_0), \varphi_1|(\omega^\omega\cap s_1), \ldots \varphi_{i-1}|(\omega^\omega\cap s_{i-1}), \varphi_i|(\omega^\omega\cap s_i), \varphi_{i+1}|(\omega^\omega\cap s_{i+1}), \ldots$$ is $\omega$-Borel-separable, i.e. $s \in C$.

\smallskip  Using the Principle of Bar Induction $\mathbf{BI}$, we conclude: $\langle\;\rangle \in C$, i.e. $\#^\mathfrak{Borel}_n\varphi_n|\omega^\omega$. 
 \end{proof}
 
 \subsection{Lusin's representation Theorem}
\begin{definition}\label{D:lusinregular}We define: $\mathcal{X}\subseteq \omega^\omega$ is \emph{regular in Lusin's sense} if and if there exists a spread $\mathcal{F}\subseteq \omega^\omega$ and a strongly injective function $\varphi:\mathcal{F}\rightarrowtail \omega^\omega$ such that $\varphi|\mathcal{F}=\mathcal{X}$.\end{definition}
\begin{theorem}[One half of Lusin's Regular Representation Theorem]\label{T:lusinrephalf}
\hfill

 For all $\mathcal{X}\subseteq \omega^\omega$, if $\mathcal{X}$ is regular in Lusin's sense, then $\mathcal{X}$ is  positively Borel.
\end{theorem}

\begin{proof} Let $\beta,\varphi$ be given such that $Spr(\beta)$ and $\varphi:\mathcal{F}_\beta\rightarrowtail \omega^\omega$.

 Note: for all $s, t$, if   $\beta(s)=\beta(t)=0$ and $s\perp t$, then $\varphi|(\mathcal{F}_\beta\cap s)\;\#\;\varphi|(\mathcal{F}_\beta\cap t)$.
 
 Using Theorem \ref{T:lusinseparation}, find for all $s,t$ such that $\beta(s)=\beta(t)=0$ and $s<_{lex}t$ a positively disjoint pair $(\mathcal{B}_{s,t,0},\mathcal{B}_{s,t,1})$ of Borel sets such that \\$\varphi|(\mathcal{F}_\beta \cap s)\subseteq \mathcal{B}_{s,t,0}$ and $\varphi|(\mathcal{F}_\beta \cap t)\subseteq \mathcal{B}_{s,t,1}$.

 Define, for each $s$ such that $\beta(s) =0$, $$\mathcal{D}_s:=\bigcap_{\beta(t)=0, 
 s<_{lex}t}\mathcal{B}_{s,t,0} \;\cap \bigcap_{\beta(t)=0,t<_{lex}s}\mathcal{B}_{t,s,1}.$$
 
 Note: for all $s$, if  $\beta(s)=0$, then $\mathcal{D}_s$ is (positively) Borel and  $\varphi|(\mathcal{F}_\beta\cap s)\subseteq \mathcal{D}_s$. \\Also note:   for all $s,t$, if $\beta(s)=\beta(t)=0$ and $s\perp t$ then  $\mathcal{D}_s\;\#\;\mathcal{D}_t$.
 \\Note: $\forall \gamma \in \mathcal{F}\forall n[\varphi|\gamma \in \mathcal{D}_{\overline \gamma n}]$ and $\forall \alpha \forall s[\bigl(\beta(s) =0 \;\wedge\;\alpha \in \mathcal{D}_s\bigr)\rightarrow \varphi|s \sqsubset \alpha]$.

\smallskip Now define, for each $n$, $$\mathcal{H}_n=\bigcup \{\mathcal{D}_s\mid \beta(s)=0\;\wedge\;s\in \omega^n\},$$ and note: $\forall n[\varphi|\mathcal{F}_\beta\subseteq \mathcal{H}_n]$. \\We thus see: $\varphi|\mathcal{F}_\beta \subseteq \bigcap_n \mathcal{H}_n$ and now prove: $\bigcap_n \mathcal{H}_n \subseteq \varphi|\mathcal{F}_\beta$.

Assume: $\alpha \in \bigcap_n\mathcal{H}_n$. Find  $\delta$ such that, for each $n$, $\delta(n)\in \omega^n$ and  $\alpha \in \mathcal{D}_{\delta(n)}$. \\Note: $\forall n[\delta(n)\sqsubset \delta(n+1)]$ and find $\gamma$ such that $\forall n[\delta(n)\sqsubset \gamma]$. \\Note: $\gamma \in \mathcal{F}_\beta$ and $\forall n[\varphi|\delta(n)\sqsubset \alpha]$.  Conclude: $\varphi|\gamma =\alpha$ and: $\alpha \in \varphi|\mathcal{F}_\beta$. 

We thus see: $\varphi|\mathcal{F}=\bigcap_n\mathcal{H}_n$ is (positively) Borel.
\end{proof}

Theorem \ref{T:lusinrephalf} shows: if $\mathcal{X}\subseteq\omega^\omega$ is regular in Lusin's sense, then $\mathcal{X}$ is (positively) Borel. The converse, a famous result in classical descriptive set theory, can not be true intuitionistically, as every $\mathcal{X}\subseteq\omega^\omega$ that is regular in Lusin's sense is strictly analytic, and, as we know from theorem \ref{T:ansan}(v), it is not even true that every closed $\mathcal{X}\subseteq\omega^\omega$ is strictly analytic. The next result shows that the converse 
of Theorem \ref{T:lusinrephalf} is also not true for strictly analytic sets.

\begin{theorem}\label{T:d2a1notregularlusin}\hfill \begin{enumerate}[\upshape (i)]\item Let $\mathcal{F}\subseteq \omega^\omega$ be a spread and let $\varphi:\mathcal{F}\twoheadrightarrow \mathbb{D}^2(\mathcal{A}_1)=\{\gamma\mid\gamma^0=\underline 0 \;\vee\;\gamma^1=\underline 0\}$ be surjective. There exist $\alpha, \gamma$ in $\mathcal{F}$ such that $\alpha\;\#\;\gamma$ and $\varphi|\alpha =\varphi|\gamma=\underline 0$. 
 \item $\mathbb{D}^2(\mathcal{A}_1)=\{\gamma\mid\gamma^0=\underline 0 \;\vee\;\gamma^1=\underline 0\}$ is strictly analytic and not regular in Lusin's sense.
 \item $\mathcal{A}_1$, $\mathcal{E}_1$, $\mathcal{A}_2$ are regular in Lusin's sense and $\mathcal{E}_2$ is not. \end{enumerate}
\end{theorem}

\begin{proof} (i) Define, for  both $i<2$, $\mathcal{P}_i:= \{\gamma\mid \gamma^i=\underline 0\}$. \\Note: $\mathbb{D}^2(\mathcal{A}_1)=\mathcal{P}_0\cup\mathcal{P}_1$ and $\mathcal{P}_0,\mathcal{P}_1$ are spreads.

\smallskip Assume:  $Spr(\beta)$  and  $\varphi:\mathcal{F}_\beta\twoheadrightarrow \mathbb{D}^2(\mathcal{A}_1)=\{\gamma\mid\gamma^0=\underline 0 \;\vee\;\gamma^1=\underline 0\}$ is surjective.

Find $\alpha$ in $\mathcal{F}$ such that $\varphi|\alpha=\underline 0$.

Note: $\forall \gamma \in \mathcal{F}\exists i <2[(\varphi|\gamma)^i=\underline 0]$.

Applying Brouwer's Continuity Principle $\mathbf{BCP}$, find $m$ and $i<2$ such that \\$\forall \gamma \in \mathcal{F}\cap\overline\alpha m[(\varphi|\gamma)^i=\underline 0]$. 

Again applying $\mathbf{BCP}$, find $n, s$ such that $s\in \omega^m$ and $\beta(s)=0$ and 
\\$\forall \delta \in \mathcal{P}_{1-i}\cap\overline{\underline 0}n \exists \gamma \in \mathcal{F}_\beta\cap s[ \varphi|\gamma = \delta]$.

Now distinguish two cases. 

\smallskip \textit{Case (a)}. $s\sqsubset \alpha$. 

Define $\delta$ in $\mathcal{P}_{1-i}\cap \underline{\overline 0}n$ such that $\delta^i\;\#\;\underline 0$.
Find $\gamma$ in $\mathcal{F}_\beta\cap s$ such that $\varphi|\gamma =\delta$. \\Conclude: $\overline \alpha m\sqsubset \gamma$ and $\delta^i=(\varphi|\gamma)^i=\underline 0$. 
 Contradiction. 

Conclude: Case (a) can not occur.

\smallskip \textit{Case (b)}. $s\perp \alpha$.

 Now find $\gamma$ in $\mathcal{F}_\beta\cap s$ such that $\varphi|\gamma = \underline 0$ and note: $\alpha\;\#\;\gamma$ and $\varphi|\alpha = \varphi|\gamma=\underline 0$. 

\medskip (ii) As we saw in (i),  $\mathbb{D}^2(\mathcal{A}_1)=\mathcal{P}_0\cup\mathcal{P}_1$ and $\mathcal{P}_0,\mathcal{P}_1$ are spreads. Conclude, using Theorem \ref{T:ansan}(v) and (viii): $\mathbb{D}^2(\mathcal{A}_1)$ is strictly analytic. It also follows from (i) that $\mathbb{D}^2(\mathcal{A}_1)$ is not regular in Lusin's sense.

\smallskip (iii) Note $\mathcal{A}_1$ is a spread, and every spread is  regular in Lusin's sense, for obvious reasons.

Define $\varphi:\omega^\omega\rightarrow\omega^\omega$ such that $\forall \alpha [\varphi|\alpha=\overline{\underline 0}\alpha(0)\ast\langle \alpha(1)+1\rangle\ast S\circ S\circ\alpha]$ and note: $\varphi:\omega^\omega\rightarrowtail \omega^\omega$ and $\varphi|\omega^\omega=\mathcal{E}_1$, so $\mathcal{E}_1$ is regular in Lusin's sense. 

Define $\psi:\omega^\omega\rightarrow\omega^\omega$ such that $\forall \alpha\forall n[(\psi|\alpha)^n = \varphi|(\alpha^n)]$ and note: $\psi:\omega^\omega\rightarrowtail \omega^\omega$ and $\psi|\omega^\omega=\mathcal{A}_2$, so $\mathcal{A}_2$ is regular in Lusin's sense. 

Assume: $\mathcal{F}\subseteq \omega^\omega$ is a spread, and $\varphi:\mathcal{F}\twoheadrightarrow \mathcal{E}_2$ is surjective.  

Slightly adapting the argument given in (i), the reader may find $\alpha, \gamma$ in $\mathcal{F}$ such that $\alpha \;\#\;\gamma$ and $\varphi|\alpha =\varphi|\gamma = \underline 0$. 

 Conclude: $\mathcal{E}_2$ is \textit{not} regular in Lusin's sense. \end{proof}

Theorem \ref{T:d2a1notregularlusin} shows that it is not so easy, for a strictly analytic (positively) Borel set, to be regular in Lusin's sense. The set $\mathcal{E}_2!$, to be discussed in the next Section, see Theorem 6.4, is an example of a set that is positively Borel and strictly analytic and also regular in Lusin's sense, but, like the set $\mathbb{D}^2(\mathcal{A}_1)$, fails to be co-analytic. It is not true, therefore, that positively Borel sets  regular in Lusin's sense must be co-analytic.

\medskip

 Lusin would perhaps have been disappointed that there is no satisfying intuitionistic counterpart to the other half of Lusin's Theorem.
 He once observed that his representation theorem may help one to believe, in spite of possible qualms about generalized inductive definitions, that,
 after all,  the collection of all positively Borel subsets of $\omega^\omega$ is a \textit{well-defined set}, see \cite{lusin}, pp. 38-39, and \cite{suslin}.

\section{Co-analytic sets}
\subsection{The class $\mathbf{\Pi}^1_1$}
\hfill

Some relevant definitions may be found in Subsubsection \ref{SSS:projective}.

\begin{definition}\label{D:co-analytic} $\mathcal{X}\subseteq \omega^\omega$ is \emph{ co-analytic} or $\mathbf{\Pi}^1_1$ if and only if  there exists $\beta$ such that $\mathcal{X}=\mathcal{UG}_\beta:=Un(\mathcal{G}_\beta)=\{\alpha\mid \forall \gamma[\ulcorner \alpha, \gamma\urcorner \in \mathcal{G}_\beta\}$. \end{definition} 

$\mathcal{X}\subseteq \omega^\omega$ thus is co-analytic if $\mathcal{X}$ is the co-projection of an open subset of $\omega^\omega$.

The next Theorem shows that the class $\mathbf{\Pi}^1_1$ behaves not so nicely as the class $\mathbf{\Sigma}^1_1$. The class $\mathbf{\Pi}^1_1$ is closed under the operation of countable  intersection but not under the operation of finite union. Most (positively) Borel subsets of $\omega^\omega$ are not co-analytic. Fortunately, every set reducing to a co-analytic set is itself co-analytic. The class $\mathbf{\Pi}^1_1$ is also closed under co-projection.

\begin{theorem}\label{T:coan} \hfill
 \begin{enumerate}[\upshape (i)]

  \item  $\mathcal{UP}^1_1:=\{\alpha\mid \alpha_{II}\in \mathcal{UG}_{\alpha_I}\}$ is  $\mathbf{\Pi}^1_1$-universal.\item $\mathcal{A}_1^1:=\{\alpha\mid \forall \gamma\exists n[\alpha(\overline \gamma n)\neq 0]\}$ is $\mathbf{\Pi}^1_1$-complete.
 \item For every infinite sequence $\mathcal{X}_0,\mathcal{X}_1,\ldots$ in $\mathbf{\Pi}^1_1$,  $\bigcap_n \mathcal{X}_n\in \mathbf{\Pi}^1_1$, i.e. \\$\forall \beta\exists \gamma[\bigcap_n \mathcal{UG}_{\beta^n}=\mathcal{UG}_\gamma]$.
\item  $\mathbb{D}^2(\mathcal{A}_1)\notin\mathbf{\Pi}^1_1$. \item $\mathbf{\Pi}^0_2 \subseteq \mathbf{\Pi}^1_1$ and $\mathbf{\Sigma}^0_2\nsubseteq \mathbf{\Pi}^1_1$. 
\item For all $\mathcal{X}\subseteq\omega^\omega$, if $\mathcal{X}\in\mathbf{\Pi}^1_1$, then $Un(\mathcal{X})\in \mathbf{\Pi}^1_1$, i.e. $\forall \beta\exists \gamma[Un(\mathcal{UG}_\beta)=\mathcal{UG}_\gamma]$.

\item For all $\mathcal{X},\mathcal{Y}\subseteq\omega^\omega$, if  $\mathcal{Y}\preceq\mathcal{X}\in \mathbf{\Pi}^1_1$, then $\mathcal{X} \in \mathbf{\Pi}^1_1$, i.e. \\$\forall\beta\forall\varphi :\omega^\omega\rightarrow \omega^\omega\exists \gamma[\{\alpha\mid\varphi|\alpha \in \mathcal{UG}_\beta\}=\mathcal{UG}_\gamma]$.
\end{enumerate}
\end{theorem}

\begin{proof}  (i) For each $\alpha$, $\alpha \in \mathcal{UP}^1_1\leftrightarrow \alpha_{II}\in \mathcal{UG}_{\alpha_I}\leftrightarrow\forall\gamma[\ulcorner \alpha_{II}, \gamma\urcorner \in \mathcal{G}_{\alpha_I}]\leftrightarrow\\ \forall \gamma \exists n [\alpha_I(\overline{\ulcorner \alpha_{II}, \gamma\urcorner}n)\neq 0]$.

Define $\beta$ such that, for all $n$, for all $a,c$ in $\omega^n$,  \\$\beta(\ulcorner  a,c\urcorner)) \neq 0$ if and only if, for some $m<n$, $\overline{\ulcorner a_{II}, c\urcorner}m<n$ and  $a_{I}(\overline{\ulcorner a_{II}, c\urcorner}m)\neq 0$.

Then, for each $\alpha$, $\alpha \in \mathcal{UG}_\beta$ if and only if $\forall \gamma[\ulcorner \alpha, \gamma\urcorner  \in \mathcal{G}_\beta]$ \\if and only if $\forall \gamma\exists n[\beta(\overline{\ulcorner \alpha, \gamma\urcorner} n)\neq 0]$ if and only if $\forall \gamma \exists n[\alpha_I(\overline{\ulcorner\alpha_{II}, \gamma\urcorner}n)\neq0]$ \\if and only if $\alpha_{II} \in \mathcal{UG}_{\alpha_I}$ if and only if $\alpha \in \mathcal{UP}^1_1$.

 Conclude: $\mathcal{UP}^1_1=\mathcal{UG}_\beta\in\mathbf{\Pi}^1_1$.

\smallskip Also: for each $\varepsilon$, $\mathcal{UG}_\varepsilon=\mathcal{UP}^1_1\upharpoonright \varepsilon$. Conclude: $\mathcal{US}^1_1$ is $\mathbf{\Sigma}^1_1$-universal.

\smallskip (ii)  For each $\alpha$, $\alpha\in \mathcal{A}_1^1\leftrightarrow \forall \gamma\exists n[\alpha(\overline \gamma n)\neq 0]$.

Define $\mathcal{G}:=\{\alpha\mid\exists n[\alpha_I(\overline{\alpha_{II}}n)\neq 0]\}$ and note $\mathcal{A}_1^1=Un(\mathcal{G})$. 

Define $\beta$ such that $\forall a[\beta(a)\neq 0\leftrightarrow \exists n[\overline{a_{II}}n<length(a_I)\;\wedge\; a_I(\overline{a_{II}}n)\neq 0]]$ and note: $\mathcal{G}=\mathcal{G}_\beta$. We thus see: $\mathcal{E}_1^1\in\mathbf{\Pi}^1_1$.

\smallskip Let $\varepsilon$ be given. Note: $\forall \alpha[\alpha \in \mathcal{UG}_\varepsilon \leftrightarrow \forall \gamma \exists n [\varepsilon(\overline{\ulcorner \alpha, \gamma\urcorner}n)\neq 0]]$. 

Define $\varphi:\omega^\omega\rightarrow\omega^\omega$ such that 
$\forall \alpha\forall k\forall c\in \omega^k[ (\varphi|\alpha)(c)=\varepsilon(\ulcorner\overline \alpha k ,c\urcorner)]$. \\ Note: $\varphi$ reduces $\mathcal{UG}_\varepsilon$ to $\mathcal{A}^1_1$.
Conclude: $\mathcal{A}_1^1$ is $\mathbf{\Pi}^1_1$-complete.

\smallskip
(iii) Let $\mathcal{X}_0,\mathcal{X}_1,\ldots$ be an infinite sequence of co-analytic subsets of $\omega^\omega$. Using $\mathbf{AC}_{0,1}$, find $\beta$ such that $\forall n[\mathcal{X}_n=\mathcal{UG}_{\beta^n}]$. Define $\mathcal{V}_0:= \{\alpha\mid\exists m[\beta^{\alpha_{II}(0)}(\overline{\ulcorner\alpha_I, \alpha_{II}\circ S\urcorner}m)\neq 0]\}$. Then: $\mathcal{V}_0 \in \mathbf{\Sigma}^0_1$ and, for all $\alpha$, $\alpha \in \bigcap_n \mathcal{X}_n\leftrightarrow \forall n\forall \gamma[\ulcorner\alpha,\gamma\urcorner\in\mathcal{G}_{\beta^n}]\leftrightarrow \alpha \in Un(\mathcal{V}_0).$ 

Conclude: $\bigcap_n\mathcal{X}_n \in \mathbf{\Pi}^1_1$. 

\smallskip (iv) Assume $\mathbb{D}^2(\mathcal{A}_1)\in\mathbf{\Pi}^1_1$. Using (ii), find $\varphi:\omega^\omega\rightarrow\omega^\omega$ reducing $\mathbb{D}^2(\mathcal{A}_1)$ to $\mathcal{A}^1_1$. Assume: $\alpha\in\overline{\mathbb{D}^2(\mathcal{A}_1)}$. \\Define $\alpha_0,\alpha_1$ such that $\forall i<2[(\alpha_i)^i =\underline 0 \;\wedge\;\forall m[\neg\exists j[m=\langle i\rangle\ast j \rightarrow \alpha_i(m)=\alpha(m)]$. \\Note: $\forall i<2[\alpha_i\in \mathbb{D}^2(\mathcal{A}_1)\;\wedge\; (\alpha \;\#\;\alpha_i\rightarrow \alpha =\alpha_{1-i})]$. \\Let $\gamma$ be given. \\Find $m,n$ such that $(\varphi|\alpha_0)(\overline \gamma n) \neq 0$ and $\forall \beta[\overline{\alpha_0}m\sqsubset \beta\rightarrow (\varphi|\alpha_0)(\overline \gamma n)=(\varphi|\beta)(\overline \gamma n)]$. \\Now distinguish two cases. \textit{Either}: $\overline{\alpha_0}m\sqsubset \alpha$ and: $(\varphi|\alpha)(\overline \gamma n)\neq 0$ \\\textit{or}: $\alpha \;\#\;\alpha_0$ and $\alpha =\alpha_1$ and $\exists p[(\varphi|\alpha)(\overline \gamma p)=(\varphi|\alpha_1)(\overline \gamma p)\neq 0]$. \\In both cases: $\exists p[(\varphi|\alpha)(\overline \gamma p)\neq 0]$. \\Conclude: $\forall \gamma \exists p[(\varphi|\alpha)(\overline \gamma p)\neq 0]$, that is: $\varphi|\alpha \in \mathcal{A}^1_1$ and: $\alpha \in \mathbb{D}^2(\mathcal{A}_1)$. \\Conclude: $\forall \alpha \in \overline{\mathbb{D}^2(\mathcal{A}_1)}[\alpha 
\in\mathbb{D}^2(\mathcal{A}_1)]$, a contradiction, according to Theorem \ref{T:disja1} in Subsubsection \ref{SSS:disjunction}.\\
Conclude: $\mathbb{D}^2(\mathcal{A}_1)\notin \mathbf{\Pi}^1_1$.

\smallskip (v) Assume: $\mathcal{G}\in\mathbf{\Sigma}^0_1$. Define $\mathcal{V}:=\{\alpha\mid\alpha_I \in \mathcal{G}\}$. Then $\mathcal{V}\in\mathbf{\Sigma}^0_1$ and $\mathcal{G}=Un(\mathcal{V})\in \mathbf{\Pi}^1_1$. Conclude: $\mathbf{\Sigma}^0_1\subseteq\mathbf{\Pi}^1_1$ and, using (iii): $\mathbf{\Pi}^0_2\subseteq\mathbf{\Pi}^1_1$. 

Note: $\mathbb{D}^2(\mathcal{A}_1) \in \mathbf{\Sigma}^0_2$ and conclude, using (iv): $\neg(\mathbf{\Sigma}^0_2\subseteq\mathbf{\Pi}^1_1)$.

\smallskip (vi) Assume: $\mathcal{X}\in\mathbf{\Pi}^1_1$ and $\varphi:\omega^\omega\rightarrow\omega^\omega$ and define: $\mathcal{Y}:=\{\alpha\mid\varphi|\alpha\in\mathcal{X}\}$. \\Find $\mathcal{G}$ in $\mathbf{\Sigma}^0_1$ such that $\mathcal{X}=Un(\mathcal{G})$. \\Then, for every $\alpha$, $\alpha\in \mathcal{Y}\leftrightarrow \varphi|\alpha \in \mathcal{X}\leftrightarrow \forall \beta[\ulcorner \varphi|\alpha,\beta \urcorner \in \mathcal{G}]$. \\Define $\mathcal{V}:=\{\alpha\mid\ulcorner \varphi|\alpha_I,\alpha_{II}\urcorner \in \mathcal{G}\}$. \\Conclude: $\mathcal{V} \in \mathbf{\Sigma}^0_1$ and $\mathcal{Y}=Un(\mathcal{V}) \in \mathbf{\Pi}^1_1$. 
\end{proof}

\subsection{The set  $\mathcal{WF}$ }
\hfill
\begin{definition}\label{D:wf} We define $\mathcal{WF}:= \{\alpha\mid \forall \beta \in (T_\alpha)^\omega\exists n[\beta(n)<_{KB}\beta(n+1)]\}$.\end{definition}

$\mathcal{WF}$ 
is the set of all $\alpha$ such that the tree $T_\alpha:=\{s\mid\forall t\sqsubset s[\alpha(t)=0]\}$ is  \textit{well-founded} with respect to the Kleene-Brouwer-ordering $<_{KB}$, see Definition \ref{D:kleenebrouwer} in Subsection \ref{SS:if}.

The following Theorem is a counterpart to Theorem \ref{T:pif}. Note that Theorem \ref{T:pif} is the statement that $\mathcal{E}^1_1$ does \emph{not} coincide with $\mathcal{IF}$. Note that both $(\mathcal{E}^1_1, \mathcal{A}^1_1)$ and $(\mathcal{IF}, \mathcal{WF})$ are  complementary $(\mathbf{\Sigma}^1_1, \mathbf{\Pi}^1_1)$-pairs, see Subsubsection \ref{SSS:projective}.

\begin{theorem} \label{T:wf}\hfill

$\mathcal{WF}=\mathcal{A}^1_1$. 
\end{theorem}

\begin{proof}  We first prove that $\mathcal{WF}$ is a subset of $\mathcal{A}^1_1$.

 Assume: $\alpha \in \mathcal{WF}$. Let $\gamma$ be given. \\Define $\beta$ such that $\beta(0)=\langle\;\rangle$ and, for each $n$, \\if $\overline \gamma (n+1) \in T_\alpha$, then $\beta(n+1)=\overline \gamma(n+1)$, and, if not, then $\beta(n+1)=\beta(n)$. \\ Note $\forall n[\beta(n)\in T_\alpha]$ and find $n$ such that $\beta(n)\le_{KB}\beta(n+1)$. \\Conclude: $\beta(n+1)\neq \overline\gamma(n+1)$ and: $\exists i\le n[\alpha(\overline \gamma i)\neq 0]$. \\We thus see: $\forall \gamma  \exists i[\alpha(\overline \gamma i)\neq 0]$, that is: $\alpha \in \mathcal{A}^1_1$. 

Conclude: $\mathcal{WF}\subseteq \mathcal{A}_1^1$.

\smallskip We now prove that $\mathcal{A}_1^1$ is a subset of $\mathcal{WF}$. This proof is more difficult and we have to use the principle of Bar Induction $\mathbf{BI}$, see Subsubsection \ref{SSS:barinduction}.

Assume: $\alpha \in \mathcal{A}^1_1$. \\Define $B:=\omega\setminus T_\alpha=\{s\mid\exists t\sqsubset s[\alpha(t) \neq 0]\}$ and note: $B$  is a bar in $\omega^\omega$. 
\\Define $C:=\{s\mid\forall\beta\in (T_\alpha)^\omega[\forall i[s\sqsubseteq\beta(i)]\rightarrow\exists j[\beta(j)\le_{KB}\beta(j+1)]]\}$ and note: $B\subseteq C$, as, for each $s$ in $B$, for each $u$ such that $s\sqsubseteq u$, $u\notin T_\alpha$.\\Also note: $C$ is monotone, that is: $\forall s\forall m[s\in C\rightarrow s\ast\langle m \rangle \in C]$.  

We now will prove that $C$ is inductive.\\Let $s$ be given such that $\forall m[s\ast\langle m\rangle \in C]$. We want to prove: $s\in C$.

Define, for each $m$, \\$P(m):= \forall \beta\in (T_\alpha)^\omega[\bigl(\forall i[s\sqsubseteq \beta(i)]\;\wedge\;s\ast\langle m \rangle \sqsubseteq \beta(0)\bigr) \rightarrow \exists j[\beta(j)\le_{KB}\beta(j+1)]]$. \\Before proving: `$s\in C$', we first prove the {\it auxiliary statement}: $\forall m[P(m)]$. \\ We use induction. \\Let $m$ be given such that $\forall k<m[P(k)]$. \\Let $\beta$ in $(T_\alpha)^\omega$  be given such that $\forall i[s\sqsubseteq \beta(i)]$ and $s\ast \langle m \rangle \sqsubseteq \beta(0)$. \\We intend to prove: $\exists j[\beta(j)\le_{KB}\beta(j+1)]$. \\Define $\beta^\ast$ such that $\beta^\ast(0)=\beta(0)$ and, for each $n$, \textit{if} $\forall i\le n+1[s\ast\langle m \rangle \sqsubseteq \beta(i)]$, then $\beta^\ast(n+1)=\beta(n+1)$ and, \textit{if not}, then $\beta^\ast(n+1)=\beta^\ast(n)$. \\Note: $\forall n[s\ast\langle m\rangle\sqsubseteq\beta^\ast(n)]$
and: $s\ast\langle m\rangle \in C$,  and find $j$ such that $\beta^\ast(j)\le_{KB}\beta^\ast(j+1)$. \\\textit{If} $\beta^\ast(j)=\beta(j)$ and $\beta^\ast(j+1)=\beta(j+1)$, conclude: $\beta(j)\le_{KB}\beta(j+1)$: we are done. \\\textit{If not}, define $i_0:=\mu i[\neg\bigl(s\ast\langle m \rangle \sqsubseteq \beta(i)\bigr)]$ and distinguish two cases.

\textit{Case (a)}. $\beta(i_0)=s$. Note: $i_0>0$ and $\beta(i_0-1)\le_{KB}\beta(i_0)$: we are done.

\textit{Case (b)}. $s\sqsubset \beta(i_0)$. 
Find $k$ such that $s\ast\langle k\rangle \sqsubseteq \beta(i_0)$. \\Note: $k\neq m$ and distinguish two cases.

\textit{Case (b1)}. $m<k$. Note: $i_0>0$ and $\beta(0)<_{KB}\beta(i_0)$ and $\exists j<i_0[\beta(j)<_{KB}\beta(j+1)]$: we are done. 

\textit{Case (b2)}. $k<m$. Define $\beta^\dag$ such that $\forall n[\beta^\dag(n)=\beta(i_0+n)]$. \\Note: $s\ast\langle k\rangle \sqsubseteq \beta^\dag(0)$ and apply $P(k)$. Find $l$ such that
$\beta^\dag(l)\le_{KB}\beta^\dag(l+1)$ and, therefore: $\beta(i_0+l)\le_{KB}\beta(i_0+l+1)$: again, we are done. 
\\We conclude: $P(m)$. \\This completes the proof of the auxiliary statement: $\forall m[P(m)]$.

We now are ready to prove: $s\in C$. Let $\beta$ in $(T_\alpha)^\omega$ be given such that $\forall i[s\sqsubseteq \beta(i)]$.  Consider $\beta(0)$ and $\beta(1)$.  \textit{Either}: we find $m$ such that $s\ast\langle  m\rangle \sqsubset \beta(0)]$ or $s\ast \langle m \rangle\sqsubseteq \beta(1)$, and, considering $\beta$ or $\beta\circ S$ and using $P(m)$, we conclude: $\exists j[\beta(j)\le_{KB}\beta(j+1)]$ \textit{or}: $\beta(0) =\beta(1)=s$ and $\beta(0)\le_{KB}\beta(1)$. \\Conclude: $\forall \beta \in (T_\alpha)^\omega[\forall i[s\sqsubseteq \beta(i)]\rightarrow \exists j[ \beta(j)\le_{KB}\beta(j+1)]]$, i.e. $s\in C$.

Using $\mathbf{BI}$, we conclude: $\langle\;\rangle\in C$,
 i.e. $\forall \beta \in (T_\alpha)^\omega\exists j[\beta(j)\le_{KB}\beta(j+1)]$, i.e. $\alpha\in \mathcal{WF}$.

We thus see: $\mathcal{A}^1_1\subseteq\mathcal{WF}$ and: $\mathcal{A}^1_1=\mathcal{WF}$.
\end{proof}

 \bigskip
The statement $\mathcal{A}^1_1=\mathcal{WF}$ is, in the formal context of Basic Intuitionistic Mathematics $\mathsf{BIM}$, an equivalent of $\mathbf{OI}(2^\omega)$, the Principle of Open induction on Cantor space $2^\omega$, see  \cite{veldman15}. 
\subsection{ $\mathsf{Sink}^\ast(\mathcal{FIN})$ and $\mathsf{Sink}^\ast(\mathcal{ALMOST}^\ast\mathcal{FIN})$.}\label{SS:finite}
\hfill

\begin{definition}We define: $\mathcal{FIN}:=\{\alpha\mid\exists m\forall n>m[\alpha(n)=0]\}$.\end{definition}

$\mathcal{FIN}$ is the set of all $\alpha$ such that $D_\alpha:=\{n\mid\alpha(n)\neq 0\}$ is a \textit{finite} subset of $\omega$.

We want to remind the reader of a fact proven in \cite[Theorem 3.3.(iii)]{veldman05}. 
\begin{theorem}\label{T:finite}\begin{enumerate}[\upshape (i)]\item $\mathbb{D}^2(\mathcal{A}_1)\npreceq \mathcal{FIN}$. \item $\mathcal{FIN}$ is $\mathbf{\Sigma}^0_2$ but not $\mathbf{\Sigma}^0_2$-complete. \item $\mathcal{FIN}$ is not $\mathbf{\Pi}^1_1$. \end{enumerate}\end{theorem}

\begin{proof} (i) Assume: $\varphi:\omega^\omega\rightarrow\omega^\omega$ reduces $\mathbb{D}^2(\mathcal{A}_1)=\{\alpha\mid\alpha^0 = \underline 0\;\vee\;\alpha^1=\underline 0\}$ to $\mathcal{FIN}$. 

We  prove that $\varphi$ maps the closure $\overline{\mathbb{D}^2(\mathcal{A}_1)}$ of $\mathbb{D}^2(\mathcal{A}_1)$ into $\mathcal{FIN}$ and thus obtain a contradiction.
\\Let $\alpha$ in $\overline{\mathbb{D}^2(\mathcal{A}_1)}$ be given. \\Define $\alpha_0,\alpha_1$ such that $\forall i<2[(\alpha_i)^i=\underline 0\;\wedge\;\forall j[\neg\exists n[j=\langle i\rangle \ast n]\rightarrow \alpha_i(j)=\alpha(j)]]$. \\Note: $\forall i<2[\alpha_i \in \mathbb{D}^2(\mathcal{A}_1)]$ and: $\neg(\alpha\;\#\;\alpha_0\;\wedge\;\alpha\;\#\;\alpha_1)$. \\Find $m_0,m_1$ such that $\forall i<2\forall n>m_i[(\varphi|\alpha_i)(n)=0]$. Define $m=\max(m_0,m_1)$. Suppose: $n>m$ and $(\varphi|\alpha)(n)\neq 0$. Then: $\alpha\;\#\;\alpha_0$ and $\alpha\;\#\;\alpha_1$, a contradiction. Conclude: $\forall n>m[(\varphi|\alpha)(n)=0]$ and: $\varphi|\alpha \in \mathcal{FIN}$ and, therefore: $\alpha \in \mathbb{D}^2(\mathcal{A}_1)$.  \\We thus see: $\overline{\mathbb{D}^2(\mathcal{A}_1)}\subseteq\mathbb{D}^2(\mathcal{A}_1)$ and, according to Theorem \ref{T:disja1} in Subsubsection \ref{SSS:disjunction}, obtain a contradiction.

Conclude: $\mathbb{D}^2(\mathcal{A}_1)\npreceq\mathcal{FIN}$. 

\smallskip (ii) $\mathcal{FIN}=\bigcup_m\{\alpha\mid\forall n>m[\alpha(n)=0]\}$ clearly is $\mathbf{\Sigma}^0_2$, but, as $\mathbb{D}^2(\mathcal{A}_1)$ is $\mathbf{\Sigma}^0_2$ and, according to (i), $\mathbb{D}^2(\mathcal{A}_1)\npreceq\mathcal{FIN}$, $\mathcal{FIN}$ is not $\mathbf{\Sigma}^0_2$-complete.

\smallskip (iii) Assume: $\varphi:\omega^\omega\rightarrow\omega^\omega$ reduces $\mathbb{D}^2(\mathcal{A}_1)$ to $\mathcal{A}^1_1$.

We  prove that $\varphi$ maps the closure $\overline{\mathbb{D}^2(\mathcal{A}_1)}$ of $\mathbb{D}^2(\mathcal{A}_1)$ into $\mathcal{A}^1_1$ and thus obtain a contradiction. 

The proof is similar to the proof of (i). \\Let $\alpha$ in $\overline{\mathbb{D}^2(\mathcal{A}_1)}$ be given.\\ Define $\alpha_0,\alpha_1$ such that $\forall i<2[(\alpha_i)^i=\underline 0\;\wedge\;\forall j[\neg\exists n[j=\langle i\rangle \ast n]\rightarrow \alpha_i(j)=\alpha(j)]]$. \\Let $\gamma$ be given. Find $n_0, n_1$ such that $\forall i<2[(\varphi|\alpha_i)(\overline \gamma n_i)\neq 0]$.\\ Note: if $\forall i<2[(\varphi|\alpha)(\overline \gamma n_i)\neq (\varphi|\alpha_i)(\overline \gamma n_i)]$, then $\forall i<2[\alpha\;\#\;\alpha_i]$, a contradiction. Conclude: $\exists i<2[(\varphi|\alpha)(\overline \gamma n_i)\neq 0]$. \\We thus see: $\forall \gamma\exists n[(\varphi|\alpha)(\overline \gamma n)\neq 0]$, that is: $\varphi|\alpha \in \mathcal{A}^1_1$, and conclude: $\alpha \in \mathbb{D}^2(\mathcal{A}_1)$. 
\\We thus see: $\overline{\mathbb{D}^2(\mathcal{A}_1)}\subseteq\mathbb{D}^2(\mathcal{A}_1)$ and, according to Theorem \ref{T:disja1}, obtain a contradiction.

Conclude: $\mathbb{D}^2(\mathcal{A}_1)\npreceq \mathcal{A}^1_1$ and: $\mathbb{D}^2(\mathcal{A}_1)\notin \mathbf{\Pi}^1_1$.
\end{proof}

 \begin{definition}\label{D:almostfinite} We define $\mathcal{ALMOST}^\ast\mathcal{FIN}:=\{\alpha\mid\forall \zeta \in [\omega]^\omega\exists n[\alpha\circ\zeta(n)=0]\}$. \end{definition} $\mathcal{ALMOST}^\ast\mathcal{FIN}$ is the set of all $\alpha$ such that $D_\alpha$ is an \textit{almost-finite} subset of $\omega$.

 \begin{lemma}\label{L:almostfincoanalytic} $\mathcal{ALMOST}^\ast\mathcal{FIN}$ is $\mathbf{\Pi}^1_1$. \end{lemma}
 
 \begin{proof} We shall prove that,  for each $\alpha$, 
 
 $\forall \zeta \in[\omega]^\omega\exists n[\alpha\circ\zeta(n)=0]$ if and only if $\forall \zeta \exists n[\alpha\circ\zeta(n)=0]\;\vee\;\zeta(n+1)\le\zeta(n)]$.
 
 The desired conclusion then follows easily.
 
 Let $\tau$ be the canonical retraction\footnote{see  Subsubsection \ref{SSS:continuousfunctions}} of $\omega^\omega$ onto the spread $[\omega]^\omega$. \\The function $\tau$ satisfies the conditions:\\  $\forall \zeta\in[\omega]^\omega[\tau|\zeta=\zeta]$ and $\forall \zeta[\zeta\;\#\;\tau|\zeta\rightarrow \exists n[\zeta(n+1)\le\zeta(n)]]$. 
 
 Let $\alpha$ be given. First assume $\forall \zeta \in[\omega]^\omega\exists n[\alpha\circ\zeta(n)=0]$. \\Let $\zeta$ be given.  Find $n$ such that $\alpha\circ(\tau|\zeta)(n)=0$. \\\textit{Either}: $(\tau|\zeta)(n)=\zeta(n)$ and $\alpha\circ\zeta(n)=0$, \textit{or}: $(\tau|\zeta)(n)\neq\zeta(n)$ and $\exists i \le n[\zeta(i+1)\le \zeta(i)]$. We thus see:  $\forall \zeta \exists n[\alpha\circ\zeta(n)=0\;\vee\;\zeta(n+1)\le\zeta(n)]$. 
 
 Now assume: $\forall \zeta \exists n[\alpha\circ\zeta(n)=0\;\vee\;\zeta(n+1)\le\zeta(n)]$. \\Let $\zeta$ in $[\omega]^\omega$ be given. Find $n$ such that $\alpha\circ\zeta(n)=0\;\vee\;\zeta(n+1)\le\zeta(n)$.  \\Conclude: $\alpha\circ\zeta(n)=0$. \\We thus see: $\forall \zeta \in [\omega]^\omega\exists n[\alpha\circ\zeta(n)=0]$. 
 \end{proof}
 
 The set $\mathcal{ALMOST}^\ast \mathcal{FIN}$ has been studied in \cite[Section 3]{veldman05}. It has been shown there that $\mathcal{ALMOST}^\ast \mathcal{FIN}$ is  not (positively) Borel, see \cite[Section 0.9.2(ii) and Theorem 3.17(iii)]{veldman05}. In particular, $\mathcal{FIN}$ is proper subset\footnote{For all $\mathcal{X}, \mathcal{Y}\subseteq \omega^\omega$, $\mathcal{X}$ is a \emph{proper subset} of $\mathcal{Y}$ if and only if $\mathcal{X}\subseteq \mathcal {Y}$ and not: $\mathcal{Y}\subseteq \mathcal{X}$, i.e. the assumption `$\mathcal{Y}\subseteq \mathcal{X}$' leads to a contradiction.}of $\mathcal{ALMOST}^\ast \mathcal{FIN}$. \\It has also been shown in \cite{veldman05} that $\mathcal{ALMOST}^\ast \mathcal{FIN}$ is the best $\mathbf{\Pi}^1_1$-approximation of $\mathcal{FIN}$, i.e.,  for every $\mathcal{Z}$ in $\mathbf{\Pi}^1_1$,  if $\mathcal{FIN} \subseteq \mathcal{Z}$, then $\mathcal{ALMOST}^\ast \mathcal{FIN} \subseteq \mathcal{Z}$, see  \cite[Theorem 3.21(v)]{veldman05}.  
  \\As one might expect, $\mathcal{ALMOST}^\ast \mathcal{FIN}$ is not $\mathbf{\Pi}^1_1$-complete, see \cite[Theorem 3.24(iii)]{veldman05}.
 
\smallskip In the following definition we introduce a new word for a well-known concept.
 \smallskip
 \begin{definition}\label{D:sink}For all $\mathcal{X},\mathcal{Y}\subseteq \omega^\omega$, we define: $\mathcal{X}$ \emph{sinks into} $\mathcal{Y}$ if and only if $\mathcal{X}\subseteq \mathcal{Y}$. 
 
 For each $\mathcal{X}\subseteq \omega^\omega$,  we define \\$\mathsf{Sink}(\mathcal{X}):=\{\beta\mid\mathcal{F}_\beta \subseteq \mathcal{X}\}$  and $\mathsf{Sink}^\ast(\mathcal{X}):=\{\beta\in\mathsf{Sink}(\mathcal{X})\mid Spr(\beta)\}$  \end{definition}
 $\mathsf{Sink}(\mathcal{X})$ is the set of the codes of all closed subsets of $\omega^\omega$ that sink into (i.e. are a subset of) $\mathcal{X}$ and $\mathsf{Sink}^\ast(\mathcal{X})$ is the set of the codes of all spreads, i.e. all \emph{closed and located} subsets of $\omega^\omega$, that sink into (i.e. are a subset of) $\mathcal{X}$.

 We now want to treat some results that, together, are a counterpart\footnote{Note that, from a classical point of view, the sets  $\mathsf{Share}(\mathcal{INF})$, $\mathsf{Sink}(\mathcal{FIN})$, for instance, are each other's complement.}  to Theorem \ref{T:shareinf}. The moral of the story is that, in order to obtain a satisfying such counterpart, one should work with $\mathcal{ALMOST}^\ast\mathcal{FIN}$ rather than with $\mathcal{FIN}$.

 \smallskip Recall: for all $\mathcal{X},\mathcal{Y}\subseteq\omega^\omega$: $\mathcal{X}\sim\mathcal{Y}$ (\textit{$\mathcal{X}$, $\mathcal{Y}$ reduce to each other / are Wadge-equivalent}), if and only if both $\mathcal{X}\preceq\mathcal{Y}$ and $\mathcal{Y}\preceq\mathcal{X}$.
 \begin{theorem}\label{T:sinkalmostfin}
  \hfill

\begin{enumerate}[\upshape(i)]

\item $\mathsf{Sink}^\ast(\mathcal{FIN}\cap2^\omega)\sim\mathcal{FIN}$. 
\item $\mathsf{Sink}^\ast(\mathcal{FIN}\cap2^\omega)\npreceq\mathcal{A}^1_1$ and $\mathsf{Sink}^\ast(\mathcal{FIN})\npreceq\mathcal{A}^1_1$.
\item $\mathcal{A}^1_1 \preceq \mathsf{Sink}^\ast(\mathcal{FIN})$. 
\item $\mathcal{A}^1_1 \preceq \mathsf{Sink}^\ast(\mathcal{ALMOST}^\ast\mathcal{FIN}\cap2^\omega)\preceq \mathsf{Sink}^\ast(\mathcal{ALMOST}^\ast\mathcal{FIN})$.
\item $\mathsf{Sink}^\ast(\mathcal{ALMOST}^\ast\mathcal{FIN})\preceq\mathcal{A}^1_1$. \item $\mathsf{Sink}^\ast(\mathcal{ALMOST}^\ast\mathcal{FIN})$ and $\mathsf{Sink}^\ast(\mathcal{ALMOST}^\ast\mathcal{FIN}\cap2^\omega)$ are $\mathbf{\Pi}^1_1$-complete.

\end{enumerate}
\end{theorem}

\begin{proof}
(i) Assume: $Spr(\beta)$ and $\mathcal{F}_\beta\subseteq\mathcal{FIN}\cap2^\omega$. \\Conclude: $\forall s[\beta(s) =0\rightarrow s \in Bin]$ and: $Fan(\beta)$ and: $\forall \gamma \in \mathcal{F}_\beta\exists m\forall n>m[\gamma(n)=0]$. \\Applying the First Axiom of Continuous Choice $\mathbf{AC}_{1,0}$, see Subsubsection \ref{SSS:countablechoice}, \\find $\varphi:\mathcal{F}_\beta \rightarrow \omega$ such that $\forall \gamma \in \mathcal{F}_\beta\forall n>\varphi(\gamma)[\gamma(n)=0]$. \\Applying the Fan Theorem $\mathbf{FT}$, see Subsubsection \ref{SSS:fantheorem}, \\find $p$ such that $\forall \gamma \in \mathcal{F}_\beta[\varphi(\gamma)\le p]$. \\ Note: $
\forall n>p\forall s \in Bin_n[\beta(s) =0\rightarrow s(n) =0]$. 
\\Conclude: for each $\beta$,  $\beta \in \mathsf{Sink}^\ast(\mathcal{FIN}\cap2^\omega)$ if and only if \\$Spr(\beta)$ and $\forall s[\beta(s) =0\rightarrow s \in Bin]$ and $\exists p\forall n>p\forall s\in Bin_n[\beta(s) =0\rightarrow s(n) =0]]$. 

\smallskip Define $\psi:\omega^\omega\rightarrow\omega^\omega$ such that, for all $\beta$, for all $n$, $(\psi|\beta)(n) = 0$ if and only if \\$\forall s\le n[\beta(s) =0\rightarrow s\in Bin]$ and $ \forall s\in Bin_{n+1}[\beta(s) =0\rightarrow s(n) =0]$. \\ Note: $\psi$ reduces $\mathsf{Sink}^\ast(\mathcal{FIN}\cap2^\omega)$ to $\mathcal{FIN}$.

\smallskip Define $\rho:\omega^\omega\rightarrow\omega^\omega$ such that, for all $\alpha$, for all $s$,  $(\rho|\alpha)(s) =0$ if and only if  \\$s\in Bin$ and $\forall i<length(s)[s(i)=1\leftrightarrow \alpha(i)\neq 0]$. \\ Note: $\rho$ reduces $\mathcal{FIN}$ to $\mathsf{Sink}^\ast(\mathcal{FIN}\cap2^\omega)$. 

\smallskip Conclude: $\mathsf{Sink}^\ast(\mathcal{FIN}\cap2^\omega)\sim\mathcal{FIN}$.

 \smallskip
 (ii) Use (i) and Theorem \ref{T:finite}(ii) and conclude: $\mathsf{Sink}^\ast(\mathcal{FIN}\cap2^\omega)\npreceq\mathcal{A}_1^1$.  
 \\Define $\varphi:\omega^\omega\rightarrow\omega^\omega$ such that, for all $\beta$, for all $s$, $(\varphi|\beta)(s)=0 $ if and only if $\\(s\in Bin \;\wedge\;\beta(s)=0)\;\vee\; \exists t\sqsubseteq s[t \notin Bin \;\wedge\;\beta(t)=0]$.
 \\Note: $\varphi$ reduces $\mathsf{Sink}^\ast(\mathcal{FIN}\cap2^\omega)$ to $\mathsf{Sink}^\ast(\mathcal{FIN})$.
 
 Conclude: $\mathsf{Sink}^\ast(\mathcal{FIN})\npreceq\mathcal{A}_1^1$.  
 
 \smallskip (iii) Define $\varphi:\omega^\omega\rightarrow\omega^\omega$ such that,  for all $\alpha$, for all $s$, \\$(\varphi|\alpha)(s) =0 $ if and only if $ \exists t\in T_\alpha\exists n[s=(S\circ t)\ast\underline{\overline 0}n]$.\footnote{Recall: $\mathit{length}(S\circ t)=\mathit{length}(t)$ and $\forall i<\mathit{length}(t)[(S\circ t)(i)=t(i)+1]$.} \\Note: for all $\alpha$, $Spr(\varphi|\alpha)$ and $\forall \gamma \in \mathcal{F}_{\varphi|\alpha}\forall n[\gamma(n)=0\rightarrow \gamma(n+1)=0]$.  \\We now prove that $\varphi$ reduces $\mathcal{A}_1^1$ to $\mathsf{Sink}^\ast(\mathcal{FIN})$.
 
 \smallskip First assume: $\alpha\in\mathcal{A}_1^1$. Also assume: $\gamma\in \mathcal{F}_{\varphi|\alpha}$. Find $\varepsilon$ such that, for each $n$, \\ if $\gamma(n)>0$, then $\varepsilon(n)+1=\gamma(n)$, and, if $\gamma(n)=0$, then $\varepsilon(n)=0$. \\Find $m$ such that $\alpha(\overline \varepsilon m)\neq 0$. Then: $\overline \varepsilon (m+1) \notin T_\alpha$ and: $\overline \gamma (m+1) \neq S\circ\overline \varepsilon (m+1)$. \\Find $k\le m$ such that $\gamma(k) =0$ and note: $\forall n>k[\gamma(n)=0]$ and $\gamma \in \mathcal{FIN}$. \\We thus see: $\forall \gamma \in \mathcal{F}_{\varphi|\alpha}[\gamma \in \mathcal{FIN}]$. \\ Conclude: $\mathcal{F}_{\varphi|\alpha}\subseteq \mathcal{FIN}$ and: $\varphi|\alpha \in \mathsf{Sink}^\ast(\mathcal{FIN})$. 
 
 \smallskip Now assume: $\varphi|\alpha \in \mathsf{Sink}^\ast(\mathcal{FIN})$. Then $\forall \gamma \in \mathcal{F}_{\varphi|\alpha}\exists m\forall n>m[\gamma(n)=0]$. \\Let $\varepsilon$ be given. Define $\gamma$ such that, for each $n$, \\\textit{if} $\overline \varepsilon(n+1)\in T_\alpha$, then $\gamma(n)=\varepsilon(n)+1$, and, \textit{if not}, then $\gamma(n)=0$. \\Note: $\gamma \in \mathcal{F}_{\varphi|\alpha}$ and find $m$ such that $\gamma(m)=0$. \\Conclude: $\overline \varepsilon(m+1) \notin T_\alpha$ and $\exists i\le m+1[\alpha(\overline \varepsilon i)\neq 0]$.  \\We thus see: $\forall \varepsilon\exists i[\alpha(\overline \varepsilon i)\neq 0]$, i.e. $\alpha\in\mathcal{A}^1_1$. 
  
  We thus see: $\forall \alpha[\alpha \in \mathcal{A}^1_1 \leftrightarrow \varphi|\alpha \in \mathsf{Sink}^\ast(\mathcal{FIN})]$, i.e. $\varphi$ reduces $\mathcal{A}^1_1$ to $ \mathsf{Sink}^\ast(\mathcal{FIN})$.

\smallskip
(iv) Define $\delta$ such that $\delta(0)=0$ and, for all $s$, for all $n$,  \\$\delta(s\ast\langle n\rangle)=\delta(s)\ast\underline{\overline 0}n\ast\langle 1\rangle$. Define $\varphi:\omega^\omega\rightarrow\omega^\omega$ such that for all $\alpha$, for all $s$,\\$(\varphi|\alpha)(s)=0$ if and only if $\exists t \in T_\alpha \exists n[ s = \delta(t)\ast\overline{\underline 0}n]$. \\Note: for all $\alpha$, $Spr(\varphi|\alpha)$ and $\mathcal{F}_{\varphi|\alpha}\subseteq 2^\omega$. \\We  prove that $\varphi$ reduces $\mathcal{A}^1_1$ to $\mathsf{Sink}^\ast(\mathcal{ALMOST}^\ast\mathcal{FIN}\cap 2^\omega)$.

\smallskip Assume: $\alpha\in \mathcal{A}_1^1$. Also assume: $\gamma \in \mathcal{F}_{\varphi|\alpha}$ and  $\zeta \in [\omega]^\omega$. \\Define $\gamma'$ such that $\forall n[\gamma'\circ\zeta(n) =1]$ and $\forall n[\forall i[n\neq\zeta(i)]\rightarrow \gamma'(n)=\gamma(n)]$. \\Define $\varepsilon$ such that $\varepsilon(0)=\mu p[\gamma'(p)=1]$ and $\forall n[\varepsilon(n+1) =\mu p>\varepsilon(n)[\gamma'(p)=1]]$. \\Note: $\varepsilon \in [\omega]^\omega$ and, for all $n$, $\delta(\overline\varepsilon n)\sqsubset\gamma'$. Find $n$ such that $\alpha(\overline \varepsilon n)\neq 0$ .\\ Note: $\overline\varepsilon(n+1)\notin T_\alpha$  and $(\varphi|\alpha)\bigl(\delta(\overline \varepsilon (n+1))\bigr)\neq 0$. \\Find 
$m$ such that $\overline{\gamma'}m=\delta(\overline\varepsilon (n+1))$. \\Note:  $(\varphi|\alpha)(\overline{\gamma'}m)\neq 0=(\varphi|\alpha)(\overline \gamma m)$ and conclude: $\overline{\gamma'}m\neq \overline \gamma m$. \\Find $i<m$ such that $\gamma'(i)\neq \gamma(i)$. \\Determine $j<m$ such that $i=\zeta(j)$ and conclude: $\gamma\circ\zeta(j) =0$. \\We thus see: $\forall \gamma\in\mathcal{F}_{\varphi|\alpha}\forall\zeta\in [\omega]^\omega\exists j[\gamma\circ\zeta(j)=0]$. \\Conclude: $\mathcal{F}_{\varphi|\alpha}\subseteq \mathcal{ALMOST}^\ast\mathcal{FIN}$ and: $ \varphi|\alpha \in \mathsf{Sink}^\ast(\mathcal{ALMOST}^\ast\mathcal{FIN}\cap2^\omega)$.

\smallskip Now assume: $\varphi|\alpha \in\mathsf{Sink}^\ast(\mathcal{ALMOST}^\ast\mathcal{FIN}\cap2^\omega)$. Let $\gamma$ be given. \\Find $\beta$ in $2^\omega$ such that $\forall n[\delta(\overline \gamma n)\sqsubset \beta]$. \\Define $\zeta$ such that $\zeta(0)=\gamma(0)$ and $\forall n[\zeta(n+1)=\zeta(n) +\gamma(n+1)+1]$. \\Note $\zeta \in [\omega]^\omega$ and $\forall n[\beta\circ\zeta(n)=1]$. Define $\beta^\ast$ such that, for each $n$, \\\textit{if} $\overline\beta(n+1)\in T_{\varphi|\alpha}$, then $\beta^\ast(n)=\beta(n)$, and \textit{if not}, then $\beta^\ast(n) =0$. \\Note: $\beta^\ast\in\mathcal{F}_{\varphi|\alpha}\subseteq \mathcal{ALMOST}^\ast\mathcal{FIN}$ and find $n$ such that $\beta^\ast\circ\zeta(n)=0$. \\Define $p:=\zeta(n)+1$ and conclude: $\overline \beta p\neq\overline{\beta^\ast}p$ and: $\overline \beta p \notin T_{\varphi|\alpha}$. \\Find $m$ such that $\overline \beta p\sqsubseteq \delta(\overline \gamma m)$  and note: $\overline \gamma m \notin T_\alpha$ and: $\exists i\le m[\alpha(\overline \gamma i)\neq 0]$. \\We thus see: $\forall \gamma\exists i[\alpha(\overline \gamma i)\neq 0]$ , i.e. $\alpha \in \mathcal{A}_1^1$.

\smallskip Conclude: for each $\alpha$, $\alpha \in \mathcal{A}^1_1$ if and only if  $\varphi|\alpha\in \mathsf{Sink}^\ast(\mathcal{ALMOST}^\ast\mathcal{FIN}\cap2^\omega)$, i.e. $\varphi$ reduces $\mathcal{A}^1_1$ to  $\mathsf{Sink}^\ast(\mathcal{ALMOST}^\ast\mathcal{FIN}\cap2^\omega)$.

\smallskip Finally, define $\psi:\omega^\omega\rightarrow\omega^\omega$ such that, for all $\beta$, for all $s$, $(\psi|\beta)(s) =0$ \\if and only if either $\beta(s) =0\;\wedge\;s\in Bin$ or $\exists t\sqsubseteq s[\beta(t)=0\;\wedge s\notin Bin]$. \\Note: $\psi$ reduces $\mathsf{Sink}^\ast(\mathcal{ALMOST}^\ast\mathcal{FIN}\cap 2^\omega)$ to $\mathsf{Sink}^\ast(\mathcal{ALMOST}^\ast\mathcal{FIN})$.

\smallskip
 (v) We first prove a \textit{preliminary observation}. For all $\beta$ such that $Spr(\beta)$, \\ $\forall \alpha\in\mathcal{F}_\beta\forall\zeta\in[\omega]^\omega\exists n[\alpha\circ\zeta(n)=0]$ if and only if \\$\forall \alpha\forall\zeta\exists n[\alpha\circ\zeta(n)=0\;\vee\;\zeta(n+1)\le\zeta(n)\;\vee\; \beta(\overline\alpha n)\neq 0]$. \\The argument is a small extension of the argument given for Lemma \ref{L:almostfincoanalytic}.

 \smallskip Let $\beta$ be given such that $Spr(\beta)$.

First assume $\forall \alpha\in\mathcal{F}_\beta\forall\zeta\in[\omega]^\omega\exists n[\alpha\circ\zeta(n)=0]$. Let $\rho,\tau$ be the canonical retractions\footnote{See Subsubsection \ref{SSS:continuousfunctions}} of $\omega^\omega$ onto the spreads $\mathcal{F}_\beta$ and $[\omega]^\omega$, respectively. Let $\alpha, \zeta$ be given. Find $n$ such that $\bigl((\rho|\alpha)\circ(\tau|\zeta)\bigr)(n)=0$. There are three cases to consider. \\\textit{Case (a)}. $\overline{\tau|\zeta}(n+1) \neq \overline\zeta(n+1)$. Then  $\exists i[\zeta(i+1)\le\zeta(i)]$. \\\textit{Case (b)}. $\overline{\tau|\zeta}(n+1) = \overline \zeta(n+1)$ and  $\overline{\rho|\alpha}(\zeta(n)+1)\neq \overline \alpha(\zeta(n)+1)$. Then $\exists i[\beta(\overline \alpha i)\neq 0]$.  \\\textit{Case (c)}. $\overline{\tau|\zeta}(n+1) = \overline \zeta(n+1)$ and  $\overline{\rho|\alpha}(\zeta(n)+1)= \overline \alpha(\zeta(n)+1)$. Then $\alpha\circ\zeta(n)=0$. \\Conclude: $\forall \alpha\forall\zeta\exists n[\alpha\circ\zeta(n)=0\;\vee\;\zeta(n+1)\le\zeta(n)\;\vee\; \beta(\overline\alpha n)\neq 0]$.
 
 Now assume $\forall \alpha\forall\zeta\exists n[\alpha\circ\zeta(n)=0\;\vee\;\zeta(n+1)\le\zeta(n)\;\vee\; \beta(\overline\alpha n)\neq 0]$. \\Let $\alpha$ be given in $\mathcal{F}_\beta$ and $\zeta$ in $[\omega]^\omega$. \\Find $n$ such that  $\alpha\circ\zeta(n)=0\;\vee\;\zeta(n+1)\le\zeta(n)\;\vee\; \beta(\overline\alpha n)\neq 0$ and conclude: $\alpha\circ\zeta(n)=0$. We thus see: $\forall \alpha \in \mathcal{F}_\beta\forall \zeta\in [\omega]^\omega\exists n[\alpha\circ\zeta(n)=0]$. 
 
 \smallskip Now observe: $\{\beta\mid Spr(\beta\}$ belongs to $ \mathbf{\Pi}^0_2$. Using our preliminary observation and  also  Theorem \ref{T:coan},  conclude: 
 $\mathsf{Sink}^\ast(\mathcal{ALMOST}^\ast\mathcal{FIN})=\\\{\beta\mid Spr(\beta)\;\wedge\;\forall\alpha\forall\zeta\exists n[\alpha\circ\zeta(n)=0\;\vee\;\zeta(n+1)\le\zeta(n)\;\vee\; \beta(\overline\alpha n)\neq 0]\}$ belongs to  $\mathbf{\Pi}^1_1$.
 
 \smallskip (vi) Use (iv) and (v).
 \end{proof}
   
Theorem \ref{T:sinkalmostfin}(i) seems to contradict  classical results:  its proof uses the strongly nonclassical axiom $\mathbf{AC}_{1,0}$.  Theorem \ref{T:sinkalmostfin}(iv) is a counterpart to Theorem \ref{T:shareinf}. Both Theorem \ref{T:sinkalmostfin}(vi) and Theorem \ref{T:shareinf} resemble a classical result due to Hurewicz that plays a key role in the sketch of the proof of a theorem by Solovay and Kaufman in \cite{kechris87}. The Solovay-Kaufman Theorem states that the set of the codes of closed sets of uniqueness and the  set of the codes of closed sets of extended uniqueness are $\mathbf{\Pi}^1_1$-complete. Note that we obtained the more \textit{`classical'} results of Theorem \ref{T:sinkalmostfin} by replacing $\mathcal{FIN}$ by $\mathcal{ALMOST}^\ast\mathcal{FIN}$.

   \subsection{Exactly one path}
\hfill
   
   \begin{definition}\label{D:e11!} We define $\mathcal{E}^1_1!:=\{\alpha\mid \exists \gamma[\forall n[\alpha(\overline\gamma n)=0]\;\wedge\forall \delta[\delta\;\#\;\gamma\rightarrow\exists n[\alpha(\overline \delta n)\neq 0]]]\}$\end{definition} $\mathcal{E}^1_1!$ is the set of all $\alpha$ admitting \textit{exactly one path}. In \cite[pp. 125-127]{kechris}, there is a fascinating argument, due to A.S.~Kechris, showing that, in classical descriptive set theory, $\mathcal{E}_1^1!$ is $\mathbf{\Pi}^1_1$-complete. We will see that this result does not go through in our intuitionistic context.
   
     \begin{definition}\label{D:d2!a1}  We define 
     
     $\mathbb{D}^2!(\mathcal{A}_1):=\{\alpha\mid\exists i<2[\alpha^i=\underline 0\;\wedge\;\alpha^{1-i}\;\#\;\underline 0]\}$, and 
     
     $\mathcal{E}_2!:=\{\alpha\mid\exists n[\alpha^n=\underline 0 \;\wedge\;\forall m\neq n[\alpha^n \;\#\;\underline 0]]\}$. \end{definition} Note: $\mathbb{D}^2!(\mathcal{A}_1)$ is $\mathbf{\Sigma}^0_2$ and $\mathcal{E}_2!$ is $\mathbf{\Sigma}^0_3$.

We will see that the set $\mathcal{E}_2!$ is an example of a subset of $\omega^\omega$ that is positively Borel and  regular in Lusin's sense\footnote{See Definition \ref{D:lusinregular}.}, see Theorem \ref{T:exactlyonepath}, but not $\mathbf{\Pi}^1_1$, see Theorem \ref{T:a11perhapsive}.
   
   \begin{theorem}\label{T:exactlyonepath}
   \hfill
 
   \begin{enumerate}[\upshape(i)]
  \item $\mathbb{D}^2!(\mathcal{A}_1) \preceq \mathcal{E}_2!$ and $\mathcal{E}_2!\preceq\mathcal{E}_1^1!$. 
  
  \item $\mathcal{A}_2\preceq \mathcal{E}_2!$ and $\mathcal{A}_1^1 \preceq \mathcal{E}_1^1!$.
  
  \item $\mathbb{D}^2!(\mathcal{A}_1)\preceq\mathcal{A}_2$, \item $\mathbb{D}^2(\mathcal{A}_1)\npreceq \mathcal{E}_2!$. \item $\mathcal{E}_2!$ is regular in Lusin's sense.  
   \end{enumerate}
   \end{theorem}

   \begin{proof} (i) Define $\varphi:\omega^\omega\rightarrow\omega^\omega$ such that $\forall \alpha[[\forall i<2[(\varphi|\alpha)^i=\alpha^i]\;\wedge\;\forall i\ge 2[(\varphi|\alpha)^i=\underline 1]$ and note: $\varphi$ reduces $\mathbb{D}^2!(\mathcal{A}_1)$ to $\mathcal{E}_2!$. 
   
  Define $\psi:\omega^\omega\rightarrow\omega^\omega$ such that $\forall \alpha\forall s[(\psi|\alpha)(s) =0 \leftrightarrow\exists n[s\sqsubset\underline n\;\wedge\;\overline{\alpha^n}s \sqsubset \underline 0]]$ and note: $\psi$ reduces $\mathcal{E}_2!$ to $\mathcal{E}^1_1!$. 
  
\smallskip
(ii) Define $\varphi:\omega^\omega\rightarrow\omega^\omega$ such that $\forall \alpha[(\varphi|\alpha)^0=\underline 0\;\wedge\;\forall i[(\varphi|\alpha)^{i+1} =\alpha^i]$ and note: \\$\varphi$ reduces $\mathcal{A}_2$ to $\mathcal{E}_2!$. 

Define $\psi:\omega^\omega\rightarrow\omega^\omega$ such that, for all $\alpha$,  $\forall n[(\psi|\alpha)(\overline{\underline 0}n)=0]$ and\\ $\forall m\forall n\forall t[(\psi|\alpha)(\overline{\underline 0}n\ast\langle m+1\rangle\ast t)=\alpha(t)]]$ and note: $\psi$ reduces $\mathcal{A}^1_1$ to $\mathcal{E}^1_1!$.

   \smallskip
   (iii) Define $\varphi:\omega^\omega\rightarrow\omega^\omega$ such that, for all $\alpha$, for all $n$,  $(\varphi|\alpha)^0(n)=\max \bigl(\alpha^0(n), \alpha^1(n)\bigr)$, and, for all $i$,  $(\varphi|\alpha)^{i+1}(n) \neq 0$ if and only if either $\overline{\alpha^0}i\sqsubset\underline 0$ or $\overline{\alpha^1}i\sqsubset\underline 0$. \\ Note: $\varphi$ reduces $\mathbb{D}^2!(\mathcal{A}_1)$ to $\mathcal{A}_2$. 
   
   \smallskip (iv) Assume: $\psi:\omega^\omega\rightarrow\omega^\omega$ maps $\mathbb{D}^2(\mathcal{A}_1)$ into $\mathcal{E}_1^1!$. \\We shall prove that $\psi$ also maps the closure $\overline{\mathbb{D}^2(\mathcal{A}_1)}$ of $\mathbb{D}^2(\mathcal{A}_1)$ into $\mathcal{E}_1^1!$ and thus does not reduce $\mathbb{D}^2(\mathcal{A}_1)$ to $\mathcal{E}_1^1!$.
  
  \smallskip 
   First, as in the proof of Theorem \ref{T:d2a1notregularlusin}, define, for  both $i<2$, $\mathcal{P}_i:= \{\beta\mid \beta^i=\underline 0\}$.
  \\Note: $\mathcal{P}_0, \mathcal{P}_1$ are spreads and $\mathbb{D}^2(\mathcal{A}_1) = \mathcal{P}_0\cup\mathcal{P}_1$.
  
   Assume: $\alpha\in\overline{\mathbb{D}^2(\mathcal{A}_1)}$. We are going to prove: $\psi|\alpha \in \mathcal{E}^1_1!$. \\The following notion is useful. We define, for all $s$,
    $s$ is \textit{fine for $
    \alpha$} if and only if    \\$\exists m\forall \beta \in \mathbb{D}^2(\mathcal{A}_1)[\overline \alpha m \sqsubset \beta\rightarrow \exists \gamma \in \mathcal{F}_{\varphi|\beta}[s\sqsubset \gamma]]$.\\ We will prove: for each $p$ there exists exactly one $s$ such that $length(s)=p$ and $s$ is fine for $\alpha$. \\Define $\alpha_0,\alpha_1$ such that, for both   $i<2$, $(\alpha_i)^i=\underline 0 $ and $\forall j[\neg\exists n[j=\langle i,n\rangle\rightarrow \alpha_i(j)=\alpha(j)]$. \\Define $\alpha_{01}$ such that $(\alpha_{01})^0=(\alpha_{01})^1=\underline 0$ and $\forall j[\neg \exists i<2\exists n[j=\langle i,n\rangle\rightarrow\alpha_{01}(j)=\alpha(j)]]$. \\Note: if $\alpha\;\#\;\alpha_0$, then $\alpha=\alpha_1\in \mathcal{P}_1$, and, if $\alpha\;\#\;\alpha_1$, then $\alpha=\alpha_0\in \mathcal{P}_0$, and, \\if $\alpha\;\#\;\alpha_{01}$, then either $\alpha \;\#\;\alpha_0$ or $\alpha\;\#\;\alpha_1$, and, therefore: $\alpha\in \mathcal{P}_0\cup\mathcal{P}_1=\mathbb{D}^2(\mathcal{A}_1)$.
   
   Note:   $\alpha_{01}\in \mathcal{P}_0 \cap \mathcal{P}_1$.
   
    Let $p$ be given. \\Using Brouwer's Continuity Principle $\mathbf{BCP}$, see Subsubsection \ref{SSS:bcpcontchoice}, find $s_0,s_1, m_0, m_1$ such that $length(s_0)=length(s_1)=p$ and  $\forall i<2\forall\beta\in \mathcal{P}_i\cap \overline{\alpha_{01}}m_i\exists \gamma\in\mathcal{F}_{\psi|\beta}[s_i\sqsubset\gamma]$. 
    
    Assume $s_0\perp s_1$. Then $\exists \gamma \in \mathcal{F}_{\psi|\alpha_{01}}\exists\delta \in \mathcal{F}_{\psi|\alpha_{01}}[s_0\sqsubset \gamma \;\wedge\;s_1\sqsubset\delta]$ and this contradicts $\psi|\alpha_{01}\in \mathcal{E}^1_1!$. Conclude: $s_0=s_1$. \\Define $s:=s_0$ and $m:=\max(m_0, m_1)$, and note:
 if $\overline \alpha m =\overline {\alpha_{01}}m$, then    $s$ is fine for  $\alpha$.
    \\ Now assume: $\overline \alpha m \neq \overline{\alpha_{01}} m$. Find $k<2$ such that $\alpha=\alpha_k$ and note $(\overline \alpha m)^{1-k}\perp \underline 0$.
    \\Find $s_2, m_2$ such that  $length(s_2)=p$ and $m<m_2$ and $\forall\beta\in \mathcal{P}_k \cap\overline   \alpha m_2\exists \gamma \in \mathcal{F}_{\psi|\beta}[s_2\sqsubset\gamma] $. \\As  $(\overline \alpha m_2)^{1-k} \perp \underline 0$, conclude: $\forall i<2\forall\beta\in \mathcal{P}_i \cap\overline   \alpha m_2\exists \gamma \in \mathcal{F}_{\psi|\beta}[s_2\sqsubset\gamma] $, and: \\$s_2$ is fine for $\alpha$. 
    \\Clearly then, for each $p$,
 one may find $s$ such that $length(s)=p$ and
  $s$ is fine for $\alpha$. 
 \\
 Suppose $s,t$ are given such that both $s,t$ are fine for $\alpha$. \\Find $m$ such that $\forall \beta \in \mathbb{D}^2(\mathcal{A}_1)[\overline \alpha m \sqsubset \alpha\rightarrow (\exists \gamma \in \mathcal{F}_{\varphi|\beta}[s\sqsubset \gamma]\;\wedge\;\exists \gamma \in \mathcal{F}_{\varphi|\beta}[t\sqsubset \gamma])]$.
 \\Find $k<2$ such that $\overline \alpha m\sqsubset \alpha_k$. \\Note: $\alpha_k \in \mathbb{D}^2(\mathcal{A}_1)$ and $\varphi|\alpha_k \in \mathcal{E}^1_1!$ and conclude: $s\sqsubseteq t \;\vee\; t \sqsubseteq s$.
 \\We thus see: if both $s,t$ are fine for $\alpha$, then $s\sqsubseteq t \;\vee\; t \sqsubseteq s$.
 \\We thus may define $\delta$ such that, for each $p$, $\overline \delta p$ is fine for $\alpha$. \\Conclude: $\delta \in \mathcal{F}_{\psi|\alpha}$, and: $\psi|\alpha \in \mathcal{E}^1_1$, i.e. $\psi|\alpha$ admits a path.\\
 We still have to prove that $\psi|\alpha$ admits {\it exactly one} path. \\Let $\eta$  be  given such that $\delta\;\#\;\eta$. \\ Note: $\psi|\alpha_0\in \mathcal{E}^1_1!$ and find $\lambda$ in $\mathcal{F}_{\psi|\alpha_0}$. \\Using the co-transitivity of the relation  $\#$, distinguish two cases. 
  
   \textit{Case (a)}: $\eta\;\#\;\lambda$. Find $n$ such that $(\psi|\alpha_0)(\overline \eta n)\neq 0$. \\Either: $(\psi|\alpha)(\overline \eta n)=(\psi|\alpha_0)(\overline \eta n)\neq 0$ or: $\alpha\;\#\;\alpha_0$ and $\alpha=\alpha_1$ and $\exists m[(\psi|\alpha)(\overline \eta m)\neq 0]$.
  
    \textit{Case (b)}:  $\delta\;\#\;\lambda$. Then: $\alpha\;\#\;\alpha_0$ and $\alpha=\alpha_1$ and $\exists m[(\psi|\alpha)(\overline \eta m)\neq 0]$. 
    \\We thus see: $\forall \eta[\eta\;\#\;\delta\rightarrow \exists p[(\psi|\alpha)(\overline \eta p) \neq 0]]$, and: $\psi|\alpha\in \mathcal{E}^1_1!$.
  
 Conclude: $\forall \alpha \in \overline{\mathbb{D}^2(\mathcal{A}_1)}[\psi|\alpha \in \mathcal{E}^1_1!]$.
 
 Now assume that $\psi$ reduces $\mathbb{D}^2(\mathcal{A}_1)$ to $\mathcal{E}^1_1!$. 
  
  Conclude: $\forall \alpha \in \overline{\mathbb{D}^2(\mathcal{A}_1)}[\alpha \in \mathbb{D}^2(\mathcal{A}_1)]$. \\According to Theorem \ref{T:disja1}, see Subsubsection \ref{SSS:disjunction}, we have a contradiction.

    \smallskip
(v) Define $\varphi:\omega^\omega\rightarrow\omega^\omega$ such that, for all $\alpha$,  $(\varphi|\alpha)^{\alpha(0)}=\underline 0 $ and \\$\forall n<\alpha(0)[(\varphi|\alpha)^n=\overline{\underline 0}\alpha^0(2n)\ast\langle\alpha^0(2n+1)+1\rangle\ast \alpha^{n+1}]$ and \\$\forall n>\alpha(0)[(\varphi|\alpha)^n=\overline{\underline 0}\alpha^0(2n-2)\ast\langle\alpha^0(2n-1)+1\rangle\ast \alpha^{n}]$. 

Then $\varphi:\omega^\omega\rightarrowtail\omega^\omega$ and $\varphi|\omega^\omega=\mathcal{E}_2!$. 

Conclude: $\mathcal{E}_2!$ is regular in Lusin's sense.
 \end{proof}

According to Theorem \ref{T:exactlyonepath}(iv), $\mathbb{D}^2(\mathcal{A}_1)\npreceq \mathcal{E}_2!$,  and, therefore, also $\mathcal{E}_2 \npreceq \mathcal{E}_2!$.  This is an \textit{intuitionistic} phenomenon, as, in classical descriptive theory, $\mathcal{E}_2\preceq\mathcal{E}_2!$. One may understand this classical fact by replacing $\mathcal{E}_2, \mathcal{E}_2!$ by sets that, from a constructive point of view, are extensions of them, although, classically, they would be judged to be the same. Theorem \ref{T:almoste2e2!} will make this clear. \begin{definition}\label{D:almoste2}We define \\$\mathcal{ALMOST}$-$\mathcal{E}_2:=\{\alpha\mid\alpha \;\#\;\mathcal{A}_2\}=\{\alpha\mid\forall \gamma \exists n[\alpha^n\bigl(\gamma(n)\bigr)= 0]\}$, and\\
 $\mathcal{ALMOST}$-$\mathcal{E}_2! :=\mathcal{ALMOST}$-$\mathcal{E}_2\cap\{\alpha|\forall m\forall n[m\neq n\rightarrow\exists p[\alpha^m(p)\neq 0 \;\vee\;\alpha^n(p)\neq 0]\}$. \end{definition}
 
 $\mathcal{ALMOST}$-$\mathcal{E}_2$ and $\mathcal{ALMOST}$-$\mathcal{E}_2!$ may be called $\mathbf{\Pi}^1_1$-\textit{approximations} to $\mathcal{E}_2$ and $\mathcal{E}_2!$, respectively. 
 \begin{theorem}\label{T:almoste2e2!}   $\mathcal{ALMOST}$-$\mathcal{E}_2 \preceq \mathcal{ALMOST}$-$\mathcal{E}_2!$.  \end{theorem}
 
 \begin{proof}  Define $\psi,\varphi:\omega^\omega\rightarrow\omega^\omega$ such that, for each $\alpha$, \\$(\psi|\alpha)(0) =0$ and 
    $(\varphi|\alpha)^0=\alpha^0=\alpha^{(\psi|\alpha)(0)}$, and, for each $n$, \begin{enumerate}[\upshape (1)]\item\textit{if} $\overline{\alpha^{(\psi|\alpha)(n)}}(n+1) \sqsubset \underline 0$, then $(\psi|\alpha)(n+1)=(\psi|\alpha)(n)$ and $(\varphi|\alpha)^{n+1}=\underline 1$, and,\item \textit{if $\overline{\alpha^{(\psi|\alpha)(n)}}(n+1) \perp \underline 0$}, then $(\psi|\alpha)(n+1)=(\psi|\alpha)(n)+1$ and $(\varphi|\alpha)^{n+1}= \alpha^{(\psi|\alpha)(n+1)}$.\end{enumerate}
 
 The idea behind these definitions is the following.
\\$\varphi$ will be the function reducing $\mathcal{ALMOST}$-$\mathcal{E}_2$ to  $\mathcal{ALMOST}$-$\mathcal{E}_2!$, and $\psi$ will be an {\it auxiliary function}.
  \\Given $\alpha$, we check its subsequences, $\alpha^0, \alpha^1, \ldots$ one by one.
 \\At stage $0$, we start with studying $\alpha^0$ and we define $(\varphi|\alpha)^0=\alpha^0$. 
 \\ At every stage $n+1$, if $(\psi|\alpha)(n)=k$, we consider $\alpha^k$, and we distinguish two cases. 
 \\{\it Case 1}. We discover that $\alpha^k \;\#\;\underline 0$, (as $\overline {\alpha^k}(n+1)\perp\underline 0$).   We now decide  to study  $\alpha^{k+1}$ at the next stage $n+1$, so we define $(\psi|\alpha)(n+1)=k+1$. \\We also define $(\varphi|\alpha)^{n+1}=\alpha^{k+1}$.  
 \\{\it Case 2}. We do not yet see that  $\alpha^k \;\#\;\underline 0$ (as $\overline {\alpha^k}(n+1)\sqsubset\underline 0$). We decide to continue our  study of  $\alpha^k$ at stage $n+1$, so we define $(\psi|\alpha)(n+1)=k$. We also define $(\varphi|\alpha)^{n+1}=\underline 1$. 
 
 \smallskip
Note: for each $\alpha$, for all $k$, if  if $\forall i<k[\alpha^i\;\#\;\underline 0]$, then there exists $j$ such  that $ (\psi|\alpha)(j)=k$. If $j_0$ is the least such $j$ and $\alpha^k = \underline 0$, then $(\varphi|\alpha)^{j_0}=\underline 0$ and, for all $i\neq j_0$, one has $(\varphi|\alpha)^i \;\#\;\underline 0$. 

\smallskip Also note: for all $n,m$, if $n<m$, then {\it either} $(\psi|\alpha)(n) <(\psi|\alpha)(m)$ and $(\varphi|\alpha)^n \;\#\;\underline 0$, {\it  or} $(\psi|\alpha)(n) =(\psi|\alpha)(m)$ and $(\varphi|\alpha)^m=\underline 1 \;\#\;\underline 0$. 
 \smallskip 
 
 Also note: for each $n$, $(\varphi|\alpha)^n = \underline 1$ or $(\varphi|\alpha)^n = \alpha^{(\psi|\alpha)(n)}$.

 We now prove that $\varphi$ reduces $\mathcal{ALMOST}$-$\mathcal{E}_2$ to $  \mathcal{ALMOST}$-$\mathcal{E}_2!$. 
 
 Assume: $\alpha \in \mathcal{ALMOST}$-$\mathcal{E}_2$. Let $\gamma$ be given. \\We want to find $m$ such that  $(\varphi|\alpha)^m\bigl(\gamma(m)\bigr)=0$. 
\\Define $\delta$ such that $\delta(0):=0$ and, for each $n$, \\\textit{if} $\forall   i\le n[(\varphi|\alpha)^{\delta(i)}\circ\gamma\circ\delta(i) \neq 0]$, then   $\delta(n+1):=\mu j[(\psi|\alpha)(j)=n+1]$, and, \\{\it if not}, then $\delta(n+1):= \delta(n)$.
  \\Note: for each $n$, if $\forall i< n[(\varphi|\alpha)^{\delta(i)}\circ\gamma\circ\delta(i)\neq 0]$, then $\forall i\le n[(\varphi|\alpha)^{\delta(i)} = \alpha^i]]$.  
 \\Define $n:=\mu k[\alpha^k\circ\gamma\circ\delta(k)=0]$.  \\Conclude: $(\varphi|\alpha)^{\delta(n)}=\alpha^n$ and: $(\varphi|\alpha)^{\delta(n)}\circ\gamma\circ \delta(n)=0$ and $\exists m[(\varphi|\alpha)^m\circ\gamma(m)=0]$. 
 \\We thus see: $\forall\gamma\exists m[(\varphi|\alpha)^m\circ\gamma(m)=0]$, that is: $\varphi|\alpha \in \mathcal{ALMOST}$-$\mathcal{E}_2$. 
 \\ As we observed already: for all $m,n$, if $m\neq n$, then either $
 (\varphi|\alpha)^m\;\#\;\underline 0$ or $(\varphi|\alpha)^n\;\#\;\underline 0$. \\Conclude: $\varphi|\alpha \in \mathcal{ALMOST}$-$\mathcal{E}_2!$.
 
\smallskip Now assume:  $\varphi|\alpha \in \mathcal{ALMOST}$-$\mathcal{E}_2!$. Let $\gamma$ be given. \\We want to find $m$ such that  $\alpha
^m\bigl(\gamma(m)\bigr)=0$. \\Define $\delta$ such that, for each $n$, $\delta(n)=\gamma\bigl((\psi|\alpha)(n)\bigr) $. 
\\Find $n$ such that $(\varphi|\alpha)^n\circ \delta(n) =0$.
\\Note: $(\varphi|\alpha)^n\;\#\;\underline 1 $ and $(\varphi|\alpha)^n=\alpha^{(\psi|\alpha)(n)}$.
\\Define $m:=(\psi|\alpha)(n)$ and note:\\ $\alpha^m(\gamma(m))= \alpha^{(\psi|\alpha)(n)}(\gamma\bigl(\psi|\alpha)(n)\bigr)=(\varphi|\alpha)^n\bigl(\delta(n)\bigr)=0$. 
\\We thus see: $\forall \gamma \exists n[\alpha^n\circ\gamma(n)=0]$, that is: $\alpha \in \mathcal{ALMOST}$-$\mathcal{E}_2$. 
 \end{proof}
     
   The following Definition has been given already in  Subsubsection \ref{SSS:perhaps}.  
  \begin{definition}\label{D:perhaps}   For every $\mathcal{X}\subseteq \omega^\omega$, we define \\$\mathsf{Perhaps}(\mathcal{X}) =\{\alpha\mid\exists \beta \in \mathcal{X}[\alpha \;\#\;\beta\rightarrow \alpha \in \mathcal{X}]\}$. \\$\mathcal{X}\subseteq\omega^\omega$ is called \emph{perhapsive} if and only if $\mathsf{Perhaps}(\mathcal{X})=\mathcal{X}$.\end{definition}
     
     The first item of the next Theorem extends Theorem \ref{T:perhaps}(iii).
     
     \begin{theorem}\label{T:a11perhapsive} \begin{enumerate}[\upshape (i)] \item  $\mathcal{A}^1_1$ is perhapsive. \item $\mathcal{E}_2!$ is  not perhapsive. \item $\mathcal{E}_2!$ and $\mathcal{E}^1_1!$ are not $\mathbf{\Pi}^1_1$. \end{enumerate}\end{theorem}
     
     \begin{proof} (i) Let $\alpha, \beta$ be given such that $\beta \in \mathcal{A}_1^1$ and $\alpha\;\#\;\beta\rightarrow \alpha \in \mathcal{A}_1^1$. \\Let $\gamma$ be given. Find $m$ such that $\beta(\overline \gamma m)\neq 0$. \\\textit{Either}: $\alpha(\overline \gamma m)=\beta(\overline \gamma(m) \neq 0$, \textit{or} $\alpha\;\#\;\beta$ and $\alpha \in \mathcal{A}^1_1$ and $\exists p[\alpha(\overline \gamma p)\neq 0]$. \\We thus see: $\forall \gamma \exists p[\alpha(\overline \gamma p)\neq 0]$, that is: $\alpha \in \mathcal{A}_1^1$. 
     
     Conclude: $\forall \alpha[\exists \beta\in \mathcal{A}_1^1[\alpha\;\#\;\beta\rightarrow \alpha\in\mathcal{A}^1_1]\rightarrow \alpha\in\mathcal{A}_1^1]$, that is: $\mathcal{A}_1^1$ is perhapsive.
     
     \smallskip (ii) Let $\mathcal{X}$ be the set of all $\alpha$ such that $\alpha(0)=0$ and, for all $n$, \\if  $n=\mu p[\alpha^0(p)\neq 0]$, then $\alpha^{n+1}=\underline 0$ and, if $n\neq\mu p[\alpha^0(p)\neq 0]$, then $ \alpha^{n+1}=\underline 1$.
     \\We shall prove that $\mathcal{X}$ is a subset of $\mathsf{Perhaps}(\mathcal{E}_2!)$ but not of $\mathcal{E}_2!$ itself. \\It then follows that $\mathcal{E}_2!$ is not perhapsive. 
\\Define $\zeta$ such that $\zeta(0)=0$ and $\zeta^0=\underline 0$ and $\forall  n[\zeta^{n+1}= \underline 1]$. Note: $\zeta \in \mathcal{X}\cap\mathcal{E}_2!$. \\Assume: $\alpha \in \mathcal{X}$ and: $\alpha \;\#\;\zeta$.  Find $i,n$ such that $\alpha^i(n)\neq \zeta^i(n)$. \\\textit{Either}: $i=0$ and $\alpha^0(n)\neq 0$, \textit{or}: $i>0$ and $\alpha^i(n)\neq \zeta^i(n)=1$ and $\alpha^0(i-1)\neq 0$. \\In both cases: $\alpha^0\;\#\;\underline 0$ and $\alpha\in \mathcal{E}_2!$.\\We thus see: $\forall \alpha \in \mathcal{X}[\alpha\;\#\;\zeta \rightarrow \alpha \in \mathcal{E}_2!]$ and conclude: $\mathcal{X}\subseteq\mathsf{Perhaps}(\mathcal{E}_2!)$. 
 
 Assume: $\mathcal{X}\subseteq \mathcal{E}_2!$.  Note that $\mathcal{X}$ is a spread containing $\zeta$. \\ Using $\mathbf{BCP}$, find $m,n$ such that $\forall \alpha \in \mathcal{X}[\overline \zeta m\sqsubset \alpha\rightarrow \alpha^n= \underline 0]$.  \\In particular: $\zeta^n =\underline 0$, and: $n=0$. But $\exists \alpha \in \mathcal{X}[\overline \zeta m\sqsubset \alpha\;\wedge\; \alpha^0 \;\#\;\underline 0]$. Contradiction.
 
   Conclude:  $\mathcal{X}\nsubseteq\mathcal{E}_2!$ while  $\mathcal{X}\subseteq\mathsf{Perhaps}(\mathcal{E}_2!)$, \\so: $\mathsf{Perhaps}(\mathcal{E}_2!)\nsubseteq\mathcal{E}_2!$ and: $\mathcal{E}_2!$ is not perhapsive.
   
   \smallskip (iii) Use (i), (ii), and Theorems \ref{T:perhaps}(i),  \ref{T:coan}(ii) and  \ref{T:exactlyonepath}(i).
     \end{proof}

     \section{$\mathcal{A}_1^1$ and $\mathcal{E}_1^1$}\label{S:a11e11}
     
     In this Section, we study the sets $\mathcal{A}^1_1=\mathcal{BAR}:=\{\alpha\mid\forall \gamma \exists n[\alpha(\overline \gamma n)\neq 0]\}$ and \\$\mathcal{E}^1_1=\mathcal{PATH}:=\{\alpha\mid\exists \gamma \forall n[\alpha(\overline \gamma n)= 0]\}$. \\We have seen that $\mathcal{A}^1_1$ is $\mathbf{\Pi}^1_1$-complete and that $\mathcal{E}^1_1$ is $\mathbf{\Sigma}^1_1$-complete, see Theorems \ref{T:coan}(ii) and \ref{T:analytic}(ii).

 \subsection{$\mathcal{A}^1_1$ positively fails to be strictly analytic}\hfill

  The following definitions have been given already  in Subsubsection \ref{SSS:infseq}.
 \begin{definition}For each $\alpha$, $T_\alpha:=\{s\mid \forall t\sqsubset s[\alpha(t)=0]\}$. 

\smallskip
For all $\alpha,\beta$, for all $\gamma$, we define:\\
$\gamma:\alpha \le^\ast \beta \leftrightarrow$
$\bigl(\forall s[s\in T_\alpha\rightarrow \gamma(s) \in T_\beta]\;\wedge\;\forall s\forall t[s\sqsubset t \rightarrow \gamma(s) \sqsubset \gamma(t)]\bigr)$, and: \\
$\gamma:\alpha <^\ast \beta \leftrightarrow$
$\bigl(\forall s[s\in T_\alpha\rightarrow \gamma(s) \in T_\beta]\;\wedge\;\forall s\forall t[s\sqsubset t \rightarrow \gamma(s) \sqsubset \gamma(t)]\;\wedge\;\gamma(\langle\;\rangle)\neq\langle\;\rangle\bigr)$.

 For all $\alpha, \beta$,  we define:
$\alpha<^\ast\beta \leftrightarrow \exists \gamma[\gamma:\alpha<^\ast \beta]$, and: $\alpha\le^\ast\beta \leftrightarrow \exists \gamma[\gamma:\alpha\le^\ast \beta]$.
\end{definition}

$T_\alpha$ is called the \textit{tree determined by $\alpha$}. Note: $\forall \alpha[0=\langle\;\rangle\in T_\alpha]$.

$\alpha \le^\ast \beta$ if and only if there exists a $\sqsubset$-preserving embedding of $T_\alpha$ into $T_\beta$.

$\alpha <^\ast \beta$ if and only if there exists $n$ in $\omega$ and a $\sqsubset$-preserving embedding of $T_\alpha$ into  $\{s\in T_\beta\mid \langle n \rangle \sqsubseteq s\}$. 
\begin{lemma}\label{L:a11}\hfill

\begin{enumerate}[\upshape (i)]

\item For all  $\alpha,\beta,\gamma$, $\alpha\le^\ast\alpha$, $(\alpha\le^\ast\beta\;\wedge\;\beta\le^\ast\gamma)\rightarrow \alpha\le^\ast\gamma)$ and: $\alpha<^\ast\beta\rightarrow\alpha\le^\ast \beta$ and: 
$(\alpha<^\ast\beta\;\wedge\;\beta\le^\ast\gamma)\rightarrow \alpha<^\ast\gamma$ and: $(\alpha\le^\ast\beta\;\wedge\;\beta<^\ast\gamma)\rightarrow \alpha<^\ast\gamma$.
\item $\forall \alpha \in \mathcal{A}_1^1\forall\beta\in\mathcal{A}_1^1[\alpha<^\ast\beta\rightarrow\alpha\;\#\;\beta]$.

\end{enumerate}

\end{lemma}

\begin{proof}

 (i) Note: for all $\alpha,\beta,\gamma,\delta,\varepsilon$, if $\delta:\alpha\le^\ast\beta$ and $\varepsilon:\beta\le^\ast\gamma$, then $\varepsilon\circ\delta:\alpha\le^\ast\gamma$. Conclude: if $\alpha\le^\ast\beta$ and $\beta\le^\ast\gamma$, then $\alpha\le^\ast\gamma$.
 
 The proofs of the other statements are also straightforward. 

\smallskip
 (ii) Let $\alpha,\beta$ in $\mathcal{A}^1_1$ be given such that $\alpha <^\ast\beta$. \\Find $\gamma$ such that $\forall s \in T_\alpha[\gamma(s)\in T_\beta]$ and $\forall s \forall t[s\sqsubset t\rightarrow \gamma(s) \sqsubset\gamma(t)]$ and $\gamma(\langle\;\rangle)\neq\langle\;\rangle$. \\Define $\varepsilon$ such that $\varepsilon(0) =\gamma(0)=\gamma
 (\langle\;\rangle)$ and, for each $n$, $\varepsilon(n+1)= \gamma\circ\varepsilon(n)$.\\ Note: for all $n$, $\varepsilon(n)\sqsubset \varepsilon(n+1)$, and, if $\varepsilon(n)\in T_\alpha$, then $\varepsilon(n+1)\in T_\beta$.   
 \\Find $\delta$ such that $\forall n[\varepsilon(n)\sqsubset \delta]$ and note:  $\exists n[\overline \delta n \notin T_\alpha]$. \\Conclude: $\exists m[\varepsilon(m)\notin T_\alpha]$ and define $p:=\mu m[\varepsilon(m)\notin T_\alpha]$. \\Note: $p>0$ and find $q$  such that $p=q+1$. \\Conclude: $\varepsilon(q)\in T_\alpha$ and $\varepsilon(p)\in T_\beta\setminus T_\alpha$ and: $\alpha\;\#\;\beta$. 
\end{proof}

 The next Theorem, Theorem \ref{T:cdiag},  shows that $\mathcal{A}_1^1$ \textit{positively fails to be strictly analytic or $\mathbf{\Sigma}^{1\ast}_1$} in the following sense: given a (continuous) function from $\omega^\omega$ into $\mathcal{A}_1^1$ one may construct an element of $\mathcal{A}_1^1$ that does not occur in the range of $\varphi$. 
\begin{theorem}\label{T:cdiag} \hfill

 \begin{enumerate}[\upshape (i)]
  \item \emph{Cantor's  diagonal argument:} 
$\forall \varphi:\omega^\omega\rightarrow\mathcal{A}_1^1\exists \alpha \in \mathcal{A}_1^1\forall \beta[\alpha\;\#\;\varphi|\beta]$. 
\item \emph{The Boundedness Theorem:} $\forall \varphi:\omega^\omega\rightarrow\mathcal{A}_1^1\exists \alpha \in \mathcal{A}_1^1\forall \beta[\varphi|\beta\le^\ast \alpha]$. 
\end{enumerate}

\end{theorem}

\begin{proof}
(i)  Assume: $\varphi:\omega^\omega\rightarrow \mathcal{A}^1_1$. \\Define $\alpha:\omega^\omega\rightarrow \omega$ such that $\forall \beta[\alpha(\beta) =(\varphi|\beta)(\beta) +1]$.  
\\Note: $\alpha \in \mathcal{A}^1_1$ and $\forall \beta[\alpha\;\#\; \varphi|\beta]$. 

\smallskip
 (ii). Assume: $\varphi:\omega^\omega\rightarrow \mathcal{A}^1_1$. 
 \\ Note: $\forall \beta\forall\delta \exists n[(\varphi|\beta)(\overline \delta n)\neq 0]$, and:
 $\forall \beta \forall \delta \exists n \exists m[\varphi^{\overline \delta n}(\overline \beta m)>1 \;\wedge\; \forall i<m[\varphi^{\overline\delta n}(\overline \beta i)=0]]$. 
 \\Define $\alpha$ such that $\forall s[\alpha(s)\neq 0\leftrightarrow \exists t\sqsubseteq s_I\exists u\sqsubseteq s_{II}[\varphi^t(u) >1\;\wedge\;\forall v\sqsubset u[\varphi^t(v)=0]]]$.
\\Note: $\alpha\in\mathcal{A}_1^1$. 
\\Let $\beta$ be given. Define $\varepsilon$ such that $\forall d\forall n[n=length(d)\rightarrow \varepsilon(d)=\ulcorner d, \overline \beta n\urcorner]$. 
\\Note: $\varepsilon: \varphi|\beta\le^\ast\alpha$. \\We thus see: $\forall \beta[\varphi|\beta\le^\ast\alpha]$.
\end{proof}

Using Lemma \ref{L:a11}, one may obtain Theorem \ref{T:cdiag}(i) from Theorem \ref{T:cdiag}(ii), as follows. Assume $\alpha\in\mathcal{A}_1^1$ and $\forall \beta[\varphi|\beta\le^\ast\alpha]$. \\Note\footnote{For each $\alpha$, $S^\ast(\alpha)$ is the element $\beta$ of $\omega^\omega$ such that $\beta(0)=0$ and $\forall n[\beta^n=\alpha]$,  see Subsubsection \ref{SSS:stumps}. If $\alpha \in \mathcal{A}^1_1$, then also $S^\ast(\alpha) \in \mathcal{A}^1_1$. $S^\ast(\alpha)$ is called the {\it successor} of $\alpha$.}: $S^\ast(\alpha)\in \mathcal{A}^1_1$ and $\forall \beta[\varphi|\beta<^\ast S^\ast(\alpha)]$ and thus, according to  Theorem \ref{T:cdiag}(i), $\forall \beta[\varphi|\beta\;\#\; S^\ast(\alpha)]$.

\subsection{$\mathcal{E}^1_1$ positively fails to be $\mathbf{\Pi}^1_1$}\hfill

The next Theorem, Theorem \ref{T:e2notpi02a2notsigma2}, should prepare the reader for Theorem  \ref{T:e11notpi11}. \\The  proof of Theorem \ref{T:e2notpi02a2notsigma2} is elementary in the sense that no use is made of intuitionistic principles like Brouwer's Continuity Principle $\mathbf{BCP}$ or the Fan Theorem  $\mathbf{FT}$.  The proof of Theorem \ref{T:e2notpi02a2notsigma2}(i)  has been given  in \cite[Section 5.4]{veldman08}. 
Theorem \ref{T:e2notpi02a2notsigma2}(iii) is a rather weak statement if one compares it to the result of the Borel Hierarchy Theorem, Theorem \ref{T:borelhiertheorem} in Subsubsection \ref{SSS:borelgeneral}. One should compare Theorem \ref{T:e2notpi02a2notsigma2}(iii) to 
Theorem \ref{T:a11notsigma11}(i). 

\begin{theorem}\label{T:e2notpi02a2notsigma2} \hfill

\begin{enumerate}[\upshape (i)] \item \emph{$\mathcal{E}_2$ positively fails to be $\mathbf{\Pi}^0_2$:} if a continuous function maps $\mathcal{E}_2$ into $\mathcal{A}_2$, it also maps some element of $\mathcal{A}_2$ into $\mathcal{A}_2$:

$\forall \varphi:\omega^\omega\rightarrow\omega^\omega[\forall\alpha \in\mathcal{E}_2 [ \varphi|\alpha \in \mathcal{A}_2]\rightarrow \exists \alpha \in \mathcal{A}_2 [\varphi|\alpha\in\mathcal{A}_2]]$. 

\item If $\mathcal{E}_2$ is contained in a set $\mathcal{X}$ that is a countable intersection of open sets, also some element of $\mathcal{A}_2$ is in $\mathcal{X}$: $\forall \beta[\mathcal{E}_2\subseteq \mathcal{F}^2_\beta\rightarrow\exists \alpha[\alpha \in  \mathcal{A}_2\cap\mathcal{F}^2_\beta]]$.

\item The assumption that $\mathcal{A}_2$ is a countable union of \emph{spreads} leads to a contradiction: $\neg \exists \beta[\forall n[Spr(\beta^n)]\;\wedge\;\mathcal{A}_2=\bigcup_n\mathcal{F}_{\beta^n}]$. 

\end{enumerate} \end{theorem}
\begin{proof} (i) Assume $\varphi:\omega^\omega\rightarrow\omega^\omega$ and $\forall \alpha \in \mathcal{E}_2[\varphi|\alpha \in \mathcal{A}_2]$.  \\Now define $\alpha$ such that, for all $n$, for all $t$, $\alpha^n(t)\neq 0$ if and only if $(\varphi|t)^n\perp\underline 0$. \\Note: for all $n$, $\alpha^n\;\#\;\underline 0$ if and only if $(\varphi|\alpha)^n\;\#\;\underline 0$.

We now prove: for all $n$, both $\alpha^n$ and $(\varphi|\alpha)^n$ are in $\mathcal{E}_1$.   
\\Let $n$ be given. Define $\alpha_n$ such that  $(\alpha_n)^n =\underline 0$ and $\forall j[\neg \exists t[j=\langle n \rangle \ast t]\rightarrow \alpha_n(j)=\alpha(j)]$. Note: $\alpha_n\in \mathcal{E}_2$ and $\varphi|\alpha_n\in \mathcal{A}_2$ and: $(\varphi|\alpha_n)^n \perp \underline  0$. \\Find $t\sqsubset \alpha_n$ such that $(\varphi|t)^n\perp \underline 0$ and  distinguish two cases. \\\textit{Either}: $t\sqsubset \alpha$ and $(\varphi|\alpha)^n\;\#\;\underline 0$ and also $\alpha^n\;\#\;\underline 0$, \\\textit{or}: $t\perp\alpha$ and $\alpha_n\perp\alpha$ and $\alpha^n\;\#\;\underline 0$ and also $(\varphi|\alpha)^n\;\#\;\underline 0$.
\\We thus see: for all $n$, $\alpha^n\;\#\;\underline 0$ and $(\varphi|\alpha)^n\;\#\;\underline 0$, i.e. $\alpha\in\mathcal{A}_2$ and $\varphi|\alpha\in \mathcal{A}_2$. 

\smallskip (ii) Let $\beta$ given such that $\mathcal{E}_2 \subseteq \mathcal{F}^2_\beta$.  Find $\varphi:\omega^\omega\rightarrow\omega^\omega$ reducing $\mathcal{F}^2_\beta$ to $\mathcal{A}_2$. \\ Note: $\forall \alpha \in \mathcal{E}_2[\varphi|\alpha \in \mathcal{A}_2]$. \\Applying (i), find $\alpha$ in $\mathcal{A}_2$ such that $\varphi|\alpha \in \mathcal{A}_2$, so $\alpha \in \mathcal{A}_2\cap \mathcal{F}^2_\beta$. 

\smallskip (iii) Let $\beta$ be given such that $\forall n[Spr(\beta^n)]$ and $\mathcal{A}_2=\mathcal{G}^2_\beta=\bigcup_n\mathcal{F}_{\beta_n}$. \\Find $\rho$ such that, for each $n$, $\rho^n:\omega^\omega\rightarrow\mathcal{F}_{\beta^n}$ is the canonical retraction of $\omega^\omega$ onto $\mathcal{F}_{\beta^n}$. Assume: $\alpha\in \mathcal{E}_2$. Note: $\forall \delta \in \mathcal{A}_2[\alpha\;\#\;\delta]$ and: $\forall n\forall \delta \in \mathcal{F}_{\beta^n}[\alpha\;\#\;\delta]$ and: $\forall n[\alpha\;\#\; \rho^n|\alpha]$ and: $\forall n\exists m[\beta^n(\overline \alpha m)\neq 0]$ and: $\forall n [\alpha \in \mathcal{G}_{\beta^n}]$ and: $\alpha\in \mathcal{F}^2_\beta$.  We thus see: $\forall \alpha \in \mathcal{E}_2[\alpha\in \mathcal{F}^2_\beta]$, that is: $\mathcal{E}_2\subseteq \mathcal{F}^2_\beta$. Applying (ii), we find $\alpha \in \mathcal{A}_2 \cap \mathcal{F}^2_\beta=\mathcal{G}^2_\beta\cap\mathcal{F}^2_\beta=\emptyset$. Contradiction.  \end{proof}
The proof of the next Theorem, Theorem \ref{T:e11notpi11}, is also elementary.

\begin{theorem}\label{T:e11notpi11}\hfill

\begin{enumerate}[\upshape (i)] \item
 \emph{$\mathcal{E}_1^1$ positively fails to be $\mathbf{\Pi}^1_1$:} If a continuous function from $\omega^\omega$ to $\omega^\omega$ maps $\mathcal{E}^1_1$ into $\mathcal{A}^1_1$, it also maps some element of $\mathcal{A}^1_1$ into $\mathcal{A}^1_1$:

$\forall \varphi:\omega^\omega\rightarrow\omega^\omega[\forall\alpha \in\mathcal{E}_1^1 [ \varphi|\alpha \in \mathcal{A}_1^1]\rightarrow \exists \alpha \in \mathcal{A}_1^1 [\varphi|\alpha\in\mathcal{A}_1^1]]$. 

\item If $\mathcal{E}^1_1$ is contained in a $\mathbf{\Pi}^1_1$ set $\mathcal{X}$, also some element of $\mathcal{A}^1_1$ is in $\mathcal{X}$: 

$\forall \beta[\mathcal{E}_1^1\subseteq \mathcal{UG}_\beta\rightarrow\exists \alpha[\alpha \in  \mathcal{A}_1^1\cap\mathcal{UG}_\beta]]$. \end{enumerate}
\end{theorem}

\begin{proof} (i) Assume $\varphi:\omega^\omega\rightarrow\omega^\omega$ and $\forall \alpha \in  \mathcal{E}^1_1[\varphi| \alpha\in\mathcal{A}^1_1]$. \\Now define $\alpha$ such that, for all $t$,  $\alpha(t)\neq 0$ if and only if $ \exists s \sqsubseteq t[(\varphi| \overline \alpha t)(s)\neq 0]$. \\ Note: for all $\gamma$, $\exists n[\alpha(\overline \gamma n)\neq 0]$ if and only if $ \exists n[(\varphi|\alpha)(\overline \gamma n)\neq 0]$.

We now prove: for all $\gamma$, $\exists n[\alpha(\overline \gamma n)\neq 0]$ and $ \exists n[(\varphi|\alpha)(\overline \gamma n)\neq 0]$.

Let $\gamma$ be given. Define $\alpha_\gamma$ such that  $\forall n[\alpha(\overline \gamma n)=0]$ and $\forall t[t\perp \gamma \rightarrow   \alpha_\gamma(t)=\alpha(t)]$. Note: $\alpha_\gamma\in \mathcal{E}_1^1$ and: $\varphi|\alpha_\gamma\in \mathcal{A}_1^1$. Find $m$ such that $(\varphi|\alpha_\gamma)(\overline \gamma m)\neq 0$.  \\Find $t\sqsubset \alpha_\gamma$ such that $(\varphi|t)(\overline \gamma m)\neq 0$ and distinguish two cases. \\\textit{Either}: $t\sqsubset \alpha$ and $(\varphi|\alpha)(\overline \gamma m)\neq 0$ and: $\exists n\le m[\alpha(\overline \gamma n)\neq 0]$, \\\textit{or}: $t\perp\alpha$ and $\alpha\perp\alpha_\gamma$ and  $\exists n[\alpha(\overline \gamma n)\neq 0]$  and: $\exists  n[(\varphi|\alpha)(\overline \gamma n)\neq 0]$.

We thus see: for all $\gamma$,  $\exists n[\alpha(\overline \gamma n)\neq 0]$ and $\exists n[(\varphi|\alpha(\overline \gamma n)\neq 0]$, i.e. \\$\alpha\in\mathcal{A}_1^1$ and $\varphi|\alpha\in \mathcal{A}_1^1$. 

\smallskip (ii) Let $\beta$ given such that $\mathcal{E}_1^1 \subseteq \mathcal{UG}_\beta$.  Find $\varphi:\omega^\omega\rightarrow\omega^\omega$ reducing $\mathcal{UG}_\beta$ to $\mathcal{A}_1^1$. \\ Note: $\forall \alpha \in \mathcal{E}_1^1[\varphi|\alpha \in \mathcal{A}_1^1]$. \\Applying (i), find $\alpha$ in $\mathcal{A}_1^1$ such that $\varphi|\alpha \in \mathcal{A}_1^1$, so $\alpha \in \mathcal{A}_1^1\cap \mathcal{UG}_\beta$. 
\end{proof}

\subsection{May one prove: `$\mathcal{A}^1_1$ is not analytic'?}\hfill

The following  Theorem should be compared to \cite[Theorem 5.2(iv)]{veldman08}.
\begin{theorem}\label{T:a11notsigma11}
\hfill
\begin{enumerate}[\upshape (i)] \item If $\mathcal{A}_2$ is a countable union of \emph{closed} sets, there exists $\alpha$ not in either $\mathcal{A}_2$ or $\mathcal{E}_2$:  $\mathcal{A}_2 \preceq\mathcal{E}_2\rightarrow\exists\alpha[\alpha \notin \mathcal{E}_2\;\wedge\; \alpha \notin \mathcal{A}_2]$. \item If $\mathcal{A}_1^1$ is analytic, there exists $\alpha$ not in either $\mathcal{A}_1^1$ or $\mathcal{E}_1^1$:\\
$\mathcal{A}^1_1\preceq\mathcal{E}_1^1\rightarrow\exists\alpha[\alpha \notin \mathcal{E}^1_1\;\wedge\; \alpha \notin \mathcal{A}^1_1]$.
\end{enumerate}
\end{theorem}

\begin{proof} (i)  Let $\varphi:\omega^\omega\rightarrow\omega^\omega$ be given.  \\Define $\alpha$ such that, for all $n$, for all $t$, $\alpha^n(t)\neq 0$ if and only if $ \exists s \sqsubseteq t[(\varphi| \overline \alpha t)^n(s)\neq 0]$.  \\Note: for all $n$,  $ \exists m[\alpha^n(m)\neq 0]$ if and only if $ \exists m[(\varphi|\alpha)^n(m)\neq 0]]$, so \\$\alpha^n \in \mathcal{E}_1$ if and only if $(\varphi|\alpha)^n \in \mathcal{E}_1$ and: $\alpha^n \in \mathcal{A}_1$ if and only if $(\varphi|\alpha)^n \in \mathcal{A}_1$.\\Conclude:  $\alpha \in \mathcal{E}_2$ if and only if $ \varphi|\alpha \in \mathcal{E}_2$ and: $\alpha \in \mathcal{A}_2$ if and only if $ \varphi|\alpha \in \mathcal{A}_2$.

Now assume, in addition: $\varphi$ reduces $\mathcal{A}_2$ to $\mathcal{E}_2$. \\If $\alpha \in \mathcal{A}_2$, then both $\varphi|\alpha \in \mathcal{E}_2$ and $\varphi|\alpha \in \mathcal{A}_2$:  contradiction. \\If $\alpha \in \mathcal{E}_2$, then both $\varphi|\alpha \in \mathcal{E}_2$ and $\alpha \in \mathcal{A}_2$:  contradiction. \\We thus see: $\alpha \notin \mathcal{A}_2$ and $\alpha \notin \mathcal{E}_2$.

\smallskip (ii) Let $\varphi:\omega^\omega\rightarrow\omega^\omega$ be given.  \\Define $\alpha$ such that, for all $t$,  $\alpha(t)\neq 0$ if and only if $\exists s \sqsubseteq t[(\varphi| \overline \alpha t)(s)\neq 0]$.  \\Note: for each $\gamma$, $\exists n[\alpha(\overline \gamma n)\neq 0]$ if and only if $\exists n[(\varphi|\alpha)(\overline \gamma n)\neq 0]$. \\Conclude: $\alpha \in \mathcal{E}_1^1$ if and only if $\varphi|\alpha \in \mathcal{E}^1_1$ and: $\alpha \in \mathcal{A}_1^1$ if and only if $\varphi|\alpha \in \mathcal{A}^1_1$.

Now assume, in addition: $\varphi$ reduces $\mathcal{A}_1^1$ to $\mathcal{E}^1_1$. \\If $\alpha \in \mathcal{A}^1_1$, then both $\varphi|\alpha \in \mathcal{E}_1^1$ and $\varphi|\alpha \in \mathcal{A}^1_1$:  contradiction. \\If $\alpha \in \mathcal{E}^1_1$, then both $\varphi|\alpha \in \mathcal{E}_1^1$ and $\alpha \in \mathcal{A}^1_1$:  contradiction. \\We thus see: $\alpha \notin \mathcal{A}^1_1$ and $\alpha \notin \mathcal{E}_1^1$. 
\end{proof}

\textit{Markov's Principle} $\mathbf{MP}$, in our view a dubious assumption, see Subsubsection 
\ref{SSS:lpomp}, proves: \\$\alpha \notin \mathcal{E}_2\Rightarrow \neg \exists n \forall m[\alpha^n(m)=0]\Rightarrow \forall  n\neg\neg\exists m[\alpha^n(m)\neq 0]\Rightarrow \forall n \exists m[\alpha^n(m)\neq 0]\Rightarrow \alpha \in \mathcal{A}_2$,\\ and thus, together with Theorem \ref{T:a11notsigma11}(i):  $\mathcal{A}_2\npreceq\mathcal{E}_2$.

 $\mathbf{MP}$ proves also the following: \\$\alpha \notin \mathcal{E}^1_1\Rightarrow \neg \exists \gamma \forall n[\alpha(\overline \gamma n)=0]\Rightarrow \forall \gamma \neg\neg\exists n[\alpha(\overline \gamma n)\neq 0]\Rightarrow \forall \gamma \exists n[\alpha(\overline \gamma n)\neq 0]\Rightarrow \alpha \in \mathcal{A}^1_1$,\\ and thus, together with Theorem \ref{T:a11notsigma11}(ii):  $\mathcal{A}_1^1\npreceq\mathcal{E}^1_1$.

Intuitionistically, one obtains the  conclusion: $\mathcal{A}_2\npreceq\mathcal{E}_2$ as a corollary of a stronger statement proven from Brouwer's Continuity Principle $\mathbf{BCP}$, see Theorem \ref{T:borelhiertheorem} in Subsubsection \ref{SSS:borelgeneral}. No such argument seems to be available for the conclusion: $\mathcal{A}^1_1\npreceq\mathcal{E}^1_1$.

One may prove: $\mathcal{A}^1_1\npreceq\mathcal{E}^1_1$, avoiding $\mathbf{MP}$,  but  using $\mathbf{KS}$,  see Subsubsection \ref{SSS:creasubj}. One may argue that $\mathcal{A}_1^1$ is \textit{definite}, and therefore, if analytic, also strictly analytic,  see  Theorem \ref{T:kripkesan}
 in Subsection \ref{SS:strictlyanalytic}.  \\We have seen that $\mathcal{A}_1^1$ is not strictly analytic, see  Theorem \ref{T:cdiag}.

\subsection{$\mathcal{E}^1_1$ and $\mathcal{A}^1_1$ positively fail to be (positively) Borel}\hfill

In classical descriptive set theory, the following statement holds: \begin{quote}{\it   A continuous function $\varphi:\omega^\omega\rightarrow\omega^\omega$ reducing $\mathcal{X}\subseteq \omega^\omega$ to $\mathcal{E}^1_1$  \\ reduces $\omega^\omega\setminus \mathcal{X}$ to $\mathcal{A}^1_1$}. \end{quote} So, if one has seen that every Borel $\mathcal{X}\subseteq \omega^\omega$ is $\mathbf{\Sigma}^1_1$ and reduces to $\mathcal{E}^1_1$, one may conclude that every Borel $\mathcal{X}\subseteq \omega^\omega$ reduces to $\mathcal{A}^1_1$ and is $\mathbf{\Pi}^1_1$.  In our constructive context, this conclusion is wrong, see Theorems \ref{T:analytic}(iv) and  Theorem \ref{T:coan}(iv).

The following subtle Lemma \ref{L:redpairsborela11e11} replaces the just mentioned statement. 

\begin{lemma}\label{L:redpairsborela11e11} 
For every complementary pair $(\mathcal{X}, \mathcal{Y})$ of positively Borel sets there exists $\varphi:\omega^\omega\rightarrow\omega^\omega$ reducing $\mathcal{X}$ to $\mathcal{E}^1_1$ and mapping $\mathcal{Y}$ into $\mathcal{A}^1_1$.

\end{lemma}

\begin{proof} We use induction on the class of complementary pairs of Borel sets and distinguish three cases.

\textit{Case 1.} Let $\beta$ be given such that \\$\mathcal{X}=\mathcal{G}_\beta=\{\alpha\mid \exists n[\beta(\overline \alpha n)\neq 0]\}$ and $\mathcal{Y}=\mathcal{F}_\beta=\{\alpha\mid \forall n[\beta(\overline \alpha n)= 0]\}$.

 Define $\varphi:\omega^\omega\rightarrow\omega^\omega$ such that \\$\forall \alpha[(\varphi|\alpha)(0)=0\;\wedge\;\forall s>0[(\varphi|\alpha)(s)=0\leftrightarrow \beta\bigl(\overline \alpha(s(0))\bigr)\neq 0]]$.

Note that $\varphi$ simultaneously reduces $\mathcal{X}$ to $\mathcal{E}_1^1$ and $\mathcal{Y}$ to $\mathcal{A}_1^1$.

Because: for each $\alpha$,  $\alpha\in \mathcal{G}_\beta \leftrightarrow \exists n[\beta(\overline \alpha n)\neq 0]\leftrightarrow \exists \gamma[\beta(\overline \alpha\gamma(0))\neq 0] \leftrightarrow \\\exists \gamma[(\varphi|\alpha)(\langle \gamma(0)\rangle)= 0]\leftrightarrow \exists \gamma\forall n[(\varphi|\alpha)(\overline \gamma n)=0]\leftrightarrow \varphi|\alpha \in \mathcal{E}^1_1$.

\smallskip And:  for each $\alpha$,  $\alpha\in \mathcal{F}_\beta \leftrightarrow \forall  n[\beta(\overline \alpha n)= 0]\leftrightarrow \forall \gamma[\beta(\overline \alpha\gamma(0))= 0] \leftrightarrow \\\forall \gamma[(\varphi|\alpha)(\langle \gamma(0)\rangle)\neq 0]\leftrightarrow\forall \gamma\exists n[(\varphi|\alpha)(\overline \gamma n)\neq 0]\leftrightarrow \varphi|\alpha \in \mathcal{A}^1_1$.

\smallskip
\textit{Case 2.} Let $\beta$ be given such that $\mathcal{X}=\mathcal{F}_\beta$ and $\mathcal{Y}=\mathcal{G}_\beta$.

Define $\varphi:\omega^\omega\rightarrow\omega^\omega$ such that $\forall \alpha\forall s[(\varphi|\alpha)(s)=0\leftrightarrow \forall j\le s[ \beta(\overline \alpha j)=0]]$.

Note that $\varphi$ simultaneously reduces $\mathcal{X}$ to $\mathcal{E}_1^1$ and $\mathcal{Y}$ to $\mathcal{A}_1^1$.

\smallskip Because: for each $\alpha$, $\alpha \in \mathcal{F}_\beta \leftrightarrow \forall n [\beta(\overline \alpha n)=0]\leftrightarrow \forall s[(\varphi|\alpha)(s)=0]\leftrightarrow\\ \forall \gamma \forall n [(\varphi|\alpha)(\overline \gamma n)=0] \leftrightarrow \exists \gamma \forall n[(\varphi|\alpha)(\overline \gamma n)=0]\leftrightarrow \varphi|\alpha \in \mathcal{E}^1_1$. 

\smallskip And: for each $\alpha$, $\alpha \in \mathcal{G}_\beta \leftrightarrow \exists n [\beta(\overline \alpha n)\neq 0]\leftrightarrow \exists s\forall t\ge s [(\varphi|\alpha)(t)\neq 0]\leftrightarrow\\ \forall \gamma \exists n [(\varphi|\alpha)(\overline \gamma n)\neq 0] \leftrightarrow \varphi|\alpha \in \mathcal{A}^1_1$.

\smallskip \textit{Case 3.} Let $(\mathcal{X}_0,\mathcal{Y}_0), (\mathcal{X}_1,\mathcal{Y}_1),\ldots$ be an infinite sequence of complementary pairs of (positively) Borel sets and let $\varphi$ be given such that, for each $n$, $\varphi^n:\omega^\omega\rightarrow\omega^\omega$ reduces $\mathcal{X}_n$ to $\mathcal{E}_1^1$ and maps $\mathcal{Y}_n$ into $\mathcal{A}_1^1$.

\smallskip\textit{Case 3a.} Define $\mathcal{X}= \bigcup_n\mathcal{X}_n$ and $\mathcal{Y}:=\bigcap_n\mathcal{Y}_n$. 

Define $\psi:\omega^\omega\rightarrow\omega^\omega$ such that $\forall\alpha[(\psi|\alpha)(0)=0 \;\wedge\; \forall n\forall s[(\psi|\alpha)(\langle n\rangle \ast s)=(\varphi^n|\alpha)(s)]]$.

Note that $\psi$  reduces $\mathcal{X}$ to $\mathcal{E}_1^1$ and maps $\mathcal{Y}$ into $\mathcal{A}_1^1$.

\smallskip Because: for each $\alpha$, $\alpha \in \mathcal{X}\leftrightarrow \exists n[\alpha \in \mathcal{X}_n] \leftrightarrow 
\exists n[\varphi^n|\alpha\in \mathcal{E}^1_1] \leftrightarrow 
\\\exists n \exists \gamma \forall m[(\varphi^n|\alpha)(\overline \gamma m)=0]
\leftrightarrow  \exists \gamma\forall m[(\psi|\alpha)(\overline \gamma m)=0] \leftrightarrow \psi|\alpha\in \mathcal{E}^1_1$.

\smallskip And: for each $\alpha$, $\alpha\in \mathcal{Y} \leftrightarrow \forall n[\alpha \in \mathcal{Y}_n] \rightarrow \forall n[\varphi^n|\alpha \in \mathcal{A}^1_1]\leftrightarrow \\\forall n \forall \gamma \exists m[(\varphi^n|\alpha)(\overline \gamma m)\neq 0]\leftrightarrow\forall \gamma \exists m[(\psi|\alpha)(\overline \gamma m)\neq 0]$, so $\alpha \in \mathcal{Y} \rightarrow \psi|\alpha \in \mathcal{A}^1_1$. 

\smallskip\textit{Case 3b.} Define $\mathcal{X}= \bigcap_n\mathcal{X}_n$ and $\mathcal{Y}:=\bigcup_n\mathcal{Y}_n$. 

Define $\psi:\omega^\omega\rightarrow\omega^\omega$ such that $\forall\alpha\forall s[(\psi|\alpha)(s)=0 \leftrightarrow \forall n\le   s\forall t\sqsubseteq s ^n[(\varphi^n|\alpha)(t)=0]]$.

Note that $\psi$  reduces $\mathcal{X}$ to $\mathcal{E}_1^1$ and maps $\mathcal{Y}$ into $\mathcal{A}_1^1$.

\smallskip Because: for each $\alpha$, $\alpha \in \mathcal{X}\leftrightarrow \forall n[\alpha \in \mathcal{X}_n]\leftrightarrow \forall n[\varphi^n|\alpha \in \mathcal{E}_1^1]\leftrightarrow \\\forall n \exists \gamma\forall m[(\varphi^n|\alpha)(\overline \gamma m)=0]\leftrightarrow\footnote{We are applying the Second Axiom of Countable Choice, $\mathbf{AC}_{0,1}$:  $\forall m\exists \gamma [m\mathcal{R}\gamma]\rightarrow \exists \gamma \forall m [m\mathcal{R}\gamma^m]$, see Subsubsection \ref{SSS:countablechoice}.} \exists \gamma \forall n \forall m[(\varphi^n|\alpha)(\overline {\gamma^n} m)=0]\leftrightarrow \exists \gamma\forall m[(\psi|\alpha)(\overline \gamma m)=0]\leftrightarrow \psi|\alpha \in \mathcal{E}^1_1$.

\smallskip And: for each $\alpha$, $\alpha \in \mathcal{Y}\leftrightarrow \exists n[\alpha \in \mathcal{Y}_n]\rightarrow \exists n[\varphi^n|\alpha \in \mathcal{A}_1^1]\leftrightarrow \\\exists n \forall \gamma\exists m[(\varphi^n|\alpha)(\overline \gamma m)\neq 0]\rightarrow\footnote{The contraposition of $\mathbf{AC}_{0,1}$: $\forall \gamma \exists m [m\mathcal{R}\gamma^m]\rightarrow \exists m \forall \gamma[m\mathcal{R}\gamma]$,  is not constructively valid, and, therefore, we have here a single arrow only.} \forall \gamma \exists n \exists m[(\varphi^n|\alpha)(\overline {\gamma^n} m)\neq 0]\leftrightarrow \forall \gamma\exists m[(\psi|\alpha)(\overline \gamma m)\neq 0]\leftrightarrow \psi|\alpha \in \mathcal{A}^1_1$, so $\alpha\in \mathcal{Y}\rightarrow \psi|\alpha \in \mathcal{A}^1_1$. 
\end{proof}

\begin{theorem}[$\mathcal{E}_1^1$ and $\mathcal{A}^1_1$ positively fail to be (positively) Borel]\label{T:e11nonborel}\hfill

\begin{enumerate}[\upshape (i)]

\item For every $\sigma$ in $\mathcal{HRS}$,  for every $\varphi:\omega^\omega\rightarrow\omega^\omega$, \\if $\varphi| \mathcal{E}^1_1\subseteq \mathcal{E}_\sigma$, then $\exists\alpha\in\mathcal{A}_1^1[\varphi|\alpha\in\mathcal{E}_\sigma]$.

\item For every $\mathcal{X}$ in $\mathfrak{Borel}$, if $\mathcal{E}^1_1\subseteq\mathcal{X}$, then $\exists\alpha\in\mathcal{A}_1^1[\alpha \in\mathcal{X}]$. 

\item For every $\sigma$ in $\mathcal{HRS}$,  for every $\varphi:\omega^\omega\rightarrow\omega^\omega$, \\if $\varphi| \mathcal{A}^1_1\subseteq \mathcal{E}_\sigma$, then $\exists\alpha\in\mathcal{E}_1^1[\varphi|\alpha\in\mathcal{E}_\sigma]$.

\item For every $\mathcal{X}$ in $\mathfrak{Borel}$, if $\mathcal{A}^1_1\subseteq\mathcal{X}$, then $\exists\alpha\in\mathcal{E}_1^1[\alpha \in\mathcal{X}]$. 
\end{enumerate}

\end{theorem}

 \begin{proof} 

(i) Let $\sigma,\varphi$ be given such that $\sigma \in \mathcal{HRS}$ and $\varphi:\omega^\omega\rightarrow\omega^\omega$ and $\varphi|\mathcal{E}^1_1\subseteq\mathcal{E}_\sigma$. \\Using Lemma \ref{L:redpairsborela11e11}, find $\psi:\omega^\omega\rightarrow\omega^\omega$ reducing $\mathcal{A}_\sigma$ to $\mathcal{E}^1_1$ and mapping $\mathcal{E}_\sigma$ into $\mathcal{A}_1^1$. \\Note that $\varphi\star\psi$\footnote{For all $\varphi, \psi:\omega^\omega\rightarrow \omega^\omega$, also $\varphi\ast\psi:\omega^\omega\rightarrow \omega^\omega$ and, for all $\alpha$, $\varphi\ast\psi|\alpha = \varphi|(\psi|\alpha)$, see Subsubsection \ref{SSS:continuousfunctions}.} maps $\mathcal{A}_\sigma$ into $\mathcal{E}_\sigma$. \\Applying the Borel Hierarchy Theorem, Theorem \ref{T:borelhiertheorem}, \\find $\beta$ in $\mathcal{E}_\sigma$ such that $(\varphi\star\psi)|\beta \in \mathcal{E}_\sigma$. \\Define $\alpha:=\psi|\beta$ and note: 
$\alpha \in \mathcal{A}_1^1$ and $\varphi|\alpha \in \mathcal{E}_\sigma$. 

\smallskip (ii) Let $\mathcal{X}$ in $\mathfrak{Borel}$ be given such that $\mathcal{E}^1_1\subseteq \mathcal{X}$. 
\\ Find $\sigma$ in $\mathcal{HRS}$ and $\varphi:\omega^\omega\rightarrow\omega^\omega$ reducing $\mathcal{X}$ to $\mathcal{E}_\sigma$. Note $\varphi|\mathcal{E}^1_1\subseteq \mathcal{E}_\sigma$. \\Applying (i), find $\alpha$ in $\mathcal{A}^1_1$ such that $\varphi|\alpha \in \mathcal{E}_\sigma$ and, therefore, $\alpha\in \mathcal{X}$. 

\smallskip (iii)  Let $\sigma,\varphi$ be given such that $\sigma \in \mathcal{HRS}$ and $\varphi:\omega^\omega\rightarrow\omega^\omega$ and $\varphi|\mathcal{A}^1_1\subseteq\mathcal{E}_\sigma$. \\Using Lemma \ref{L:redpairsborela11e11}, find $\psi:\omega^\omega\rightarrow\omega^\omega$ reducing $\mathcal{E}_\sigma$ to $\mathcal{E}^1_1$ and mapping $\mathcal{A}_\sigma$ into $\mathcal{A}_1^1$. \\Note that $\varphi\star\psi$ maps $\mathcal{A}_\sigma$ into $\mathcal{E}_\sigma$. \\Applying the Borel Hierarchy Theorem, Theorem \ref{T:borelhiertheorem}, \\find $\beta$ in $\mathcal{E}_\sigma$ such that $(\varphi\star\psi)|\beta \in \mathcal{E}_\sigma$. \\Define $\alpha:=\psi|\beta$ and note: 
$\alpha \in \mathcal{E}_1^1$ and $\varphi|\alpha \in \mathcal{E}_\sigma$. 

\smallskip (ii) Let $\mathcal{X}$ in $\mathfrak{Borel}$ be given such that $\mathcal{A}^1_1\subseteq \mathcal{X}$. 
\\ Find $\sigma$ in $\mathcal{HRS}$ and $\varphi:\omega^\omega\rightarrow\omega^\omega$ reducing $\mathcal{X}$ to $\mathcal{E}_\sigma$. Note $\varphi|\mathcal{A}^1_1\subseteq \mathcal{E}_\sigma$. \\Applying (iii), find $\alpha$ in $\mathcal{E}^1_1$ such that $\varphi|\alpha \in \mathcal{E}_\sigma$ and, therefore, $\alpha\in \mathcal{X}$.
\end{proof}

\subsection{Other results showing that $\mathcal{E}_1^1$ and $\mathcal{A}^1_1$ are not (positively) Borel}\hfill

   $\mathcal{MONPATH}:=\{\alpha\mid\exists\gamma\in\mathcal{F}_\alpha\forall n[\gamma(n)\le\gamma(n+1)\le 1]\}$ is what might be called  a \textit{simple} $\mathbf{\Sigma}^1_1$ set, as, from a classical point of view, $\mathcal{MONPATH}$ is $\mathbf{\Pi}^0_1$.  \\The assumption that $\mathcal{MONPATH}$ is (positively) Borel leads to a contradiction,  see  \cite[Theorem 2.23(vi)]{veldman05}. \\It follows that $\mathcal{E}^1_1$ is not positively Borel, but the statement of Theorem \ref{T:e11nonborel}(ii) is a stronger conclusion.
   
     As we mentioned in Subsection \ref{SS:finite},   \\$\mathcal{ALMOST^\ast FIN}:=\{\alpha\mid\forall \zeta\in[\omega]^\omega\exists n[\alpha\circ\zeta(n)=0]\}$ is  $\mathbf{\Pi}^1_1$ but \textit{not} (positively)  Borel. $\mathcal{ALMOST^\ast FIN}$ might be called a \textit{simple} $\mathbf{\Pi}^1_1$ set, as, from a classical point of view, $\mathcal{ALMOST^\ast FIN}$ is $\mathbf{\Sigma}^0_2$.  \\It follows that also   $\mathcal{A}_1^1$ is \textit{not} (positively) Borel, but the statement of Theorem \ref{T:e11nonborel}(iv) is a stronger conclusion.
   
   \smallskip 
     
     As one might expect, the results about $\mathcal{MONPATH}$ and   $\mathcal{ALMOST^\ast FIN}$ strongly use Brouwer's Continuity Principle $\mathbf{BCP}$. 
\subsection{One half of Souslin's Theorem}\hfill

\begin{theorem}\label{T:souslin}

\hfill

\begin{enumerate}[\upshape (i)]
\item For every $\sigma$ in $\mathcal{STP}$, $\{\alpha\mid\alpha\le^\ast\sigma\}\in \mathfrak{Borel}$. 
\item  Every $\mathcal{X}\subseteq \mathcal{N}$ that is both \emph{strictly} analytic and co-analytic is (positively) Borel: $\mathbf{\Sigma}^{1\ast}_1\cap\mathbf{\Pi}^1_1\subseteq\mathfrak{Borel}$. 

\end{enumerate}
\end{theorem}

\begin{proof}

 (i) Note:  $\forall \alpha[\alpha\le^\ast 1^\ast\leftrightarrow \alpha(0)\neq 0]$. \\Also note:  for all $\sigma\neq 1^\ast$  in $\mathcal{STP}$, $\forall \alpha[\alpha\le^\ast\sigma\leftrightarrow \forall m\exists n[\alpha^m\le^\ast\sigma^n]]$. \\Now use induction on $\mathcal{STP}$. 
 
 \smallskip
(ii) Assume: $\mathcal{X}\in \mathbf{\Sigma}^{1\ast}_1\cap\mathbf{\Pi}^1_1$. \\If $\mathcal{X} =\emptyset$, clearly $\mathcal{X}\in\mathfrak{Borel}$. \\Assume $\mathcal{X}$ is inhabited. Find $\varphi:\omega^\omega\rightarrow\omega^\omega$ such that $\mathcal{X}=\varphi|\omega^\omega$. 
\\Find $\psi:\omega^\omega\rightarrow\omega^\omega$ reducing $\mathcal{X}$ to $\mathcal{A}^1_1$.   \\Using 
Theorem \ref{T:cdiag}(ii), find $\beta$ in $\mathcal{A}^1_1$ such that $\forall\alpha[(\psi\star\varphi)(\alpha)\le^\ast\beta]$.
\\Note that $D_\beta$ is a bar in $\omega^\omega$.
 \\Using  Brouwer's Thesis on bars in $\omega^\omega$ $\mathbf{BT}$, see Subsubsection \ref{SSS:barinduction}, \\find a stump $\sigma$ such that $D_\beta \cap
 T_\sigma$ is bar in $\omega^\omega$.  \\Conclude: $\forall \alpha[\alpha\le^\ast \beta\rightarrow \alpha\le^\ast \sigma]$. \\Conclude, using (i): $\mathcal{X}=\{\gamma\mid  \psi|\gamma\le^\ast\sigma\}\in\mathfrak{Borel}$. \end{proof}

Theorem  \ref{T:souslin}(ii) is of limited application as every $\mathbf{\Pi}^1_1$ subset
of $\omega^\omega$ is perhapsive, see Theorem \ref{T:a11perhapsive}(i), and ``most'' positively Borel sets are not. Therefore, there are not ``many'' positively Borel sets that are both co-analytic and strictly
analytic. The converse of Theorem  \ref{T:souslin}(ii), although classically a well-known fact,  is far from true.

\section{Countable and almost-countable spreads}

\subsection{Countable spreads}\label{SS:countablespread}
Countable closed subsets of the set of the real numbers were among the first objects studied by Cantor. Ome might say that this study led him to discover set theory.

In our constructive context we study {\it located} and closed subsets of $\omega^\omega$, i.e. {\it spreads},  and ask ourselves what could be a useful notion of countability. 

\begin{definition}\label{D:countablespreads}

For each $\delta $, we define $En_\delta=\{\delta^n\mid n\in \omega\}$.
 
We also define: $\mathcal{COUNT}:=\{\beta\mid  Spr(\beta)\;\wedge\;\exists\delta[\mathcal{F}_\beta\subseteq En_\delta]\}.$
\end{definition}

$En_\delta$ is called the subset of $\omega^\omega$ {\it enumerated by $\delta$}, see Subsubsection \ref{SSS:infseq}. 

If $\beta\in \mathcal{COUNT}$, we call $\mathcal{F}_\beta$ an {\it (at most) countable spread}.

\begin{definition} $\mathcal{X}\subseteq \mathcal{N}$ is called \emph{discrete} if and only if $\forall \alpha \in \mathcal{X}\forall \beta \in \mathcal{X}[\alpha \;\#\;\beta \vee \alpha = \beta]$.  \end{definition}

Recall that $\mathcal{FIN}$ is the set of all $\alpha$ such that $\exists n\forall m\ge n[\alpha(m)=0]$, i.e.\\$D_\alpha:=\{m\mid \alpha(m)\neq 0\}$ is a {\it finite} subset of $\omega$. 

\smallskip Like Theorem \ref{T:uncsigma11complete}, the following Theorem \ref{T:countspread}  should be compared to a classical result due to W.~Hurewicz, see \cite[Theorem 27.5]{kechris}.  
\begin{theorem}\label{T:countspread}\hfill
\begin{enumerate}[\upshape (i)]
\item \emph{For every spread $\mathcal{F}\subseteq \mathcal{N}$: $\mathcal{F}$ is (at most) countable if and only if $\mathcal{F}$ is discrete:}\\$\forall \beta[\beta \in \mathcal{COUNT}\leftrightarrow (Spr(\beta)\;\wedge\;\forall \gamma_0 \in \mathcal{F}_\beta\forall\gamma_1\in\mathcal{F}_\beta[\gamma_0\;\#\;\gamma_1 \;\vee\; \gamma_0 =\gamma_1])]$.
\item $\mathcal{FIN}\preceq \mathcal{COUNT}$. \item $\mathcal{A}_1^1\preceq\mathcal{COUNT}$. \item  \emph{$\mathcal{COUNT}$ is not the co-projection of a closed subset of $\mathcal{N}$ but it is the co-projection of a (positively) Borel subset of $\mathcal{N}$:} $\mathcal{COUNT}$ is not $\mathbf{\Pi}^1_1$ but  $\mathcal{COUNT}$ is  $
\mathbf{\Pi}^{1+}_1$.
\end{enumerate}
\end{theorem}

\begin{proof} (i) Assume $\beta \in \mathcal{COUNT}$, i.e. $Spr(\beta)$ and $\mathcal{F}_\beta$ is (at most) countable.  \\Note:  if $\beta(0)\neq 0$ then  $\mathcal{F}_\beta=\emptyset$ is discrete. \\Assume: $\beta(0)=0$ and find $\delta$ such that $\mathcal{F}_\beta\subseteq En_\delta$. \\ Then $\forall \gamma \in \mathcal{F}_\beta\exists n[\gamma =\delta^n]$. \\Let $\gamma_0,\gamma_1$ in $\mathcal{F}_\beta$ be given.  Using Brouwer's Continuity Principle $\mathbf{BCP}$, see Subsubsection \ref{SSS:bcpcontchoice},   find $n_0,m_0$ such that   $\forall \gamma \in\mathcal{F}_\beta[\overline{\gamma_0}m_0\sqsubset \gamma\rightarrow \gamma =\delta^{n_0}]$. and find $n_1, m_1$ such that $\forall \gamma \in\mathcal{F}_\beta[\overline{\gamma_1}m_1\sqsubset \gamma\rightarrow \gamma =\delta^{n_1}]$.\\ Note: \textit{if} $\overline{\gamma_0}m_0\perp \overline{\gamma_1}m_1$, then $\gamma_0\;\#\;\gamma_1$, and, \textit{if not}, then $\gamma_0=\delta^{n_0}=\gamma_1$. \\We thus see: $\forall \gamma_0\in \mathcal{F}_\beta\forall\gamma_1\in\mathcal{F}_\beta[\gamma_0\;\#\;\gamma_1\;\vee\;\gamma_0=\gamma_1]$, i.e. $\mathcal{F}_\beta$ is discrete.

\smallskip Now assume  $Spr(\beta)$ and $\mathcal{F}_\beta$ is discrete. \\ We  may assume: $\beta(0)=0$, i.e. $\mathcal{F}_\beta$ is inhabited. \\Define $\varepsilon$ such that, for all $s$, $\varepsilon(s)=0$ if and only if $\beta(s_I)=\beta(s_{II})=0$. \\Note: $Spr(\varepsilon)$ and for all $\gamma$, $ \gamma\in \mathcal{F}_\varepsilon$ if and only if both $\gamma_I$ and $\gamma_{II}$ are in $\mathcal{F}_\beta$.  
\\Conclude: $\forall\gamma \in \mathcal{F}_\varepsilon[\gamma_I\;\#\;\gamma_{II}\;\vee\;\gamma_I=\gamma_{II}]$. \\Using the First Axiom of Continuous Choice $\mathbf{AC}_{1,0}$, see Subsubsection \ref{SSS:bcpcontchoice}, \\find $\varphi:\mathcal{F}_\varepsilon \rightarrow \omega$ such that  $\forall \gamma\in \mathcal{F}_\varepsilon[\bigl(\varphi(\gamma)= 0\rightarrow \gamma_I\;\#\;\gamma_{II}\bigr)\;\wedge\;\bigl(\varphi(\gamma)>0\rightarrow \gamma_I=\gamma_{II}\bigr)]$. 
\\Note: $\forall \gamma \in \mathcal{F}_\beta[\varphi(\ulcorner \gamma,\gamma\urcorner)>0]$ and, for all $n$, if $\beta(n)=0$ and $\varphi|\ulcorner n, n\urcorner\perp \langle 0 \rangle$, then \\there exists exactly one $\gamma \in \mathcal{F}_\beta$ such that $n\sqsubset \gamma$. 
\\Find $\delta$ such that, for each $n$, if $\beta(n)=0$ and  $\varphi|\ulcorner n, n\urcorner\perp \langle 0 \rangle$, then $n\sqsubset \delta^n$ and $\delta^n\in\mathcal{F}_\beta$, and note: $\mathcal{F}_\beta\subseteq En_\delta$. 
\\We thus see: $\mathcal{F}_\beta$ is (at most) countable. 

\smallskip (ii) Define $\varphi:\omega^\omega\rightarrow\omega^\omega$ such that $\forall \alpha\forall s[(\varphi|\alpha)(s) =0 \leftrightarrow \exists m\exists k[s=\overline\alpha m\ast\overline{\underline 0}k]]$. \\ We shall prove that $\varphi$ reduces $\mathcal{FIN}$ to $\mathcal{COUNT}$.
\\Note: for every $\alpha$, $Spr(\varphi|\alpha)$ and $\alpha \in \mathcal{F}_{\varphi|\alpha}$.  

First, let $\alpha$ in $\mathcal{FIN}$ be given. Find  $m$  such that $\forall n\ge m[\alpha(n)=0]$. \\Note: $\forall\gamma[\gamma \in \mathcal{F}_{\varphi|\alpha}\leftrightarrow\exists k\le m[\gamma = \overline \alpha k\ast\underline 0]]$. \\Define $\delta$ such that  $\forall k\le m[\delta^{k} = \overline \alpha k\ast\underline 0]$.  \\Note: $\mathcal{F}_{\varphi|\alpha}\subseteq En_\delta$ and: $\varphi|\alpha\in \mathcal{COUNT}$.\\Clearly, for every $\alpha$, if $\alpha\in \mathcal{FIN}$, then $\varphi|\alpha \in \mathcal{COUNT}$. 

 Now let $\alpha$ be given such that $\varphi|\alpha \in \mathcal{COUNT}$.\\According to (i): $\mathcal{F}_{\varphi|\alpha}$ is discrete. \\Note: $\alpha \in \mathcal{F}_{\varphi|\alpha}$.  Using Brouwer's Continuity Principle $\mathbf{BCP}$, see Subsubsection \ref{SSS:bcpcontchoice}, find $m$ such that $\forall \gamma\in \mathcal{F}_{\varphi|\alpha}[\overline{\alpha}m\sqsubset \gamma\rightarrow \alpha=\gamma]$. \\Conclude:   $\forall n\ge m[\alpha(n)=0]$ and: $\alpha \in \mathcal{FIN}$. \\Clearly, for every $\alpha$, if $\varphi|\alpha\in \mathcal{COUNT}$, then $\alpha \in \mathcal{FIN}$. 
 
 We thus see that $\varphi$ reduces $\mathcal{FIN}$ to $\mathcal{COUNT}$. 

\smallskip (iii)  Recall that we defined, for each $\alpha$, $T_\alpha = \{t \mid \forall u\sqsubset t[\alpha(u)=0]\}$.

Define $\varphi:\omega^\omega\rightarrow\omega^\omega$ such that, for all $\alpha$, for all $s$,  \\$(\varphi|\alpha)(s) =0$ if and only if  $\exists t\in T_\alpha\exists k[s=t\ast\overline{\underline 0}k]$.

We shall prove that $\varphi$ reduces $\mathcal{A}^1_1$ to $\mathcal{COUNT}$. 
\\Note: $\forall \alpha[Spr(\varphi|\alpha)]$. 

 First, assume: $\alpha\in\mathcal{A}^1_1$. Let $\gamma_0,\gamma_1$ in $\mathcal{F}_{\varphi|\alpha}$ be given.\\ Find $n_0:=\mu n[\alpha(\overline{\gamma_0} n)\neq 0]$ and  $n_1:=\mu n[\alpha(\overline{\gamma_1}n)\neq 0]$.\\Note: $\overline{\gamma_0}n_0 \in T_\alpha$ and $\overline{\gamma_0}(n_0+1)\notin T_\alpha$ and $\gamma_0=\overline{ \gamma_0} n_0\ast\underline 0$. \\Similarly,   $\gamma_1=\overline{ \gamma_1} n_1\ast\underline 0$. \\\textit{If} $\overline{\gamma_0}n_0 \perp \overline{\gamma_1}n_1$, then $\gamma_0\;\#\;\gamma_1$ and, \textit{if not}, then $\gamma_0=\gamma_1$. \\We thus see: $\forall \gamma_0\in\mathcal{F}_{\varphi|\alpha}\forall \gamma_1\in\mathcal{F}_{\varphi|\alpha}[\gamma_0\;\#\;\gamma_1\;\vee\;\gamma_0=\gamma_1]$, i.e. $\mathcal{F}_{\varphi|\alpha}$ is discrete. \\
  Using (i), conclude: $\varphi|\alpha \in \mathcal{COUNT}$. 
\\Clearly, for each $\alpha$, if $\alpha \in \mathcal{A}^1_1$, then $\varphi|\alpha \in \mathcal{COUNT}$.

Now let $\alpha$ be given such that $\varphi|\alpha\in \mathcal{COUNT}$.   \\ Let $\gamma$ be given. Define $\gamma^\ast$ such that, \\for each $n$, \textit{if} $\overline \gamma(n+1)\in T_\alpha$,  then $\gamma^\ast(n)=\gamma(n)$, and, \textit{if not}, then $\gamma^\ast(n)=0$. \\Note: $\gamma^\ast \in \mathcal{F}_{\varphi|\alpha}$. \\According to (i), $\mathcal{F}_{\varphi|\alpha}$ is discrete. Using Brouwer's Continuity Principle $\mathbf{BCP}$, \\find $n$ such that $\forall \delta \in \mathcal{F}_{\varphi|\alpha}[\overline{\gamma^\ast}n\sqsubset\delta\rightarrow\gamma^\ast =\delta]$. \\  Suppose: $\forall m\le n[\alpha(\overline{\gamma^\ast}m)=0]$. Then $\forall p[\overline{\gamma^\ast}n\ast\langle p\rangle \in T_\alpha$ and $(\varphi|\alpha)(\overline{\gamma^\ast}n\ast\langle p\rangle)=0]$. \\Conclude: $\exists m\le n[\alpha(\overline{\gamma^\ast}m)\neq 0]$, and: $\exists m\le n[\alpha(\overline{\gamma}m)\neq 0]$. \\We thus see: $\forall \gamma\exists m[\alpha(\overline{\gamma}m)\neq 0]$, i.e. $\alpha\in \mathcal{A}^1_1$.
\\Clearly, for each $\alpha$, if $\varphi|\alpha \in \mathcal{COUNT}$, then $\alpha \in \mathcal{A}^1_1$.

We thus see that $\varphi$ reduces $\mathcal{A}_1^1$ to $\mathcal{COUNT}$. 

\smallskip (iv) As $\mathcal{FIN}$ reduces to $\mathcal{COUNT}$, see (ii),  and $\mathcal{FIN}$ is not $\mathbf{\Pi}^1_1$, see Theorem \ref{T:finite}(iii),  also $\mathcal{COUNT}$ is not $\mathbf{\Pi}^1_1$. 

 Note, considering the proof of (i): for all $\beta$, $\beta\in\mathcal{COUNT}$ if and only if \\$Spr(\beta)$ and $\mathcal{F}_\beta$ is discrete, i.e. \\$\forall \gamma\in \mathcal{F}_\beta\exists n\forall s\forall t[(\beta(s) =\beta(t)=0\;\wedge\;\overline \gamma n\sqsubseteq s\;\wedge\;\overline\gamma n \sqsubseteq t)\rightarrow (s\sqsubseteq t \;\vee\; t\sqsubseteq s)]$.  
\\Conclude, using the last observation of Subsubsection \ref{SSS:continuousfunctions}: \\for all $\beta$, $\beta \in \mathcal{COUNT}$ if and only if $Spr(\beta)$ and \\$\forall \gamma\exists n\forall s\forall t[(\beta(s) =\beta(t)=0\;\wedge\;\overline \gamma n\sqsubseteq s\;\wedge\;\overline\gamma n \sqsubseteq t)\rightarrow (s\sqsubseteq t \;\vee\; t\sqsubseteq s)]$. 
\\ Let $\mathcal{X}$ be the set of all $\beta$ such that $Spr(\beta_I)$ and  either $\exists n[\beta_I(\overline{\beta_{II}}n)\neq 0]$ or  \\$ \exists n\forall s\forall t[(\beta_I(s) =\beta_I(t)=0\;\wedge\;\overline{\beta_{II}} n\sqsubseteq s\;\wedge\;\overline{\beta_{II}} n \sqsubseteq t)\rightarrow (s\sqsubseteq t \;\vee\; t\sqsubseteq s)]$  and note: $\mathcal{X}\in\mathbf{\Pi}^0_3$ and: $\mathcal{COUNT}=Un(\mathcal{X})$ and: $\mathcal{COUNT}$ is $\mathbf{\Pi}^{1+}_1$. 
\end{proof}
\subsection{Almost-countable spreads}\label{SS:almostcountablespreads}\hfill

One might feel that the notion of a {\it countable spread} as introduced in Subsection \ref{SS:countablespread} is perhaps too strong. We therefore introduce a weaker notion. 

\smallskip

Note: for each $\delta$,  for each $\gamma$, if $\forall n[\gamma \;\#\;\delta^n]$, one may define $\alpha$ such that, for each $n$, $\alpha(n)=\mu(p)[\overline \gamma p \perp \delta^n]$.

Conclude: $\forall n[\gamma \;\#\;\delta^n]$ if and only if $ \exists \alpha\forall n[\overline \gamma \alpha(n)\perp\delta^n]$. \\One may consider $\alpha$ such that $\forall n[\overline \gamma \alpha(n)\perp\delta^n]$ as \textit{evidence} for the fact: $\forall n[\gamma \;\#\;\delta^n]$. 

\begin{definition}\label{D:almostbelongsto} For all $\gamma,\delta$, we define: \\$\gamma$ \emph{almost belongs to} $En_\delta=\{\delta^n\mid n\in \omega\}$ if and only if $\forall \alpha\exists n[\overline \gamma \alpha(n)\perp \delta^n]$ \end{definition}

So $\gamma$ almost belongs to $En_\delta$ if every attempt to give evidence   that $\gamma$ is apart from every element of $En_\delta$ fails in finitely many steps.
\begin{lemma}\label{L:en1}  $\;$  \begin{enumerate}[\upshape (i)] \item For all $\gamma, \delta, \varepsilon$, if $En_\delta\subseteq En_\varepsilon$ and  $\gamma$ almost belongs to $En_\delta$, \\then  $\gamma$ almost belongs to $En_\varepsilon$.  \item For all $\gamma, \delta, \varepsilon$, if $En_\delta= En_\varepsilon$, then  \\ $\gamma$ almost belongs to $En_\delta$ if and only if $\gamma$ almost belongs to $En_\varepsilon$. \end{enumerate}\end{lemma}

\begin{proof} (i) Let $\delta,\varepsilon$ be given such that $En_\delta\subseteq En_\varepsilon$. \\Let $\gamma$ be given such that $\forall \alpha \exists n[\overline\gamma \alpha(n)\sqsubset \delta^n]$. \\Using  the First Axiom of Countable Choice $\mathbf{AC}_{0,0}$, see Subsubsection \ref{SSS:countablechoice}, \\find $\zeta$ such that $\forall n[\delta^n = \varepsilon^{\zeta(n)}]$.  \\Let $\alpha$ be given. \\Find $n$ such that $\overline\gamma\alpha\circ\zeta(n)\sqsubset\delta^n=\varepsilon^{\zeta(n)}$ and conclude: $\exists m[\overline\gamma \alpha(m)\sqsubset\varepsilon^m]$. \\Conclude: $\forall \alpha \exists n[\overline\gamma \alpha(n)\sqsubset \varepsilon^n]$. 

\smallskip (ii) immediately follows from (i). \end{proof}
Define $\delta$ such that $\forall n[\delta^n=n\ast\underline 0]$ and note: $\mathcal{FIN}=\{\delta^n\mid n\in\omega\}=En_\delta$.

Recall: $\mathcal{ALMOST}^\ast\mathcal{FIN}=\{\gamma\mid\forall \zeta \in [\omega]^\omega\exists n[\gamma\circ\zeta(n)=0]\}$, see Definition \ref{D:almostfinite}.

\begin{lemma}\label{L:en2} $\;$\\
For each $\gamma$,  $\gamma \in \mathcal{ALMOST}^\ast\mathcal{FIN}$ if and only if  $\gamma$ almost belongs to $\mathcal{FIN}$. \end{lemma}

\begin{proof} Let $\gamma$ in $\mathcal{ALMOST}^\ast\mathcal{FIN}$ be given.   \\We want to prove that $\gamma$ almost belongs to $\mathcal{FIN}=\{n\ast\underline 0\mid n \in \omega\}$ . \\ 
Let $\alpha$ be given.  We want to prove: $\exists n[\overline \gamma \alpha(n)\sqsubset n\ast\underline 0]$.
\\To this end, we define $\zeta$ in $[\omega]^\omega$, step by step.  \\ \textit{If} $\overline\gamma\alpha(0)\perp\underline 0$, define $\zeta(0)=\mu i <\alpha(m)[\gamma(i)\neq 0]$, and, \textit{if not}, define $\zeta(0)=0$. \\
Now assume $p>0$ and we defined $\zeta(0), \zeta(1), \ldots,\zeta(p-1)$. Define  $m:=\overline\gamma\bigl(\zeta(p-1)+1\bigr)$.  \\
 \textit{If} $\overline\gamma\alpha(m)\perp m\ast\underline 0$, i.e. $\overline\gamma\alpha(\overline\gamma\bigl(\zeta(p-1)+1\bigr))\perp \overline\gamma\bigl(\zeta(p-1)+1\bigr)\ast\underline 0$, \\
 define $\zeta(p)=\mu i <\alpha(m)[i>\zeta(p-1)\;\wedge\gamma(i)\neq 0]$, and, \\\textit{if not}, define $\zeta(p)=\zeta(p-1)+1$. 
 \\Now find $n$ such that $\gamma\circ\zeta(n)=0$ and conclude: for some $p\le n$ we must have seen $ \overline \gamma\alpha(m)\sqsubset m\ast\underline 0$, where $m=\overline\gamma\bigl(\zeta(p-1)+1\bigr)$.\\
 We thus see that $\gamma$ almost belongs to $\mathcal{FIN}$.  

\smallskip Conversely, let $\gamma$ be given such that $\gamma$ almost belongs to $\mathcal{FIN}$, i.e. \\$\forall \alpha\exists n[\overline\gamma\alpha(n)\sqsubset n\ast\underline 0]$. \\Assume: $\zeta \in[\omega]^\omega$. 
\\Find $\eta$ in $[\omega]^\omega$ such that $\forall n[\zeta\circ\eta(n)> length(n)]$. \\
Define $\alpha$ such that, for each $n$, $\alpha(n) = \zeta\circ \eta (n)+1$. 
\\Find $n$ such that $\overline\gamma \alpha(n)\sqsubset n\ast\underline 0$ and conclude: $\gamma\circ\zeta\circ\eta(n)=0$.\\ We thus see: $\forall \zeta \in[\omega]^\omega\exists n[\gamma\circ\zeta(n)=0]$, i.e. $\gamma \in \mathcal{ALMOST^\ast FIN}$. \end{proof}
\begin{definition}\label{D:almost-countablespreads}
For each $\delta$, we  let $\mathcal{ALMOST}^\ast(En_\delta)$ be the set of all $\gamma$ that almost belong to $En_\delta$, i.e. such that $\forall \alpha\exists n [\overline \gamma \alpha(n) \sqsubset \delta^n]$. 

\smallskip
We also define:  $\mathcal{ALMOST^\ast COUNT}:=\{\beta\mid  Spr(\beta)\;\wedge\;\exists\delta[\mathcal{F}_\beta\subseteq \mathcal{ALMOST}^\ast({En_\delta})]\}$.

\end{definition}

If $\beta\in \mathcal{ALMOST^\ast COUNT}$, we call $\mathcal{F}_\beta$ an {\it almost-countable spread}.

\begin{lemma}\label{L:almostcountablespread} For each $\beta$, if $\mathcal{F}_\beta$ is an almost-countable spread, then there exists $\varepsilon$ in $(\mathcal{F}_\beta)^\omega$ such that $\mathcal{F}_\beta \subseteq \mathcal{ALMOST}^\ast(En_\varepsilon)$. \end{lemma}
                             
 \begin{proof} Let $\beta, \delta$ be given such that $Spr(\beta$ and  $\mathcal{F}_\beta \subseteq \mathcal{ALMOST}^\ast(En_\delta)$.  \\Assume: $\beta(0)=0$ and let $\rho$ be the retraction of $\omega^\omega$ onto $\mathcal{F}_\beta$. \\Define $\varepsilon$ such that $\forall n[\varepsilon^n = \rho|\delta^n]$ and note: $\forall n[\varepsilon^n \in \mathcal{F}_\beta]$. 
 \\ We now prove: $\mathcal{F}_\beta \subseteq \mathcal{ALMOST}^\ast(En_\varepsilon)$.
 \\Assume $\gamma \in \mathcal{F}_\beta$ and let $\alpha$ be given.  Find $n$ such that $\overline\gamma \alpha(n)\sqsubset \delta^n$. \\Conclude: 	$\beta\bigl(\overline\gamma\alpha(n)\bigr)=0$ and: $\overline\gamma \alpha(n)\sqsubset \varepsilon^n$.
 \\ We thus see:  $\mathcal{F}_\beta \subseteq \mathcal{ALMOST}^\ast(En_\varepsilon)$.\end{proof}

                                                                                                                                                                                                                                                                                                                                                                                                                                                                                                                                                                                                                                         \begin{lemma}\label{L:imagealmostcountable} If $\mathcal{F}, \mathcal{H}$ are spreads and $\mathcal{F}$ maps onto $\mathcal{H}$ and $\mathcal{F}$ is almost-countable, also $\mathcal{H}$ is almost-countable. \end{lemma} 
                                                                                                                                                                                                                                                                                                                                                                                                                                                                                                                                                                                                                                         
\begin{proof}  Let $\beta_0, \beta_1$ be spread-laws such that $\mathcal{F}_{\beta_0}$ is at-most-countable. \\Assume $\varphi:\mathcal{F}_{\beta_0}\rightarrow \mathcal{F}_{\beta_1}$ is surjective. 
\\Find $\delta$ such that $\mathcal{F}_{\beta_0}\subseteq Almost^\ast(En_\delta)$. Define $\varepsilon$ such that $\forall n[\varepsilon^n=\varphi|\delta^n]$. 
\\Assume: $\zeta \in \mathcal{F}_{\beta_1}$ and find $\gamma$ in $\mathcal{F}_{\beta_0}$ such that $\varphi|\gamma=\zeta$. \\Let $\alpha$ be given. Find $\eta$ such that $\forall n[\overline{ \varepsilon^n}\alpha (n)\sqsubseteq \varphi|\overline{\delta^n}\eta(n)]$. \\Find $n$ such that $\overline \gamma \eta(n)\sqsubset \delta^n$ and conclude: $\overline\zeta \alpha(n)=\overline{\varphi|\gamma}\alpha(n)\sqsubset \varphi|\delta^n=\varepsilon^n$. \\We thus see: $\forall \zeta \in \mathcal{F}_{\beta_1}\forall \alpha \exists n[\overline{\zeta}\alpha(n)\sqsubset \varepsilon^n]$, i.e. $\mathcal{F}_{\beta^1}\subseteq Almost^\ast(En_\varepsilon)$ and: \\$\mathcal{F}_{\beta_1}$ is almost-countable. 
\end{proof}

      \begin{theorem}\label{T:countandfin}\begin{enumerate}[\upshape (i)] \item For each $\beta$ such that $Spr(\beta)$, \\$\mathcal{F}_\beta$ is a countable spread if and only if  $\mathcal{F}_\beta$ embeds into $\mathcal{FIN}$. \item                                     For each $\beta$ such that $Spr(\beta)$,  \\ if $\mathcal{F}_\beta$ is an almost-countable spread, then  $\mathcal{F}_\beta$ embeds into $\mathcal{ALMOST^\ast FIN}$. \end{enumerate}\end{theorem}
      
      \begin{proof}
      (i) Assume: $Spr(\beta)$ and $\mathcal{F}_\beta$ is an inhabited countable spread. \\Find $\delta$ in $(\mathcal{F}_\beta)^\omega$ such that $\mathcal{F}_\beta =En_\delta$, i.e. $\forall \gamma\in\mathcal{F}_\beta\exists n[\gamma=\delta^n]$. 
      \\Using the First Axiom of Continuous Choice $\mathbf{AC}_{1,0}$, see Subsubsection \ref{SSS:bcpcontchoice}, \\find $\varphi:\mathcal{F}_\beta\rightarrow\omega$ such that $\forall \gamma\in\mathcal{F}_\beta[\gamma=\delta^{\varphi(\gamma)}]$.  
      \\Define $\psi:\mathcal{F}_\beta\rightarrow \omega^\omega$ such that $\forall \gamma \in\mathcal{F}_\beta[\psi|\gamma=\overline{\underline 1}\varphi(\gamma)\ast\underline 0]$ and note: $\psi:\mathcal{F}_\beta\rightarrowtail\mathcal{FIN}$.
      
  \smallskip    Conversely, assume: $Spr(\beta)$ and: $\mathcal{F}_\beta$ embeds into $\mathcal{FIN}$. \\Find $\varphi$ such that $\varphi:\mathcal{F}_\beta\rightarrowtail \mathcal{FIN}$.\\ Note: $\mathcal{FIN}$ is discrete, i.e. for all $\delta_0, \delta_1$ in $\mathcal{FIN}$, either $\delta_0=\delta_1$ or $\delta_0\;\#\;\delta_1$. \\ Conclude: for all $\gamma_0, \gamma_1$ in $\mathcal{F}_\beta$, either $\varphi|\gamma_0 = \varphi|\gamma_1$ or $\varphi|\gamma_0\;\#\;\varphi|\gamma_1$, and, \\
  
  therefore,  either $\gamma_0=\gamma_1$ or $\gamma_0\;\#\;\gamma_1$, i.e. $\mathcal{F}_\beta$ is discrete.\\Using Theorem \ref{T:countspread}(i), conclude:   $\mathcal{F}_\beta$ is a countable spread.
      
      \smallskip (ii) Assume: $Spr(\beta)$ and $\mathcal{F}_\beta$ is an inhabited almost-countable spread. 
      \\Using Lemma \ref{L:almostcountablespread}, find $\delta$ in $(\mathcal{F}_\beta)^\omega$ such that $\mathcal{F}_\beta =\mathcal{ALMOST}^\ast(En_\delta)$. 
      
      We first prove the following observation: \\for all $s$ such that $\beta(s)=0$ there exists $n$ such that $s\sqsubset\delta^n$.
\\Let $s$ be given such that $\beta(s)=0$. Find $\gamma$ in $\mathcal{F}_\beta$ such that $s\sqsubset\gamma$. \\
Then find $n$ such that $\overline\gamma length(s)\sqsubset\delta^n$ and conclude: $s\sqsubset\delta^n$. 

\smallskip
      
      Now define $\varphi:\mathcal{F}_\beta\rightarrow\omega^\omega$ such that, for all $ \gamma$ in $\mathcal{F}_\beta$, for all $n$, \\\textit{if} $\mu p[\overline \gamma n\sqsubset\delta^p]<\mu p[\overline \gamma(n+1)\sqsubset\delta^p]]$, then $(\varphi|\gamma)(n)= \mu p[\overline \gamma(n+1)\sqsubset\delta^p]$, and, \\\textit{if $\mu p[\overline \gamma n\sqsubset\delta^p]=\mu p[\overline \gamma(n+1)\sqsubset\delta^p]]$}, then $(\varphi|\gamma)(n)=0$. 
      
      We prove that $\varphi$ is a strongly injective function from $\mathcal{F}_\beta$ into $\omega^\omega$. \\
      Let $\gamma_0,\gamma_1$ in $\mathcal{F}_\beta$ be given such that $\gamma_0\;\#\;\gamma_1$. Find $n$ such that $\overline{\gamma_0}n\neq\overline{\gamma_1}n$. \\
      Note: $\mu p[\overline{\gamma_0}n\sqsubset\delta^p]\neq \mu p[\overline{\gamma_1}n\sqsubset \delta^p]$.
      \\ Conclude: $\exists i\le n[(\varphi|\gamma_0)(i)\neq(\varphi|\gamma_1)(i)]$ and: $\varphi|\gamma_0\;\#\;\varphi|\gamma_1$. 
      
      We  prove that $\varphi$ maps $\mathcal{F}_\beta$ into $\mathcal{ALMOST^\ast FIN}$.\\
      Let $\gamma$ in $\mathcal{F}_\beta$ be given and consider $\varphi|\gamma$. Let $\zeta$ in $[\omega]^\omega$ be given. 
   \\ Find $n$ such that $\overline\gamma\bigl(\zeta(n)+1\bigr)\sqsubset \delta^n$.   Assume: $\forall i\le n [(\varphi|\gamma)\bigl(\zeta(i)\bigr)\neq 0]$. \\
   Conclude: $\forall i< n[0<(\varphi|\gamma)\bigl(\zeta(i)\bigr)<(\varphi|\gamma)\bigl(\zeta(i+1)\bigr)]$, and: $(\varphi|\gamma)(\zeta(n)\bigr)
\ge n+1$. Conclude: $\mu p[\overline \gamma\bigl( \zeta(n)+1\bigr)\sqsubset\delta^p]\ge n+1$ and also: $\overline\gamma\bigl(\zeta(n)+1\bigr)\sqsubset\delta^n$. Contradiction. Conclude: $\exists i \le n[(\varphi|\gamma)\bigl(\zeta(i)\bigr)=0]$. \\ 
Clearly,  $\forall \zeta\in[\omega]^\omega\exists i[(\varphi|\gamma)\bigl(\zeta(i)\bigr)=0]$, i.e.  $\varphi|\gamma\in\mathcal{ALMOST^\ast FIN}$. 
 \end{proof}                                                                                                                                                                                                                                                                                                                                                                                                                                                                                                                                                                                                                                        
   We did not succeed in proving the converse of Theorem \ref{T:countandfin}(ii).

       \subsection{Cantor-Bendixson sets} 
\begin{definition} Let $\varepsilon, \beta$ be given. We define $\nu =\mathsf{CB}(\varepsilon, \beta)$ in $2^\omega$ as follows.

For each $s$, $\nu(s)=0$ if and only if either $s\sqsubset \varepsilon$ or \\there exist $m,n,t$ such that $s=\overline \varepsilon m\ast \langle n \rangle\ast t$ and $ \varepsilon(m)\neq n$ and    $\beta^{J(m,n)}(t)=0$.
\end{definition} 

\begin{lemma}\label{L:cb} Let $\varepsilon, \beta$ be given. \begin{enumerate}[\upshape (i)] \item If, for all $n$, $Spr(\beta^n)$, then $Spr\bigl(\mathsf{CB}(\varepsilon, \beta)\bigr)$. \item If, for all $n$, $\beta^n \in \mathcal{ALMOST^\ast COUNT}$, then $\mathsf{CB}(\varepsilon, \beta) \in \mathcal{ALMOST^\ast COUNT}$. \end{enumerate} \end{lemma}

\begin{proof} (i) The proof is straightforward and left to the reader. \\If, for all $n$, $Spr(\beta^n)$, and  $\nu = \mathsf{CB}(\varepsilon, \beta)$, we call $\varepsilon$ the \emph{spine} of the spread $\mathcal{F}_\nu$. 

\smallskip

(ii) Assume: for all $n$, $\beta^n \in \mathcal{ALMOST^\ast COUNT}$. 
\\Using the Second Axiom of Countable Choice $\mathbf{AC}_{0,1}$, see Subsubsection \ref{SSS:countablechoice},  \\find $\delta$ such that, for all $n$,  $\mathcal{F}_{\beta^n}\subseteq\mathcal{ALMOST^\ast}(En_{\delta^n})$. \\Define $\eta$ such that  $\eta^0=\varepsilon$ and, for all $m, n, p$, \\if $\varepsilon(m)\neq n$, then   $\eta^{J(J(m,n),p)+1}=\overline \varepsilon m\ast\langle n \rangle \ast\delta^{J(m,n),p}$. 
 
\smallskip Define $\nu:=\mathsf{CB}(\varepsilon, \beta)$. We prove that $\mathcal{F}_\nu$ is a subset of $\mathcal{ALMOST}^\ast(En_\eta)$. 
\\Asume: $\gamma \in \mathcal{F}_\nu$. Let $\alpha$ be given. We want to prove: $\exists n[\overline \gamma\alpha(n)\sqsubset \eta^n]$. 
\\If $\overline \gamma \alpha(0)\sqsubset \eta^o=\varepsilon$, we are done. 
\\Now assume: $\overline \gamma\alpha(0) \perp \eta^0=\varepsilon$.  Define  $m:=\mu p[\gamma(p)\neq\varepsilon(p)]$ and $n:= \gamma(m)$.\\  
Define $k:=J(m,n)$ and $s:=\overline \varepsilon m\ast\langle n \rangle$.  \\
Note: $s\sqsubset \gamma$ and find $\mu$ such that $\gamma=s\ast\mu$. \\Note:  $\mu \in \mathcal{F}_{\beta^{k}}$.
\\Find $p$ such that $\overline \mu \bigl(\alpha(J(k,p)+1)\bigr) \sqsubset \delta^{k, p}$. \\
 Conclude: $\overline{\gamma} \alpha\bigl(J(k,p)+1\bigr)\sqsubset s\ast \overline{\mu} \alpha(\bigl (J(k,p)+1\bigr)\sqsubset  s\ast \delta^{k, p}=\eta^{J(k,p)+1}$. 
\\Conclude: $\forall \alpha \exists n[\overline \gamma \alpha( n)\sqsubset \eta^n]$, that is: $\gamma \in \mathcal{ALMOST}^\ast(En_\eta)$. \\We thus see: $\mathcal{F}_\nu \subseteq \mathcal{ALMOST^\ast}(En_\eta)$ and: $\nu=\mathsf{CB}(\varepsilon, \beta)\in\mathcal{ALMOST^\ast COUNT}$. 
 \end{proof}
 \begin{definition}\label{D:cantorbendixsonsets}We introduce a subset $\mathcal{CB}$ of $\omega^\omega$ by means of the following inductive definition. 

\begin{enumerate}[\upshape (i)] \item For all $\beta$, if $Spr(\beta)$ and $\beta(0)\neq 0$, (so $\mathcal{F}_\beta=\emptyset$),  then $\beta \in \mathcal{CB}$,  and, \item  for all $\varepsilon$, for all $\beta$, if for all $n$, $\beta^n \in \mathcal{CB}$, then $\mathsf{CB}(\varepsilon, \beta) \in \mathcal{CB}$, and, \item all members of $\mathcal{CB}$ are given by (i), (ii). \end{enumerate}\end{definition}

The following theorem may be compared to Cantor's result  \cite[Theorem C]{cantor1} in \cite[page 220]{cantor2}, and to a related intuitionistic result: \cite[Theorems 9.1 and  9.2]{veldman18}.
\begin{theorem}\label{T:cbcount}  $\mathcal{ALMOST^\ast COUNT}=\mathcal{CB}$. \end{theorem}

\begin{proof} 

Using Lemma \ref{L:cb} and induction, we conclude: $\mathcal{CB}\subseteq \mathcal{ALMOST^\ast COUNT}$.

\smallskip We now prove  that $\mathcal{ALMOST^\ast COUNT}$ is a subset of $\mathcal{CB}$.

Let $\beta$ in $\mathcal{ALMOST^\ast COUNT}$ be given. One may assume: $\beta(0)=0$.  \\ Using Lemma \ref{L:almostcountablespread},  find $\delta$ in $                   (\mathcal{F}_\beta)^\omega$                                                                                                                                                                                                                                                                      such that $\mathcal{F}_\beta\subseteq \mathcal{ALMOST}^\ast(En_\delta)$. \\ Now define $\beta^+$ in $2^\omega$ such that, for all $c$, $\beta^+(c)=0$ if and only if \\$\forall i<length(c)[\beta\bigl(c(i)\bigr)=0]\;\wedge\;\bigl(i+1<length(c)\rightarrow c(i)\sqsubset c(i+1)\bigr)]$.\\ Note: $Spr(\beta^+)$. \\Define $B:=\{c\mid\exists i<length(c)[c(i)\sqsubset\delta^i]\}$.  We now prove: $B$ is a bar in $\mathcal{F}_{\beta^+}$. \\Let $\gamma$ in $\mathcal{F}_{\beta^+}$ be given. Find $\zeta$ in $\mathcal{F}_\beta$ such that $\forall n[\gamma(n)\sqsubset \zeta]$.\\ Find $\alpha$ such that $\forall n[\gamma(n)=\overline\zeta\alpha(n)]$. \\Find $n$ such that $\overline\zeta\alpha(n)\sqsubset\delta^n$ and, therefore: $\gamma(n)\sqsubset\delta^n$ and: $\overline \gamma(n+1) \in B$.\\ We thus see: $Bar_{\mathcal{F}_{\beta^+}}(B)$. \\ We define: $\tilde{\langle\;\rangle}=\langle\;\rangle$, and, for each $n>0$, for each  $c$ in $\omega^n$, $\tilde c:= c(n-1)$.
\\Define 
$C:=\bigcup_n\{c\in \omega^n\mid \beta^+(c)=0\;\wedge\;\bigl(\forall i< n[c(i)\perp\delta^i] \rightarrow \;^{\tilde c}\beta \in \mathcal{CB}]\bigr)\}$. \\Note: $B\subseteq C$ and  $C$ is monotone in $\{s\mid\beta^+(s)=0\}$. \\Let $c,n$ be given such that $c\in \omega^n$ and $\beta^+(c)=0$ and $\forall t
[\beta^+(c\ast\langle t \rangle )=0\rightarrow c\ast\langle t \rangle\in C]$. \\ Assume: $\forall i<n[c(i)\perp\delta^i]$.  
\\Note: for all $t$, if $\tilde c\sqsubset t$ and $\beta(t)=0$ and $ t\perp \delta^n$, then $c\ast\langle t \rangle \in C$ and $ ^t\beta \in \mathcal{CB}$. 
\\Find $\varepsilon$ such that $\tilde c\ast\varepsilon \in \mathcal{F}_\beta$, and, if $\tilde c\sqsubset\delta^n$, then $\delta^n = \tilde c\ast \varepsilon $. 
\\Define $\nu:=\;^{\tilde c}\beta$ and note: $\varepsilon \in \mathcal{F}_\nu$ and, for all $s$, if $\nu(s)=0$ and $\varepsilon\perp s$, then $^s\nu \in \mathcal{CB}$. \\ In particular, for all $m,n,s$, if $s=\overline \varepsilon m\ast \langle n \rangle$ and $\varepsilon(m) \neq n$, then $^s\nu \in \mathcal{CB}$.
\\Conclude: $\nu =\;^{\tilde c}\beta \in \mathcal{CB}$, and: $c\in C$.
\\We thus see: $C$ is inductive in $\{s\mid\beta^+(s)=0\}$. 
\\ Using the Principle of Bar Induction $\mathbf{BI}
$, see Subsubsection \ref{SSS:barinduction}, we conclude: \\$\langle\;\rangle \in C$, i.e.  $\beta \in \mathcal{CB}$. 
 
 We thus see: $\mathcal{ALMOST^\ast COUNT}\subseteq\mathcal{CB}$. \end{proof}                                                                                                                                                                                                                                                                                                                                                                                                                                                                                                                                                                                                              
      
\subsection{Reducible spreads}

\begin{definition}\label{D:reduciblespreads} For each $\sigma$ in $\mathcal{STP}$, we define the collection $\mathcal{RED}_\sigma$ of \textit{codes of  $\sigma$-reducible spreads}, as follows, by induction. \begin{enumerate}[\upshape (i)] \item   $\mathcal{RED}_{1^\ast}=\mathcal{RED}_{\underline 1}:=\{\underline 1\}$, and, \item for every $\sigma\neq 1^\ast$ in $\mathcal{STP}$, \\$\mathcal{RED}_\sigma$ is the set of all $\beta$ in $2^\omega$ such that $ Spr(\beta)$ and, for some $\varepsilon$ in $\mathcal{F}_\beta$,  \\$\forall m \forall n[\bigl(\varepsilon(n)\neq m \;\wedge\; \beta(\overline\varepsilon n\ast\langle m \rangle) =0  \bigr)\rightarrow \exists p[^{\overline\varepsilon n\ast\langle m \rangle}\beta\in \mathcal{RED}_{\sigma^p}]]$.  \end{enumerate} We also define $\mathcal{RED}:=\bigcup_{\sigma \in \mathcal{STP}} \mathcal{RED}_\sigma$. \end{definition}

If $\beta \in \mathcal{RED}\sigma$, then $\mathcal{F}_\beta$ is called a \textit{$\sigma$-reducible spread}.

 If $\beta \in \mathcal{RED}$, then $\mathcal{F}_\beta$ is called a \textit{reducible spread}.

 \smallskip The notion of a reducible spread  goes back to Cantor. We here introduce this notion without bringing up the operation of taking the derivative of a given $\mathcal{X}\subseteq \omega^\omega$.  Cantor defined a closed set to be {\it reducible} if one, by repeating the operation of taking the derivative, if needed transfinitely many times, ends up with the empty set. 

Note that, for all $\sigma$ in $\mathcal{STP}$, for all $\beta$ such that $Spr(\beta)$, $\mathcal{F}_\beta$ is $\sigma$-reducible if and only if $s\ast\mathcal{F}_\beta$ is $\sigma$-reducible.

Also note  that, for all $\beta_0,\beta_1$ such that $\forall i<2[Spr(\beta_i)]$ and $\mathcal{F}_{\beta_0}\subseteq\mathcal{F}_{\beta_1}$, for all $\sigma$ in $\mathcal{STP}$, if  $\mathcal{F}_{\beta_1}$ is $\sigma$-reducible, then  $\mathcal{F}_{\beta_0}$ is $\sigma$-reducible.

\begin{theorem}\label{T:cantorreducible} $\mathcal{CB}=\mathcal{RED}$. \end{theorem}
\begin{proof}
We first prove: $\mathcal{CB}\subseteq \mathcal{RED}$, using induction on $\mathcal{CB}$.  

(1) For all $\beta$,  if $Spr(\beta)$ and $\beta(0)\neq 0$, then $\mathcal{F}_\beta=\emptyset$ and  $\beta \in \mathcal{RED}_{1^\ast}$.  

(2) Let $\beta, \varepsilon$ be given such that $Spr(\beta)$ and $\varepsilon \in \mathcal{F}_\beta$ and $\forall n  \exists \sigma \in \mathcal{STP}[\beta^n\in\mathcal{DER}_\sigma]$. Using $\mathbf{AC}_{0,1}$, find $\tau$ in $\mathcal{STP}$ such that $\tau(0)=0$ and $\forall n[\beta^n \in\mathcal{DER}_{\tau^n}]$. \\ Conclude: $\mathsf{CB}(\varepsilon, \beta)\in\mathcal{RED}_\tau$. 

(3) Using induction on $\mathcal{CB}$, conclude: $\mathcal{CB}\subseteq\bigcup_{\sigma\in\mathcal{STP}}\mathcal{RED}_\sigma$.

\smallskip We now prove: 
$ \mathcal{RED}\subseteq \mathcal{CB}$, using induction on $\mathcal{STP}$.

(1) For all $\sigma$ in $\mathcal{STP}$, for all $\beta$ if if $\sigma(0)\neq 0$ and $\beta\in \mathcal{RED}_\sigma$, then $\mathcal{F}_\beta=\emptyset$ and $\beta\in \mathcal{CB}$. 

(2) Let $\sigma$ in $\mathcal{STP}$ be given such that $\sigma(0)=0$ and $\forall n[\mathcal{RED}_{\sigma^n}\subseteq\mathcal{CB}]$. 

Let $\beta$ in $\mathcal{RED_\sigma}$ be given. \\Find $\varepsilon$ in $\mathcal{F}_\beta$ such that $\forall s[(\beta(s)=0\;\wedge\;s\perp\varepsilon)\rightarrow \exists n[^s\beta\in \mathcal{RED}_{\sigma^n}]]$. \\Conclude $\forall s[(\beta(s)=0\;\wedge\;s\perp\varepsilon)\rightarrow \;^s\beta\in \mathcal{CB}]$, and: $\beta\in\mathcal{CB}$. \\ Define $\gamma$ such that, for all $m,n$, \\if $\varepsilon(m)\neq n$, then  $\gamma^{J(m,n)}=\;^{\overline \varepsilon m\ast\langle n \rangle}\beta$ and, if $\varepsilon(m) =n$, then $\gamma^{J(m,n)}=\underline 1$. 

Note: for all $n$, $\gamma^n \in \mathcal{CB}$ and 
  $\beta=\mathsf{CB}(\varepsilon, \gamma) \in \mathcal{CB}$. 

  We thus see: $\mathcal{RED}_\sigma\subseteq \mathcal{CB}$.

(3) Using induction on $\mathcal{STP}$, we conclude: $\forall \sigma \in \mathcal{STP}[\mathcal{RED}_\sigma\subseteq \mathcal{CB}]$. 
\end{proof}
\subsection{Perhaps$_\sigma$-countable spread}\label{SS:perhapscountable}
In this Subsection we will see that there are many notions of countability for spreads in between the notion of a countable spread, see Subsection \ref{SS:countablespread}, and the notion fo an almost-countable spread, see Subsection \ref{SS:almostcountablespreads}.  \begin{definition}\label{D:perhapscountable} For each inhabited $\mathcal{X}\subseteq\omega^\omega$, for each $\sigma$ in $\mathcal{STP}$, we define $\mathbb{P}(\sigma,\mathcal{X})\subseteq\omega^\omega$, the \textit{$\sigma$-th perhapsive extension of $\mathcal{X}$}, as follows, by induction. For every $\sigma$ in $\mathcal{STP}$,\begin{enumerate}[\upshape (i)] \item if $\sigma(0) \neq 0$, then  $\mathbb{P}(\sigma,\mathcal{X})=\mathcal{X}$, and, \item if $\sigma(0)=0$, then  $\mathbb{P}(\sigma,\mathcal{X})=\{\alpha\mid\exists \beta\in\mathcal{X}[\alpha\;\#\;\beta\rightarrow \exists n[\alpha \in  \mathbb{P}(\sigma^n,\mathcal{X})]]\}$. \end{enumerate}\end{definition}
In  \cite[Theorem 3.19]{veldman09}, one may find the straightforward proof that, for all inhabited $\mathcal{X},\mathcal{Y}\subseteq\omega^\omega$,  for all  $\sigma, \tau$ in $\mathcal{STP}$,  if $\mathcal{X}\subseteq\mathcal{Y}$ and $\sigma \le \tau$, then $\mathbb{P}(\sigma, \mathcal{X})\subseteq \mathbb{P}(\tau, \mathcal{Y})$.

\begin{definition}\label{D:sigmacountablespread} Let $\beta, \sigma$ be given such that $Spr(\beta)$ and $\sigma \in \mathcal{STP}$.  \\The spread $\mathcal{F}_\beta$ is called \emph{perhaps$_\sigma$-countable} if and only if $\exists \delta[\mathcal{F}_\beta\subseteq\mathbb{P}(\sigma, En_\delta)]$.\end{definition}
 
 \smallskip The proof of the third item of the next Theorem, Theorem \ref{T:cbperhapscount}, resembles the proof of: `$\mathcal{ALMOST^\ast COUNT}\subseteq \mathcal{CB}$', see Theorem \ref{T:cbcount}.
\begin{theorem}\label{T:cbperhapscount} \begin{enumerate}[\upshape (i)] \item $\forall\delta[ \mathcal{ALMOST}^\ast(En_\delta)=\bigcup_{\sigma\in\mathcal{STP}}\mathbb{P}(\sigma, En_\delta)]]$. \item $\mathcal{ALMOST^\ast FIN} =\bigcup_{\sigma\in\mathcal{STP}}\mathbb{P}(\sigma,\mathcal{FIN})$.  \item For all $\beta, \delta, \varphi$, if $Spr(\beta)$ and $\varphi:\mathcal{F}_\beta\rightarrow\mathcal{ALMOST}^\ast (En_\delta)$, \\then $\exists \sigma\in\mathcal{STP}[\varphi:\mathcal{F}_\beta \rightarrow\mathbb{P}(\sigma,En_\delta)]$. \item For all $\beta, \varphi$, if $Spr(\beta)$ and $\varphi:\mathcal{F}_\beta\rightarrow\mathcal{ALMOST}^\ast (\mathcal{FIN})$, \\then $\exists \sigma\in\mathcal{STP}[\varphi:\mathcal{F}_\beta \rightarrow\mathbb{P}(\sigma,\mathcal{FIN})]$. \item $\forall \beta \in \mathcal{CB}\exists \sigma\in\mathcal{STP}\exists\varphi[\varphi: \mathcal{F}_\beta\rightarrowtail\mathbb{P}(\sigma, \mathcal{FIN})]$. \end{enumerate}\end{theorem}

\begin{proof} (i) Let $\delta$ be given.\\ We first prove: $\bigcup_{\sigma\in\mathcal{STP}}\mathbb{P}(\sigma, En_\delta) \subseteq \mathcal{ALMOST}^\ast(En_\delta)$, using induction on $\mathcal{STP}$.
\\First note:  $\mathbb{P}(1^\ast,En_\delta)=En_\delta\subseteq \mathcal{ALMOST}^\ast(En_\delta)$. \\Now let $\sigma$ in $\mathcal{STP}$ be given such that $\sigma\neq 1^\ast$ and  $\forall n[\mathbb{P}(\sigma^n, En_\delta)\subseteq \mathcal{ALMOST}^\ast(En_\delta)]$. \\Assume: $\gamma\in\mathbb{P}(\sigma, En_\delta)$. Find $n$ such that $\gamma\;\#\;\delta^n\rightarrow \exists m[\alpha \in \mathbb{P}(\sigma^m, En_\delta)]$. \\Let $\alpha$ be given and distinguish two cases. \\\textit{Case (a)}: $\overline\gamma \alpha(n)\sqsubset \delta^n$. \\\textit{Case (b)}: $\overline\gamma\alpha(n)\perp \delta^n$. Find $m$ such that $\gamma \in \mathbb{P}(\sigma^m,En_\delta)$. \\Conclude: $\gamma \in \mathcal{ALMOST}^\ast(En_\delta)$ and: $\exists p[\overline \gamma \alpha (p) \sqsubset \delta^p]$. \\We thus see, in both cases: $\exists p[\overline \gamma \alpha (p) \sqsubset \delta^p]$. \\Conclude: $\forall \gamma \in\mathbb{P}(\sigma,En_\delta)\forall \alpha \exists p[\overline \gamma \alpha p \sqsubset \delta^p]$, that is: $\mathbb{P}(\sigma,En_\delta)\subseteq \mathcal{ALMOST}^\ast(En_\delta)$. Using induction on $\mathcal{STP}$, conclude: $\bigcup_{\sigma\in\mathcal{STP}}\mathbb{P}(\sigma, En_\delta)\subseteq \mathcal{ALMOST}^\ast(En_\delta)$.

\smallskip We now prove: $\mathcal{ALMOST}^\ast(En_\delta)\subseteq \bigcup_{\sigma\in\mathcal{STP}}\mathbb{P}(\sigma, En_\delta) $. \\ Let $\gamma$ in $\mathcal{ALMOST}^\ast(En_\delta)$ be given. \\Define $B:=\bigcup_p\{a\in \omega^p\mid\exists i<p[\overline\gamma a(i) \sqsubset\delta^i]\}$ and note: $B$ is a bar in $\omega^\omega.$ \\ Define $C:=\bigcup_p\{a\in \omega^p\mid\forall i< p[\overline\gamma a(i)\perp\delta^i]\rightarrow \exists\sigma\in \mathcal{STP}[\gamma \in \mathbb{P}(\sigma, En_\delta)]\}$. \\Note: $B\subseteq C$ and $C$ is monotone. We now prove that $C$ is inductive.
\\ Let $a$ be given such that $\forall n[a\ast\langle n \rangle \in C]$. Define $p:=length(a)$. \\Assume $\forall i< p[\overline\gamma a(i) \perp\delta^i]$.   \\Using the Second Axiom of Countable Choice $\mathbf{AC}_{0,1}$, see Subsubsection \ref{SSS:countablechoice}, \\find $\tau$ in $\mathcal{STP}$ such that $\forall b[\overline \gamma b \perp \delta^p\rightarrow \gamma\in\mathbb{P}(\tau^b, En_\delta)]$.  \\Conclude: if  $\gamma \perp \delta^p$, then $ \exists b[\gamma \in \mathbb{P}(\tau^b,En_\delta)]$, i.e. $\gamma \in \mathbb{P}(\tau, En_\delta)$. \\We thus see:  if $\forall i< length(a)[\overline\gamma a(i) \perp\delta^i$, then $\exists\tau[\gamma \in \mathbb{P}(\tau, En_\delta)]$, that is: $a\in C$. \\
Conclude: $C$ is inductive. \\Using the Principle of Bar Induction $\mathbf{BI}$, see Subsubsection \ref{SSS:barinduction}, we find:\\ $\langle\;\rangle \in C$, i.e. $\exists\tau[\gamma \in \mathbb{P}(\tau, En_\delta)]$.

We thus see: $\mathcal{ALMOST}^\ast(En_\delta)\subseteq\bigcup_{\sigma \in \mathcal{STP}}\mathbb{P}(\sigma, En_\delta)$.

\smallskip (ii) This follows from (i) and Lemma \ref{L:en2}.

\smallskip (iii) Let $\beta,\delta, \varphi$ be given such that $Spr(\beta)$ and $\varphi:\mathcal{F}_\beta\rightarrow
 \mathcal{ALMOST}(En_\delta)$. \\ Note: $\forall \gamma \in \mathcal{F}_\beta \forall \alpha\exists n[\overline{\varphi|\gamma}\alpha(n)\sqsubset \delta^n]$.  \\Define $\beta^+$ such that, for each $c$,  $\beta^+(c)=0$ if and only if \\$\forall i[i+1<length(c)\rightarrow c(i)\sqsubset c(i+1)]$ and \\$\forall i<length(c)[\beta(c_I(i))=0\;\wedge\; length\bigl(\varphi|c_I(i))\ge c_{II}(i)])]$. Note: $Spr(\beta^+)$. \\Define $B:=\bigcup_p\{c\in \omega^p\mid\exists i<p[\overline{\varphi|c_I(i)}c_{II}(i)\sqsubset\delta^i]\}$.  \\We now prove that $B$ is a bar in $\mathcal{F}_{\beta^+}$. \\Let $\gamma$ in $\mathcal{F}_{\beta^+}$ be given. Find $\zeta$ in $\mathcal{F}_\beta$ such that $\forall n[\gamma_I(n)\sqsubset \zeta]$.  \\Find $n$ such that $\overline\zeta\gamma_{II}(n)\sqsubset\delta^n$ and, therefore:  $\overline \gamma(n+1) \in B$. \\Conclude: $Bar_{\mathcal{F}_{\beta^+}}(B)$. 

For each $c$ such that $\beta^+(c)=0$ we define $\tilde c$ as follows. $\tilde{0}=0$ and, for each $c$, for all $n$, if $n=length(c)>0$, then $\tilde c:=c_I(n-1)$.  Let $C$ be the set of all $c$ such that 
$ \beta^+(c)=0$ and, if $\forall i< length(c)[\overline{c_I(i)}c_{II}(i)\perp\delta^i]$, then $\exists \sigma \in \mathcal{STP}[\varphi:\mathcal{F}_\beta\cap\tilde c\rightarrow\mathbb{P}(\sigma,En_\delta)]$. \\Note: $B\subseteq C$ and $C$ is monotone in $\{s\mid\beta^+(s)=0\}$.  
\\We now prove that $C$ is inductive in $\{s\mid\beta^+(s)=0\}$.   \\Let $c$ be given such that $\beta^+(c)=0$ and $\forall t
[\beta^+(c\ast\langle t \rangle )=0\rightarrow c\ast\langle t \rangle\in C]$. \\Find $n:=length(c)$. Assume: $\forall i<n[\overline{c_I(i)}c_{II}(i)\perp\delta^i]$.  \\Note: $\forall t[(\beta^+ (c\ast\langle t \rangle)=0\;\wedge\; \overline{t_I}c_{II}(n)\perp\delta^n)\rightarrow \exists\sigma\in\mathcal{STP}[\varphi:\mathcal{F}_\beta\cap K(t)\rightarrow\mathbb{P}(\sigma,En_\delta)]]$. Using the Second Axiom of Countable Choice $\mathbf{AC}_{0,1}$, see Subsubsection \ref{SSS:countablechoice}, \\find $\tau$ in $\mathcal{STP}$ such that, for all $t$, \\
if $\beta^+(c\ast\langle t \rangle)=0$ and  $\overline{t_I}c_{II}(n)\perp\delta^n$, then $\varphi:\mathcal{F}_\beta\cap K(t)\rightarrow\mathbb{P}(\tau^t,En_\delta)$. \\Clearly, $\forall \gamma\in\mathcal{F}_\beta\cap\tilde c[\varphi|\gamma\;\#\;\delta^n\rightarrow \exists t[\varphi|\gamma \in\mathbb{P}(\tau^t,En_\delta)]]$ and: $\varphi:\mathcal{F}_\beta\cap\tilde c\rightarrow \mathbb{P}(\tau,En_\delta)$.
We thus see: $C$ is inductive in $\{s\mid\beta^+(s)=0\}$. 

 Using the Principle of Bar Induction $\mathbf{BI}$, see Subsubsection \ref{SSS:barinduction}, we conclude: \\$\langle\;\rangle \in C$, i.e. $\exists \sigma \in \mathcal{STP}[\varphi:\mathcal{F}_\beta\rightarrow \mathbb{P}(\sigma, En_\delta)]$.

 \smallskip (iv) This is an immediate consequence of (iii), as $\exists \delta[\mathcal{FIN}=En_\delta]$.

 \smallskip (v) This follows from (iii) and Theorem \ref{T:countandfin}(ii).
\end{proof}

 \subsection{Special and very special Cantor-Bendixson sets }
\begin{definition}\label{D:specialcb}We define a function $\sigma\mapsto cb_\sigma$ from $\mathcal{STP}$ to $\omega^\omega$, as follows.
\begin{enumerate}[\upshape(i)]\item  $cb_{1^\ast}=\underline 1$, and,\item for all $\sigma\neq 1^\ast $ in $\mathcal{STP}$, $cb_\sigma$ satisfies: $\forall m[cb_\sigma(\overline{\underline 0}m)=0]$ and  \\$\forall m\forall n\forall s[cb_\sigma(\overline{\underline 0}m\ast\langle n+1\rangle\ast s)=cb_{\sigma^n}(s)]$. \end{enumerate}

\emph{Note: if $\sigma\neq 1^\ast$, then $ cb_\sigma=\mathsf{CB}(\underline 0, \beta)$, where, for all $m,n$, $\beta^{J(m,n+1)}=cb_{\sigma^n}$}.

\smallskip

We also define a function $\sigma\mapsto cb^\lozenge_\sigma$ from $\mathcal{STP}$ to $\omega^\omega$, as follows.
\begin{enumerate}[\upshape(i)]\item  $cb^\lozenge_{1^\ast}=\underline 1$, and,\item for all $\sigma\neq 1^\ast$ in $\mathcal{STP}$,  $cb^\lozenge_\sigma$ satisfies: $\forall m[cb_\sigma^\lozenge(\overline{\underline 0}m)=0]$ and \\$\forall m\forall s[cb^\lozenge_\sigma(\overline{\underline 0}m\ast\langle 1\rangle\ast s)=cb^\lozenge_{\sigma^{L(m)}}(s)]$ and $\forall m\forall n\forall s[cb_\sigma^\lozenge(\overline{\underline 0}m\ast\langle n+2 \rangle\ast s)= 1]$. \end{enumerate}

\emph{Note: if $\sigma\neq 1^\ast$, then $ cb^\lozenge_\sigma=\mathsf{CB}(\underline 0, \beta)$, where, \\for all $m$, $\beta^{J(m,1)}=cb^\lozenge_{\sigma^{L(m)}}$ and, for all $m,n$,  $\beta^{J(m,n+2)}=\underline 1$.}\end{definition}

\smallskip

Note: for each $\sigma$ in $\mathcal{STP}$, $cb_\sigma$ is a spread-law and $cb^\lozenge_\sigma$ is a fan-law and $\mathcal{F}^\lozenge_{cb_\sigma}\subseteq 2^\omega$.

Note: for each $\sigma$ in $\mathcal{STP}$, for each $n$, \\$\mathcal{F}_{cb_\sigma}$ embeds into $\mathcal{F}_{cb_\sigma}\cap\overline{\underline 0}n$, and $\mathcal{F}_{cb^\lozenge_\sigma}$ embeds into $\mathcal{F}^\lozenge_{cb_\sigma}\cap\overline{\underline 0}n$.

The sets $\mathcal{F}_{cb_\sigma}$, where $\sigma \in\mathcal{STP}$, are called: \textit{special Cantor-Bendixson sets}.

The sets $\mathcal{F}^\lozenge_{cb_\sigma}$, where $\sigma \in\mathcal{STP}$, are called: \textit{very special Cantor-Bendixson sets}.

The latter sets occur  in \cite{veldman05} and \cite{veldman09}.

\begin{lemma}\label{L:embspecialveryspecial} For all $\sigma$ in $\mathcal{STP}$, $\mathcal{F}_{cb_\sigma}$ embeds into $\mathcal{F}_{cb^\lozenge_\sigma}$. \end{lemma}\begin{proof} We use induction on $\mathcal{STP}$. \\First note: $\mathcal{F}_{cb_{1^\ast}}=\mathcal{F}_{cb^\lozenge_{1^\ast}}=\emptyset$, so, for $\sigma = 1^\ast$,   the statement is trivial. \\Let $\sigma\neq 1^\ast$ in $\mathcal{STP}$ be given such that, for  all $n$,  $ \mathcal{F}_{cb_{\sigma^n}}$ embeds into $\mathcal{F}_{cb^\lozenge_{\sigma^n}}$. \\Using $\mathbf{AC}_{0,1}$, find $\varphi$ such that, for all $n$, $\varphi^n$ embeds $ \mathcal{F}_{cb_{\sigma^n}}$ into $\mathcal{F}_{cb^\lozenge_{\sigma^n}}$. \\Define $\psi:\mathcal{F}_{cb_\sigma}\rightarrow \omega^\omega$ such that $\psi|\underline 0=\underline 0$ and for all $m,n$, \\for all $\alpha$ in $\mathcal{F}_{cb_{\sigma^n}}$, $\psi|\overline{\underline 0}m\ast\langle n +1\rangle\ast\alpha=\underline{\overline 0}J(n,m)\ast \langle 1 \rangle\ast\varphi^n|\alpha$. \\Then  $\psi$ embeds $\mathcal{F}_{cb_\sigma}$ into $\mathcal{F}_{cb^\lozenge_\sigma}$.
\end{proof}
 The proof of the following lemma does not use the Fan Theorem. 
\begin{lemma}[The Fan Theorem for very special Cantor-Bendixson sets]\label{L:fantheoremcbs}\hfill

 For every $\sigma$ in $\mathcal{STP}$,  for every $B\subseteq \omega$, every bar in 
$\mathcal{F}_{cb^\lozenge_\sigma}$ has a finite subbar.\end{lemma}

\begin{proof} We use induction on $\mathcal{STP}$. \\Assume $ \sigma \in \mathcal{STP}$. If $\sigma=1^\ast$, there is nothing to prove. \\So assume $\sigma\neq 1^\ast$ and, for each $n$, every bar in $\mathcal{F}_{cb^\lozenge_{\sigma^n}}$ has a finite subbar. \\Now assume $B\subseteq\omega$ is a bar in $\mathcal{F}_{cb^\lozenge_\sigma}$. Find $n$ such that $\underline{\overline 0}n\in B$. \\Using the induction hypothesis, find finite subsets $B_0, B_1, \ldots,B_{n-1}$ of $B$ such that, \\for each $i<n$, $B_i$ is bar in $\mathcal{F}_{cb^\lozenge_\sigma}\cap \overline{\underline 0}i\ast\langle 1 \rangle$. \\Note: the finite set $\{\underline{\overline 0}n\}\cup\bigcup_{i<n}B_i$ is bar in $\mathcal{F}_{cb_\sigma}$. 
\end{proof}
The next Theorem shows that every Cantor-Bendixson set is, in a certain sense, equinumerous to a special Cantor-Bendixson set.
\begin{theorem}\label{T:cardinalityofspreads} For every Cantor-Bendixson set $\mathcal{F}$  there exists a special \\Cantor-Bendixson set $\mathcal{H}$ such that $\mathcal{H}$ maps onto $\mathcal{F}$ and $\mathcal{F}$ embeds into $\mathcal{H}$:  \\$\forall \beta \in \mathcal{CB}\exists \sigma\in \mathcal{STP}[\exists \varphi[\varphi: \mathcal{F}_{cb_\sigma}\twoheadrightarrow\mathcal{F}_\beta]\;\wedge\; \exists \psi[\psi:\mathcal{F}_\beta \rightarrowtail\mathcal{F}_{cb_\sigma}]]$.  \end{theorem}

\begin{proof} We use induction on $\mathcal{CB}$.
\\If $\beta(0) \neq 0$, so $\mathcal{F}_\beta =\emptyset$, one may take $\sigma=1^\ast$, as also $\mathcal{F}_{cb_\sigma}=\emptyset$. 
\\Now let $\beta,\varepsilon$ be given such that $Spr(\beta)$ and $\varepsilon\in\mathcal{F}_\beta$ and, for all $m,n,s$, if $\varepsilon(n)\neq m$ and $s=\overline\varepsilon n\ast\langle m \rangle$, then there exist $\sigma$ in $\mathcal{STP}$  such that $  \mathcal{F}_{cb_\sigma}$ maps onto  $\mathcal{F}_{^s\beta}$ and $\mathcal{F}_{^s\beta}$ embeds into $\mathcal{F}_{cb_\sigma}$. 

Using the Second Axiom of Countable Choice $\mathbf{AC}_{0,1}$, see Subsubsection \ref{SSS:countablechoice},\\ find $\tau, \varphi,\psi$ such that  $ 1^\ast\neq \tau\in\mathcal{STP}$ and, for all $m,n,s$, if $\varepsilon(m)\neq n$ and $s=\overline\varepsilon m \ast\langle n \rangle$, then $\varphi^s: \mathcal{F}_{cb_{\tau^s}}\twoheadrightarrow \mathcal{F}_{^s\beta}$ and $\psi^s:\mathcal{F}_{^s\beta}\rightarrowtail\mathcal{F}_{cb_{\tau^s}}$. 

Define $C:=\{s\mid \beta(s)=0\;\wedge\;\exists m\exists n[s=\overline \varepsilon m \ast\langle n\rangle\;\wedge\; \varepsilon(m)\neq n]\}$.

Define $\rho:\mathcal{F}_{cb_\tau}\rightarrow \omega^\omega$ such that $\rho|\underline 0=\varepsilon$ and, \\ for all $s$, if $s\in C$, then, for all $\gamma$  in $\mathcal{F}_{cb_{\tau^s}}$, $\rho|(\overline{\underline 0}s\ast\langle s+1\rangle \ast\gamma)=s\ast\varphi^s|\gamma$ and, \\for each $\delta$ in $\mathcal{F}_{cb_\tau}$, if there is no $s$ in $C$ such that $\overline{\underline 0}s\ast\langle s +1\rangle \sqsubset \delta$, then $\rho|\gamma=\varepsilon$. \\Clearly, $\rho$ maps $\mathcal{F}_{cb_\tau}$ onto $\mathcal{F}_\beta$.

 Define $\chi:\mathcal{F}_\beta \rightarrow \omega^\omega$ such that $\chi|\varepsilon=\underline 0$ and, for all $s$ in $C$,  for all $\gamma\in\mathcal{F}_{^s\beta}$,\\ $\chi|(s\ast\gamma)=\overline{\underline 0} s\ast\langle s+1\rangle\ast\psi^s|\gamma$. \\Clearly, $\chi$ embeds $\mathcal{F}_\beta$ into $\mathcal{F}_{cb_\tau} $.
\end{proof}

 The next result, Theorem \ref{T:cardinalityoffans}, gives a refinement of  Theorem \ref{T:cardinalityofspreads}: every {\it finitary} Cantor-Bendixson set is, what one might call, equinumerous to a {\it very} special Cantor-Bendixson set.  
 
\begin{theorem}\label{T:cardinalityoffans} For every Cantor-Bendixson-set $\mathcal{F}$ that is a fan  there exists a very special Cantor-Bendixson-set $\mathcal{H}$ such that $\mathcal{H}$ maps onto $\mathcal{F}$ and $\mathcal{F}$ embeds into $\mathcal{H}$:

$\forall \beta \in \mathcal{CB}[Fan(\beta)\rightarrow\exists \sigma\in \mathcal{STP}[\exists \varphi[\varphi: \mathcal{F}_{cb^\lozenge_\sigma}\twoheadrightarrow\mathcal{F}_\beta]\;\wedge\; \exists \psi[\psi:\mathcal{F}_\beta \rightarrowtail\mathcal{F}_{cb^\lozenge_\sigma}]]$.  \end{theorem}

\begin{proof} We use induction on $\mathcal{CB}$.
\\If $\beta(0) \neq 0$,  take $\sigma=1^\ast$, and note: $\mathcal{F}_\beta= \mathcal{F}_{cb_{1^\ast}}=\emptyset$. 
\\Now let $\beta,\varepsilon$ be given such that $Fan(\beta)$ and $\varepsilon\in\mathcal{F}_\beta$ and for all $m,n,s$, \\if $\varepsilon(m)\neq n$ and $s=\overline\varepsilon m\ast\langle n \rangle$, \\then there exist $\sigma$ in $\mathcal{STP}$  such that $  \mathcal{F}^\lozenge_{cb_\sigma}$ maps onto  $\mathcal{F}_{^s\beta}$ and $\mathcal{F}_{^s\beta}$ embeds into $\mathcal{F}^\lozenge_{cb_\sigma}$. 
 
 Using the Second Axiom of Countable Choice $\mathbf{AC}_{0,1}$, see Subsubsection \ref{SSS:countablechoice},\\ find $\tau, \varphi,\psi$ such that  $ 1^\ast\neq \tau\in\mathcal{STP}$ and, for all $m,n,s$, if $\varepsilon(m)\neq n$ and $s=\overline\varepsilon m\ast\langle n \rangle$, then $\varphi^s: \mathcal{F}_{cb_{\tau^s}}\twoheadrightarrow \mathcal{F}^\lozenge_{^s\beta}$ and $\psi^s:\mathcal{F}_{^s\beta}\rightarrowtail\mathcal{F}^\lozenge_{cb_{\tau^s}}$.

Define $C:=\{s\mid \beta(s)=0\;\wedge\;\exists m\exists n[s=\overline \varepsilon m \ast\langle n\rangle\;\wedge\; \varepsilon(m)\neq n]\}$. \\Note:  $Fan(\beta)$, and thus: $\forall m\exists p\forall s\ge p[s\in C\rightarrow\mathit{length}(s)\ge m]$. 

Using the First Axiom of Countable Choice $\mathbf{AC}_{0,0}$, see Subsubsection \ref{SSS:countablechoice},\\ find $\zeta$ such that $\forall m\forall s\ge \zeta(m)[s\in C\rightarrow\mathit{length}(s)\ge m]$.

Define $\rho:\mathcal{F}^\lozenge_{cb_\tau}\rightarrow \omega^\omega$ such that $\rho|\underline 0=\varepsilon$ and, \\for all $s$ in $C$, for all $ \gamma\in \mathcal{F}^\lozenge_{cb_{\tau^n}}$, $\rho|(\overline{\underline 0}J(s,0)\ast\langle 1\rangle\ast \gamma) = s\ast\varphi^s|\gamma$ and, \\for all $\delta$ in $\mathcal{F}_{cb_\tau}$, if there is no $s$ in $C$ such that  $\overline{\underline 0}J(s,0)\ast\langle 1\rangle\sqsubset \delta$, then $\rho|\delta=\varepsilon$.\\ Note: $\rho$ is well-defined and: $\forall m\forall \gamma \in \mathcal{F}_{cb_\tau}[\overline{\underline  0} J\bigl(\zeta(m),0\bigr)\sqsubset\gamma\rightarrow \overline\varepsilon m\sqsubset\rho|\gamma]$. \\Clearly, $\rho:\mathcal{F}_{cb_\tau}\twoheadrightarrow\mathcal{F}_\beta$.

 Define $\chi:\mathcal{F}_\beta \rightarrow \omega^\omega$ such that $\chi|\varepsilon=\underline 0$ and, for all $s$ in $C$, for all  $\gamma$ in $\mathcal{F}_{^s\beta}$, \\$\chi|(s\ast\gamma)=\overline{\underline 0} J(s,0)\ast\langle 1\rangle\ast\psi^s|\gamma$. \\Clearly, $\chi:\mathcal{F}_\beta\rightarrowtail\mathcal{F}_{cb_\tau} $.
\end{proof}

\begin{corollary}\label{COR:cardinalityofspreads}Let $\beta$ be given such that $Spr(\beta)$. 

$
\mathcal{F}_\beta$ is almost-countable if and only if $\exists \sigma \in \mathcal{STP}\exists \varphi[\varphi: \mathcal{F}_{cb_\sigma}\twoheadrightarrow \mathcal{F}_\beta]$.  \end{corollary}
\begin{proof}
Use Theorems \ref{T:cbcount} and \ref{T:cardinalityofspreads} and Lemma \ref{L:imagealmostcountable}.
\end{proof}

The second item of the following Theorem seems to be of some interest in itself. It is an extension of Theorem \ref{T:perfspr}(iii).
\begin{theorem}\label{T:injectiveandsurjective}\hfill

\begin{enumerate}[\upshape (i)] \item For all $\beta$, if $\forall i<2[Spr(\beta^i)]$  and $\exists \varphi[\varphi:\mathcal{F}_{\beta^0}\twoheadrightarrow\mathcal{F}_{\beta^1}]$, then $\exists\psi[\psi:\mathcal{F}_{\beta^1}\rightarrowtail\mathcal{F}_{\beta^0}]$. \item For all $\beta$, if $Spr(\beta^0)$ and $Fan(\beta^1)$ and $\exists \psi[\psi:\mathcal{F}_{\beta^1}\rightarrowtail\mathcal{F}_{\beta^0}]$, \\then $\exists\varphi[\varphi:\mathcal{F}_{\beta^0}\twoheadrightarrow\mathcal{F}_{\beta^1}]$. \end{enumerate}\end{theorem}

\begin{proof}
(i) Let $\beta,\varphi$ be given such that $\varphi:\mathcal{F}_{\beta^0}\twoheadrightarrow\mathcal{F}_{\beta^1}$, and, therefore: \\$\forall \gamma\in\mathcal{F}_{\beta^1}\exists\alpha\in\mathcal{F}_{\beta^0}[\varphi|\alpha=\gamma]$. 
\\Using the Second Axiom of Continuous Choice $\mathbf{AC}_{1,1}$,  see Subsubsection \ref{SSS:bcpcontchoice},\\ find $\psi:\mathcal{F}_{\beta^1}\rightarrow\mathcal{F}_{\beta^0}$ such that $\forall \gamma \in \mathcal{F}_{\beta^1}[\varphi|(\psi|\gamma)=\gamma]$.  
\\We prove that $\psi$ is strongly injective. \\Let $\gamma, \delta$ in $\mathcal{F}_{\beta^1}$ be given such that $\gamma\;\#\;\delta$. Find $n$ such that $\overline\gamma n\perp \delta$. \\Find $m$ such that $\forall \alpha\in\mathcal{F}_{\beta^0}[\overline{\psi|\gamma}m=\overline  \alpha m \rightarrow \overline{\varphi|(\psi|\gamma)}n=\overline { \varphi|\alpha}n]$. \\Consider $\alpha:=\psi|\delta$ and conclude:  $\overline{\psi|\gamma} m \neq \overline{\psi|\delta}m$.\\We thus see: $\forall \gamma \in \mathcal{F}_{\beta^1}\forall \delta\in \mathcal{F}_{\beta^1}[\gamma \;\#\;\delta\rightarrow\psi|\gamma \;\#\;\psi|\delta]$, that is: $\psi:\mathcal{F}_{\beta^1}\rightarrowtail\mathcal{F}_{\beta^0}$. 

\smallskip 
 (ii) Let $\beta,\psi$ be given such that  $Spr(\beta^0)$ and $Fan(\beta^1)$ and $\psi:\mathcal{F}_{\beta^1}\rightarrowtail\mathcal{F}_{\beta^0}$.

  We first define $\delta$ such that 
  $\forall s[\delta(s)=0 \leftrightarrow \exists \alpha\in\mathcal{F}_{\beta^1}[s\sqsubset\psi|\alpha]]$.
  
  Let $s$ be given. Note $\forall \alpha \in \mathcal{F}_{\beta^1}\exists m[s\sqsubset \psi|\overline \alpha m\;\vee\;s\perp \psi|\overline \alpha m]$. 
  
  Using the Fan Theorem $\mathbf{FT}$, see Subsubsection \ref{SSS:fantheorem}, find $m$ such that \\$\forall \alpha \in \mathcal{F}_{\beta^1}[s\sqsubset \psi|\overline \alpha m\;\vee\;s\perp \psi|\overline \alpha m]$, i.e. $\forall t \in \omega^m[\beta^1(t)=0\rightarrow s \sqsubset \psi|t\;\vee\; s\perp \psi|t]$.
  \\Define $\delta(s):=0$ if $\exists t\in \omega^m[\beta^1(t)=0 \;\wedge\; s\sqsubset \psi|t]$ and \\$\delta(s):=1$ if $\forall t \in \omega^m[\beta^1(t)=0\rightarrow s\perp  \psi|t]$. 
 \\ Conclude: $\forall s[\delta(s)=0 \leftrightarrow \exists \alpha\in \mathcal{F}_{\beta^1}[s\sqsubset \psi|\alpha]]$. 
  \\Note: $Spr(\delta)$. 
   \\Also note, using $\mathbf{FT}$ again:  for each $m$, the set $\{\overline {\psi|\alpha}m\mid \alpha \in \mathcal{F}_{\beta^1}\}$ is finite. \\Conclude:  $Fan(\delta)$.

 We now construct  $ \tau:\mathcal{F}_\delta\rightarrow\mathcal{F}_{\beta^1}$ such that $\forall \varepsilon \in \mathcal{F}_\delta [\psi|(\tau|\varepsilon)=\varepsilon]$.

\smallskip Let $\varepsilon$ in $F_\delta$ be given. \\We claim: for all $s,t$ if $\beta^1(s)=\beta^1(t)=0$ and $s\perp t$, then there exists $n$\\ such that either $\forall \alpha \in  \mathcal{F}_{\beta^1}\cap s)[\psi|\overline\alpha n\perp\overline\varepsilon n]$ or $\forall \alpha \in  \mathcal{F}_{\beta^1}\cap t)[\psi|\overline\alpha n\perp \overline \varepsilon n]$.

\smallskip We prove this claim as follows. \\Let $s,t$ be given such that $\beta^1(s) = \beta^1(t) =0$ and $s \perp t$. \\Note: $\forall \alpha\in\mathcal{F}_{\beta^1}\cap s\forall\gamma\in \mathcal{F}_{\beta^1}\cap t[\psi|\alpha\perp\psi|\gamma]$. \\Conclude:   $\forall \alpha\in\mathcal{F}_{\beta^1}\cap s\forall\gamma\in \mathcal{F}_{\beta^1}\cap t \exists n[\psi|\overline\alpha n\perp \overline \varepsilon n\;\vee\;\psi|\overline\gamma n\perp\overline \varepsilon n]$. \\Using  the Fan Theorem $\mathbf{FT}$,  find $n$ such that $n\ge length(s)$ and $n\ge length(t)$ and \\$\forall \alpha\in\mathcal{F}_{\beta^1}\cap s\forall\gamma\in \mathcal{F}_{\beta^1}\cap t [\psi|\overline\alpha n\perp \overline \varepsilon n\;\vee\;\psi|\overline\gamma n\perp\overline \varepsilon n]$. \\Define $A:=\{u\in \omega^n \mid \beta^1(u)=0\;\wedge\ s\sqsubseteq u\}$ and $B:=\{u\in \omega^n \mid \beta^1(u)=0\;\wedge\ t\sqsubseteq u\}$. \\  Note: $\forall u\in A\forall v\in B[\psi|u\perp\overline \varepsilon n\;\vee\;\psi|v\perp\overline \varepsilon n]$. Note that $A,B$ are finite sets. \\Conclude, using Lemma \ref{L:logicalhelp}, either
$\forall u\in A[\psi|u\perp\overline \varepsilon n]$ or $\forall v\in B[\psi|v\perp\overline   \varepsilon n]$, i.e. \\either $\forall \alpha \in \mathcal{F}_{\beta^1}\cap s[\psi|\overline \alpha n\perp \overline \varepsilon n]$ or  $\forall \alpha\in \mathcal{F}_{\beta^1}\cap t[\psi|\overline\alpha n\perp \overline \varepsilon n]$.

\smallskip

Using the above  fact repeatedly and keeping in mind that $\{k\mid\beta^1(\langle k \rangle)=0\}$ is a finite set,  conclude: 
$\exists k\exists n[\beta^1(\langle k\rangle)=0\;\wedge\; \forall \alpha \in \mathcal{F}_{\beta^1}[\alpha(0)\neq k\rightarrow \psi|\overline\alpha n\perp \overline \varepsilon n]]$.

We now define the promised $\tau$, inductively, first specifying $\tau^0$, then $\tau^1$, and so on. 
 \\We start with $\tau^0$. Let $s$ be given and define $n:=length(s)$. \\Find out if there exists $k$ such that $\beta^1(\langle k\rangle)=0$ and \\$\forall j[\bigl(j\neq k \;\wedge\;\beta^1(a\ast\langle j \rangle)=0\bigr) \rightarrow  \forall \alpha\in\mathcal{F}_{\beta^1}\cap a\ast\langle j \rangle [ \psi|\overline\alpha n\perp s]]$. \\If so, find such $k$ and define $\tau^{m+1}(s)=k+1$, and,  if not, define $\tau^{m+1}(s)=0$.

\smallskip \noindent Assume: $m>0$ is given and $\tau^0, \tau^1,\ldots\tau^{m-1} $ have been defined. \\We  define $\tau^{m}$ as follows. \\Let $s$ be given. If $\delta(s)\neq 0$ or $\exists i< m\neg\exists j<length(s)[\tau^i(\overline s j)>0]$, define $\tau^{m+1}(s)=0$. \\Assume $\delta(s)= 0$ and $\forall i< m\exists j<length(s)[\tau^i(\overline s j)>0]$. \\Find $a$ such that $length(a)=m$  and $\forall i< m\exists j<length(s)[\tau^i(\overline s j)=a(i)+1]$.  \\(One might say: $a:=\overline{\tau|s}m$, although this is a little previous, \\as $\tau$ is still under construction.) \\Note that $\{k\mid\beta^1(a\ast\langle k \rangle)=0\}$ is a finite set. Define $n:=length(s)$.  \\Again using the claim we proved a moment ago, find out if there exists $k$ \\such that $\beta^1(a\ast \langle k\rangle)=0$ and $\forall j[\bigl(j\neq k \;\wedge\;\beta^1(a\ast\langle j \rangle)=0\bigr) \rightarrow  \forall \alpha\in\mathcal{F}_{\beta^1}\cap a\ast\langle j \rangle [ \psi|\overline\alpha n\perp s]]$. \\If so, find such $k$ and define $\tau^{m+1}(s)=k+1$, and,  if not, define $\tau^{m+1}(s)=0$.

\smallskip \noindent Note: $\tau:\mathcal{F}_\delta\rightarrow\omega^\omega$ and $\forall \varepsilon \in \mathcal{F}_\delta[\tau|\varepsilon \in \mathcal{F}_{\beta^1}\;\wedge\;\forall \alpha\in\mathcal{F}_{\beta_1}[\alpha\perp(\tau|\varepsilon)\rightarrow \psi|\alpha\perp \varepsilon] ]$. \\In particular: $\forall \alpha\in\mathcal{F}_{\beta^1}[\alpha\perp\bigl(\tau|(\psi|\alpha)\bigr)\rightarrow \psi|\alpha\perp\psi|\alpha]$. \\Conclude: $\forall \alpha\in\mathcal{F}_{\beta^1}[\tau|(\psi|\alpha)=\alpha]$ and $\tau:\mathcal{F}_\delta\twoheadrightarrow\mathcal{F}_{\beta^1}$. 

\smallskip \noindent Assume: $\varepsilon \in \mathcal{F}_\delta$ and  $\psi|(\tau|\varepsilon)\perp\varepsilon$. Find $m$ such that $\psi|(\tau|\overline\varepsilon m) \perp\varepsilon$. \\Note: $\forall \alpha\in\mathcal{F}_{\beta^1}[(\tau|\overline\varepsilon m) \sqsubset \alpha\rightarrow \psi|\alpha\perp\varepsilon]$. \\Conclude  $\forall \alpha \in \mathcal{F}_{\beta^1}[\psi|\alpha\perp\varepsilon]$ and: $\forall \alpha \in \mathcal{F}_{\beta^1}\exists n[\psi|\overline \alpha n\perp\overline \varepsilon n]$. \\Using  $\mathbf{FT}$ again,   find $n$ such that $\forall \alpha \in \mathcal{F}_{\beta^1}[\psi|\overline \alpha n\perp\overline \varepsilon n]$, and we have to conclude: $\delta(\overline \varepsilon n)\neq 0$ and $\varepsilon \notin \mathcal{F}_\delta$.  Contradiction.  \\Conclude: $\forall \varepsilon \in \mathcal{F}_\delta[\psi|(\tau|\varepsilon)=\varepsilon]$.
 
 \smallskip \noindent
Let $\rho:\omega^\omega\rightarrow\mathcal{F}_{\delta}$ be the canonical retraction of $\omega^\omega$ onto $\mathcal{F}_\delta$. \\Define $\varphi:\mathcal{F}_{\beta^0} \rightarrow\mathcal{F}_{\beta^1}$ such that $\forall \gamma\in\mathcal{F}_{\beta^0}[\varphi|\gamma=\tau|(\rho|\gamma)]$. \\ Note: $\forall \alpha\in \mathcal{F}_{\beta^1}[\varphi|(\psi|\alpha)=\alpha]$ and $\varphi: \mathcal{F}_{\beta^0}\twoheadrightarrow \mathcal{F}_{\beta^1}$. 
\end{proof}\begin{corollary}\label{COR:cardinalityofspreads2}Let $\beta$ be given such that $Fan(\beta)$. \\$\mathcal{F}_\beta$ is almost-countable if and only if $\exists \sigma\in\mathcal{STP}\exists \varphi[\varphi:\mathcal{F}_\beta\rightarrowtail\mathcal{F}_{cb_\sigma}]$. 

  \end{corollary}\begin{proof} Every almost-countable spread $\mathcal{F}_\beta$ embeds into some $\mathcal{F}_{cb_\sigma}$, see Theorem \ref{T:cardinalityofspreads}. Conversely, if $Fan(\beta)$ and $\mathcal{F}_\beta$ embeds into some $\mathcal{F}_{cb_\sigma}$, then, according to Theorem \ref{T:injectiveandsurjective},   $\exists \psi[\psi:\mathcal{F}_{cb_\sigma} \twoheadrightarrow \mathcal{F}_\beta]$, and, according to Lemma \ref{L:imagealmostcountable},  $\mathcal{F}_\beta$ is almost-countable. \end{proof}
 
\subsubsection{A comment} G. Ronzitti, on  page 63 of her Ph.D. dissertation \cite{ronzitti} and in the last definition of her paper \cite{ronzitti2}, suggested\footnote{We describe her suggestion in the language of this paper.} to call a spread $\mathcal{F}_\beta$  \textit{countable} if and only if   $\exists \sigma \in \mathcal{STP}\exists\varphi[\varphi:\mathcal{F}_{cb^\lozenge_\sigma}\twoheadrightarrow \mathcal{F}_\beta]$.  Unfortunately, following this suggestion, one would have to call the set $\{\underline n\mid n \in \omega\}$ a not-countable set. Corollary \ref{COR:cardinalityofspreads} shows the suggestion makes sense if one uses  the non-compact Cantor-Bendixson sets given by the function $\sigma\mapsto cb_\sigma$. The suggestion is also a good suggestion if one restricts oneself to fans, rather than spreads, see Theorem \ref{T:cardinalityoffans} and Lemma \ref{L:imagealmostcountable}.

\subsection{The Cantor-Bendixson Hierarchy}
\begin{lemma}\label{L:cbh} For all  $\sigma$ in $\mathcal{STP}$, for all $\delta$, if  $\mathcal{F}_{cb_\sigma}$ embeds into $En_\delta$, then $\sigma\le S^\ast (1^\ast)$.\end{lemma}\begin{proof} Let $\sigma,\delta$ be given such that $\sigma\in\mathcal{STP}$ and $\mathcal{F}_{cb_\sigma}$ embeds into $En_\delta$.  Then, according to Theorem \ref{T:countspread}(i), $\forall \gamma_0 \in \mathcal{F}_{cb_\sigma}\forall \gamma_1\in\mathcal{F}_{cb_\sigma}[\gamma_0 = \gamma_1\;\vee\;\gamma_0\;\#\;\gamma_1]$. Using $\mathbf{BCP}$, find $m$ such that $\forall \gamma \in \mathcal{F}_{cb_\sigma}[\overline{\underline 0}m\sqsubset \gamma\rightarrow \underline 0 =\gamma]$. Conclude: $\forall n[\mathcal{F}_{cb_{\sigma^n}} =\emptyset]$ and: $\forall n[\sigma^n \le 1^\ast]$ and: $\sigma\le S^\ast(1^\ast)$.
\end{proof} 
\begin{theorem}[The Cantor-Bendixson Hierarchy Theorem]\label{T:cbhierarchy}
 \hfill
\begin{enumerate}[\upshape(i)]
 \item For all $\sigma,\tau$ in $\mathcal{STP}$, if $\mathcal{F}_{cb_\sigma}$ is $\tau$-reducible, i.e. $cb_\sigma \in \mathcal{RED}_\tau$ \footnote{See Definition \ref{D:reduciblespreads}.}, then $\sigma\le \tau$. 
\item For all $\sigma,\tau$ in $\mathcal{STP}$, for all $\delta$, if $\mathcal{F}_{cb_\sigma}$ embeds into $\mathbb{P}(\tau, En_\delta)$, then $\sigma \le S^\ast(\tau)$.\item For all $\sigma,\tau$ in $\mathcal{STP}$, if $\mathcal{F}_{cb_\sigma}$ embeds into $\mathbb{P}(\tau, \mathcal{FIN})$, then $\sigma \le S^\ast(\tau)$.
\item For all $\sigma,\tau$ in $\mathcal{STP}$, for all $\delta$ in $(\mathcal{F}_{cb_\sigma})^\omega$, if $\mathcal{F}_{cb_\sigma}\subseteq\mathbb{P}(\tau, En_\delta)$, then $\sigma\le S^\ast(\tau)$.

\end{enumerate}

\end{theorem}
\begin{proof} (i) We use induction on $\mathcal{STP}$. \\First, note that,  for each $\sigma$ in $\mathcal{STP}$, $\mathcal{F}_{cb_\sigma}$ is $1^\ast$-reducible if and only if $\mathcal{F}_{cb_\sigma} =\emptyset$ if and only if $\sigma = 1^\ast$ if and only if $\sigma\le 1^\ast$. \\Next, assume that we are given $\tau\neq 1^\ast$ in $\mathcal{STP}$ such that, for each $n$, for each $\sigma$ in $\mathcal{STP}$, if $\mathcal{F}_{cb_\sigma}$ is $\tau^n$-reducible, then $\sigma\le \tau^n$. 

Assume that we are given $\sigma$ such that $\mathcal{F}_{cb_\sigma}$ is $\tau$-reducible. \\Find $\varepsilon$ in $\mathcal{F}_{cb_\sigma}$,  such that for all $m,n$, \\if  $\varepsilon(m)\neq n$ and $\beta(\overline\varepsilon m\ast\langle n \rangle) =0$, then, for some $p$, $\mathcal{F}_{cb_{\sigma}}\cap\overline\varepsilon m\ast\langle n\rangle$ is  $\tau^p$-reducible.  \\  
Let $p$ be given. \\Consider $s:=\langle p+1\rangle$ and  $t:=\langle 0,p+1\rangle$ and note: \textit{either} $s\perp\varepsilon$ \textit{or} $t\perp\varepsilon$. \\Find $m$ such that \textit{either} $\mathcal{F}_{cb_\sigma}\cap \langle p+1\rangle =\langle p+1\rangle\ast\mathcal{F}_{\sigma^p}$ is $\tau^m$-reducible, \\\textit{or}  $\mathcal{F}_{cb_\sigma}\cap \langle 0, p+1\rangle =\langle 0, p+1\rangle\ast\mathcal{F}_{\sigma^p}$ is $\tau^m$-reducible. \\Conclude: $\mathcal{F}_{\sigma^p}$ is $\tau^m$-reducible and: $\sigma^p \le \tau^m$. \\Conclude: $\forall p\exists m[\sigma^p\le\tau^m]$ and: $\sigma \le \tau$.

\smallskip (ii)  We use induction on $\mathcal{STP}$. \\ By Lemma \ref{L:cbh}, for each $\sigma$ in $\mathcal{STP}$, for each $\delta$,  if $\mathcal{F}_{cb_\sigma}$ embeds into $\mathbb{P}(1^\ast, En_\delta)=En_\delta$, then  $\sigma \le S^\ast(1^\ast)$.  \\Next, assume that we are given $\tau\neq 1^\ast$ in $\mathcal{STP}$ such that, for each $n$, for each $\sigma$ in $\mathcal{STP}$, for each $\delta$, if $\mathcal{F}_{cb_\sigma}$ embeds into  $\mathbb{P}(\tau^n, En_\delta)$, then $\sigma\le S^\ast(\tau^n)$. \\Further assume that we are given $\sigma,\delta$ such that $\sigma\in\mathcal{STP}$ and $\mathcal{F}_{cb_\sigma}$ embeds into $\mathbb{P}(\tau, En_\delta)$. Find $\varphi$ embedding $\mathcal{F}_{cb_\sigma}$ into $\mathbb{P}(\tau, En_\delta)$.   \\ Note: $\forall \gamma \in \mathcal{F}_{cb_\sigma}\exists p[\varphi|\gamma \;\#\;\delta^p\rightarrow \exists n[\varphi|\gamma \in \mathbb{P}(\tau^n, En_\delta)]]$. \\Using Brouwer's Continuity Principle $\mathbf{BCP}$, see Subsubsection \ref{SSS:bcpcontchoice}, \\find $m,p$ such that $\forall \gamma \in \mathcal{F}_{cb_\sigma}[(\overline{\underline 0}m\sqsubset\gamma\;\wedge\;\varphi|\gamma \;\#\;\delta^p)\rightarrow \exists n[\varphi|\gamma \in \mathbb{P}(\tau^n, En_\delta)]]$. \\ Consider $\gamma_0:=\overline{\underline 0}m\ast\langle p+1\rangle\ast\underline 0$ and  $\gamma_1:=\overline{\underline 0}(m+1)\ast\langle p+1\rangle\ast\underline 0$. \\Note $\varphi|\gamma_0\;\#\;\varphi|\gamma_1$  and find $i<2$ such that $\varphi|\gamma^i\;\#\;\delta^p$. Find $j,n$ such that $\varphi|\overline{\gamma_i}j\perp\overline{\delta^p}n$. Note: $\forall \gamma \in \mathcal{F}_{cb_\sigma}\cap \overline{\gamma^i}j \exists i[\varphi|\gamma\in \mathbb{P}(\tau^i, En_\delta)]$. Using $\mathbf{BCP}$ again, \\find $k,l$ such that $k>j$ and $\forall \gamma \in\mathcal{F}_{cb_\sigma}[\overline{\gamma_i}k\sqsubset\gamma\rightarrow \varphi|\gamma\in\mathbb{P}(\tau^l,En_\delta)]$. \\Note: $\mathcal{F}_{cb_{\sigma^p}}$ embeds into $\mathcal{F}_{cb_\sigma}\cap\overline{\gamma_i}k$ and $\varphi$ embeds $\mathcal{F}_{cb_\sigma}\cap\overline{\gamma_i}k$ into $\mathbb{P}(\tau^l, En_\delta)$.  \\Conclude: $\mathcal{F}_{cb_{\sigma^p}}$ embeds into $\mathbb{P}(\tau^l, En_\delta)$, and: $\sigma^p \le S^\ast(\tau^l)$. \\Conclude: $\forall p\exists l[\sigma^p\le S^\ast(\tau^l)\le \tau = (S^\ast(\tau)^l]$ and: $\sigma \le S^\ast(\tau)$.

\smallskip (iii) Note: $\exists\delta[\mathcal{FIN}=En_\delta]$ and apply (ii).

\smallskip (iv) This is an immediate consequence of (ii).
\end{proof}

\section{The second level and the collapse of the projective hierarchy}
\subsection{The classes $\mathbf{\Sigma}^1_2$ and $\mathbf{\Pi}^1_2$}\label{SSS:classeslevel2}\hfill

Some relevant  definitions may be found in Subsubsection \ref{SSS:projective}.

\begin{definition} $\mathcal{X}\subseteq \omega^\omega$ is $\mathbf{\Sigma}^1_2$ if and only if there exists
 $\beta$ such that \\$\mathcal{X}=\mathcal{E}\mathcal{U}\mathcal{G}_\beta :=Ex\bigl(Un(\mathcal{G}_\beta)\bigr)=\{\alpha\mid \exists\delta \forall \gamma[\ulcorner\ulcorner \alpha, \gamma\urcorner, \delta\urcorner \in \mathcal{G}_\beta\}$.   

$\mathcal{X}\subseteq \omega^\omega$ is $\mathbf{\Pi}^1_2$ if and only if there exists $\beta$ such that \\$\mathcal{X}=\mathcal{UEF}_\beta:=Un\bigl(Ex(\mathcal{F}_\beta)\bigr)=\{\alpha\mid \forall\delta \exists \gamma[\ulcorner\ulcorner \alpha, \gamma\urcorner, \delta\urcorner \in \mathcal{F}_\beta\}$.

\end{definition}

Let $\beta,\varepsilon,\zeta$ be given such that $\varepsilon \in \mathcal{EUG}_\beta$ and $\zeta \in \mathcal{UEF}_\beta$.  \\Find $\delta$ such that $\forall \gamma \exists n[\beta(\overline{\ulcorner\ulcorner\varepsilon, \gamma\urcorner,\delta\urcorner})n\neq 0]$. \\Find $\gamma$ such that 
$\forall n[\beta(\overline{\ulcorner\ulcorner\zeta, \gamma\urcorner,\delta\urcorner}n)=0]$. \\Find $n$ such that $\beta(\overline{\ulcorner\ulcorner\varepsilon, \gamma\urcorner,\delta\urcorner}n)\neq 0]$ and conclude: $\overline \varepsilon n \neq \overline \zeta n$ and $\varepsilon\;\#\;\zeta$.

We thus see that, for each $\beta$,  $\mathcal{EUG}_\beta\;\#\;\mathcal{UEF}_\beta$. 

\medskip The next Theorem shows some properties of the classes $\mathbf{\Sigma}^1_2$ and $\mathbf{\Pi}^1_2$. 
Note that we do  not prove that the class $\mathbf{\Pi}^1_2$ is closed under the operation of countable union or even under the operation of finite union.

\smallskip\begin{theorem}\label{T:secondlevel} \begin{enumerate}[\upshape (i)]\hfill

 \item $\mathcal{US}^1_2:=\{\alpha\mid \alpha_{II}\in \mathcal{EUG}_{\alpha_I}\} $ is $\mathbf{\Sigma}^1_2$-universal and \\$\mathcal{UP}^1_2:=\{\alpha\mid \alpha_{II}\in \mathcal{UEG}_{\alpha_I})\} $ is $\mathbf{\Pi}^1_2$-universal.
\item $\mathcal{E}^1_2:=\{\alpha\mid \exists \delta\forall \gamma\exists n[\alpha(\overline{\ulcorner  \gamma, \delta\urcorner}n)\neq 0]\}$ is $\mathbf{\Sigma}^1_2$-complete and\\ $\mathcal{A}^1_2:=\{\alpha\mid \forall \delta\exists \gamma\forall n[\alpha(\overline{\ulcorner  \gamma, \delta  \urcorner}n)= 0]\}$ is $\mathbf{\Pi}^1_2$-complete.

\item $\mathbf{\Sigma}^1_2$ is closed under the operations of countable union and countable intersection and  $\mathbf{\Pi}^1_2$ is closed under the operation of countable intersection:

$\forall\beta\exists\varepsilon\exists\zeta[\bigcup_m\mathcal{EUG}_{\beta^m} =\mathcal{EUG}_\varepsilon \;\wedge\;\bigcap_m\mathcal{UEF}_{\beta^m} =\mathcal{UEF}_\varepsilon \;\wedge\; 
\bigcap_m\mathcal{EUG}_{\beta^m} =\mathcal{EUG}_\zeta]$.

\item For all $\mathcal{X}\subseteq \omega^\omega$, if $\mathcal{X}\in \mathbf{\Sigma}^1_2$, then $Ex(\mathcal{X})\in \mathbf{\Sigma}^1_2$, and, if $\mathcal{X}\in \mathbf{\Pi}^1_2$, then $Un(\mathcal{X})\in \mathbf{\Pi}^1_2$:  

$\forall\beta\exists\eta[Ex(\mathcal{EUG}_\beta)=\mathcal{EUG}_\eta\;\wedge\;Un(\mathcal{UEF}_\beta)=\mathcal{UEF}_\eta]$.

\item For all $\mathcal{X},\mathcal{Y}\subseteq\omega^\omega$ such that  $\mathcal{X}\preceq\mathcal{Y}$, \\if $\mathcal{Y}\in \mathbf{\Sigma}^1_2$, then $\mathcal{X} \in \mathbf{\Sigma}^1_2$, and, if $\mathcal{Y}\in \mathbf{\Pi}^1_2$, then $\mathcal{X} \in \mathbf{\Pi}^1_2$:

$\forall \beta\forall \varphi:\omega^\omega\rightarrow \omega^\omega\exists\gamma[\{\alpha\mid\varphi|\alpha\in\mathcal{EUG}_\beta\}=\mathcal{EUG}_\gamma\;\wedge\;\{\alpha\mid\varphi|\alpha\in\mathcal{UEF}_\beta\}=\mathcal{UEF}_\gamma]$.

\item $\mathbf{\Sigma}^1_1\cup\mathbf{\Pi}^1_1\subseteq\mathbf{\Sigma}^1_2\cap\mathbf{\Pi}^1_2$. 

\end{enumerate}

\end{theorem}

\begin{proof}
(i) Note: for each $\alpha$, $\alpha \in \mathcal{US}^1_2\leftrightarrow \alpha_{II}\in \mathcal{EUG}_{\alpha_I}\leftrightarrow \exists \delta[\ulcorner\alpha_{II}, \delta\urcorner\in \mathcal{UG}_{\alpha_I}]\leftrightarrow \exists\delta\forall\gamma[\ulcorner\ulcorner\alpha_{II},\delta\urcorner,\gamma\urcorner \in \mathcal{G}_{\alpha_I}]\leftrightarrow \exists \delta\forall\gamma \exists n[\alpha_I(\overline{\ulcorner\ulcorner\alpha_{II},\delta\urcorner,\gamma\urcorner}n)\neq 0]$. Define $\beta$ such that, 
  for all $a,c,d$, if $length(a)=length(d)=length(c)$, then $\beta(\ulcorner\ulcorner a, d\urcorner,c\urcorner)=a_I(\ulcorner\ulcorner a_{II},d\urcorner,c\urcorner)$. \\Note: for all $\alpha$, $\alpha \in \mathcal{US}^1_2\leftrightarrow \exists\delta\forall\gamma\exists n[\beta(\overline{\ulcorner\ulcorner \alpha, \delta\urcorner,\gamma\urcorner}n)\neq 0]\leftrightarrow \alpha \in \mathcal{EUG}_\beta$. \\Conclude: $\mathcal{US}^1_2\in \mathbf{\Sigma}^1_2$. \\Also note: for each $\varepsilon$, $\mathcal{EUG}_\varepsilon=\mathcal{US}^1_2\upharpoonright\varepsilon$ We thus see: $\mathcal{US}^1_2$ is $\mathbf{\Sigma}^1_2$-universal. 

Similarly, for each $\alpha$, $\alpha \in \mathcal{UP}^1_2\leftrightarrow \forall \delta\exists\gamma \forall n[\alpha_I(\overline{\ulcorner\ulcorner\alpha_{II},\delta\urcorner,\gamma\urcorner}n)= 0]]$. \\Define $\beta$ as above and conclude: $\mathcal{UP}^1_2=\mathcal{UEF}_\beta \in \mathbf{\Pi}^0_2$.\\ Note: for each $\varepsilon$, $\mathcal{UEF}_\varepsilon=\mathcal{UP}^1_2\upharpoonright\varepsilon$. We thus see: $\mathcal{UP}^1_2$ is $\mathbf{\Sigma}^1_2$-universal.

\smallskip (ii) Define $\beta$ such that, for all  $a,c, d$,  $ \beta(a,c, d)\neq 0$ if and only if  \\$length(a)=length(c)=length(d)>0$ and $\exists i<length(a)[a(\ulcorner \overline c i, \overline d i\urcorner)\neq 0]$.  \\Note: for each $\alpha$, $\exists \delta \forall \gamma \exists n[\alpha(\ulcorner\overline\gamma n, \overline \delta n\urcorner)\neq 0]$ if and only if $ \exists \delta \forall \gamma\exists n[\beta(\overline{\ulcorner \gamma, \delta\urcorner}n)\neq0]$, \\and:   $\forall \delta \exists \gamma \forall n[\alpha(\ulcorner\overline\gamma n, \overline \delta n\urcorner)=0]$ if and only if $\forall \delta \exists \gamma \forall n[\beta(\overline{\ulcorner \gamma, \delta\urcorner}n)=0]$.\\ Conclude:  $\mathcal{E}^1_2=\mathcal{EUG}_\beta\in\mathbf{\Sigma}^1_2$ and $\mathcal{A}^1_2=\mathcal{UEF}_\beta\in\mathbf{\Pi}^1_2$. 

Let $\varepsilon$ be given. Define $\varphi:\omega^\omega\rightarrow\omega^\omega$ such that, for all $\alpha$, for all $c,d$, \\if  $length(c)=length(d)$, then $ (\varphi|\alpha)(\ulcorner c, d \urcorner) = \varepsilon(\ulcorner\overline\alpha n, c\urcorner, d\urcorner)]$.\\Note: for all $\alpha$, $\exists\gamma\forall\delta\exists n[\varepsilon(\overline{\ulcorner\alpha, \gamma\urcorner,\delta\urcorner}n)\neq 0]$ if and only if $\exists\gamma\forall\delta\exists n[(\varphi|\alpha)(\ulcorner\overline \gamma n, \overline\delta n\urcorner)\neq 0]$, \\i.e. $\alpha \in \mathcal{EUG}_\varepsilon$ if and only if $\varphi|\alpha \in \mathcal{E}^1_2$, and:
\\$\forall\gamma\exists\delta\forall n[\varepsilon(\overline{\ulcorner\alpha, \gamma\urcorner,\delta\urcorner}n)= 0]$ if and only if $\forall\gamma\exists\delta\forall n[(\varphi|\alpha)(\ulcorner\overline \gamma n, \overline\delta n\urcorner)= 0]$,\\ i.e. $\alpha \in \mathcal{UEF}_\varepsilon$ if and only if $ \varphi|\alpha \in \mathcal{A}^1_2$. \\We thus see that $\varphi$ reduces the pair $(\mathcal{EUG}_\varepsilon, \mathcal{UEF}_\varepsilon)$ to the pair $(\mathcal{E}^1_2, \mathcal{A}^1_2)$. 

We may conclude that $\mathcal{E}^1_2$ is $\mathbf{\Sigma}^1_2$-complete and that $\mathcal{A}^1_2$ is $\mathbf{\Sigma}^1_2$-complete.

\smallskip (iii) Let $\beta$ be given.

For each $\alpha$, $\alpha \in \bigcup_m\mathcal{UEG}_{\beta^m}$ if and only if $\exists m\exists \delta\forall\gamma\exists n[\beta^m(\overline{\ulcorner \ulcorner \alpha, \gamma\urcorner,\delta\urcorner}n)\neq 0]$. \\Define $\varepsilon$ such that,  for all $m, a,c, d$, 
$\varepsilon(\ulcorner\ulcorner  a, c \urcorner, \langle m\rangle\ast d\urcorner)=\beta^m(\ulcorner\ulcorner  a, c \urcorner, d\urcorner)]$, and\\ 
$\beta(\ulcorner\ulcorner  0, 0 \urcorner, 0\urcorner)=0$. Note that, for each $m$, for all $\alpha$, $\gamma$,  $\delta$,  \\$\ulcorner \ulcorner\alpha,\gamma\urcorner, \langle m\rangle \ast \delta\urcorner  \in \mathcal{G}_\varepsilon$ if and only if $\ulcorner \ulcorner\alpha, \gamma\urcorner,\delta\urcorner \in \mathcal{G}_{\beta^m}$. \\Therefore, for each $\alpha$, $\alpha \in \mathcal{EUG}_\varepsilon$ if and only if $\exists m[\alpha \in \mathcal{EUG}_{\beta^m}]$ and:\\ $\mathcal{EUG}_\varepsilon =\bigcup_m\mathcal{EUG}_{\beta^m}$.  

Also note that, for each $m$, forall $\alpha, \gamma, \delta$, \\$\ulcorner \ulcorner \alpha, \gamma\urcorner, \langle m\rangle \ast \delta\urcorner  \in \mathcal{F}_\varepsilon$ if and only if $\ulcorner \ulcorner\alpha, \gamma\urcorner, \delta\urcorner \in \mathcal{F}_{\beta^m}$. \\Therefore, for each $\alpha$, $\alpha \in \mathcal{UEF}_\varepsilon$ if and only if  $\forall m[\alpha \in \mathcal{UEF}_{\beta^m}]$, \\ i.e.:  $\mathcal{UEF}_\varepsilon =\bigcap_m\mathcal{UEF}_{\beta^m}$.

Also, for each $\alpha$, \\$\alpha \in \bigcap_m\mathcal{UEG}_{\beta^m}$ if and only if $\forall m\exists \delta\forall\gamma\exists n[\beta^m(\overline{\ulcorner \ulcorner \ulcorner\alpha, \gamma\urcorner,\delta\urcorner}n)\neq 0]$. \\
Then, by $\mathbf{AC}_{0,1}$, $\alpha \in \bigcap_m\mathcal{UEG}_{\beta^m}$ if and only if $\exists \delta\forall m\forall\gamma\exists n[\beta^m(\overline{\ulcorner \ulcorner \alpha, \gamma\urcorner,\delta^m\urcorner}n)\neq 0]$\\ if and only if 
$\exists \delta\forall\gamma\exists n[\beta^{\gamma(0)}(\overline{\ulcorner \ulcorner \alpha, \gamma\circ S\urcorner,\delta^{\gamma(0)}\urcorner}n)\neq 0]$.\\
Define $\zeta$ such that, for all $a,c,d$,   $\zeta(\ulcorner \ulcorner a, c\urcorner,d\urcorner)\neq 0$ if and only if \\$ length(a)=length(c)=length(d)>0$ and  \\$\exists i\le length(a)[\beta^{c(0)}(\ulcorner \ulcorner \overline a i, \overline{c\circ S} i \urcorner,\overline{d^{c(0)}}i \urcorner)\neq 0]$. 
\\Note that,  for all $\alpha, \delta$, \\ $\forall \gamma \exists n [\beta^{\gamma(0)}(\overline{\ulcorner \ulcorner \alpha, \gamma\circ S\urcorner,\delta^{\gamma(0)}\urcorner}n)\neq 0]$ if and only if 
$\forall \gamma \exists n[\ulcorner\ulcorner\alpha, \gamma\urcorner, \delta\urcorner  \in \mathcal{G}_\zeta]$.  
\\Conclude: for all $\alpha$, $\alpha\in \bigcap_m \mathcal{UEG}_{\beta^m}$ if and only if $\alpha \in \mathcal{UEG}_\zeta$, \\

i.e.: $ \mathcal{EUG}_\zeta=\bigcap_m\mathcal{EUG}_{\beta^m}$.

 \smallskip (iv)  Let $\beta$ be given. Note: for all $\alpha$, \\$\alpha \in  Ex(\mathcal{EUG}_\beta) $ if and only if $\exists\varepsilon\exists\delta\forall \gamma\exists n[\beta(\overline{\ulcorner\ulcorner\ulcorner\alpha,\gamma\urcorner,\delta\urcorner,\varepsilon\urcorner}n)\neq 0]$, and:  \\$\alpha \in Un(\mathcal{UEF}_\beta)$ if and only if $ \forall\varepsilon\forall\delta \exists\gamma\forall n[\beta(\overline{\ulcorner\ulcorner\ulcorner\alpha,\gamma\urcorner,\delta\urcorner,\varepsilon\urcorner}n)= 0]$. \\Define $\eta$ such that, for all $a,c,d$, if  $length(a)=length(c)=length(d)$, then \\$\eta(\ulcorner\ulcorner a,c\urcorner, d\urcorner)=\beta(\ulcorner\ulcorner\ulcorner a, c\urcorner, d_I\urcorner,d_{II})].$  \\One easily verifies: $Ex(\mathcal{EUG}_\beta)=\mathcal{EUG}_\eta$ and: $Un(\mathcal{UEF}_\beta)=\mathcal{UEF}_\eta$. 
 
 \smallskip (v) Let $\beta,\varphi$ be given such that $\varphi:\omega^\omega\rightarrow\omega^\omega$. Note that, for each $\alpha$, \\$\varphi|\alpha \in \mathcal{EUG}_\beta $ if and only if $\exists\delta\forall\gamma\exists n[\overline{\ulcorner\ulcorner\varphi|\alpha,\gamma\urcorner,\delta\urcorner}n)\neq 0]$ and:\\
 $\varphi|\alpha \in \mathcal{UEF}_\beta $ if and only if $\forall\delta\exists\gamma\forall n[\overline{\ulcorner\ulcorner\varphi|\alpha,\gamma\urcorner,\delta\urcorner}n)= 0]$. \\Define $\varepsilon$ such that for all $a,c,d$ if $length(a)=length(c)=length(d)$, then \\$\varepsilon(\ulcorner\ulcorner a, c\urcorner, d \urcorner)\neq 0$ if and only if $\exists i[length(\varphi|a)\ge i \;\wedge\;\beta(\ulcorner \ulcorner\overline{\varphi|a}i,\overline c i\urcorner, \overline d i\urcorner)\neq 0]$. 
 
 Then: $\{\alpha\mid\varphi|\alpha \in \mathcal{UEG}_\beta\}=\mathcal{UEG}_\varepsilon$ and: $\{\alpha\mid\varphi|\alpha \in \mathcal{EUF}_\beta\}=\mathcal{EUF}_\varepsilon$. \end{proof}

\subsection{The collapse of the projective hierarchy}

\begin{theorem}\label{T:collapse} \hfill

\begin{enumerate}[\upshape(i)]
  
  \item For all $\mathcal{X}\subseteq\omega^\omega$, if $\mathcal{X}\in\mathbf{\Sigma}^1_2$, then
$Un(\mathcal{X})\in\mathbf{\Sigma}^1_2$:
$\forall \beta\exists\varepsilon[Un(\mathcal{EUG}_\beta)=\mathcal{EUG}_\varepsilon]$. 
\item $\mathbf{\Pi}^1_2\subseteq \mathbf{\Sigma}^1_2$, and for all $\mathcal{X}\subseteq \omega^\omega$ , if $\mathcal{X}$ is  (positively) projective,  then $\mathcal{X}\in\mathbf{\Sigma}^1_2$.
\end{enumerate}
\end{theorem}

\begin{proof}

(i) Let $\beta$ be given. Using $\mathbf{AC}_{1,1}$,  note: for all $\alpha$, \\$\alpha \in Un(\mathcal{EUG}_\beta)$ if and only if  $\forall \varepsilon \exists \delta\forall \gamma\exists n[\beta(\overline{\ulcorner\ulcorner\ulcorner\alpha,\gamma\urcorner, \delta\urcorner ,\varepsilon\urcorner}n)\neq 0]$ if and only if  \\$\exists\varphi[\varphi\in\mathcal{A}_1^1\;\wedge\;\varphi(0)=0\;\wedge\;\forall \varepsilon \forall \gamma\exists n[\beta(\overline{\ulcorner\ulcorner\ulcorner\alpha,\gamma\urcorner, \varphi|\varepsilon\urcorner ,\varepsilon\urcorner}n)\neq 0]] $ if and only if  \\$ \exists\varphi[\varphi\in\mathcal{A}_1^1\;\wedge\;\varphi(0)=0\;\wedge\\ \forall\varepsilon\forall \gamma\exists n\exists m[length(\varphi|\overline{\varepsilon}m)\ge n\;\wedge\;\beta(\ulcorner\ulcorner\ulcorner\overline \alpha n,\overline \gamma n\urcorner,\overline{(\varphi|\overline{\varepsilon}m)}n\urcorner, \overline{\varepsilon}n\urcorner)\neq 0]].$

Using Theorem \ref{T:secondlevel}, we conclude: $Un(\mathcal{EUG}_\beta) \in \mathbf{\Sigma}^1_2$.

(ii) This follows from (i). 
\end{proof}

\subsubsection{}  Theorem \ref{T:collapse} shows that, in intuitionistic mathematics,
$\mathbf{\Sigma}^1_2$ is the class of all positively projective sets. \\Many difficult questions remain, for instance, if
$\mathbf{\Pi}^1_2$ is a proper subclass of $\mathbf{\Sigma}^1_2$ and if the class $\mathbf{\Pi}^1_2$ is
closed under the operation of disjunction. We were unable to answer these
questions.

Note that the  projection of a positively Borel set is analytic. It is not true however, that the co-projection of
a positively Borel set is always co-analytic, for the simple reason that some positively Borel sets, like $\mathbb{D}^2(\mathcal{A}_1)$ \footnote{See Theorem \ref{T:coan}(iv).},
 are not  co-analytic.

\begin{lemma}\label{L:a12e12}
$\forall \varphi:\omega^\omega\rightarrow\omega^\omega\exists \alpha[(\alpha \in\mathcal{E}^1_2\leftrightarrow \varphi|\alpha\in \mathcal{E}^1_2)\;\wedge\;(\alpha \in\mathcal{A}^1_2\leftrightarrow \varphi|\alpha\in \mathcal{A}^1_2)]$.
\end{lemma}

\begin{proof} Let $\varphi:\omega^\omega\rightarrow\omega^\omega$ be given. Define $\alpha$ such that for all $p, c, d$, if $length(c)=length(d)$ and $ p=\ulcorner c, d\urcorner$, then $\alpha(p)\neq 0$ if and only if, for some $ m\le length(c)$, $\ulcorner\overline c m,\overline d m\urcorner < length(\varphi|\overline \alpha p)$ and $(\varphi|\overline\alpha p)(\ulcorner\overline c m, \overline d m\urcorner) \neq 0]$. 

Note that, for all $\gamma, \delta$, $\exists m[\alpha(\overline{\ulcorner\gamma,\delta\urcorner}m)\neq 0]$ if and only if $\exists m[(\varphi|\alpha)(\overline{\ulcorner\gamma,\delta\urcorner}m)\neq 0]$.

Conclude: $\exists\gamma\forall \delta\exists n[\alpha(\overline{\ulcorner\gamma,\delta\urcorner}n)\neq 0] \leftrightarrow \exists\gamma\forall \delta\exists n[(\phi|\alpha)(\overline{\ulcorner\gamma,\delta\urcorner}n)\neq 0]$, that is: $\alpha\in\mathcal{E}^1_2\leftrightarrow \varphi|\alpha\in\mathcal{E}^1_2$, and also:  $\forall \gamma\exists \delta\forall n[\alpha(\overline{\ulcorner\gamma,\delta\urcorner}n)= 0] \leftrightarrow \forall \gamma\exists \delta\forall n[(\varphi|\alpha)(\overline{\ulcorner\gamma,\delta\urcorner}n)= 0]$, that is: $\alpha\in\mathcal{A}^1_2\leftrightarrow \varphi|\alpha\in\mathcal{A}^1_2$.
\end{proof}

Note that the classical mathematician would conclude, from Lemma \ref{L:a12e12}: $\mathcal{A}^1_2\npreceq\mathcal{E}^1_2$ and $\mathcal{E}^1_2\npreceq\mathcal{A}^1_2$.

\begin{theorem}\label{T:ruin} \hfill
\begin{enumerate}[\upshape(i)]
\item $\exists\alpha[\alpha \notin \mathcal{E}^1_2\;\wedge\;\alpha \notin \mathcal{A}^1_2].$
\item  $\exists \gamma[\gamma\notin \mathcal{US}^1_2\;\wedge\;\gamma \notin \mathcal{UP}^1_2]$.

\end{enumerate}
\end{theorem}

\begin{proof}
(i) Using Theorems \ref{T:collapse}(i) and \ref{T:secondlevel}(ii), find $\varphi:\omega^\omega\rightarrow\omega^\omega$ reducing $\mathcal{A}^1_2$ to $\mathcal{E}^1_2$.  
Applying Lemma \ref{L:a12e12}, find $\alpha$ such that $ \alpha\in\mathcal{E}^1_2\leftrightarrow \varphi|\alpha\in\mathcal{E}^1_2$ and $\alpha\in\mathcal{A}^1_2\leftrightarrow \varphi|\alpha\in\mathcal{A}^1_2$. 
\\Assume: $\alpha\in\mathcal{E}^1_2$. Conclude: $\varphi|\alpha \in \mathcal{A}^1_2$ and: $\alpha \in \mathcal{A}^1_2$. Contradiction, as $\mathcal{A}^1_2\;\#\;\mathcal{E}^1_2$.
\\Conclude: $\alpha \notin \mathcal{E}^1_2$ and: $\varphi|\alpha \notin \mathcal{E}^1_2$ and: $\alpha \notin \mathcal{A}^1_2$. 

\smallskip  (ii) Define $\mathcal{DP}^1_2:=\{\alpha\mid
\ulcorner\alpha,\alpha\urcorner\in \mathcal{UP}^1_2\}$. According to Theorem \ref{T:collapse}(i), 
$\mathcal{DP}^1_2\in\mathbf{\Sigma}^1_2$.  Using Theorem \ref{T:secondlevel}(iii), find $\beta$ such that $\mathcal{DP}^1_2=\mathcal{US}^1_2{\upharpoonright}\beta$.  \\Note: for every $\alpha$, 
$\ulcorner\alpha,\alpha\urcorner\in \mathcal{UP}^1_2\leftrightarrow \alpha\in\mathcal{DP}^1_2\leftrightarrow\ulcorner\beta,\alpha\urcorner \in\mathcal{US}^1_2$. \\Define
$\gamma:=\ulcorner\beta,\beta\urcorner$.  and note: $\gamma \notin \mathcal{US}^1_2$ and $\gamma \notin \mathcal{UP}^1_2$, as $\mathcal{US}^1_2\;\#\;\mathcal{UP}^1_2$.
\end{proof}

\bigskip\noindent
  Theorem \ref{T:ruin} has some noteworthy consequences. \\Assume:   $\alpha\notin\mathcal{E}^1_2\cup\mathcal{A}^1_2$. Then: \begin{enumerate}[\upshape (i)]
 \item $\neg \exists \delta \forall \gamma \exists n [
\alpha(\ulcorner \overline \gamma n, 
\overline \delta n \urcorner) \neq 0$, and 
\item $\neg \forall \delta \exists \gamma \forall n [  \alpha (\ulcorner\overline \gamma n, \overline
\delta n \urcorner ) = 0 ]$, and 
\item $\forall \delta \forall \gamma \forall n[\ulcorner \alpha (\overline \gamma n, \overline
\delta n \urcorner ) = 0  \; \vee \; \alpha (\ulcorner\overline \gamma n, \overline
\delta n \urcorner ) \neq 0 ]$.
\end{enumerate}

  Theorem \ref{T:ruin} thus shows that,  in intuitionistic mathematics it is
possible that statements \begin{enumerate}[\upshape (i)]
\item
 $\mathsf{\neg\exists x\forall y\exists z[P(x,y,z)]}$, and
\item $\mathsf{\neg\forall x\exists y\forall z[\neg P(x,y,z)]}$,  and 
\item $\mathsf{\forall x\forall y\forall
z[P(x,y,z)\;\vee\;\neg P(x,y,z)]}$, \end{enumerate}
are simultaneously true. The example depends on $\mathbf{AC}_{1,1}$.  
Another example, depending only on $\mathbf{BCP}$, has been given in \cite[Section 5.5]{veldman08}:

\begin{enumerate}[\upshape (i)]
 \item $\neg \exists \alpha \forall n \exists m[\alpha(n) = 0 \; \wedge \;\alpha(m) \neq 0]$, and
  \item $\neg \forall \alpha \exists n \forall m[\alpha(n) \neq 0 \; \vee \; \alpha(m) = 0]$, and 
\item $\forall \alpha \forall n \forall m[\bigl( \alpha(n) = 0 \; \wedge \;\alpha(m) \neq 0 \bigr) \; \vee \; \bigl(\alpha(n) \neq 0 \; \vee \; \alpha(m)=0 \bigr)]$.
\end{enumerate}

\subsection{A parallel: the collapse of the (positive) arithmetical hierarchy}\hfill

It has been observed by J.R.~Moschovakis that, in the context of intuitionistic arithmetic,  Church's Thesis $\mathbf{CT}$ causes the collapse of the (positively) arithmetical hierarchy, just as $\mathbf{AC}_{1,1}$ causes the collapse of the (positively) projective hierarchy, see \cite{moschovakis02} and \cite{moschovakis05}. It seems useful to explain this.

 Let $T\subseteq \omega^3$ be Kleene's $T$-predicate. $T$ is a \textit{(Kalm\'ar-)elementary} subset of $\omega^3$ and, for all $e,n,z$, $T(e,n,z)$ stands for: \textit{`$z$ is the code of a succesful computation according to the  algorithm coded by $e$ at the argument $n$}'. Let $U$ be the elementary function from $\omega$ to $\omega$ extracting from each succesful computation $z$ its outcome $U(z)$. Every $e$ determines a \textit{partial} function $\varphi_e$ from $\omega$ to $\omega$ by: $$\forall n[\varphi_e(n)\simeq U(\mu z[T(e,n,z])].$$  For each $e$, $W_e:=\{n\mid\exists z[T(e,n,z)]\}$ is the domain of the partial function $\varphi_e$. 
 
For every $X\subseteq\omega$, we define the {\it projection} $Ex_0(X):=\{m\mid\exists n[\langle m, n\rangle\in X]\}$ and the {\it co-projection} $Un_0(X):=\{m\mid\forall n[\langle m, n \rangle\in X]\}$.

One defines $\Sigma^0_1:=\{W_e\mid e \in \omega\}$ and 
 $\Pi^0_1:=\{\omega\setminus W_e\mid e \in \omega\}$, and, for each $m>0$, $\Sigma^0_{m+1} :=\{Ex_0(X)\mid X\in \Pi^0_m\}$ and  $\Pi^0_{m+1} :=\{Un_0(X)\mid X\in \Sigma^0_m\}$.
 
 One may prove: $\forall m>0[\Sigma^0_m\cup\Pi^0_m\subseteq\Sigma^0_{m+1}\cap\Pi^0_{m+1}]$.
 
  Using the following strong form of \textit{Church's Thesis} $\mathbf{CT}$: for every $R\subseteq\omega\times\omega$, $$\forall m\exists n[mRn]\rightarrow\exists e\forall m\exists z[T(e,m,z)\;\wedge\; mRU(z)],$$ one may prove that,  for every $X$ in $\Sigma^0_3$, also $Un_0(X)\in \Sigma^0_3$, as follows:
 
Assume $X\in \Sigma^0_3$. Find $e$ such that $X=Ex_0\bigl(Un_0( W_e)\bigr)$. \\Consider $Y=Un_0(X)=\{m\mid \forall q[\langle m , q\rangle \in X]\} =\\\{m\mid\forall q \exists n\forall p\exists z[T(e, \langle m,p,n,q\rangle,z)]\}=\\\{ m\mid\exists f\forall q\forall p\exists u\exists z[T(f, q, u)\;\wedge\; T(e, \langle m, p, U(u),q\rangle,z)]\} \in \Sigma^0_3$.

One may conclude: $\Pi^0_3\subseteq \Sigma^0_3$ and: $\bigcup_m \Sigma^0_m\subseteq \Sigma^0_3$.

\smallskip
Find $f$ such that $\{e\mid \forall p\exists n\forall z[\neg T(e, \langle e, n,p\rangle,z)]\}=\{ m\mid \exists p\forall n\exists z[T(f, \langle m,n,p\rangle, z)]\}$, and note: 
$\forall p\exists n\forall z[\neg T(f, \langle f, n,p\rangle,z)]\leftrightarrow \exists p\forall n\exists z[T(f, \langle f,n,p\rangle, z)]$, and, therefore: 

$\neg \forall p\exists n\forall z[\neg T(f, \langle f, p, n\rangle,z)]$ and $ \neg\exists p\forall n\exists z[T(f, \langle f,p,n\rangle, z)]$.

Again, we see that  three statements of the form
\begin{enumerate}[\upshape (i)]
 \item $\mathsf{\neg\exists x\forall y\exists z[P(x,y,z)]}$, and
\item $\mathsf{\neg\forall x\exists y\forall z[\neg P(x,y,z)]}$,  and 
\item $\mathsf{\forall x\forall y\forall
z[P(x,y,z)\vee\neg P(x,y,z)]}$,\end{enumerate} 
may be true simultaneously.

\end{document}